\documentclass[11pt,fleqn]{amsart}
\usepackage[cp1250]{inputenc}
\usepackage{latexsym}
\usepackage{amsmath}
\usepackage{amsfonts}
\usepackage{amscd}
\usepackage{amsthm}
\usepackage{amssymb}
\usepackage[T1]{fontenc}
\usepackage{amsrefs}
\usepackage[sans]{dsfont}
\usepackage{graphicx}
\usepackage{subcaption}
\usepackage{color}
\author[\.Zywilla~Fechner and Robert Stelzer]{\.Zywilla~Fechner $^{1,2}$ and Robert Stelzer $^2$}
\address{$^1$ Institute of Mathematics, University of Silesia, Bankowa 14, 40-007 Katowice, Poland\\$^2$ Institute of Mathematical Finance, Ulm University, Helmholtzstra\ss e 18, 
89069 Ulm, Germany}
\email{zfechner@gmail.com, robert.stelzer@uni-ulm.de}
\title[Limit behaviour of the truncated pathwise Fourier-transformation]{Limit behaviour of the truncated pathwise Fourier-transformation of L\'evy-driven CARMA processes for non-equidistant discrete time observations}
\newtheorem{theorem}{Theorem}[section]

\newtheorem{assumption}[theorem]{Assumption}
\newtheorem{lemma}[theorem]{Lemma}
\newtheorem{prop}[theorem]{Proposition}

\theoremstyle{remark}

\newtheorem{example}[theorem]{Example}

\newcommand{\N}{\mathbb{N}}

\newcommand{\R}{\mathbb{R}}
\newcommand{\C}{\mathbb{C}}

\newcommand{\E}{\mathbb{E}}
\newcommand{\Cov}{\mathbb{C}\mathrm{ov}}
\renewcommand{\P}{\mathbb{P}}
\renewcommand{\b}{\mathbf{b}}
\newcommand{\X}{\mathbf{X}}

\newcommand{\A}{\mathbf{A}}
\newcommand{\e}{\mathbf{e}}

\newcommand{\F}{\mathcal{F}}
\newcommand{\Var}{\operatorname{Var}}
\newcommand{\ARMApq}{\operatorname{ARMA}(p,q)}

\renewcommand{\[}{\left[} \renewcommand{\]}{\right]}
\renewcommand{\proof}{\noindent\textit{Proof. }}

\newcommand{\1}{\mathds{1}}
\allowdisplaybreaks[4]
\begin{document}

\begin{abstract} 
This paper considers a continuous time analogue of the classical autoregressive moving average processes, L\'evy-driven CARMA processes. 
First we describe limiting properties of the periodogram by means of the so-called truncated Fourier transform if observations are available continuously. 
The obtained results are in accordance with their counterparts from the discrete-time case. 
Then we discuss the numerical approximation of the truncated Fourier transform  based on non-equidistant high frequency data. 
In order to ensure convergence of the numerical approximation to the true value of the truncated Fourier transform a certain control on the maximal distance between observations and 
the number of observations is needed. We obtain both convergence to the continuous time quantity and asymptotic normality under a high-frequency infinite time horizon limit.
\end{abstract}
\maketitle
\section{Introduction}
%Discrete time series are widely used in many applications like econometrics, finance, weather forecast, astronomy among others. %---------------------
%-------------------------------------------------------------------------------------------------------------------------------------------------------------------------------------------------------
The classical autoregressive moving average process $\operatorname{ARMA}$ has been broadly discussed in the literature. For a comprehensive discussion see e.g. 
the monograph by Brockwell  and Davis \cite{BrockwellDavis}
and references therein. In the discrete time models we restrict ourselves to observations at fixed equidistant points in time. In many cases these observations made at discrete times come from an underlying continuous process,
thus the natural question arises: can we model also the time series in continuous time?
One of the earliest results dealing with properties of such processes can be found in Doob \cite{Doob}. Later this problem was discussed by Brockwell in
 \cite{Brockwell2001a} for continuous time $\operatorname{ARMA}$ processes driven by Gaussian noise. 
The next step was to extend these ideas to the models with noise modelled by jump processes, so-called L\'evy-driven 
$\operatorname{CARMA}$ models introduced by Brockwell in  \cite{Brockwell2001}. In these papers time series are modelled as continuous time processes with continuous time noises (with or without jumps)
and the inference is based mainly on discrete equidistant data. One of the latest results can be found in the paper \cite{BrockwellDavisYang} of Brockwell, Davis and Yang, 
which consideres QML estimations of the $\operatorname{AR}$ and $\operatorname{MA}$ parameters based on equidistant observations.

The estimation procedure of L\'evy-driven $\operatorname{CARMA}$ processes in high-frequency settings has been discussed by Fasen and Fuchs in \cite{FuchsFasen2013}, 
where the authors deal with the limit behaviour of the
periodogram of $\operatorname{CARMA}$ processes under equidistant sampling when the sampling interval tends to $0$. The results are analogous to $\operatorname{ARMA}$ processes: the periodogram for  
$\operatorname{CARMA}$ processes is not a consistent estimator of the spectral density function, but after appropriate smoothing the consistency can be obtained. 
%Using a smoothed periodogram 
%for equidistant observations the authors discussed a Whittle-type estimation in order to find the value of the parameters of $\operatorname{CARMA}$ models. 
Some related results were discussed by Fasen and Fuchs in \cite{FuchsFasen2013a},
where asymptotic distributions of periodograms of $\operatorname{CARMA}$ processes driven by a symmetric $\alpha$-stable L\'evy noise are obtained 
and where it is shown that the vector composed of periodograms for various frequencies converges in distribution to a function of a multidimensional stable random vector. Likewise, Fasen \cite{Fasen2013} considers the behaviour of the periodogram for an equidistantly sampled continuous time moving average process when only the number of observations goes to infinity.
%The statistical inference in the frequency domain for continuous time $\operatorname{MA}$ processes  was discussed by V. Fasen in \cite{Fasen2013}. 
%The moving average process is sampled at an equally spaced time grid in the frequency domain. In the case of finite second moments of the underlying L\'evy process the asymptotic behaviour of the sampled process 
%is described: the periodogram is not a consistent estimator for the spectral density of the sampled process. 
%Moreover, different periodogram frequencies are asymptotically independent and exponentially distributed like in the discrete time $\operatorname{ARMA}$ process.
%The asymptotic behaviour of the sample autocovariance of the multivariate $\operatorname{CARMA}$ models the was investigated by V. Fasen in \cite{Fasen2014}.

The problem of statistical analysis of such processes has been studied further for example by Gillberg in his dissertation \cite{Gillberg2006}, where
different approaches to the estimation of $\operatorname{CARMA}$ processes with Gaussian noise are discussed both using equidistant and non-equidistant observations. 
The author works mainly in the frequency domain. He describes the properties of the truncated Fourier transform of a $\operatorname{CARMA}$ process with Gaussian noise on a fixed interval $[0,T]$ 
based on equidistant frequencies. In the non-equidistant case he has used a method based on splines in order to find an approximation of the spectral density.

Another approach for the estimation of a zero-mean stationary process $(Y_t)_{t\in \R}$ with finite second-order moments and continuous covariance function 
has been discussed by Lii and Masry in \cite{LiiMasry1992} and \cite{LiiMasry1994},
where they described some properties of a smoothed periodogram. Here observations are assumed to be given on a random grid $(\tau_k)$ of an interval $[0,T]$, 
where $\tau_k$ is a stationary point process on the real line which is independent of $(Y_t)_{t\in \R}$.

In the present paper we are going to describe the asymptotic behaviour of the so-called truncated Fourier 
transform of a $\operatorname{CARMA}$ process, which is a building block for an estimation of the spectral density of a $\operatorname{CARMA}$ process. 
We are going to use some of the ideas from \cite{Gillberg2006} to prove results in more general settings.

The paper is structured as follows: first we recall second order L\'evy-driven $\operatorname{CARMA}$ models and summarize the results needed later in Section \ref{sec2}. 
Then we define in Section \ref{sec3} the truncated Fourier transform of a $\operatorname{CARMA}$ process and we investigate its asymptotic properties at a fixed frequency: for a non-zero frequency we obtain that the limiting 
law of the real and  imaginary part is the two dimensional normal distribution with mean zero and the covariance matrix depending on the spectral density of the $\operatorname{CARMA}$ process. If we consider the truncated Fourier transform at zero, 
we obtain  a one dimensional normal law with mean zero and variance depending only on two parameters of the $\operatorname{CARMA}$ process. 
We show that the limiting law of the joint distribution of the squared modulus of the truncated Fourier transform at different positive frequencies converges to a vector of independent and exponentially 
distributed random variables with mean depending on the values of the spectral density.
All these results can be interpreted as the limiting behaviour of the truncated Fourier transform when the $\operatorname{CARMA}$ process is observed continuously.
The next step in Section \ref{Section:ApproximationofIntegrals} is to approximate the truncated Fourier transform when the $\operatorname{CARMA}$ process is observed on a non-equidistant deterministic grid. 
%For fixed $T>0$ we consider the grid of interval $[0,T]$ consisting of $N$ elements. The maximal distance between the elements of the grid is denoted by 
%$h_{\max}(T)$. 
In order to find a numerical approximation value of the truncated Fourier transform we apply 
the trapezoidal rule. We are interested in the convergence of the truncated Fourier transform when the length of the interval $T$ goes to infinity and the mesh of the grid to zero.
Since the interplay of the length of the interval, of the number of elements of the grid and of the maximal distance between the elements of the grid 
plays a crucial role, in order to ensure the convergence of the approximating sum to the true value of the truncated Fourier transform we have to impose some limiting conditions on these quantities.
In the last Section \ref{sec4} we look at some illustrative simulations of the truncated Fourier transform based on non-equidistant observations. We consider Ornstein-Uhlenbeck type (CAR(1)) and CARMA(2,1) processes driven by a standard Brownian motion, a Variance Gamma process and a ``two-sided Poisson process'' and we compare our simulations with the theoretical asymptotic distributions described earlier.

\subsection*{Notation}
The symbol $\N:=\{1,2,3,\dots  \}$ denotes the set of positive integers, $\N_0:=\N\cup \{ 0 \}$, $\R$ is the set of real numbers and $\C$ denotes the set of complex numbers.
The symbol $\R^{m\times n}$, resp. $\C^{m\times n}$ denotes the space of real- (resp. complex-) valued matrices with $m$ rows and $n$ columns.
For $A\in \C^{m\times n}$ the symbol $A^T$ denotes the transposed of a matrix $A$. We are working on a given filtered  probability space $(\Omega, \mathcal{F}, (\mathcal{F}_t)_{t\geq 0}, \P)$ satisfying the usual hypothesis (cf. Protter \cite{Protter}, Chapter 1).

Moreover, by $X\stackrel{d}{=}Y$ we denote that the random variables $X$ and $Y$ are equal in distribution.
%-------------------------------------------------------------------------------------
%-------------------------------------------------------------------------------------
%-------------------------------------------------------------------------------------
%-------------------------Preliminaries----------------------------------------------------
%-------------------------------------------------------------------------------------
%-------------------------------------------------------------------------------------
%-------------------------------------------------------------------------------------
%--------------------------------------------------------------------------------------------------------------
\section{Preliminaries}\label{sec2}
%--------------------------------------------------------------------------------------------------------------
We begin with the model set-up given by Brockwell (see \cite{Brockwell2001},  \cite{Brockwell2001a}). 
A second-order L\'evy-driven continuous-time $\ARMApq$ process is defined in terms of a state-space representation of the formal differential equation
\begin{equation}\label{eq: ARMAFomalEquatioab}
a(D)Y(t) = b(D)DL(t),\quad t  \geq 0.
\end{equation}
%Second Moments Assumptions CF Brockwell 2001 Def 1.1
Here, $D$ denotes differentiation with respect to $t$, non-negative integers $p,q$ satisfying $p>q$ and $(L(t))_{t\geq 0}$ is a one dimensional L\'evy process 
(i.e. a continuous time process with stationary and independent increments and $L(0)=0$ a.s.)
with $\E L(1)^2 < \infty$.  A comprehensive monograph dealing with L\'evy  processes is e.g. \cite{Applebaum2009}. The polynomials
$$
a(z) := z^p + a_1z^{p-1} + \dots + a_p,\quad 
b(z) := b_0 + b_1z + \dots + b_{p-1}z^{p-1},
$$
are called the \textit{autoregressive}- and \textit{moving average} polynomial, respectively. We assume that  $b_q\neq 0$ and $b_j = 0$ for $q < j < p$. The \textit{state-space representation}
consists of the \textit{observation} and \textit{state equations}:
\begin{equation}
Y(t)=\b^T\X (t),
\label{eq:ObservationArma}
\end{equation}
\begin{equation}
d\X(t)=\A\X(t)dt+\mathbf{e} dL(t),
\label{eq:StateRepresentationArma}
\end{equation}
where 
$$ 
\A:=\[
\begin{matrix}
0&1&0&\ldots&0\\
0&0&1&\ldots&0\\
\vdots&\vdots&\vdots&\ddots&\vdots\\
0&0&0&\ldots&1\\
-a_p&-a_{p-1}&-a_{p-2}&\ldots&-a_1\\
\end{matrix}
\],\quad \X(t):=\[
\begin{matrix}
X(t)\\
X^{(1)}(t)\\
\vdots\\
X^{(p-2)}(t)\\
X^{(p-1)}(t)\\
\end{matrix}
\],\quad 
$$
$$ \mathbf{e}:=\[ 0, \dots 0, 1 \]^T , \quad \mathbf{b}:=\[ b_0,b_1, \dots b_{p-1}\]^T,$$
i.e. $$ \A\in \R^{p\times p},\quad \X(t)\in \R^{p\times 1}\quad \mathbf{e}\in \R^p,\quad   \mathbf{b} \in \R^p.$$

If $p=1$, we set $\A=-a_1$.
\begin{assumption}\label{ass:FiniteSecondMomentsLevy}
The L\'evy process satisfies $\E L(1)=0$ and $\E |L(1)|^2=\sigma^2<\infty$.
\end{assumption}
Observe that $\E[L(s)L(t)]=\min\{ s,t \} \E |L(1)|^2$.
It was shown by Brockwell in \cite{Brockwell2009} that the solution $\X(t)$ of \eqref{eq:StateRepresentationArma} satisfies
\begin{equation}
\X(t)=e^{\mathbf{A}t}\X(0)+\int_0^t e^{\mathbf{A}(t-u)}\mathbf{e}dL(u),
\label{eq:StateRepresentationArmaSolIntegral}
\end{equation}
where the integral is defined as the $L^2$-limit of approximating Riemann sums. 
\begin{assumption}\label{ass:IndependenceLevyAndX}
$\X(0)$ is independent of  $(L(t))_{t\geq 0}$.
\end{assumption}
From now on let us assume that Assumption \ref{ass:IndependenceLevyAndX} holds. It is well-known (\cite[Proposition 2]{Brockwell2009}) that under Assumptions \ref{ass:FiniteSecondMomentsLevy}  and \ref{ass:IndependenceLevyAndX} 
the process  $\{ \X(t)\}_{t\geq 0}$ is strictly stationary and causal iff  $\X(0)$  has the same distribution as $\int_0^\infty e^{\mathbf{A}u}\mathbf{e}dL(u)$ and the $p$ (not necessarily distinct) eigenvalues $\lambda_1,\dots,\lambda_p$ of $\A$ have negative real parts, i.e.
$$ \Re (\lambda_i)<0,\quad i=1,\dots, p.$$
Now we extend the L\'evy process $(L(u))_{u\geq 0}$ to the whole line in the usual way: Let $\widetilde{L}=(\widetilde{L}(t))_{t\geq 0}$ be an independent copy 
of $(L(t))_{t\geq 0}$. For $t \in \R$ we define
$$L^*(t):=L(t)\1_{[0,\infty)}(t)+\widetilde{L}(-t-)\1_{(-\infty,0]}(t). $$
In order to get stationary solutions of \eqref{eq:StateRepresentationArma} we need the following assumptions:
\begin{assumption}\label{ass:EigenvaluesDistinctNonnegativeParts}
All eigenvalues of $\A$ have strictly negative real parts.
\end{assumption}

\begin{assumption}\label{ass:EqualityOfDistributionAtZero}
$$ \X(0)\stackrel{d}{=}\int_{-\infty}^0 e^{-\A u}\e dL^*(u)$$
\end{assumption}
In Brockwell \cite{Brockwell2009} it was shown that if Assumptions \ref{ass:EigenvaluesDistinctNonnegativeParts} and \ref{ass:EqualityOfDistributionAtZero} are satisfied the process $\{\X(t)\}_{t\in \R}$ given by
\begin{equation}\label{eq:SticlyStationaryStateRepresSolX}
\X(t)=\int_{-\infty}^t e^{\A(t-u)}\e dL^*(u)
\end{equation}
is a strictly stationary solution of \eqref{eq:StateRepresentationArma} (with $L$ replaced by $L^*$) for $t \in \R$ with corresponding $\operatorname{CARMA}$ process
\begin{equation}
Y(t)=\int_{-\infty}^t \b^Te^{\A(t-u)}\e dL^*(u).
\label{eq:CarmaSolution}
\end{equation}
For $t\geq 0$ one can rewrite it in the following form
\begin{equation}
Y(t)=\b^Te^{\A t}\X(0) +\int_{0}^t \b^Te^{\A(t-u)}\e dL(u).
\label{eq:CarmaSolution0}
\end{equation}
In the present paper the spectral density of a $\operatorname{CARMA}$ process plays a crucial role. The spectral density is the Fourier transform of the autocovariance function $\gamma_Y(h):=\Cov(Y(0),Y(h))$ for $h\in\mathbb{R}$.
The spectral density of a $\operatorname{CARMA}$  process is 
\begin{equation}
f_Y(\omega)=\frac{1}{2\pi}\int_{-\infty}^\infty \gamma_Y (h)e^{-ih\omega}dh=\frac{\sigma^2}{2\pi}\frac{|b(i\omega)|^2}{|a(i\omega)|^2},\quad \omega\in \R.
\label{eq:CARMASpectralDensity}
\end{equation}
%It is worth underlining that both the kernel function \eqref{eq:CARMAKernelDef} and the spectral density \eqref{eq:CARMASpectralDensity} contain the fraction $\frac{b(i\lambda)}{a(i\lambda)}$. V. Fasen and F. Fuchs in \cite{FuchsFasen2013} used the fraction $\frac{b(i\lambda)}{a(i\lambda)}$ to construct an estimator based on the smoothed periodogram for high-frequency equidistant data.
%We are going to present another approach based on the so called truncated Fourier transform.
%-------------------------------------------------------------------------------------
%-------------------------------------------------------------------------------------
%-------------------------------------------------------------------------------------
%-------------------------Limit behaviour of TFT----------------------------------------------------
%-------------------------------------------------------------------------------------
%-------------------------------------------------------------------------------------
%-------------------------------------------------------------------------------------
%--------------------------------------------------------------------------------------------------------------
\section{Limit behaviour of the Fourier transform}\label{sec3}
%--------------------------------------------------------------------------------------------------------------
In this section we are going to deal with the Fourier transform of the $\operatorname{CARMA}$ process assuming that the observations are given continuously on the time interval $[0,T]$.
A similar idea for Gaussian $\operatorname{CARMA}$ processes was presented in \cite{Gillberg2006} for equidistant observations. 
The truncated continuous-time Fourier transform of the process $Y$ at a fixed frequency $\omega\in \R$ is given by
$$\F_T(Y)(\omega):=\frac{1}{\sqrt{T}}\int_0^T Y(t)e^{-i\omega t}dt .$$
Observe that the norming constant $\frac{1}{\sqrt{T}}$ is taken as this ensures convergence in distribution for $T\to \infty$ as will be shown later.
 
\subsection{Properties of the Truncated Fourier Transform of a $\operatorname{CARMA}$ Process}\label{sec31}
First we derive an alternative representation.
\begin{lemma}\label{Lem:FourierTransformCarmaGeneral}
Let $\X$ and $Y$ be processes given by the state-space representation  \eqref{eq:ObservationArma} and \eqref{eq:StateRepresentationArma}. Suppose that Assumptions \ref{ass:FiniteSecondMomentsLevy}, \ref{ass:IndependenceLevyAndX} and \ref{ass:EigenvaluesDistinctNonnegativeParts} are satisfied. Then the truncated Fourier transform of the $\operatorname{CARMA}$ process $Y$ at a fixed frequency $\omega\in \R$ is of the form 
\begin{equation}
\label{eq:TruncatedFourierCarma}
\F_T(Y)(\omega)=\frac{1}{\sqrt{T}}\frac{b(i\omega)}{a(i\omega)}\int_0^T e^{-i\omega t} dL(t)+\frac{1}{\sqrt{T}}\b^T(i\omega I-A)^{-1}\left(\X(0)-e^{-i\omega T} \X(T) \right),
 \end{equation}
or equivalently
\begin{align}
\label{eq:TruncatedFourierCarmaPrim}
\F_T(Y)(\omega)=&\frac{1}{\sqrt{T}}\b^T(i\omega I-A)^{-1} \\\nonumber &\times \[ \int_0^T \left(e^{-i\omega u}- e^{-i\omega T} e^{\mathbf{A}(T-u)}\right)\mathbf{e} dL(u) +\left(I- e^{(-i\omega I+\mathbf{A})T}\right)\X(0)\].
\end{align}
\end{lemma}
\textsc{Proof}.
Let $\omega$ be an arbitrary frequency. Observe that by Corollary 3.4 from \cite[p. 51]{SchlemmStelzer2012} one has
$$\b^T(\A-i\omega I)^{-1}\e =-\frac{b(i\omega)}{a(i\omega)}. $$
Denote
\begin{align*}
F(t)&=\b^T(\A-i\omega I)^{-1}e^{(\A-i\omega I)t},\quad t\in [0,T],\\
G(t)&=\int_0^t e^{-\A u}\e dL(u)\quad t\in [0,T].
\end{align*}
Observe that $G(0)=0$ and since $F$ is continuous and of finite variation, we get $[F,G]=0$, where $[\cdot,\cdot]$ denotes the usual quadratic covariation of semimartingales (see e.g. \cite{Protter}).
Applying the (multidimensional) integration by parts formula 
\begin{align*}
 \int_0^T dF(t)G(t)&=F(T)G(T)-F(0)G(0)-\int_0^T F(t)dG(t)-[F,G]
 \\&=F(T)G(T)-\int_0^T F(t)dG(t)
\end{align*}
we obtain
\begin{align*}
&\int_0^T dF(t)G(t)=\int_0^T\b^T (\A-i\omega I)^{-1}(\A-i\omega I)  e^{(\A-i\omega I) t}\int_0^t e^{-\A u}\e dL(u)dt\\
&=\int_0^T\int_0^t \b^Te^{\A(t-u)}\e dL (u)e^{-i\omega t}dt
\\
&=  \b^T(\A-i\omega I)^{-1}e^{(\A-i\omega I) T}\int_0^T e^{-\A t}\e dL(t)-\int_0^T\b^T(\A-i\omega I)^{-1}e^{(\A-i\omega I) t} e^{-\A t} \e dL(t)
\\&= \b^T(\A-i\omega I)^{-1} e^{-i\omega T}\int_0^T e^{\A(T-t)}\e dL(t)+\frac{b(i\omega)}{a(i\omega)}\int_0^T e^{-i\omega t} dL(t).
\end{align*}
%Observe that the assumptions of the Stochastic Fubini Theorem (cf. Protter, Theorem IV. 64. \cite{Protter}) are satisfied.
Thus
\begin{eqnarray}\label{eq:IntInteDoATMinustdLt}
&\int_0^T\int_0^t \b^Te^{\A(t-u)}\e dL (u)e^{-i\omega t}dt
\\ \nonumber
&=\b^T(\A-i\omega I)^{-1} e^{-i\omega T}\int_0^T e^{\A(T-t)}\e dL(t)+\frac{b(i\omega)}{a(i\omega)}\int_0^T e^{-i\omega t} dL(t).
\end{eqnarray}

Using the form of the strictly stationary solution of \eqref{eq:StateRepresentationArma} given in \eqref{eq:StateRepresentationArmaSolIntegral} we get
\begin{equation}\label{eq:InteDoATMinustdLt}
 \int_0^T e^{\A(T-t)}\e dL(t)=\X(T)-e^{\A T}\X(0).
\end{equation}
Moreover, since $\int_0^T e^{(\A-i\omega I)t}dt= (i\omega I-\A)^{-1}(I-e^{(\A-i\omega I)T})$, we have 
\begin{equation}\label{eq:InteDoATMinustdLt11}
\int_0^T\b^Te^{(\A-i\omega I)t}\X(0)dt=\b^T (i\omega I-\A)^{-1}(I-e^{(\A-i\omega I)T})\X(0).
\end{equation}
We have
\begin{align*}
&\F_T(Y)(\omega)=\frac{1}{\sqrt{T}}\int_0^T Y(t)e^{-i\omega t}dt
\\&\stackrel{\eqref{eq:CarmaSolution0}}{=} \frac{1}{\sqrt{T}}\int_0^T\left(\b^Te^{\A t}\X(0) +\int_{0}^t \b^Te^{\A (t-u)}\e dL(u)\right)e^{-i\omega t}dt
\\&\stackrel{\eqref{eq:IntInteDoATMinustdLt},\eqref{eq:InteDoATMinustdLt}, \eqref{eq:InteDoATMinustdLt11}}{=}  
\frac{1}{\sqrt{T}}\frac{b(i\omega)}{a(i\omega)}\int_0^T e^{-i\omega u} dL(u)+\frac{1}{\sqrt{T}}\b^T(i\omega I-\A)^{-1}\left(\X(0)-e^{-i\omega T} \X(T) \right).
\end{align*}
To get the equivalent form note,
\begin{align*}
&\sqrt{T}\F(Y)(\omega)=\b^T(i\omega I-\A)^{-1}\[ \e \int_0^T e^{-i\omega u} dL(u) +\left(\X(0)-e^{-i\omega T} \X(T) \right)\]
\\&\stackrel{\eqref{eq:InteDoATMinustdLt}}{=} \b^T(i\omega I-\A)^{-1}
\[ \int_0^T \left(e^{-i\omega u}- e^{-i\omega T} e^{\mathbf{A}(T-u)}\right)\mathbf{e} dL(u) +\left(I- e^{(-i\omega I+\mathbf{A})T}\right)\X(0)\],
\end{align*}
which completes the proof of this Lemma. $\Box$

The next step is to calculate moments of the truncated Fourier transform. First, recall the so-called \textit{compensation formula}:
If $(L_t)_{t\geq 0}$ is a L\'evy process with finite first moments and $f$ is a bounded deterministic function, then 
\begin{equation}
\E\[ \int_0^T f(u)dL_u \]=\E[L_1]\int_0^T f(s)ds. 
\label{eq:CompensationFormula}
\end{equation}

%-------------------------------------------------------------------------------------------------------------
%----------------Behme----If $f$ is a bounded deterministic function then this formula holds (Behme, p. 46-47)
%-------------------------------------------------------------------------------------------------------------
Secondly, observe that the solution of the system \eqref{eq:ObservationArma} and \eqref{eq:StateRepresentationArma} is of the form \eqref{eq:StateRepresentationArmaSolIntegral}, 
where $\X$ is the process with mean $m(t)=\E[\X(t)]$ and 
$P_X(t)= \E[\X(t) \X(t)^T]$ satisfying 
\begin{align}
\nonumber
m_X(t)&=e^{\A t}m_X(0)\\
\label{eq:PX_IntegralRepresentation}
P_X(t)&=e^{\A t}P_X(0) e^{\A^Tt}+\sigma^2\int_0^t e^{\mathbf{A}(t-u)} \mathbf{e} \mathbf{e}^T e^{\mathbf{A^T}(t-u)}  du 
\end{align}
In particular, for stationary processes these solutions are constant and the so called Lyapunov equation
\begin{equation}
\A P_X+P_X\A^T+\sigma^2 \e \e^T=0
\label{eq:LyapunovStationaryCond}
\end{equation}
holds true. For L\'evy-driven CARMA processes the form of the autocovariance function in terms of solutions of Lyapunov equations is formulated 
e.g. in \cite[Proposition 3.13.]{MarquardtStelzer2007}.
%-------------------------------------------------------------------------------------------------------------
%-------------------------------------------------------------------------------------------------------------
%----------------Covariance of the truncated Fourie transform
%-------------------------------------------------------------------------------------------------------------

We are first going to show that the truncated Fourier transform of a stationary $\operatorname{CARMA}$ process is a zero-mean random variable.
Next, we find the covariance between the truncated Fourier transform at two different frequencies. As we have mentioned earlier, the spectral density function plays a central role.
\begin{theorem}\label{thm:FourierTransformCovariances}
Let $\X$ and $Y$ be processes given by the state-space representation  \eqref{eq:ObservationArma} and \eqref{eq:StateRepresentationArma}. Suppose that Assumptions \ref{ass:FiniteSecondMomentsLevy}, \ref{ass:IndependenceLevyAndX}, \ref{ass:EigenvaluesDistinctNonnegativeParts} and \ref{ass:EqualityOfDistributionAtZero} are satisfied. Then
$\E(\F_T(Y)(\omega))=0$ for all $\omega \in \R$.
For $\omega_1,\omega_2 \in \R$ we have
\begin{equation}
\E\[\F_T(Y)(\omega_1)\F_T(Y)(\omega_2)\]=\sigma^2 \frac{|b(i\omega_1)|^2}{|a(i\omega_1)|^2}+\frac{1}{T}K(T,\omega_1,-\omega_1), \quad\mathrm{if\quad } \omega_1=-\omega_2
\label{eq:IloczynTransformatFourieraomega1RownaMinusomega2}
\end{equation}
and
\begin{equation}
\E\[\F_T(Y)(\omega_1)\F_T(Y)(\omega_2)\]=\frac{1}{T}K_1(T,\omega_1,\omega_2), \quad \mathrm{if\quad }\omega_1\neq -\omega_2,
\label{eq:IloczynTransformatFourieraomega1RozneMinusomega2}
\end{equation}
where $K$ is a bounded function of $T$ given by \eqref{eq:KWzorDlugi} below and 
$$ K_1(T,\omega_1,\omega_2)
=K(T,\omega_1,\omega_2)+\b^T(i\omega_1 I-\mathbf{A})^{-1} \sigma^2 \frac{1-\exp(-Ti(\omega_1+\omega_2))}{i(\omega_1+\omega_2)} \mathbf{e} \mathbf{e}^T(i\omega_2 I-\mathbf{A}^T)^{-1}\b.$$ 
\end{theorem}
\textsc{Proof}.
For the first part it is enough to observe that by the compensation formula $ \E \left( \int_0^T e^{-i\omega u} dL(u)\right)=0 $ and $\E[\X(t)]=0$. For the second part
observe that using Lemma \ref{Lem:FourierTransformCarmaGeneral} and formula \eqref{eq:TruncatedFourierCarmaPrim} we have
\begin{align*}
&\E[\F_T(Y)(\omega_1)\F_T(Y)(\omega_2)]=\frac{1}{T}\b^T(i\omega_1 I-\A)^{-1}\times
\\&\E\[ \left(\int_0^T \left(e^{-i\omega_1 u}- e^{-i\omega_1 T} e^{\mathbf{A}(T-u)}\right)\mathbf{e} dL(u)\right.+\left.\left(I- e^{(-i\omega_1 I+\mathbf{A})T}\right)\X(0) \right)\right)  \times
\\&\left. \left(\int_0^T\mathbf{e}^T \left(e^{-i\omega_2 u}- e^{-i\omega_2 T} e^{\mathbf{A}^T(T-u)}\right) dL(u) +\X(0)^T\left(I- e^{(-i\omega_2 I+\mathbf{A}^T)T}\right)\right)\]\times
\\& (i\omega_2 I-\A^T)^{-1}\b = \frac{1}{T}\b^T(i\omega_1 I-\A)^{-1} \widetilde{I} (i\omega_2 I-\A^T)^{-1}\b,
\end{align*}
where $\widetilde{I}=I_1+I_2+I_3+I_4$ with
\begin{align*}
I_1&:=\E\[ \int_0^T \left(e^{-i\omega_1 u}- e^{-i\omega_1 T} e^{\mathbf{A}(T-u)}\right)\mathbf{e} dL(u) \cdot  \int_0^T\mathbf{e}^T \left(e^{-i\omega_2 u}- e^{-i\omega_2 T} e^{\mathbf{A}^T(T-u)}\right) dL(u) \]
\\I_2&:=\E\[ \int_0^T \left(e^{-i\omega_1 u}- e^{-i\omega_1 T} e^{\mathbf{A}(T-u)}\right)\mathbf{e} dL(u) \cdot \X(0)^T\left(I- e^{(-i\omega_2 I+\mathbf{A}^T)T}\right)\] 
\\I_3&:=\E\[\left(I- e^{(-i\omega_1 I+\mathbf{A})T}\right)\X(0) \cdot \int_0^T\mathbf{e}^T \left(e^{-i\omega_2 u}- e^{-i\omega_2 T} e^{\mathbf{A}^T(T-u)}\right) dL(u)\] 
\\I_4&:= \E\[ \left(I- e^{(-i\omega_1 I+\mathbf{A})T}\right)\X(0) \X(0)^T\left(I- e^{(-i\omega_2 I+\mathbf{A}^T)T}\right)\].
\end{align*}
%!!!!Zbadac zbieznosc calki I1!!!!!
We have that $I_2=I_3=0$ since $(L_t)_{t\geq 0}$ is independent of $\X(0)$.
Observe that by the It\^o isometry, the compensation formula and the fact that $\E[[L,L]_1]=\operatorname{Var} (L(1))=\sigma^2$ we have
\begin{align*}
I_1^1:=&\E\[ \int_0^T e^{-i\omega_1 u}\mathbf{e} dL(u)  \int_0^T \mathbf{e}^Te^{-i\omega_2 u} dL(u)\]=\E\[ \int_0^T e^{-i(\omega_1 +\omega_2) u}\mathbf{e} \mathbf{e}^T d[L,L]_u\]
\\&=\E[[L,L]_1]\int_0^T e^{-i(\omega_1 +\omega_2)u}\mathbf{e} \mathbf{e}^Tdu=\sigma^2 \int_0^T e^{-i(\omega_1 +\omega_2)u}\mathbf{e} \mathbf{e}^Tdu.
\end{align*}
Thus
\begin{equation}
I^1_1=\left\{ 
\begin{matrix}
\sigma^2 T \mathbf{e} \mathbf{e}^T,\quad \omega_1=-\omega_2,\\
\sigma^2 \frac{1-\exp(-Ti(\omega_1+\omega_2))}{i(\omega_1+\omega_2)} \mathbf{e} \mathbf{e}^T,\quad \omega_1\neq-\omega_2. \\
\end{matrix}
\right.
\label{eq:IloczynTransformatFourieraI1Value}
\end{equation}
Thus, if $\omega_1=-\omega_2$, then 
\begin{align*}
\frac{1}{T}\b^T(i\omega_1 I-\A)^{-1} I^1_1  (i\omega_2 I-\A^T)^{-1}\b&=\frac{1}{T} \cdot \sigma^2 T\b^T(i\omega_1 I-\A)^{-1} \mathbf{e} \mathbf{e}^T(i\omega_2 I-\A^T)^{-1}\b\end{align*}\begin{align*}
\\&= \sigma^2 \b^T(\A-i\omega_1 I)^{-1} \mathbf{e} \mathbf{e}^T(i\omega_1 I+\A^T)^{-1}\b
\\&= \sigma^2 \left(-\frac{b(i\omega_1)}{a(i\omega_1)}\right)\left(-\frac{b(-i\omega_1)}{a(-i\omega_1)}\right)= \sigma^2 \frac{|b(i\omega_1)|^2}{|a(i\omega_1)|^2}.
\end{align*}
Now
\begin{align*}
I_1^2:&=\E\[ \int_0^T e^{-i\omega_1 u}\mathbf{e} dL(u)  \int_0^T \mathbf{e}^Te^{-i\omega_2 T} e^{\mathbf{A}^T(T-u)} dL(u)\]
\\&=e^{-i\omega_2 T}\E\[ \int_0^T e^{-i\omega_1 u} \mathbf{e} \mathbf{e}^Te^{\mathbf{A}^T(T-u)} d [L,L]_u\]
\\&=e^{-i\omega_2 T} \E[[L,L]_1] \int_0^T e^{-i\omega_1 u} \mathbf{e} \mathbf{e}^Te^{\mathbf{A}^T(T-u)} du
\\&=e^{-i\omega_2 T} \sigma^2 \int_0^T e^{-i\omega_1 u} \mathbf{e} \mathbf{e}^Te^{\mathbf{A}^T(T-u)} du.
\end{align*}
In the same way
\begin{align*}
I_1^3:&=\E\[ \int_0^T e^{-i\omega_1 T} e^{\mathbf{A}(T-u)}\mathbf{e} dL(u)  \int_0^T\mathbf{e}^T  e^{-i\omega_2 u}dL(u)\]\\
&=e^{-i\omega_1 T} \sigma^2 \int_0^T e^{\mathbf{A}(T-u)}\mathbf{e} \mathbf{e}^T  e^{-i\omega_2 u} du.
\end{align*}

Combining these two we arrive at
\begin{align*}
I_1^2+I_1^3=&e^{-i(\omega_1+\omega_2)T}\sigma^2\left[ \mathbf{e} \mathbf{e}^T(i\omega_1 I+\A^T)^{-1}\left(e^{(i\omega_1I+\A^T)T} -I \right) \right.\\
&+\left. (i\omega_2 I+\A)^{-1} \left(e^{(i\omega_2I+\A)T} -I \right)  \mathbf{e} \mathbf{e}^T\right].
\end{align*}

Now
\begin{align*}
I_1^4:&=\E\[ \int_0^T e^{-i\omega_1 T} e^{\mathbf{A}(T-u)}\mathbf{e} dL(u)  \int_0^T\mathbf{e}^T e^{-i\omega_2 T} e^{\mathbf{A}^T(T-u)}   dL(u)\] \\
&=e^{-i(\omega_1+\omega_2) T} \sigma^2 \int_0^T e^{\mathbf{A}(T-u)} \mathbf{e} \mathbf{e}^T e^{\mathbf{A}^T (T-u)}  du.
\end{align*}
Now
\begin{align*}
I_4=&\E[\X(0) \X(0)^T]-e^{-i\omega_1 T}e^{\A T}\E[\X(0) \X(0)^T]-e^{-i\omega_2 T}\E[\X(0) \X(0)^T]e^{\A^TT}
\\&+ e^{-i(\omega_1+\omega_2) T}e^{\A T} \E[\X(0) \X(0)^T] e^{\A^TT}.
\end{align*}
By stationarity we have $$ \E[\X(0) \X(0)^T]=:P_X=P_X(0)=P_X(T), $$ where $P_X$ satisfies \eqref{eq:LyapunovStationaryCond}. Combining this with \eqref{eq:PX_IntegralRepresentation} we obtain
\begin{align*}
I_1^4+I_4= &P_X-e^{-i\omega_1 T}e^{\A T}P_X-e^{-i\omega_2 T}P_Xe^{\A^TT}+ e^{-i(\omega_1+\omega_2) T}e^{\A T} P_X e^{\A^TT}
\\&+e^{-i(\omega_1+\omega_2) T}(P_X-e^{\A T}P_X e^{\A^T T})
\\=&e^{-i\omega_1 T}P_X\left(I-e^{\A T}\right)+e^{-i\omega_2 T}\left(I-e^{\A^T}\right)P_X
\\&+P_X\left(1- e^{-i\omega_1 T}-e^{-i\omega_2 T}+e^{-i(\omega_1+\omega_2) T}\right).
\end{align*}
Since $\A$ is a stable matrix, $e^{\A T}$ is bounded. 
%-------------------------------------------------------------------------------------------------------------
%----------------References for this fact
%-------------------------------------------------------------------------------------------------------------

Thus
\begin{align}
\nonumber
K(T,\omega_1,\omega_2)=&\b^T(i\omega_1 I-\A)^{-1} \left[ 
e^{-i(\omega_1+\omega_2)T}\sigma^2\left[ \mathbf{e} \mathbf{e}^T(i\omega_1 I+\A^T)^{-1}\left(e^{(i\omega_1I+\A^T)T} -I \right) \right. \right.\\
\label{eq:KWzorDlugi}
&+\left. (i\omega_2 I+\A)^{-1} \left(e^{(i\omega_2I+\A)T} -I \right)  \mathbf{e} \mathbf{e}^T\right]\\
\nonumber
&+e^{-i\omega_1 T}P_X\left(I-e^{\A^TT}\right)+e^{-i\omega_2 T}\left(I-e^{\A T}\right)P_X
\\
\nonumber
&+P_X\left(1- e^{-i\omega_1 T}-e^{-i\omega_2 T}+e^{-i(\omega_1+\omega_2) T}\right)
\left. \right] (i\omega_2 I-\A^T)^{-1}\b
\end{align}
is bounded  in $T$ for fixed $\omega_1,\omega_2\in \R$. $\quad \Box$

Now we give the form of the covariance matrix. Put

$$\Sigma(\omega_1,\omega_2):=[\Sigma_{ij}]_{1\leq i,j \leq 4}=\E\[ 
\[
\begin{matrix}
\Re \mathcal{F}_T(Y)(\omega_1)\\
\Im \mathcal{F}_T(Y)(\omega_1)\\
\Re \mathcal{F}_T(Y)(\omega_2)\\
\Im \mathcal{F}_T(Y)(\omega_2)\\
\end{matrix}
\]
\[
\begin{matrix}
\Re \mathcal{F}_T(Y)(\omega_1)\\
\Im \mathcal{F}_T(Y)(\omega_1)\\
\Re \mathcal{F}_T(Y)(\omega_2)\\
\Im \mathcal{F}_T(Y)(\omega_2)\\
\end{matrix}
\]^T
 \] .$$
\begin{theorem} 
Let $\X$ and $Y$ be processes given by the state-space representation  \eqref{eq:ObservationArma} and \eqref{eq:StateRepresentationArma}. Suppose that Assumptions \ref{ass:FiniteSecondMomentsLevy}, \ref{ass:IndependenceLevyAndX}, \ref{ass:EigenvaluesDistinctNonnegativeParts} and \ref{ass:EqualityOfDistributionAtZero} are satisfied.
For $\omega_1\neq \omega_2$ and $\omega_1\neq - \omega_2$ there exists a bounded matrix $K_2\in \C^{4\times 4}$ such that
$$\Sigma(\omega_1,\omega_2)=\frac{1}{2} \sigma^2 \operatorname{diag}\left(\frac{|b(i\omega_1)|^2}{|a(i\omega_1)|^2}, \frac{|b(i\omega_1)|^2}{|a(i\omega_1)|^2}, \frac{|b(i\omega_2)|^2}{|a(i\omega_2)|^2}, 
\frac{|b(i\omega_2)|^2}{|a(i\omega_2)|^2}\right)
+\frac{1}{T}K_2.$$
%where  $  f_Y $ is the spectral density function of the process $Y$.
\end{theorem}
\proof
For $k,l=1,2$ let us denote $$ \Sigma_1(\omega_1,\omega_2):=\E\[\Re \mathcal{F}_T(Y)(\omega_1)\Re \mathcal{F}_T(Y)(\omega_2)\], 
\quad \Sigma_2(\omega_1,\omega_2):=\E\[\Im \mathcal{F}_T(Y)(\omega_1)\Im \mathcal{F}_T(Y)(\omega_2)\],$$
$$ \Sigma_3(\omega_1,\omega_2):=\E\[\Re \mathcal{F}_T(Y)(\omega_1)\Im \mathcal{F}_T(Y)(\omega_2)\].$$
All entries $\Sigma_{i,j}$ of the matrix $\Sigma$ are of one of the above forms. Indeed, $\Sigma_{11}$, $\Sigma_{33}$ are of the form $\Sigma_1$ for $k=l$ and $k,l\in\{ 1,2 \}$. Similarly, $\Sigma_{22}$, $\Sigma_{44}$ are of the form $\Sigma_2$ for $k=l$ and $k,l\in\{ 1,2 \}$. Moreover, 
$\Sigma_{13}$, $\Sigma_{31}$ are of the form $\Sigma_1$ for $k\neq l$ and $k,l\in\{ 1,2 \}$ and 
$\Sigma_{24}$, $\Sigma_{42}$ are of the form $\Sigma_2$ for $k\neq l$ and $k,l\in\{ 1,2 \}$. All other elements are of the form $\Sigma_3$. 

Observe that for each $\omega$ we have
$$\Re \mathcal{F}_T(Y)(\omega)=\frac{\mathcal{F}_T(Y)(\omega)+\mathcal{F}_T(Y)(-\omega)}{2},\quad \Im \mathcal{F}_T(Y)(\omega)=\frac{\mathcal{F}_T(Y)(\omega)-\mathcal{F}_T(Y)(-\omega)}{2i} .$$
Using Theorem \ref{thm:FourierTransformCovariances} we obtain
\begin{align*}
\Sigma_1(\omega_1,\omega_2):=\left\{  
\begin{matrix}
\sigma^2\frac{|b(0)|^2}{|a(0)|^2}+\frac{1}{T}K(0),&\quad\omega_1=\omega_2=0;\\
\frac{1}{2}\sigma^2\frac{|b(i\omega_1)|^2}{|a(i\omega_1)|^2}+\frac{1}{T}K_{1,1}(\omega_1),&\quad\omega_1=\omega_2;\\
\frac{1}{2}\sigma^2\frac{|b(i\omega_1)|^2}{|a(i\omega_1)|^2}+\frac{1}{T}K_{1,2}(\omega_1),&\quad \omega_1=-\omega_2;\\
\frac{1}{T}K_{1,3}(\omega_1,\omega_2),&\quad\omega_1\neq\omega_2\, \mathrm{and} \,\omega_1\neq -\omega_2,\\
\end{matrix}
\right.
\end{align*}
\begin{align*}
\Sigma_2(\omega_1,\omega_2):=\left\{  
\begin{matrix}
0,&\quad\omega_1=0 \,\,\mathrm{or}\,\,\omega_2=0;\\
\frac{1}{2}\sigma^2\frac{|b(i\omega_1)|^2}{|a(i\omega_1)|^2}+\frac{1}{T}K_{2,1}(\omega_1),&\quad\omega_1=\omega_2;\\
-\frac{1}{2}\sigma^2\frac{|b(i\omega_1)|^2}{|a(i\omega_1)|^2}-\frac{1}{T}K_{2,2}(\omega_1),&\quad \omega_1=-\omega_2;\\
\frac{1}{T}K_{2,3}(\omega_1,\omega_2),&\quad\omega_1\neq\omega_2\,\, \mathrm{and} \,\,\omega_1\neq -\omega_2,\\
\end{matrix}
\right.
\end{align*}
\begin{align*}
\Sigma_3(\omega_1,\omega_2):=\left\{  
\begin{matrix}
0,&\quad\omega_2=0;\\
\frac{1}{T}K_{3,1}(\omega_1),&\quad\omega_1=\omega_2\, \mathrm{or} \,\omega_1=-\omega_2;\\
\frac{1}{T}K_{3,2}(\omega_1,\omega_2),&\quad\omega_1\neq\omega_2\,\, \mathrm{and} \,\,\omega_1\neq -\omega_2.\\
\end{matrix}
\right.
\end{align*}
Here $K$ is given by \eqref{eq:KWzorDlugi} and $K_{i,j}$ are bounded in $T$ for $i,j=1,2,3$. $\Box$

Now we are going to investigate asymptotic properties of the truncated Fourier transform. 
First, we will show that the second summand of \eqref{eq:TruncatedFourierCarma} converges in probability to zero.
\begin{lemma}\label{lem:LimitInProbabilityFirstSummandFT}
Let $\X$ and $Y$ be processes given by the state-space representation  \eqref{eq:ObservationArma} and \eqref{eq:StateRepresentationArma}. 
Suppose that Assumptions \ref{ass:FiniteSecondMomentsLevy}, \ref{ass:IndependenceLevyAndX} and \ref{ass:EigenvaluesDistinctNonnegativeParts} are satisfied.
Let $$\tilde{Z}(T):= \mathcal{F}_T(Y)(\omega)-\frac{1}{\sqrt{T}}\frac{b(i\omega)}{a(i\omega)}\int_0^T e^{-i\omega t} dL(t).$$
Then
$$ \P-\lim_{T\to \infty} |\tilde{Z}(T)|=0. $$
\end{lemma}
\proof
Observe that 
\begin{align*}
|\tilde{Z}(T)|&= \left|\frac{1}{\sqrt{T}}\b^T(i\omega I-A)^{-1}\[\X(0)-e^{-i\omega T}\X(T)\]\right|
\\&\leq  \frac{1}{\sqrt{T}}\left|\b^T(i\omega I-A)^{-1}\X(0)\right|+\frac{1}{\sqrt{T}}\left|\b^T(i\omega I-A)^{-1}\X(T)\right|. 
\end{align*}
Obviously, $$\lim_{T\to \infty}\frac{1}{\sqrt{T}}\left|\b^T(i\omega I-A)^{-1}\X(0)\right|=0 \quad a.s. \mathrm{\,\,as\,\,} T\to \infty.$$
Because of stationarity, $\X(T)$ is bounded in probability and $\frac{1}{\sqrt{T}}$ converges to zero thus 
$$ \frac{1}{\sqrt{T}}\left|\b^T(i\omega I-A)^{-1}\X(T)\right| \to 0 \mathrm{\,\,in\,\,probability} .$$
Therefore
$$ \P-\lim_{T\to \infty} |\tilde{Z}(T)|=0.$$
This completes the proof. $\Box$

%Observe that 
%\begin{align*}
%|\tilde{Z}(T)|&\leq \left|\frac{1}{\sqrt{T}}\b^T(i\omega I-A)^{-1}\[\X(0)(1-e^{-i\omega T})\]\right|\\&+ \left|\frac{1}{\sqrt{T}}\b^T(i\omega I-A)^{-1} e^{-i\omega T}\int_0^T e^{A(T-u)}\e dL(u)\right|. 
%\end{align*}
%Observe that $\b^T(i\omega I-A)^{-1}\[\X(0)(1-e^{-i\omega T})\]$ is bounded so $$\left|\frac{1}{\sqrt{T}}\b^T(i\omega I-A)^{-1}\[\X(0)(1-e^{-i\omega T})\]\right|$$ converges to zero a.s. as $T\to \infty$,
%thus it converges to zero in probability.
%Applying \cite[Theorem  4.3.16]{Applebaum2009} for $f(u)=e^{-Au}\e $ we obtain that the process $\left(\int_{-\infty}^T e^{A(T-u)}\e dL(u)\right)_{T\geq 0}$ is strictly stationary, 
%thus it converges in distribution thus it is bounded in probability. Moreover, $\frac{1}{\sqrt{T}}\left|\b^T(i\omega I-A)^{-1} e^{-i\omega T}\right|$ is convergence almost surely so in probability to zero 
%and thus  $\left|\frac{1}{\sqrt{T}}\b^T(i\omega I-A)^{-1} e^{-i\omega T}\int_0^T e^{A(T-u)}\e dL(u)\right|$ converges in probability to zero. $\Box$

Now we will show that the first summand of formula \eqref{eq:TruncatedFourierCarma} converges in distribution. Thus, together with Lemma \ref{lem:LimitInProbabilityFirstSummandFT} 
we obtain the limit in distribution of the truncated Fourier transform. We have two cases: 
the first case is if the frequency $\omega=0$. Then the truncated Fourier transform is a real valued function. In the second case for frequencies $\omega\neq 0$ the truncated Fourier transform is a complex valued function.
In both cases we first give the description of the distribution of the truncated Fourier transform and afterwards we describe the distribution of the  squared modulus of the truncated Fourier transform.

\begin{theorem}\label{lem:RozkladZ(T)OmegaRowna0} Let $\X$ and $Y$ be processes given by the state-space representation  \eqref{eq:ObservationArma} and \eqref{eq:StateRepresentationArma}. 
Suppose that Assumptions \ref{ass:FiniteSecondMomentsLevy} and \ref{ass:EigenvaluesDistinctNonnegativeParts} are satisfied.
Let $$Z(T):=\frac{1}{\sqrt{T}}\frac{b(0)}{a(0)}\int_0^T dL(t). $$Then
$$ d-\lim_{T\to \infty} Z(T) \sim \mathcal{N}\left(0,\left(\frac{b(0)}{a(0)}\right)^2\sigma^2\right),\quad d-\lim_{T\to \infty} \frac{1}{\sigma^2}  \left|\frac{a(0)Z(T)}{b(0) }\right|^2 \sim \chi^2(1)$$
\end{theorem}
 \proof 
Observe that $\int_0^T dL(t)=L(T) $, thus $Z(T)= \frac{1}{\sqrt{T}} \frac{b(0)}{a(0)} L(T) $. By the standard Central Limit Theorem 
$ d-\lim_{T\to \infty} \frac{1}{\sqrt{T}} L(T) = \mathcal{N}(0,\sigma^2) $. 
Therefore
$ d-\lim_{T\to \infty} \frac{1}{\sqrt{T}} \frac{b(0)}{a(0)}L(T) = \mathcal{N}\left(0,\left(\frac{b(0)}{a(0)}\right)^2\sigma^2\right) $. 

Observe that for all $n\in \N$ the random variable 
$\frac{a(0) Z(n)}{b(0)\sigma} \sim \mathcal{N}\left(0,1\right)$. Then by the continuous mapping theorem we have $d-\lim_{T\to \infty}\frac{1}{\sigma^2}\left|\frac{a(0)Z(T)}{b(0)}\right|^2 \sim \chi^2(1)$.  $\Box$

In order to find the asymptotic distribution of the truncated Fourier transform we use the multivariate Central Limit Theorem. Note that we state all results for positive frequencies as the corresponding results for negative  can be obtained by taking the complex conjugate.

\begin{theorem}\label{lem:RozkladZ(T)} 
Let $\X$ and $Y$ be processes given by the state-space representation  \eqref{eq:ObservationArma} and \eqref{eq:StateRepresentationArma}. Suppose that Assumptions \ref{ass:FiniteSecondMomentsLevy}  and \ref{ass:EigenvaluesDistinctNonnegativeParts} are satisfied. Assume that $\omega > 0$.
Put $$Z(T):=\frac{1}{\sqrt{T}}\frac{b(i\omega)}{a(i\omega)}\int_0^T e^{-i\omega t} dL(t) $$
and $$ Z(T)= \[\begin{matrix}
\Re Z(T)\\
\Im Z(T)\\
\end{matrix}\].$$
Then
$$ d-\lim_{T\to \infty} Z(T) \sim \mathcal{N}(0,\Sigma),$$
where $\Sigma=\frac{\sigma^2}{2}\left|\frac{b(i\omega)}{a(i\omega)} \right|^2I_{2\times 2}$.
\end{theorem}
\proof
We firstshow that $\frac{1}{\sqrt{N}}\int_0^{ \frac{2\pi N}{\omega}} e^{-i\omega t} dL(t)$ is asymptotically normal.
For $N\in \N$ and $j\in \{0,\dots, N-1 \}$ put
$$X_{j}:=\[\begin{matrix}
X^1_{j}\\
X^2_{j}\\
\end{matrix}\]
:=
\[\begin{matrix}
\int_{2\pi j/\omega}^{2\pi (j+1)/\omega} \cos (\omega t) dL(t)\\
\int_{2\pi j/\omega}^{2\pi (j+1)/\omega} \sin (\omega t) dL(t)\\
\end{matrix}\].
$$
Observe that $X_{j}$ are independent and identically distributed random vectors with mean zero and the covariance matrix
$\widetilde{\Sigma_1}:=\frac{\sigma^2\pi}{\omega}I_{2\times 2}$. 
Therefore, $$\int_0^{ \frac{2\pi N}{\omega}} e^{-i\omega t} dL(t)=\sum_{j=0}^{N-1} X_j.$$
Applying the classical CLT we obtain
$$\sqrt{ N}\left(\frac{1}{ N} \sum_{j=0}^{N-1} X_j\right)=\frac{1}{\sqrt{N}}\int_0^{ \frac{2\pi N}{\omega}} e^{-i\omega t} dL(t) 
\to \mathcal{N}\sim \mathcal{N}(0,\widetilde{\Sigma_1}) \quad \mathrm{as}\quad N \to \infty.$$
So
$\frac{\sqrt{\omega}}{\sqrt{2\pi N}}\int_0^{ \frac{2\pi N}{\omega}} e^{-i\omega t} dL(t)\to \mathcal{N}(0,\Sigma_1) $, where $\Sigma_1=\frac{\omega}{2 \pi}\widetilde{\Sigma_1}=\frac{\sigma^2}{2}I_{2\times 2}$.
Put 
$$A:=\[ 
\begin{matrix} 
\Re \left( \frac{b(i\omega)}{a(i\omega)} \right)& \Im \left( \frac{b(i\omega)}{a(i\omega)} \right)\\
\Im \left( \frac{b(i\omega)}{a(i\omega)} \right)& -\Re \left( \frac{b(i\omega)}{a(i\omega)} \right) 
\end{matrix}
 \] .$$
Observe that
\begin{align*}
&A\cdot\[\begin{matrix}
\frac{\sqrt{\omega}}{\sqrt{2\pi N}}\int_{0}^{2\pi N/\omega} \cos (\omega t) dL(t)\\
\frac{\sqrt{\omega}}{\sqrt{2\pi N}}\int_{0}^{2\pi N/\omega} \sin (\omega t) dL(t)\\
\end{matrix}\]=\[\begin{matrix}
\Re Z(2\pi N)\\
\Im Z(2\pi N)\\
\end{matrix}\]. \\
\end{align*}
Thus $Z=A\cdot X$ is normally distributed with mean zero and the covariance matrix $\Sigma =A\Sigma_1 A^T= \frac{\sigma^2}{2}\left|\frac{b(i\omega)}{a(i\omega)} \right|^2I_{2\times 2}$. $\Box$

Now we apply this theorem to find the asymptotic distribution of the truncated Fourier transform squared.

\begin{theorem}\label{lem:RozkladZ(T)Squared}
Let $\X$ and $Y$ be processes given by the state-space representation  \eqref{eq:ObservationArma} and \eqref{eq:StateRepresentationArma}. Suppose that Assumptions \ref{ass:FiniteSecondMomentsLevy} and \ref{ass:EigenvaluesDistinctNonnegativeParts}  are satisfied.
Let $Z$ be defined as in Theorem \ref{lem:RozkladZ(T)}. Then 
%$$|Z|^2=\lambda_1(U_1^2+U_2^2),$$ 
%where $U_1$ and $U_2$ are independent standard normal variables and
%$\lambda_{1}=\frac{\sigma^2}{2}\left| \frac{b(i\omega)}{a(i\omega)} \right|^2 
%$ and thus $
$|Z|^2\sim \operatorname{Exp}\left(\sigma^2\left| \frac{b(i\omega)}{a(i\omega)} \right|^2\right)$, where $ \operatorname{Exp}(\lambda)$ denotes the exponential distribution with mean $\lambda$.
\end{theorem}
\textsc{Proof} 
We use the notation of the proof of Theorem \ref{lem:RozkladZ(T)}. 
%Eigenvalues of the matrix $\Sigma$ defined in Theorem \ref{lem:RozkladZ(T)} are given by
%$$\lambda_{1}=\lambda_{2}=\frac{\sigma^2}{2}
% \left| \frac{b(i\omega)}{a(i\omega)} \right|^2 
%$$
%Observe that $|Z|^2=Z_1^2+Z_2^2$. Thus we are interested in the sum of two \textbf{dependent} random variables. 
%Let $\mathbf{Z}=[Z_1\quad Z_2]^T$ and us consider the following quadratic form
%$$ Q(\mathbf{Z}):=\mathbf{Z}^T I \mathbf{Z}=Z_1^2+Z_2^2$$
%Using the spectral theorem one can show that $$ Q(\mathbf{Z})=\lambda_1 U_1^2+\lambda_2 U_2^2=\lambda_1 \left(U_1^2+U_2^2\right),$$
%where $\lambda_1=\lambda_2$ is positive and $U_i$ are two independent standard normally distributed random variables.
Thus $|Z|^2$ is proportional to chi-square random variables with two degrees of freedom, % (cf. \cite[Chapter 3]{MathaiProvost} )
i.e. $|Z|^2=\frac{\sigma^2}{2} \left| \frac{b(i\omega)}{a(i\omega)} \right|^2 X$, where $X\sim \chi^2\left(2\right)$. 
 Thus  $|Z|^2 \sim \operatorname{\Gamma}\left(1,\frac{\sigma^2}{2}\left| \frac{b(i\omega)}{a(i\omega)} \right|^2\right)$ so 
 $|Z|^2 \sim \operatorname{Exp}\left(\sigma^2\left| \frac{b(i\omega)}{a(i\omega)} \right|^2\right)$. $\Box$

Now we are going to give the description of the convergence of the random vector consisting of the truncated Fourier transform at different frequencies.

\begin{theorem}\label{thm:LimitTheoremSeveralFrequencies}
Let $\X$ and $Y$ be processes given by the state-space representation  \eqref{eq:ObservationArma} and \eqref{eq:StateRepresentationArma}. Suppose that Assumptions \ref{ass:FiniteSecondMomentsLevy}, \ref{ass:IndependenceLevyAndX} and  \ref{ass:EigenvaluesDistinctNonnegativeParts} are satisfied.
 Let $0<\omega_1<\dots<\omega_d$ be fixed frequencies. 
Then 
$\[\Re\left(\mathcal{F}_T(Y)(\omega_j)\right), \Im\left(\mathcal{F}_T(Y)(\omega_j) \right) \]^T_{j=1,\dots,d}$ converges to  $\mathcal{N}\left(0,\frac{\sigma^2}{2} \mathbf{B}\right)$, with 
$$\mathbf{B}=\operatorname{diag}\left(\left| \frac{b(i\omega_1)}{a(i\omega_1)} \right|^2,\left| \frac{b(i\omega_1)}{a(i\omega_1)} \right|^2,\dots,\left| \frac{b(i\omega_d)}{a(i\omega_d)} \right|^2,\left|
\frac{b(i\omega_d)}{a(i\omega_d)} \right|^2\right)$$
and  $\[|\mathcal{F}_T(\omega_j)|^2 \]^T_{j=1,\dots,d}$ 
converges to a random vector whose coordinates are independent $\operatorname{Exp}\left(\sigma^2\left| \frac{b(i\omega_j)}{a(i\omega_j)} \right|^2\right)$ distributed random variables
 for $j=1,\dots,d$.
\end{theorem}
\textsc{Proof}
For fixed $n\in \N$ and $k=1,\dots, n$, put
$$ X^{(2i-1)}_k(\omega_i):=\int_{2(k-1)\pi}^{2k\pi} \cos(\omega_i t) dL(t),\quad X^{(2i)}_k(\omega_i):=\int_{2(k-1)\pi}^{2k\pi} \sin(\omega_i t) dL(t),\quad i=1,\dots,d.$$
Let $\left(s_n^{(2i-1)}\right)^2=\sum_{k=1}^n \Var \[X^{(2i-1)}(\omega_i)\]$ and  $\left(s_n^{(2i)}\right)^2=\sum_{k=1}^n \Var \[X^{(2i)}(\omega_i)\]$. Put 
$$ Z^{(2i-1)}_n(\omega_i):=\frac{\sum_{k=1}^n X^{(2i-1)}_k(\omega_i)}{s_n^{(2i-1)}}, \quad Z^{(2i)}_n(\omega_i):=\frac{\sum_{k=1}^n X^{(2i)}_k(\omega_i)}{s_n^{(2i)}} .$$
Then we will show that by the Cramer-Wold-device the random vector $\mathbf{Z}\in \R^{2d}$ with $\mathbf{Z}=\[Z^{(2i-1)}_n(\omega_i),Z^{(2i)}_n(\omega_i) \]^T_{i=1,\dots,d}$ 
converges in distribution to 
$\mathcal{N}\left(0,I_{2d\times 2d}\right)$.

We first apply the Lindeberg-Feller Central Limit Theorem (see e.g. Billingsley \cite{Billingsley1995}) to each coordinate of the vector $\mathbf{Z}$. 
Observe that for all $i=1,\dots,d$ by the It\^o isometry we obtain
\begin{align*}
\Var\left( X^{(2i-1)}_k(\omega_i)\right)&=\Var\left( \int_{2(k-1)\pi}^{2k\pi} \cos(\omega_i t) dL(t) \right) = \sigma^2 \int_{2(k-1)\pi}^{2k\pi} \cos^2(\omega_i t) dt
\\&=\sigma^2\frac{4\pi \omega_i+\sin(4\pi \omega_i k)-\sin(4\pi \omega_i (k-1))}{4\omega_i}.
\end{align*}
Thus
$$\left(s_n^{(2i-1)}\right)^2=\sum_{k=1}^n \Var \[X_k^{(2i-1)}(\omega_i)\]=\sigma^2\frac{4n\pi \omega_i+\sin(4\pi \omega_i n)}{4\omega_i}.$$
In the same way,
$$\left(s_n^{(2i)}\right)^2=\sum_{k=1}^n \Var \[X_k^{(2i)}(\omega_i)\]=\sigma^2\frac{4n\pi \omega_i-\sin(4\pi \omega_i n)}{4\omega_i}.$$
Observe that $$\lim_{n\to \infty}  \frac{1}{n} \left(s_n^{(2i-1)}\right)^2= \lim_{n\to \infty}  \frac{1}{n} \left(s_n^{(2i)}\right)^2=\sigma^2\pi.$$
If the Lindeberg condition is satisfied, the $2i$-th, respectively $2i-1$-th coordinate of $\mathbf{Z}$ for $i=1,\dots,d$, i.e.
\begin{align*}
Z^{(2i-1)}_n(\omega_i)&=\frac{2\sqrt{\omega_i}}{\sigma\sqrt{4\pi n\omega_i+\sin(4\pi n\omega_i)}}\int_0^{2\pi n} \cos(\omega_i t)dL(t)\\
Z^{(2i)}_n(\omega_i)&=\frac{2\sqrt{\omega_i}}{\sigma\sqrt{4\pi n\omega_i-\sin(4\pi n\omega_i)}}\int_0^{2\pi n} \sin(\omega_i t)dL(t)\\
\end{align*}
converges to $\mathcal{N}(0,1)$. Taking
$$Y^{(2i-1)}_n(\omega_i)=\frac{\sigma \sqrt{4\pi n\omega_i+\sin(4\pi n\omega_i)}}{2\sqrt{2\pi n \omega_i}}, 
\quad Y^{(2i)}_n(\omega_i)=\frac{\sigma \sqrt{4\pi n\omega_i-\sin(4\pi n\omega_i)}}{2\sqrt{2\pi n \omega_i}} $$
and noting that 
$$\lim_{n\to \infty} Y^{(2i-1)}_n(\omega_i)=\frac{\sigma}{\sqrt{2}},\quad \lim_{n\to \infty} Y^{(2i)}_n(\omega_i)=\frac{\sigma}{\sqrt{2}}  $$
is constant at all frequencies, by Slutsky arguments for $i=1,\dots, d$ we get
\begin{align}
\frac{1}{\sqrt{2\pi n}}\int_0^{2\pi n} \cos(\omega_i t) dL(t)&=Z^{(2i-1)}_n(\omega_i)Y^{(2i-1)}_n(\omega_i)\to \mathcal{N}\left(0,\frac{\sigma^2}{2}\right), \label{zy1}\\
\frac{1}{\sqrt{2\pi n}}\int_0^{2\pi n} \sin(\omega_i t) dL(t)&=Z^{(2i)}_n(\omega_i)Y^{(2i)}_n(\omega_i)\to \mathcal{N}\left(0,\frac{\sigma^2}{2}\right)\label{zy2}.
\end{align}

Now we are going to prove the Lindeberg condition for odd coordinates of $\mathbf{Z}$ (for the even ones an analogous reasoning holds), i.e. for all $\epsilon>0$ it holds
$$ \lim_{n\to \infty} \frac{1}{\left(s_n^{(2i-1)}\right)^2} \sum_{k=1}^n  \E\[\left(X^{(2i-1)}_k(\omega_i)\right)^2\1_{\left\{| X^{(2i-1)}_k(\omega_i)|>\epsilon s_n^{(2i-1)}\right\}}\]=0.$$
Observe that if the random variables $\{X^{(2i-1)}_k(\omega_i)\}$ are uniformly square integrable, then they satisfy the Lindeberg condition. Indeed,
\begin{align*}
&\left(s_n^{(2i-1)}\right)^{-2}\sum_{k=1}^n  \E\[\left(X^{(2i-1)}_k(\omega_i)\right)^2\1_{\left\{\left| X^{(2i-1)}_k(\omega_i)\right|>\epsilon s_n^{(2i-1)}\right\}}\]
\\&= \left(s_n^{(2i-1)}\right)^{-2}\sum_{k=1}^n  \E\[\left(X^{(2i-1)}_k(\omega_i)\right)^2\1_{\left\{\left| X^{(2i-1)}_k(\omega_i)\right|^2>\left(\epsilon s_n^{(2i-1)}\right)^2\right\}}\]
\\&\leq \left(s_n^{(2i-1)}\right)^{-2} n \sup_{k=1,...,n} \E\[\left(X^{(2i-1)}_k(\omega_i)\right)^2\1_{\left\{\left| X^{(2i-1)}_k(\omega_i)\right|^2>\left(\epsilon s_n^{(2i-1)}\right)^2\right\}}\]
\\&= \left(\sigma^2\frac{4n\pi \omega_i+\sin(4\pi \omega_i n)}{4\omega_i}\right) ^{-1} n \sup_{k=1,...,n} \E\[\left(X^{(2i-1)}_k(\omega_i)\right)^2
\1_{\left\{\left| X^{(2i-1)}_k(\omega_i)\right|^2>\left(\epsilon s_n^{(2i-1)}\right)^2\right\}}\]
%\frac{n}{nc} \sup_{i=1,...,n}E(X_i^21_{|X_i|^2\geq \epsilon^2 C_n^2)
\end{align*}
Since $\lim_{n\to \infty}\left(\sigma^2\frac{4n\pi \omega_i+\sin(4\pi \omega_i n)}{4\omega_i}\right) ^{-1} n\to \frac{1}{\pi \sigma^2}$, uniform square integrability implies in this case the Lindeberg condition. 
It remains to show the uniform square integrability of $\left\{ X^{(2i-1)}_k(\omega_i) \right\}_{k\in \N}$.

Assume first, that our driving process $(L(t))_{t\geq 0}$ is of bounded variation. 
Then
$$M_k= \left| \int_{2(k-1)\pi}^{2k\pi} \cos(\omega_i t) dL_t\right|\leq  \int_{2(k-1)\pi}^{2k\pi} |\cos(\omega_i t)| d|L_t|\leq  \int_{2(k-1)\pi}^{2k\pi}  d|L_t|,$$
where $|\cdot|$ denotes the total variation of the process.
But $\int_{2(k-1)\pi}^{2k\pi}  d|L_t|\stackrel{d}{=}\int_{0}^{2\pi}  d|L_t|$. 
We have
\begin{align*}
\E\[ |M_k|^2\1_{\left\{|M_k|>K\right\}} \]&\leq \E\[ \left| \int_{2(k-1)\pi}^{2k\pi} d|L_t|\right|^2\1_{\left\{\left| \int_{2(k-1)\pi}^{2k\pi} d|L_t|\right|>K\right\}} \]
\\&= \E\[ \left| \int_{0}^{2\pi} d|L_t|\right|^2\1_{\left\{\left| \int_{0}^{2\pi} d|L_t|\right|>K\right\}} \].
\end{align*}
By the square integrability of $\int_{0}^{2\pi} d|L_t|$, which is implied by the square integrability of $(L(t))_{t\geq 0}$ we obtain the uniform integrability of $(M_k)_{k\in \N}$.

Now we assume that $(L(t))$ is a square integrable martingale with finite moments of all orders. Observe that $X_k$ is square integrable for all $k\in \N$. By the Burkholder-Davis-Gundy Inequality (see e.g. Protter \cite{Protter})
 for each $p\geq 1$ there exists a positive constant $C_p$ such that
$$\E\[\left( \int_{2(k-1)\pi}^{2k\pi} \cos(\omega_i t) dL(t)  \right)^p  \]\leq C_p\E\[\[ \int_{2(k-1)\pi}^{2k\pi} \cos(\omega_i t) dL(t), \int_{2(k-1)\pi}^{2k\pi} \cos(\omega_i t) dL(t)\]^{p/2}\].$$
Since $$ \[ \int_{2(k-1)\pi}^{2k\pi} \cos(\omega_i t) dL(t), \int_{2(k-1)\pi}^{2k\pi} \cos(\omega_i t) dL(t)\]=  \int_{2(k-1)\pi}^{2k\pi} \cos^2(\omega_i t)d[L,L]_t$$
using the above inequality for $p=4$ we obtain
\begin{align*}\E\[\left( \int_{2(k-1)\pi}^{2k\pi} \cos(\omega_i t) dL(t)  \right)^4  \]&\leq C_4\E\[ \int_{2(k-1)\pi}^{2k\pi} \cos^2(\omega_i t)d[L,L]_t\]\\&\leq
\sigma^2  \int_{2(k-1)\pi}^{2k\pi} \cos^2(\omega_i t)dt <C\end{align*}
for some constant $C$. Since $\{X^{(2i-1)}_k(\omega_i)\}$ are square integrable and $\{X^{(2i-1)}_k(\omega_i)\}$ are bounded in $L^4(\Omega,\mathcal{F},\P)$ they are uniformly square integrable.

As any L\'evy process is by the L\'evy-It\^o decomposition the sum of a finite variation L\'evy process and an independent square integrable martingale with moments of all orders, we obtain the claimed uniform squre integrability for all driving L\'evy processes.

Likewise one shows that $\theta^T Z$ converges in distribution to $\mathcal{N}\left(0,\frac{\sigma^2}{2}\theta^T \theta\right)$ for all $\theta \in \R^{2d}$. So the Cramer-Wold device concludes.

Therefore $\[Z^{(2i-1)}_n(\omega_i)Y^{(2i-1)}_n(\omega_i),Z^{(2i)}_n(\omega_i)Y^{(2i)}_n(\omega_i) \]^T_{i=1,\dots,d}$ converges in distribution to 
$\mathcal{N}\left(0,\frac{\sigma^2}{2} I_{2d\times 2d}\right)$ and thus using Lemma \ref{lem:LimitInProbabilityFirstSummandFT} and equations \eqref{zy1}, \eqref{zy2} $\[\mathcal{F}_T(\omega_j) \]^T_{j=1,\dots,d}$ converges to  $\mathcal{N}\left(0,\frac{\sigma^2}{2} \mathbf{B}\right)$, 
where $\mathbf{B}$ is defined above.
Repeating the reasoning from the proof of Theorem  \ref{lem:RozkladZ(T)Squared} we obtain that $\[|\mathcal{F}_T(\omega_j)|^2 \]^T_{j=1,\dots,d}$ 
converges to a vector of independent, exponentially distributed random variables with 
$\operatorname{Exp}\left(\sigma^2\left| \frac{b(i\omega_j)}{a(i\omega_j)} \right|^2\right)$ for $j=1,\dots,d$. $\Box$

Note that Theorem \ref{lem:RozkladZ(T)} is basically a special case of Theorem \ref{thm:LimitTheoremSeveralFrequencies} . However, the proof in the case of several frequencies is much more complicated and a more elementary reasoning was also presented. 

The limiting result is the analogue of the one for discrete time $\operatorname{ARMA}$ models. (See e.g.  \cite{BrockwellDavis} Chapter 10.)
%!!!!!!!!!!!! Is it possible that $\lambda_{1,2}=0$ - Give the interpretation !!!!!!!!!!!!!!!!!!!!!!!!!
%-------------------------------------------------------------------------------------
%-------------------------------------------------------------------------------------
%-------------------------------------------------------------------------------------
%-------------Approximation of Integral-----------------------------------------------
%-------------------------------------------------------------------------------------
%-------------------------------------------------------------------------------------
%-------------------------------------------------------------------------------------
%----------------------------------------------
\subsection{Numerical Approximation of Integrals and Limiting Behaviour of the Truncated Pathwise Fourier Transform Based on Non-equidistant Discrete Grids}\label{Section:ApproximationofIntegrals}
%-----------------------------------------------
In this section we deal with the numerical approximation of the integral
\begin{equation}
\F_T(Y)(\omega):=\frac{1}{\sqrt{T}}\int_0^T Y(t)e^{-i\omega t}dt.
\label{eq:FT_numerics}
\end{equation}
Our aim is to describe conditions under which we are able to calculate numerically the truncated Fourier transform of a $\operatorname{CARMA}$ process based on non-equidistant observations.
The main result is the following:
\begin{theorem}\label{thm:ErrorEstimationIntegralTR}
Let $\X$ and $Y$ be processes given by the state-space representation  \eqref{eq:ObservationArma} and \eqref{eq:StateRepresentationArma}. 
Suppose that Assumptions \ref{ass:FiniteSecondMomentsLevy}, \ref{ass:IndependenceLevyAndX}, \ref{ass:EigenvaluesDistinctNonnegativeParts} and \ref{ass:EqualityOfDistributionAtZero} are satisfied. Assume that $F\colon \R \to \R^d$ be a twice continuously differentiable function with $\| F''\|_\infty<\infty$.
Let $\left(x_i^{(T)}\right)_{i=0,\dots,N(T)-2}$ be a partition of the interval $[a,b]$ with $x^{(T)}_0=a$ and $x^{(T)}_{N(T)-1}=b$ and let $h_{\max}(T)=\max_{j=0,\dots, N(T)-1}\left(x^{(T)}_{j+1}-x^{(T)}_{j}\right)$.
Put 
\begin{equation}
\alpha_0^{(N(T))}=\frac{x^{(T)}_1-x^{(T)}_0}{2}F\left(x^{(T)}_0\right),\quad \alpha_{N(T)-1}^{(N(T))}=\frac{x^{(T)}_{N(T)-1}-x^{(T)}_{N(T)-2}}{2}F\left(x^{(T)}_{N(T)-1}\right),
\label{eq:TNFHWithAlphasDef0nminus1}
\end{equation}
\begin{equation}
\alpha_j^{(N(T))}=\frac{x^{(T)}_{j+1}-x^{(T)}_{j-1}}{2}F\left(x^{(T)}_j\right),\quad  j=1,\dots,N(T)-2.
\label{eq:TNFHWithAlphasDef1nminus2}
\end{equation}

Then there exist positive constants $C_1,C_2$ such that 
$$\E\[\left \| \sum_{j=0}^{N(T)-1}\alpha_j^{(N(T))}Y\left(x_j^{(T)}\right)-\int_a^bY(t)F(t)dt\right \|^2\]\leq C_1 (C_2+T) N(T)^2 h^6_{\max}(T)   $$
and thus if  $\lim_{T\to \infty}T N(T)^2 h^6_{\max}(T) =0$, then 
$$\lim_{T\to \infty}\left\|\sum_{j=0}^{N(T)-1}\alpha_j^{(N(T))}F\left(x_j^{(T)}\right)-\int_a^bY(t)F(t)dt\right\|_{L^2}=0 .$$ 

\end{theorem}

We begin by establishing an error bound of the trapezoidal method for non-equidistant data.
For a very accessible approach of quadrature rules procedures we refer to \cite{TalvilaWiersma2012}.
Recall the basic properties of the trapezoidal rule:
\begin{lemma}\label{lem:TrapezoidalRuleEquistantError}
Let $f\colon [a,b]\to \R$ be a twice continuously differentiable function. Write 
\begin{equation}
\int_a^b f(x)dx=\frac{b-a}{2}[f(a)+f(b)]+E^T(f).
\label{eq:TROneStepApproximation}
\end{equation}
Then 
\begin{equation}
|E^T(f)|\leq \frac{(b-a)^3}{12}\sup_{x\in[a,b]}| f''(x)|.
%\label{eq:TROneStepApproximation}
\end{equation}
For the composite trapezoidal rule for an equidistant grid $a< a+(b-a)\frac{1}{n}<\dots, a+(b-a)\frac{i}{n} <\dots< b$ we have
\begin{equation}
\int_a^b f(x)dx=\frac{b-a}{2n} \[ f(a)+2\sum_{i=1}^{n-1}f\left(a+(b-a)\frac{i}{n}\right) +f(b)\] +E_n^T(f).
\label{eq:TROneStepApproximationComposite}
\end{equation}
Then 
\begin{equation}
|E_n^T(f)|\leq \frac{(b-a)^3}{12 n^2}\sup_{x\in[a,b]} |f''(x)|.
%\label{eq:TROneStepApproximationComposite}
\end{equation}
\end{lemma}
A proof can be found e.g. in \cite{TalvilaWiersma2012}. 

Now we are going to formulate a version of the trapezoidal rule for non-equidistant points. We assume that we have some control on the maximal distance between observations.

\begin{lemma}\label{lem:TRNonEquidistant}
Let $a=x_0<x_1<\dots<x_{N-1}<x_{N}=b$ be an arbitrary partition of the interval $[a,b]$ and assume that  $f\colon [a,b]\to \R$ is a twice continuously differentiable function. 
Put $h_{\max}=\max_{j=0,\dots, N-2}(x_{j+1}-x_{j})$.  
%If $$ N h_{\max}(T) ^3 \to 0,\quad N\to \infty, $$
Then
$$ \int_a^b f(x)dx = \sum_{i=0}^{N-1}\frac{x_{j+1}-x_{j}}{2}[f(x_{j})+f(x_{j+1})]+E^T(f),$$
where $|E^T(f)|\leq  N \|f''\|_\infty \frac{h_{\max} ^3}{12}$.
\end{lemma}
\textsc{Proof}

Let us write
$$ [a,b]=\bigcup_{j=0}^{N-1} [x_j,x_{j+1}], \quad I_j:=[x_j,x_{j+1}]$$
and apply Lemma \ref{lem:TrapezoidalRuleEquistantError} for each interval $I_j$. Therefore
$$ \int_{x_j}^{x_{j+1}} f(x)dx=\frac{x_{j+1}-x_{j}}{2}[f(x_{j})+f(x_{j+1})]+E_j^T(f),$$
with
$$ |E_j^T(f)|\leq \frac{|x_{j+1}-x_{j}|^3}{12 }\sup_{x\in [x_j,x_{j+1}]} |f''(x)|. $$
For each $i=0,1,\dots, N-1$ we have
$$ \sup_{x\in [x_j,x_{j+1}]}| f''(x) |\leq \sup_{x\in [a,b]}| f''(x)|=:\|f''\|_\infty.$$
Therefore
$$ \int_a^b f(x)dx = \sum_{i=0}^{N-1}\frac{x_{j+1}-x_{j}}{2}[f(x_{j})+f(x_{j+1})]+E^T(f),$$
where
\begin{align*}
|E^T(f)|&= \left|\sum_{i=0}^{N-1} E_j^T(f)\right|\leq  \|f''\|_\infty \sum_{i=0}^{N-1} \frac{(x_{j+1}-x_{j})^3}{12 }
\\&\leq \|f''\|_\infty \sum_{i=0}^{N-1} \frac{h_{\max}^3}{12} = N \|f''\|_\infty \frac{h_{\max} ^3}{12}.
\end{align*}
This completes the proof. $\Box$

We use some results and ideas from \cite{BrockwellSchlemm}. The aim is to find an approximation similar to Proposition 5.4 of \cite{BrockwellSchlemm} of the integral appearing in the truncated Fourier transform 
in the case that the observations of the process $Y$ are given on a non-equidistant grid.
Let 
$$T^N_{[0,T]}f=\sum_{j=0}^{N-1}\frac{x_{j+1}-x_{j}}{2}[f(x_{j})+f(x_{j+1})] $$
be the trapezoidal rule discussed in Lemma \ref{lem:TRNonEquidistant}. 
Recall first the Fubini type theorem for stochastic integrals from \cite{BrockwellSchlemm}.
\begin{lemma}\cite[Theorem 2.4]{BrockwellSchlemm}\label{lem:FubiniUnboundedDomain}
Let $[a,b]\subset \R$ be a bounded interval and $(L(t))_{t\geq 0}$ be a L\'evy process with finite second moments. Assume that $F\colon [a,b]\times \R \to \R^{d}$ is a bounded function $\mathcal{B}([a,b])\otimes \mathcal{B}([-s,t])$-measurable for all $s,t\in (0,\infty) $ and the family
$\{  u\mapsto F(s,u)\}_{u\in [a,b]}$ is uniformly absolutely integrable and uniformly converges to zero as $|u|\rightarrow 0$. Then
\begin{equation}
\int_a^b \int_\R F(s,u)dL(u)ds= \int_\R \int_a^bF(s,u)ds dL(u) \quad a.s.
\label{eq:StochasticFubiniUnbounded}
\end{equation}
\end{lemma}

In the paper \cite{BrockwellSchlemm} the assumption about measurability in the statement of the theorem is not explicitely stated. However, an inspection of their proof combined with results from \cite{Veraar2012} shows that the precise statement has to be in the above form.

Secondly, note that for non-equidistant data the corresponding error estimation \cite[Proposition A.6]{BrockwellSchlemm} has the following form:

\begin{prop}\label{prop:OszacowaniaNormy}
Let $[a,b]\subset \R$ be a compact interval and use the notation of Lemma \ref{lem:TRNonEquidistant}. 
\begin{enumerate}
	\item If $f\colon [a,b]\to \R$ is twice continuously differentiable, then 
	$$\left|\int_a^b f(s)ds-T^N_{[a,b]}f\right|\leq N \| f'' \|_\infty \frac{h^3_{\max}}{12}. $$
		\item If $F\colon [a,b]\to \R^d$ is twice continuously differentiable, then 
	$$\left\|\int_a^b F(s)ds-T^N_{[a,b]}F\right\|\leq \sqrt{d} N \| F'' \|_\infty \frac{h^3_{\max}}{12}, $$
	where $\| \cdot \|$ denotes the Euclidean norm in $\R^d$.
			%\item If $F\colon [a,b]\to \R^{d\times d}$ is twice differentiable, then 
	%$$\left\|\int_a^b F(s)ds-T^N_{[a,b]}F\right\|_2\leq \sqrt{d^3} N \| F'' \|_\infty \frac{h^3_{\max}}{12}, $$
%	where $\| \cdot \|_2$ denotes the operator  norm in $\R^{d\times d}$.
\end{enumerate}
Here $\| F'' \|_\infty:=\max_{i=1,\dots,d}\sup_{t_i\in[a,b]}\| F''(t_i) \|$. 
\end{prop}

Put
$$ E^{T,N}_{fg}:=T^N_{[0,T]}f(\cdot)g(\cdot)-\int_0^Tg(s)f(s)ds.$$

\textsc{Proof of Theorem \ref{thm:ErrorEstimationIntegralTR}}. Assume that we have observed the process $Y$ on the grid 
$0=x^{(T)}_0<x^{(T)}_1<\dots<x^{(T)}_{N(T)-1}=T$. We have
\begin{equation}
T^N_{[0,T]}FY=\sum_{j=0}^{N(T)-1}\alpha_j^{(N(T))}Y(x^{(T)}_j),
\label{eq:TNFHWithAlphas}
\end{equation}
where $\alpha_j^{(N(T))}$ ($j=0,\dots,N(T)-1$) are the coefficients given by \eqref{eq:TNFHWithAlphasDef0nminus1} and \eqref{eq:TNFHWithAlphasDef1nminus2}.

Observe that  $f(a)=\int f(s)\delta_a(s)ds$. Moreover, for all $j=0,\dots, N(T)-1$ we know that $x^{(T)}_j\in [0,T]$, therefore for all $u\in [0,T]$ and for all $j=0,\dots, N(T)-1$
we have
$$ \1_{[u,T]}\left(x^{(T)}_j\right)=  \1_{\[0,x^{(T)}_j\]}(u).$$
Thus by \eqref{eq:TNFHWithAlphas} we have
\begin{align*}
&T^N_{[0,T]}FY=\sum_{j=0}^{N(T)-1}\alpha_j^{(N(T))}Y\left(x_j^{(T)}\right)
=\sum_{j=0}^{N(T)-1}\alpha_j^{(N(T))}\int_{-\infty}^{x^{(T)}_j} \b^T e^{\A\left(x^{(T)}_j-u\right)}\e dL^*(u)
\\&=\sum_{j=0}^{N(T)-1}\alpha_j^{(N(T))}\left(\int_{-\infty}^{0} \b^T e^{\A\left(x^{(T)}_j-u\right)}\e dL^*(u)+\int_{0}^{x^{(T)}_j} \b^T e^{\A\left(x^{(T)}_j-u\right)}\e dL^*(u)\right)
\\&=\int_{-\infty}^{0} \sum_{j=0}^{N(T)-1}\alpha_j^{(N(T))}\b^T e^{\A\left(x^{(T)}_j-u\right)}\e dL^*(u)\\&+\int_{0}^{T} \sum_{j=0}^{N(T)-1}\alpha_j^{(N(T))}\b^T e^{\A\left(x^{(T)}_j-u\right)}\e \1_{\[0,x^{(T)}_j\]}(u)dL^*(u)
\\&=\int_{-\infty}^{0} \sum_{j=0}^{N(T)-1}\alpha_j^{(N(T))}\b^T e^{\A\left(x^{(T)}_j-u\right)}\e dL^*(u)\\&+\int_{0}^{T} \sum_{j=0}^{N(T)-1}\alpha_j^{(N(T))}\b^T e^{\A\left(x^{(T)}_j-u\right)}\e \1_{[u,T]}\left(x^{(T)}_j\right)dL^*(u)
\\&=\int_{-\infty}^{0}\int_0^T \sum_{j=0}^{N(T)-1}\alpha_j^{(N(T))} \delta_{x^{(T)}_j}(s)\b^T e^{\A(s-u)}\e ds dL^*(u)
\\&+\int_{0}^{T} \int_u^T\sum_{j=0}^{N(T)-1}\alpha_j^{(N(T))} \delta_{x_j^{(T)}}(s)\b^T e^{\A(s-u)}\e ds dL^*(u)
\\&=\int_{-\infty}^T  \int_{\max\{0,u\}}^T \sum_{j=0}^{N(T)-1}\alpha_j^{(N(T))} \delta_{x_j^{(T)}}(s) \b^Te^{\A(s-u)}\e ds dL^*(u).
\end{align*}
Thus using the representation \eqref{eq:CarmaSolution} and the Fubini-type Theorem \ref{lem:FubiniUnboundedDomain}  we have
\begin{align*}
\int_0^TF(s)Y(s)ds&=\int_0^TF(s)\int_{-\infty}^s \b^Te^{\A(s-u)}\e dL^*(u) ds
\\&=\int_{-\infty}^T  \int_{\max\{0,u\}}^T F(s) \b^Te^{\A(s-u)}\e ds dL^*(u)
\\&=\int_{-\infty}^0  \int_{0}^T F(s) \b^Te^{\A(s-u)}\e ds dL^*(u)
\\& +\int_{0}^T  \int_{u}^T F(s) \b^Te^{\A(s-u)}\e ds dL^*(u).
\end{align*}
Thus
\begin{align*}
E^{T,N}_{FY}&=T^N_{[0,T]}FY-\int_0^TF(s)Y(s)ds= \int_0^T\left(\sum_{j=0}^{N(T)-1}\alpha_j^{(N(T))}\delta_{x_j^{(T)}}(s)-F(s)\right)Y(s) ds
\\&=\int_{-\infty}^0  \int_{0}^T\left(\sum_{j=0}^{N(T)-1}\alpha_j^{(N(T))}\delta_{x_j^{(T)}}(s)-F(s)\right) \b^Te^{\A(s-u)}\e ds dL^*(u)
\\& +\int_{0}^T  \int_{u}^T\left(\sum_{j=0}^{N(T)-1}\alpha_j^{(N(T))}\delta_{x_j^{(T)}}(s)-F(s)\right) \b^Te^{\A(s-u)}\e ds dL^*(u).
\end{align*}
Let us denote
\begin{align}
\label{eq:gammaNDfTR}
\Gamma^{(N)}(u):&= \int_{0}^T\left(\sum_{j=0}^{N(T)-1}\alpha_j^{(N(T))}\delta_{x_j^{(T)}}(s)-F(s)\right) \b^Te^{\A(s-u)} \e ds, \quad u\leq 0,
\\
%\label{eq:gammaNDfTR}
G^{(N)}(u):&= \int_u^T \left(\sum_{j=0}^{N(T)-1}\alpha_j^{(N(T))}\delta_{x_j^{(T)}}(s)-F(s)\right) \b^Te^{\A(s-u)}  \e ds,\quad u\in [0,T].
\end{align}
By Assumption \ref{ass:EigenvaluesDistinctNonnegativeParts} we know that there exist positive constants $\alpha$, $\beta$ such that
\begin{equation}
\|\exp(\A t)\|\leq  \beta \exp(-\alpha t).
\label{eq:OszacowanieNormyWyklHurwiz}
\end{equation}
%!!!Submultiplicativity!!!
Note that by Lemma \ref{lem:TRNonEquidistant} and Proposition \ref{prop:OszacowaniaNormy} we have
\begin{align*}
&\left\|\int_{u_0}^T\left(\sum_{j=0}^{N(T)-1}\alpha_j^{(N(T))}\delta_{x_j^{(T)}}(s)-F(s)\right) \b^Te^{\A(s-u)} \e ds\right\|_{\R^d}
\\&= \left\|\sum_{j=0}^{N(T)-1}\alpha_j^{(N(T))}\delta_{x_j^{(T)}}(s)\b^Te^{\A(x_j-u)} \e-\int_{u_0}^T F(s)\b^Te^{\A(s-u)} \e ds\right\|_{\R^d}
 \\&\leq \sqrt{d} N(T) \| \widetilde {F}''(u) \|_\infty \frac{h^3_{\max}(T)}{12}
\end{align*}
with $\widetilde {F}(s)=F(s)\b^Te^{\A(s-u)} \e$ and $u_0\in [0,T]$. If $u_0=0$, then 
there exist $\widetilde{\alpha}>0$ and $D>0$ such that $\| \widetilde {F}''(u) \|_\infty\leq D \exp(\widetilde{\alpha}u)$ for $u\leq 0$.
Therefore there exists a constant $D_1>0$ such that
$$\| \Gamma^{(N)}(u)\|_{\R^d}\leq \sqrt{d} N(T) \| \widetilde {F}|_{[0,T]}'' \|_\infty \frac{h^3_{\max}(T)}{12} \leq D_1 N(T) h^3_{\max}(T) \exp(\widetilde{\alpha}u), \quad u\leq 0.$$

If now $u_0=u$ is any element of $[0,T]$, then there exists $D>0$ such that $\| \widetilde {F}''|_{[u,T]} \|_\infty\leq D $ for $u\in[0,T]$.
Therefore there exists a constant $D_2>0$ such that for $u\in[0,T]$ we have
$$\| G^{(N)}(u)\|_{\R^d}\leq \sqrt{d} N(T) \| \widetilde {F}|_{[u,T]}'' \|_\infty \frac{h^3_{\max}(T)}{12} \leq D_2 N(T) h^3_{\max}(T).$$

By the It\^o isometry 
\begin{align*}
\left\|\int^0_{-\infty} \Gamma^{(N)}(u)dL^*(u)\right\|^2_{L^2}&= \E\[\left(\int^0_{-\infty}  \Gamma^{(N)}(u)dL^*(u)\right)^T\left(\int^0_{-\infty}  \Gamma^{(N)}(u)dL^*(u)\right)\]
\\&= \E\[\left(\int^0_{-\infty}  \Gamma^{(N)}(u)dL^*(u)\right)^T\left(\int^0_{-\infty}  \Gamma^{(N)}(u)dL^*(u)\right)\]
\\&=\sigma^2 \left(\int^0_{-\infty}  [\Gamma^{(N)}(u)]^T \Gamma^{(N)}(u)du\right)
\\&\leq \sigma^2 \int^0_{-\infty}  \| \Gamma^{(N)}(u) \|_{\R^d}^2 du
\\&\leq \sigma^2 \int^0_{-\infty} \left(D_1 N(T) h^3_{\max}(T) \exp(\widetilde{\alpha}u)\right)^2du
\\&=\sigma^2N(T)^2 h^6_{\max}(T)D_1^2 \int^0_{-\infty}\exp(2\widetilde{\alpha}u)du
\\&=D_{\Gamma} N(T)^2 h^6_{\max}(T),
\end{align*}
where $D_\Gamma>0$ is a constant. In a similar way we obtain 
\begin{align*}
\left\|\int_0^T G^{(N)}(u)dL(u)\right\|^2_{L^2}&= \E\[\left(\int_0^T G^{(N)}(u)dL(u)\right)^T\left(\int_0^T G^{(N)}(u)dL(u)\right)\]
\\&= \E\[\left(\int_0^T G^{(N)}(u)dL(u)\right)^T\left(\int_0^T G^{(N)}(u)dL(u)\right)\]
\\&=\sigma^2 \left(\int_0^T [G^{(N)}(u)]^T G^{(N)}(u)du\right)
\\&\leq \sigma^2 \int_0^T \| G^{(N)}(u) \|_{\R^d}^2 du\leq \sigma^2 \int_0^T \left( D_2 N(T) h^3_{\max}(T) \right)^2du
\\&= D_G N(T)^2 T h^6_{\max}(T) 
\end{align*}
for some constant $D_G>0$. Therefore
$$ \| E_{FY}^{T,N} \|_{L^2}^2\leq 2\[ \left\|\int^0_{-\infty} \Gamma^{(N)}(u)dL^*(u)\right\|^2_{L^2}+\left\|\int_0^T G^{(N)}(u)dL(u)\right\|^2_{L^2}\],$$
thus
$$ \| E_{FY}^{T,N} \|_{L^2}^2\leq C_1 \left((C_2+T)N(T)^2 h^6_{\max}(T)\right)^2,$$
where $C_1,C_2$ are positive constants. If %$\lim_{T\to \infty}N^2 h^6_{\max}(T) =0$ and 
$\lim_{T\to \infty}T N(T)^2 h^6_{\max}(T)=0$, then $\lim_{T\to \infty}\| E_{FY}^{T,N} \|_{L^2}^2=0$. 
This completes the proof. $\Box$

Now we are going to apply Theorem \ref{thm:ErrorEstimationIntegralTR} to find a numerical approximation of the truncated Fourier transform. Using the notation of Theorem \ref{thm:ErrorEstimationIntegralTR} we denote the trapezoidal approximation of 
$$\mathcal{F}_T(Y)(\omega)=\frac{1}{\sqrt{T}}\int_0^T Y(t)e^{-i\omega t}dt $$
by $\mathcal{T}_T(Y)(\omega)$, i.e.
$$ \mathcal{T}_T(Y)(\omega) =\frac{1}{\sqrt{T}}\sum_{j=0}^{N-1}\alpha_j^{(N)} Y\left( x_j^{(N)} \right),$$
where the grid points $ \left(x_j^{(N)}\right)_{j=0,\dots,N(T)-1}$ are given as in Theorem \ref{thm:ErrorEstimationIntegralTR} and \begin{equation*}
\alpha_0^{(N(T))}=\frac{x^{(T)}_1-x^{(T)}_0}{2}F\left(x^{(T)}_0\right),\quad \alpha_{N(T)-1}^{(N(T))}=\frac{x^{(T)}_{N(T)-1}-x^{(T)}_{N(T)-2}}{2}F\left(x^{(T)}_{N(T)-1}\right),
\end{equation*}
\begin{equation*}
\alpha_j^{(N(T))}=\frac{x^{(T)}_{j+1}-x^{(T)}_{j-1}}{2}F\left(x^{(T)}_j\right),\quad  j=1,\dots,N(T)-2.
\end{equation*} with $F(x)=e^{-i\omega x}$.

\begin{theorem}\label{thm:ApproxTFTDoubleLimit}
Let $\X$ and $Y$ be processes given by the state-space representation  \eqref{eq:ObservationArma} and \eqref{eq:StateRepresentationArma}. 
Suppose that Assumptions \ref{ass:FiniteSecondMomentsLevy}, \ref{ass:IndependenceLevyAndX}, \ref{ass:EigenvaluesDistinctNonnegativeParts} and \ref{ass:EqualityOfDistributionAtZero} are satisfied and that the process $Y$ is observed at not  necessarily equidistant points $0=x^{(T)}_0<x^{(T)}_1<\dots<x^{(T)}_{N(T)-1}=T$. Let $h_{\max}(T):=\max_{j=0,\dots, N(T)-2}\left(x^{(T)}_{j+1}-x^{(T)}_j\right)$. 
If 
 $$ \lim_{T\to \infty}N(T)h^3_{\max}(T)=0,$$
then
$$  \lim_{T\to \infty}\| \mathcal{T}_T(Y)(\omega)- \mathcal{F}_T(Y)(\omega)\|_{L^2}=0  $$
and thus also 
$$ \P- \lim_{T\to \infty}\[ \mathcal{T}_T(Y)(\omega)- \mathcal{F}_T(Y)(\omega)\]=0.$$
\end{theorem}
\textsc{Proof} We identify $\C$ with $\R^2$ in the canonical way.
Applying Theorem \ref{thm:ErrorEstimationIntegralTR} for $d=2$, $F(t)=[\cos(\omega t),-\sin(\omega t)]^T$ we get 
$$\E\[\left \| \sum_{j=0}^{N(T)-1}\alpha_j^{(N(T))} Y\left( x_j^{(T)} \right)-\int_0^TY(t) 
     \[\begin{matrix}
       \cos(\omega t)\\
       -\sin(\omega t)\\
        \end{matrix}\]
dt\right \|^2\]\leq C_1 (C_2+T) N(T)^2 h^6_{\max}(T)  .$$
Dividing both sides by $T>0$ we obtain
$$\E\[\left \| \sum_{j=0}^{N(T)-1}\frac{\alpha_j^{(N(T))}}{\sqrt{T}} Y\left( x_j^{(T)} \right)-\frac{1}{\sqrt{T}}\int_0^TY(t) 
     \[\begin{matrix}
       \cos(\omega t)\\
       -\sin(\omega t)\\
        \end{matrix}\]
dt\right \|^2\]\leq \frac{C_1 C_2}{T}+ C_1N(T)^2 h^6_{\max}(T)  .$$

Passing to the limit with $T\to \infty$ and using the assumption $ \lim_{T\to \infty}N(T)h^3_{\max}(T)=0$  we get the assertion. $\Box$

Now we are going to state the central limit theorem for the truncated Fourier transform:
\begin{theorem}\label{thm:ApproxTFTDoubleLimitDistribution}
Let $\X$ and $Y$ be processes given by the state-space representation  \eqref{eq:ObservationArma} and \eqref{eq:StateRepresentationArma}. 
Suppose that Assumptions \ref{ass:FiniteSecondMomentsLevy}, \ref{ass:IndependenceLevyAndX}, \ref{ass:EigenvaluesDistinctNonnegativeParts} and \ref{ass:EqualityOfDistributionAtZero} are satisfied and the process $Y$ is observed at not necessarily equidistant points $0=x^{(T)}_0<x^{(T)}_1<\dots<x^{(T)}_{N(T)-1}=T$. Let $h_{\max}(T):=\max_{j=0,\dots, N(T)-2}(x^{(T)}_{j+1}-x^{(T)}_j)$. 
Let $\alpha_j^{(N)}$ be defined as in Theorem \ref{thm:ErrorEstimationIntegralTR}.
Assume that 
 $$ \lim_{T\to \infty}N(T)h^3_{\max}(T)=0.$$
Put $\Sigma=\frac{\sigma^2}{2}\left| \frac{b(i\omega)}{a(i\omega)} \right|^2 I_{2\times 2}$.
If $\omega\neq 0$, then 
\begin{align*}
d-\lim_{T\to \infty}\[\begin{matrix}
\Re(\mathcal{T}_TY(\omega)) \\ \Im(\mathcal{T}_TY(\omega))
\end{matrix}\] &= \mathcal{N}\left(0,\Sigma\right),\\
d-\lim_{T\to \infty}\left( \Re(\mathcal{T}_TY(\omega))^2+\Im(\mathcal{T}_TY(\omega))^2  \right)&=\operatorname{Exp} \left(\sigma^2\left| \frac{b(i\omega)}{a(i\omega)} \right|^2\right).
\end{align*}
If $\omega=0$, then $$ d- \lim_{T\to \infty} \mathcal{T}_TY(0) = \mathcal{N}\left(0,\left(\frac{b(0)}{a(0)}\right)^2\sigma^2\right)$$
and $$ d-\lim_{T\to \infty}  \frac{1}{\sigma^2}\left|\frac{a(0) \mathcal{T}_TY(0)}{b(0) }\right|^2 \sim \chi^2(1).$$
\end{theorem}
Clearly, an analogous statement using Theorem \ref{thm:LimitTheoremSeveralFrequencies} holds for the joint distribution when the truncated Fourier transform is taken at different frequencies.

\textsc{Proof of Theorem \ref{thm:ApproxTFTDoubleLimitDistribution}}.
Put
$$ Z(T):=\frac{1}{\sqrt{T}}\frac{b(i\omega)}{a(i\omega)}\int_0^T e^{-i\omega t}dL(t) $$
and consider the following two-dimensional random vectors:
$$ \mathbf{Z}_n :=\[\begin{matrix}
\Re(Z(T)) \\ \Im(Z(T))
\end{matrix}\], \quad  \mathbf{U}_n :=\[\begin{matrix}
\Re(\mathcal{F}_TY(\omega)) \\ \Im(\mathcal{F}_TY(\omega))
\end{matrix}\],\quad \mathbf{V}_n :=\[\begin{matrix}
\Re(\mathcal{T}_TY(\omega)) \\ \Im(\mathcal{T}_TY(\omega))
\end{matrix}\].$$
%Here we understand $n=\lfloor T \rfloor$
Observe that it is enough to consider the above limits for $T=n$.
By Lemma \ref{lem:LimitInProbabilityFirstSummandFT} we know that 
$$\P-\lim_{n\to \infty} \|\mathbf{U}_n-\mathbf{Z}_n\|
=\mathbf{0}.$$
From Theorem \ref{lem:RozkladZ(T)} we get
$$d-\lim_{n\to \infty} \mathbf{Z}_n= \mathcal{N}(0,\Sigma).$$
Therefore 
$$d-\lim_{n\to \infty} \mathbf{U}_n= \mathcal{N}(0,\Sigma).$$
By Theorem \ref{thm:ApproxTFTDoubleLimit} we have
$$\P-\lim_{n\to \infty} (\mathbf{V}_n-\mathbf{U}_n)= 0.$$
Therefore, 
$$d-\lim_{n\to \infty} \mathbf{V}_n= \mathcal{N}(0,\Sigma).$$
In the same way we obtain
$$ d-\lim_{n\to \infty} |\mathbf{Z}|^2=\operatorname{Exp} \left(\sigma^2\left| \frac{b(i\omega)}{a(i\omega)} \right|^2\right). $$
In order to obtain the assertion for $\omega=0$ we repeat the above resonings applying Theorem \ref{lem:RozkladZ(T)OmegaRowna0} instead of Theorem \ref{lem:LimitInProbabilityFirstSummandFT}. $\Box$
%-------------------------------------------------------------------------------------
%-------------------------------------------------------------------------------------
%-------------------------------------------------------------------------------------
%-------------------------Simulations-------------------------------------------------
%-------------------------------------------------------------------------------------
%-------------------------------------------------------------------------------------
%-------------------------------------------------------------------------------------
\section{Illustrative simulations}\label{sec4}
We now turn to a numerical illustration of the theoretical convergence results given in Section \ref{Section:ApproximationofIntegrals}. We are looking at simulations of CARMA processes and their numerically approximated truncated Fourier transform over different time horizons and maximal grid widths. To illustrate the convergence to the asymptotic normal distribution we shall look at several frequencies and different driving L\'evy processes, standard Brownian motion, a Variance Gamma process and a ``two sided Poisson process''. Of course,  the truncated Fourier transform of $(Y_t)$ is obtained using the trapezoidal rule based on non-equidistant observations of the CARMA process
$(Y_t)$ given on the interval $[0,T]$. On an interval $[0,T]$ we generate a  
non-equidistant grid in the following way: we fix the maximal distance $h_{\max}(T)$ between elements of the grid and from each interval 
$\[i\cdot \frac{1}{2} h_{\max}(T), (i+1)\cdot \frac{1}{2} h_{\max}(T)\right)$ for $i=0,1,\dots, N-1$ we draw a number according to the uniform distribution. 
This results in a non-equidistant grid with the number of points being $N(T)=2T/h_{\max} +1$.

For our simulations we used the R Project for Statistical Computing. For the simulation we first generate the non-equidistant grid by the above procedure and then join it with a regular grid of mesh 0.001, which is still on average five times finer than the non-equidistant grid of the largest time horizon considered. On this joint grid the CARMA process $Y$ is simulated with a standard Euler scheme for the state space representation. Afterwards only the simulated values at the times of the original non-equidistant grid are used to compute the approximation of the truncated Fourier transform with the trapezoidal rule.
In all cases we simulate 2000 independent paths of the CARMA process and compute the associated values of the truncated Fourier transform at the following frequencies:
$$ [\omega_1,\omega_2,\omega_3,\omega_4] =[0,\, 0.1,\, 1,\,10].$$
 For the non-zero frequencies real and imaginary part have to be considered separately. However, in the following we look only at the real parts as the behaviour of the imaginary parts is most similar.  Mainly, the results are presented via QQ-plots where the theoretical values follow the (limiting) law described in Theorem \ref{thm:ApproxTFTDoubleLimitDistribution}.

We are going to consider CARMA processes with the following autoregressive and moving average orders: $(p,q)=(1,0)$, i.e. an Ornstein-Uhlenbeck type process, and $(p,q)=(2,1)$. For the time horizon $T$ and the maximum distance of the non-equidistant observation times we consider the pairs  $(T=10,\,h_{\max}=0.1)$, $(T=50,\,h_{\max}=0.05)$ and $(T=100,\,h_{\max}=0.01)$. 

For each case we consider three different driving  L\'evy noises: standard Brownian Motion, a Variance Gamma process and a ``two sided Poisson process''.
For the definition and properties of the Variance Gamma process we refer to \cite{MadanCarrChang1998} and references therein. 
We construct the process in the following way: $V_t=G_t^1-G_t^2$, where $G_t^1$ and $G_t^2$ are independent Gamma processes with shape parameter $1$ and scale parameter $4$. Likewise the  
``two sided Poisson process'' is the difference of two independent Poisson processes with rate $10$, i.e. a compound Poisson process with rate $20$ and jumps $+1$ and $-1$ both with probability $1/2$.

\begin{figure}[tp]    
  \begin{subfigure}[b]{0.32\linewidth}
    \centering
    \includegraphics[width=\linewidth]{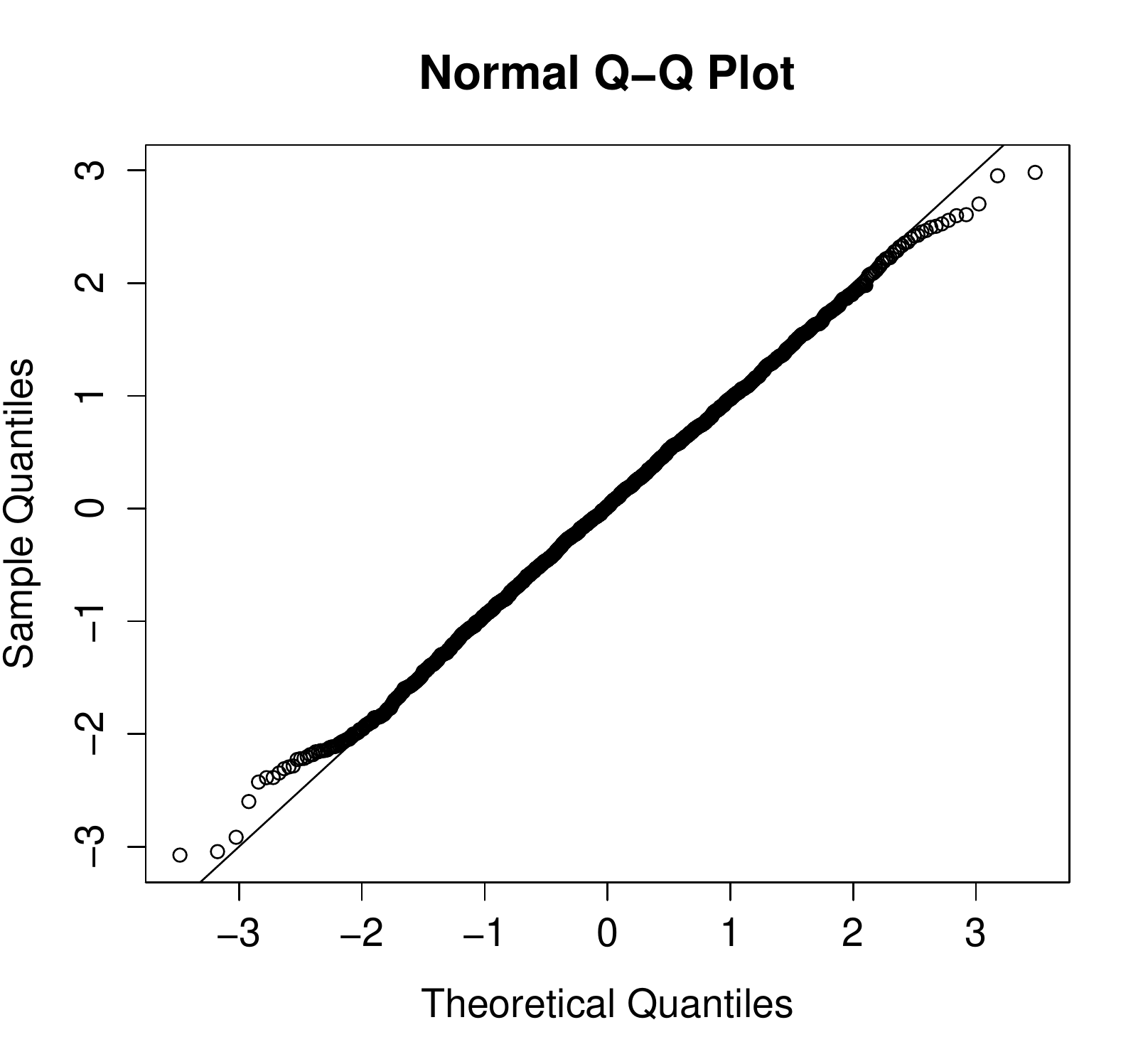} 
    \caption{ $\omega=0, T=10$} 
  
    %\vspace{4ex}
  \end{subfigure}%% 
 \begin{subfigure}[b]{0.32\linewidth}
    \centering
    \includegraphics[width=\linewidth]{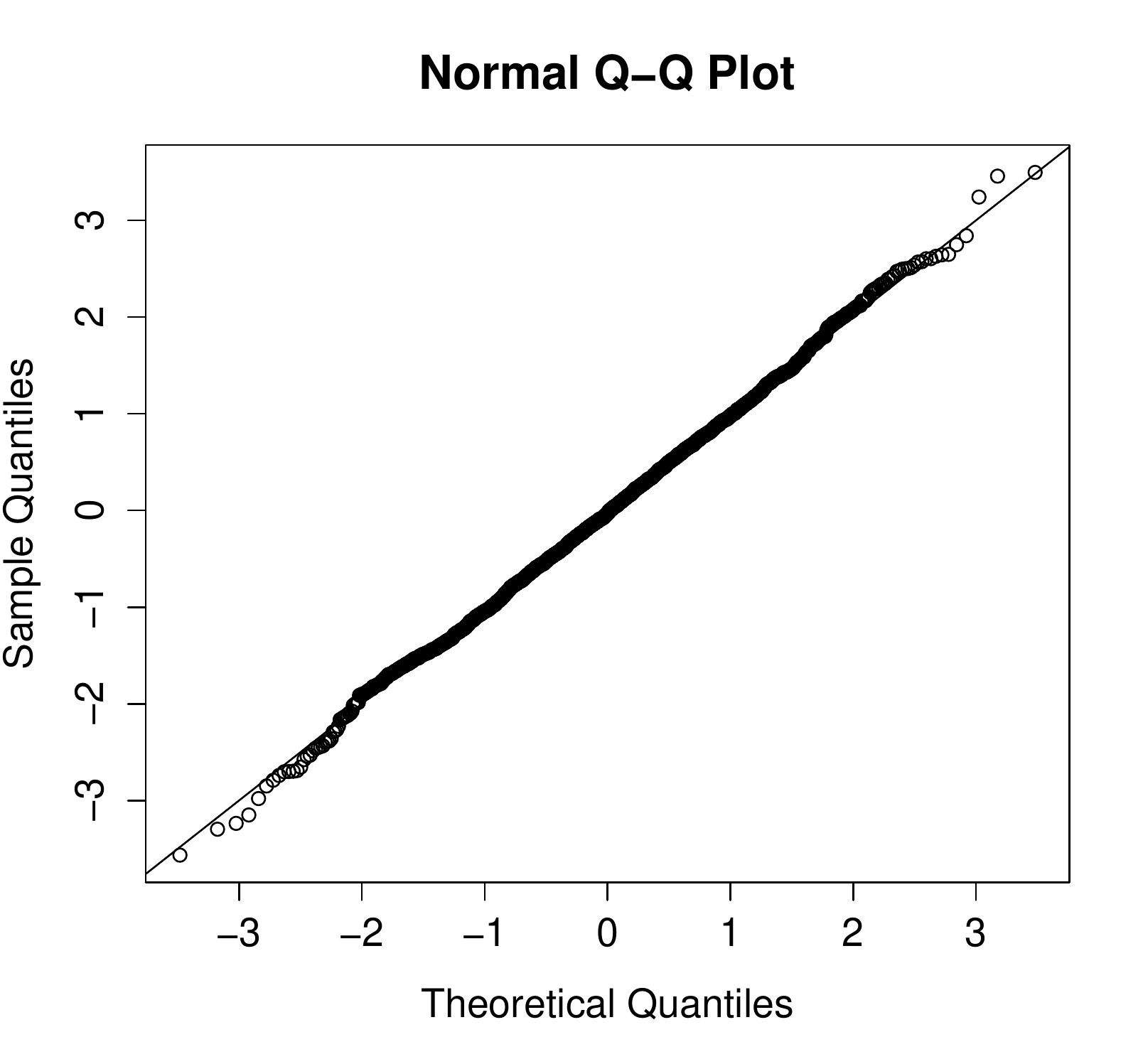} 
    \caption{$\omega=0, T=50$} 
    %\vspace{4ex}
  \end{subfigure}
 \begin{subfigure}[b]{0.32\linewidth}
    \centering
    \includegraphics[width=\linewidth]{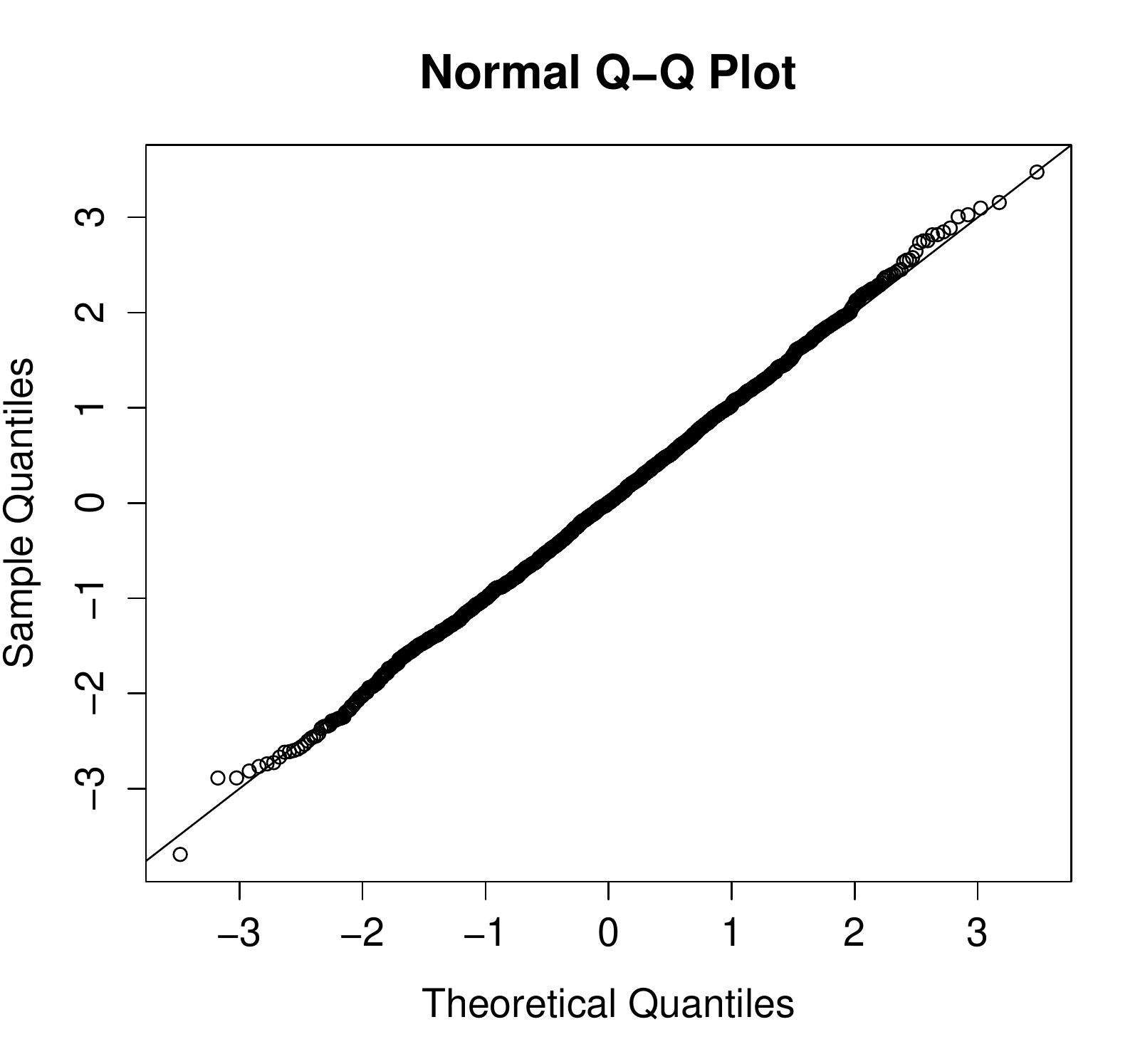} 
    \caption{$\omega=0, T=100$} 
    %\vspace{4ex}
  \end{subfigure}
  \begin{subfigure}[b]{0.32\linewidth}
    \centering
    \includegraphics[width=\linewidth]{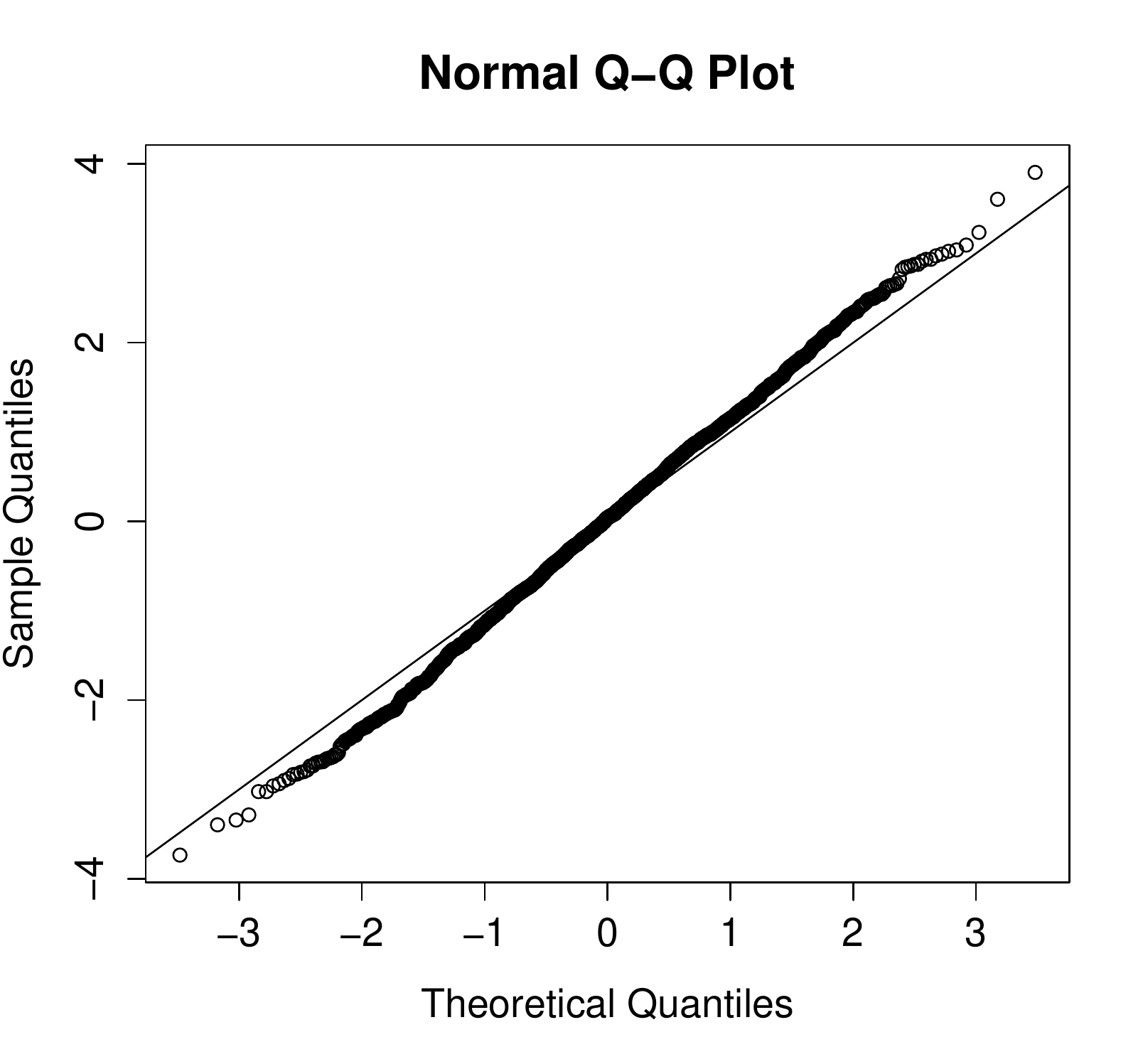} 
    \caption{$\omega=0.1, T=10$} 
    %\vspace{4ex}
  \end{subfigure}%% 
 \begin{subfigure}[b]{0.32\linewidth}
    \centering
    \includegraphics[width=\linewidth]{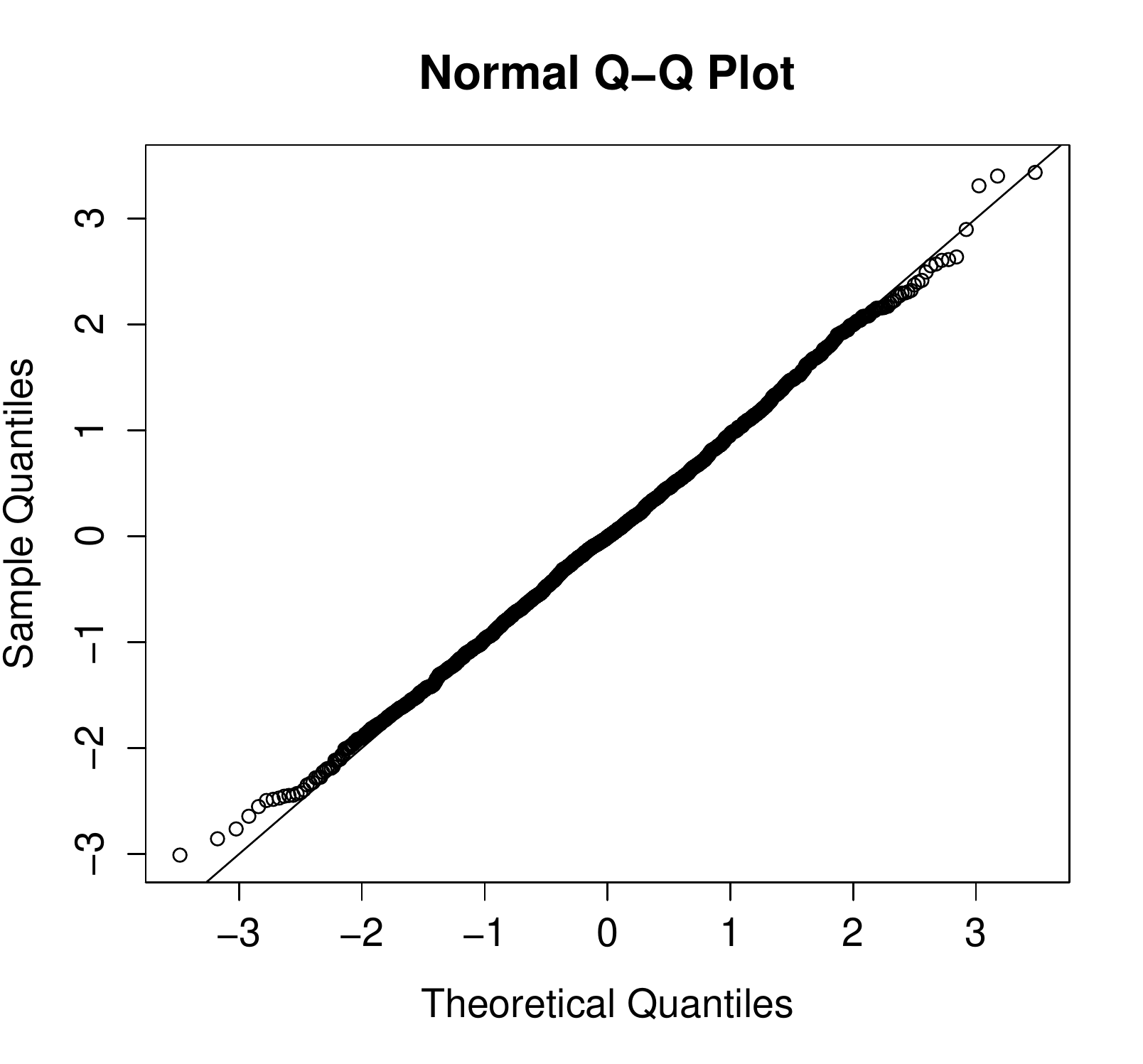} 
    \caption{$\omega=0.1, T=50$} 
    %\vspace{4ex}
  \end{subfigure}
 \begin{subfigure}[b]{0.32\linewidth}
    \centering
    \includegraphics[width=\linewidth]{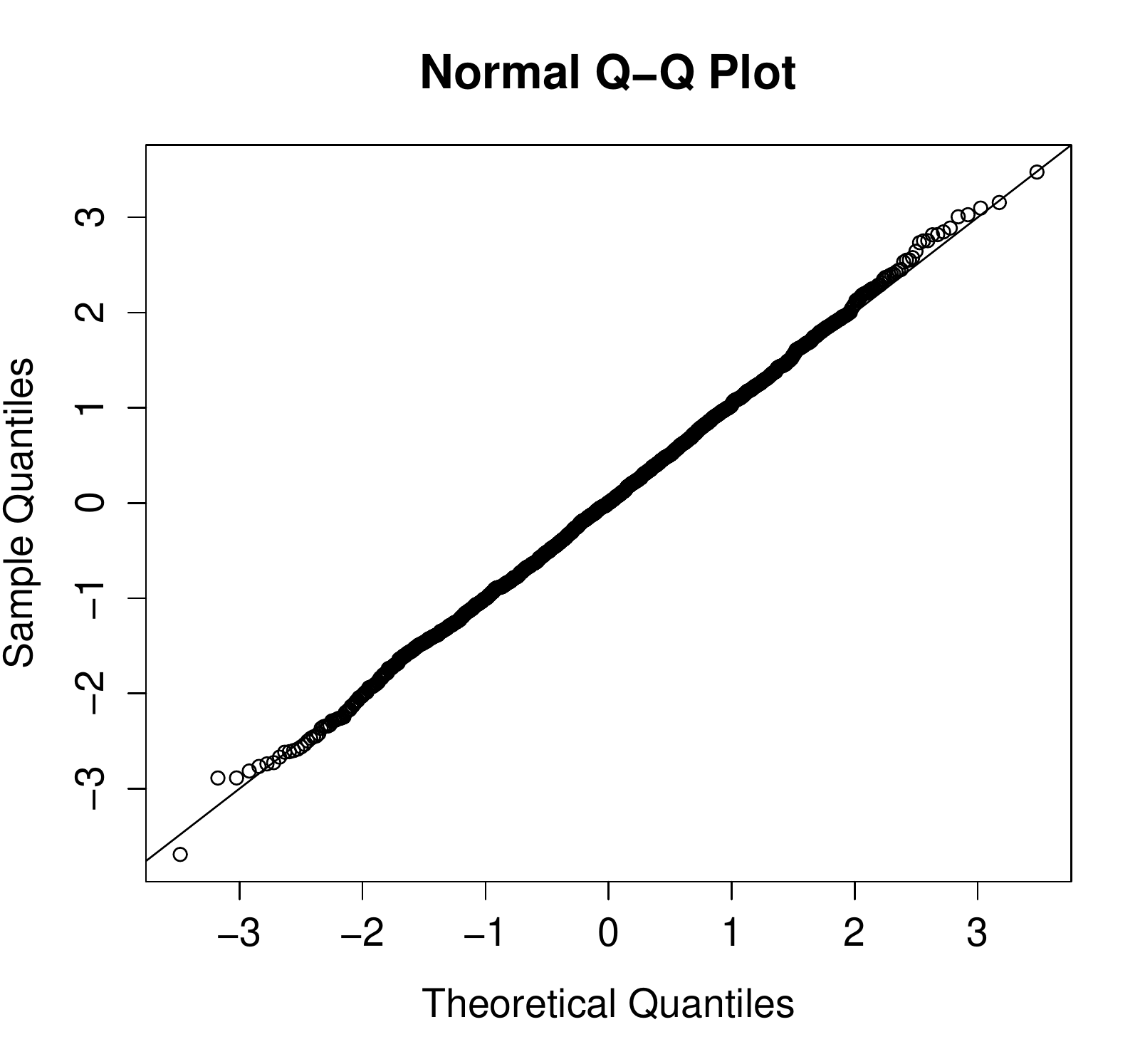} 
    \caption{$\omega=0.1, T=100$} 
    %\vspace{4ex}
  \end{subfigure}
  \begin{subfigure}[b]{0.32\linewidth}
    \centering
    \includegraphics[width=\linewidth]{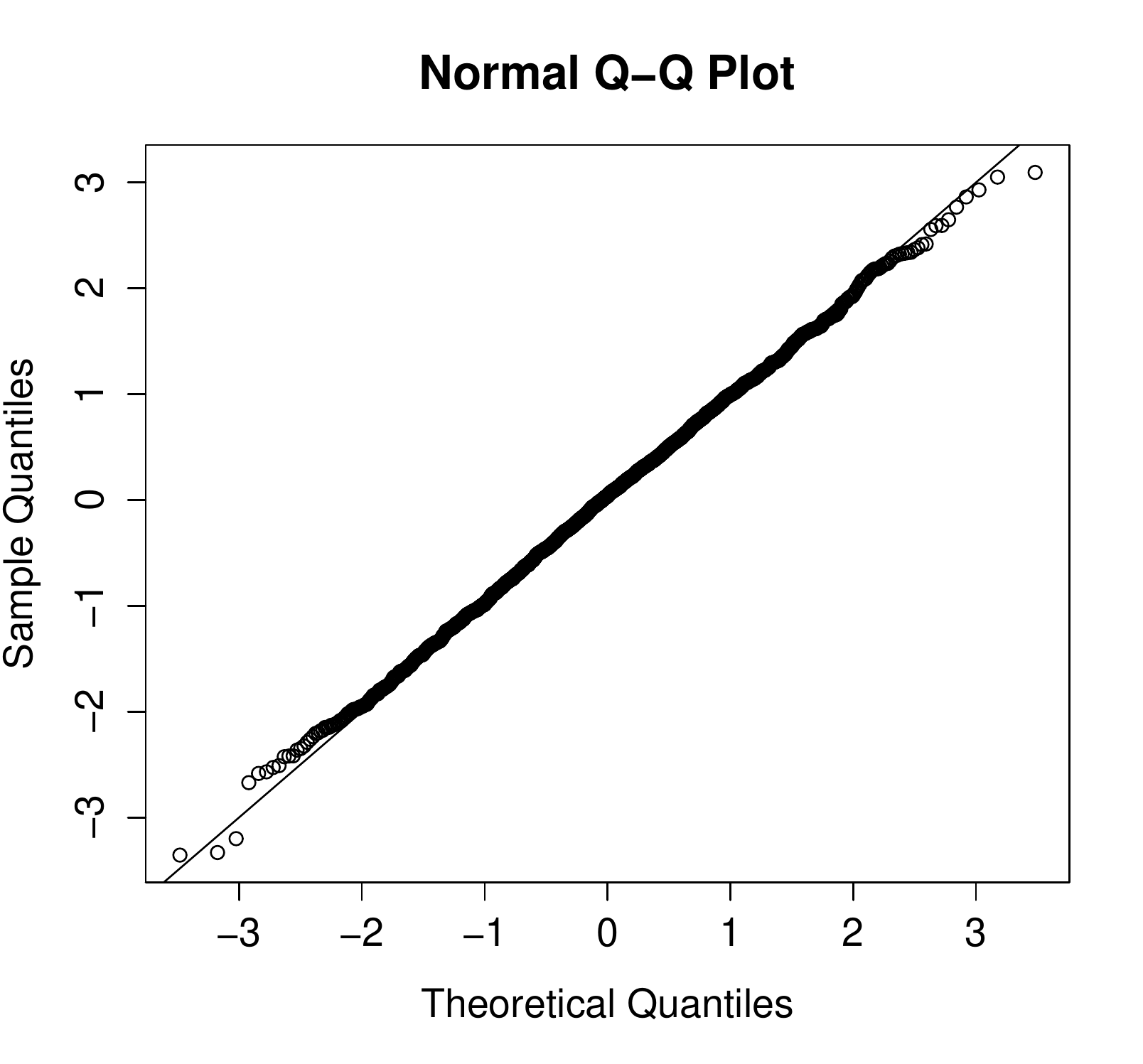} 
    \caption{$\omega=1, T=10$} 
    %\vspace{4ex}
  \end{subfigure}%% 
 \begin{subfigure}[b]{0.32\linewidth}
    \centering
    \includegraphics[width=\linewidth]{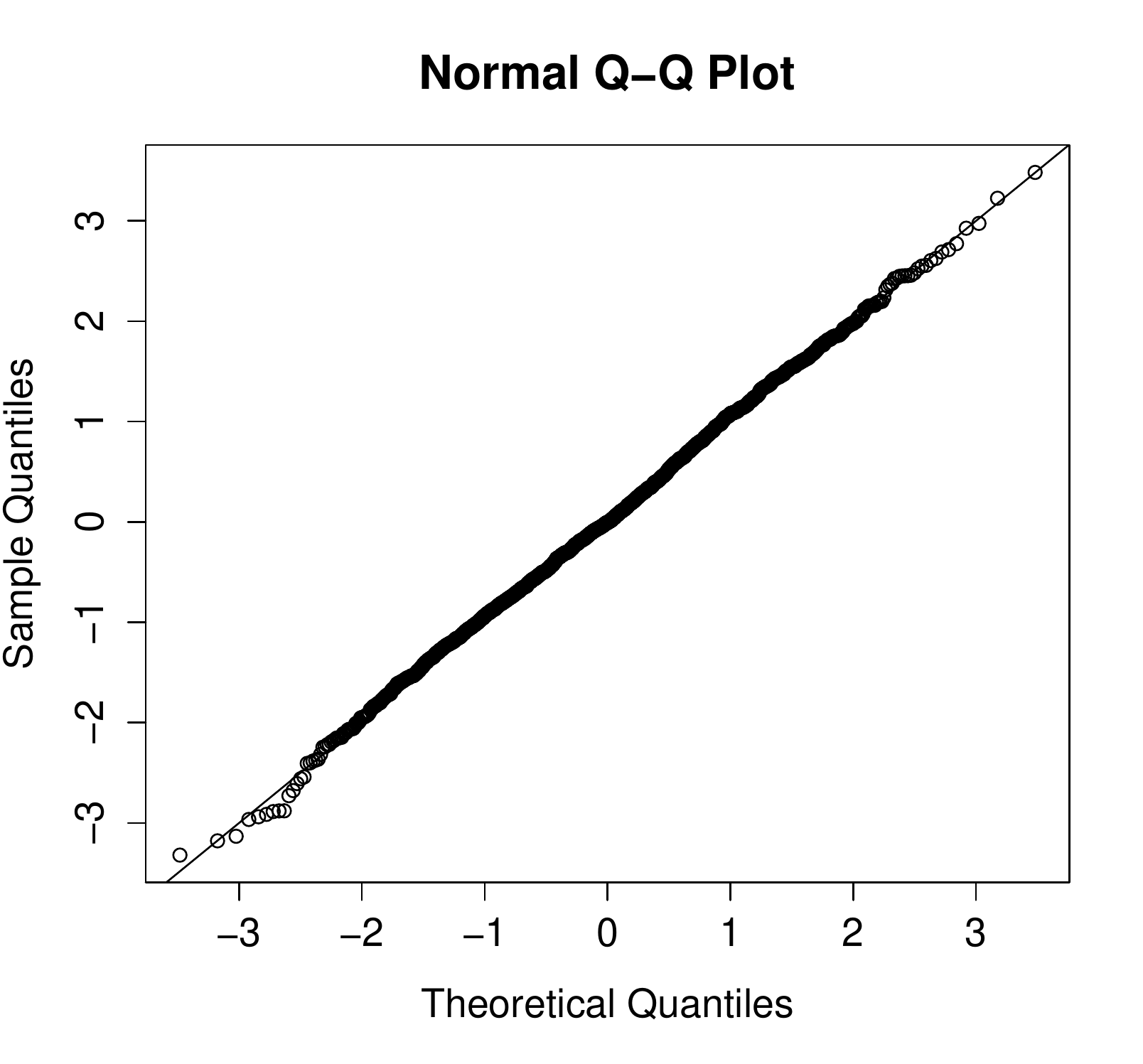} 
    \caption{$\omega=1, T=50$} 
    %\vspace{4ex}
  \end{subfigure}
 \begin{subfigure}[b]{0.32\linewidth}
    \centering
    \includegraphics[width=\linewidth]{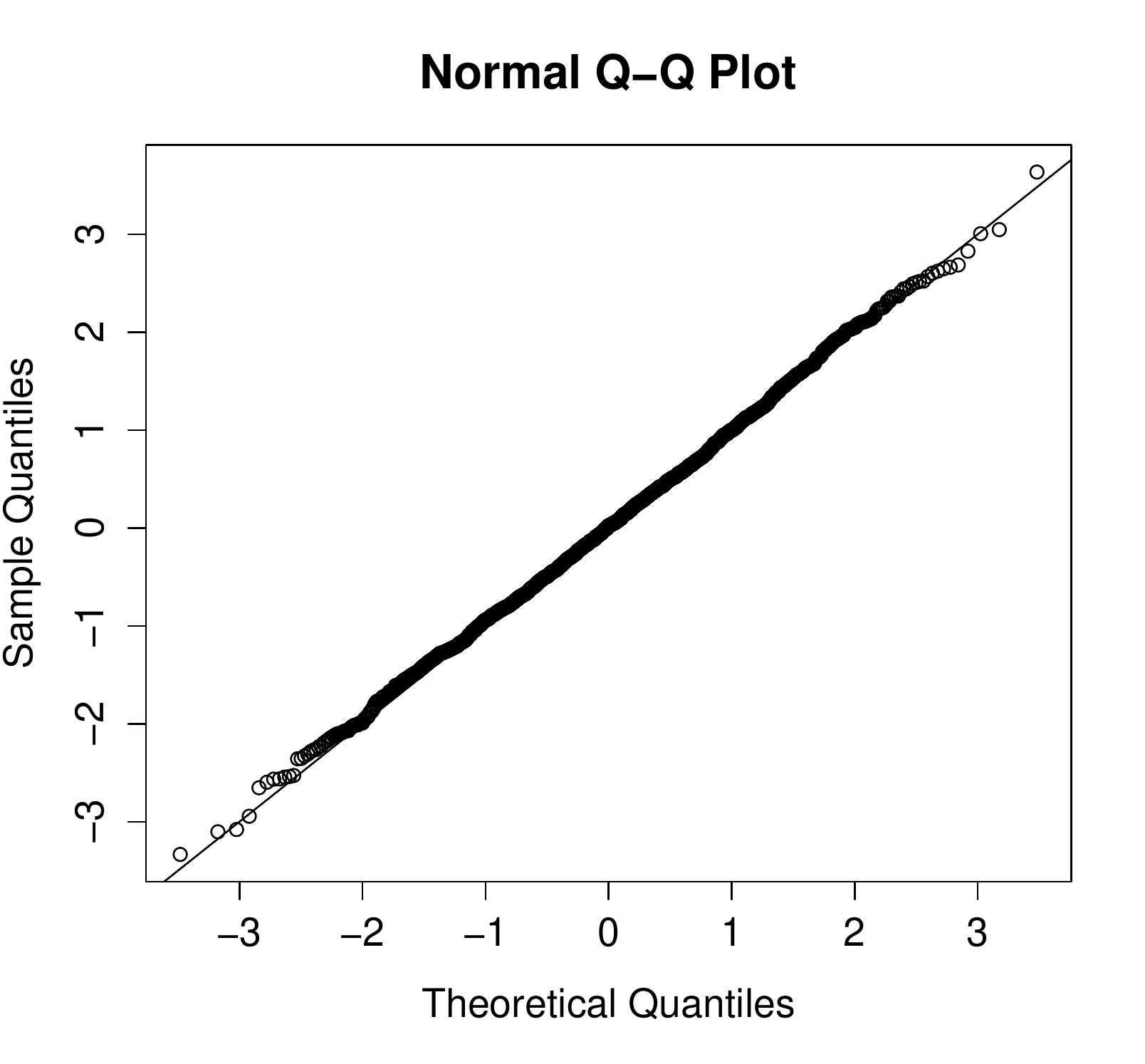} 
    \caption{$\omega=1, T=100$} 
    %\vspace{4ex}
  \end{subfigure}
  \begin{subfigure}[b]{0.32\linewidth}
    \centering
    \includegraphics[width=\linewidth]{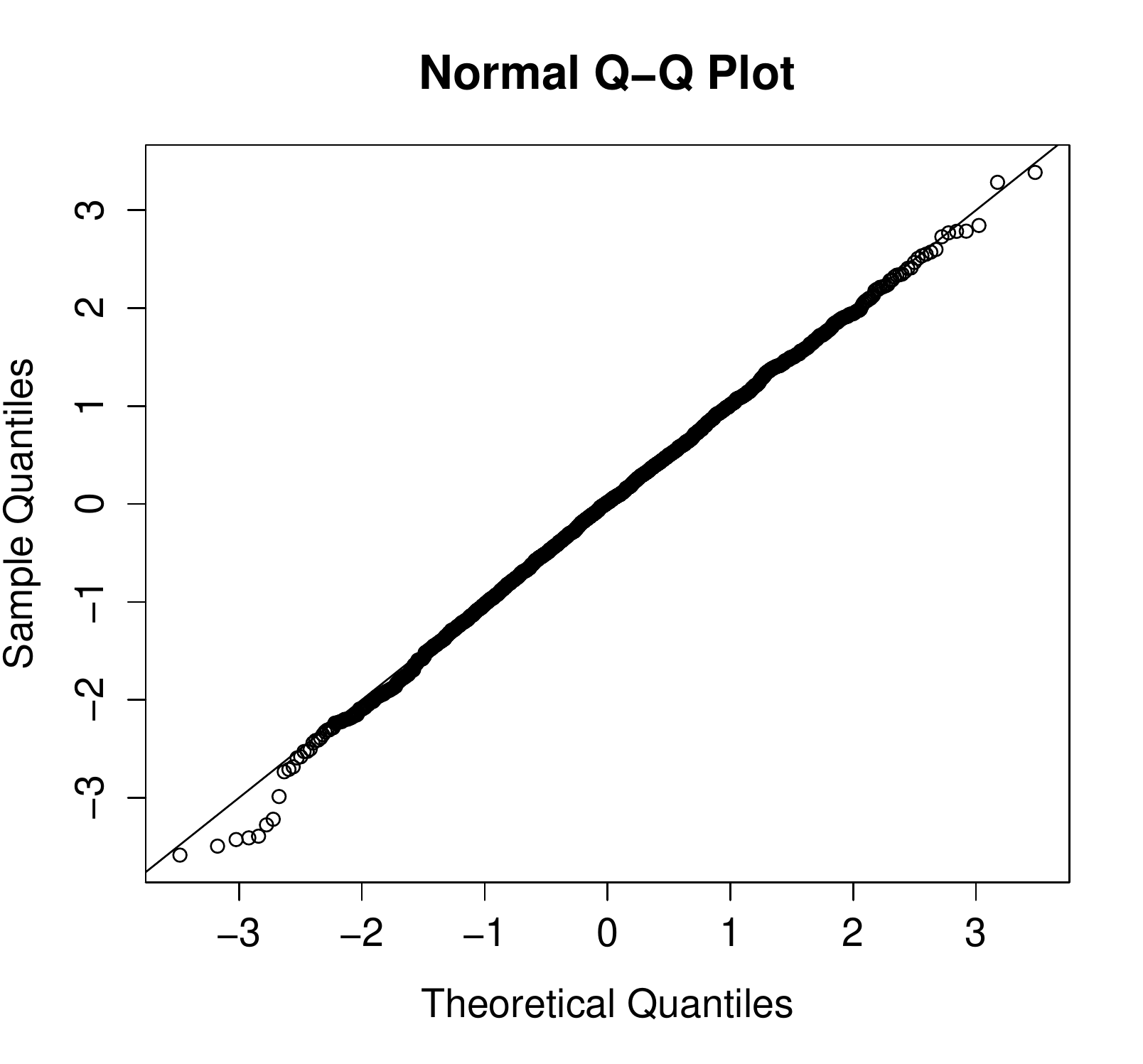} 
    \caption{$\omega=10, T=10$}  
    %\vspace{4ex}
  \end{subfigure}%% 
 \begin{subfigure}[b]{0.32\linewidth}
    \centering
    \includegraphics[width=\linewidth]{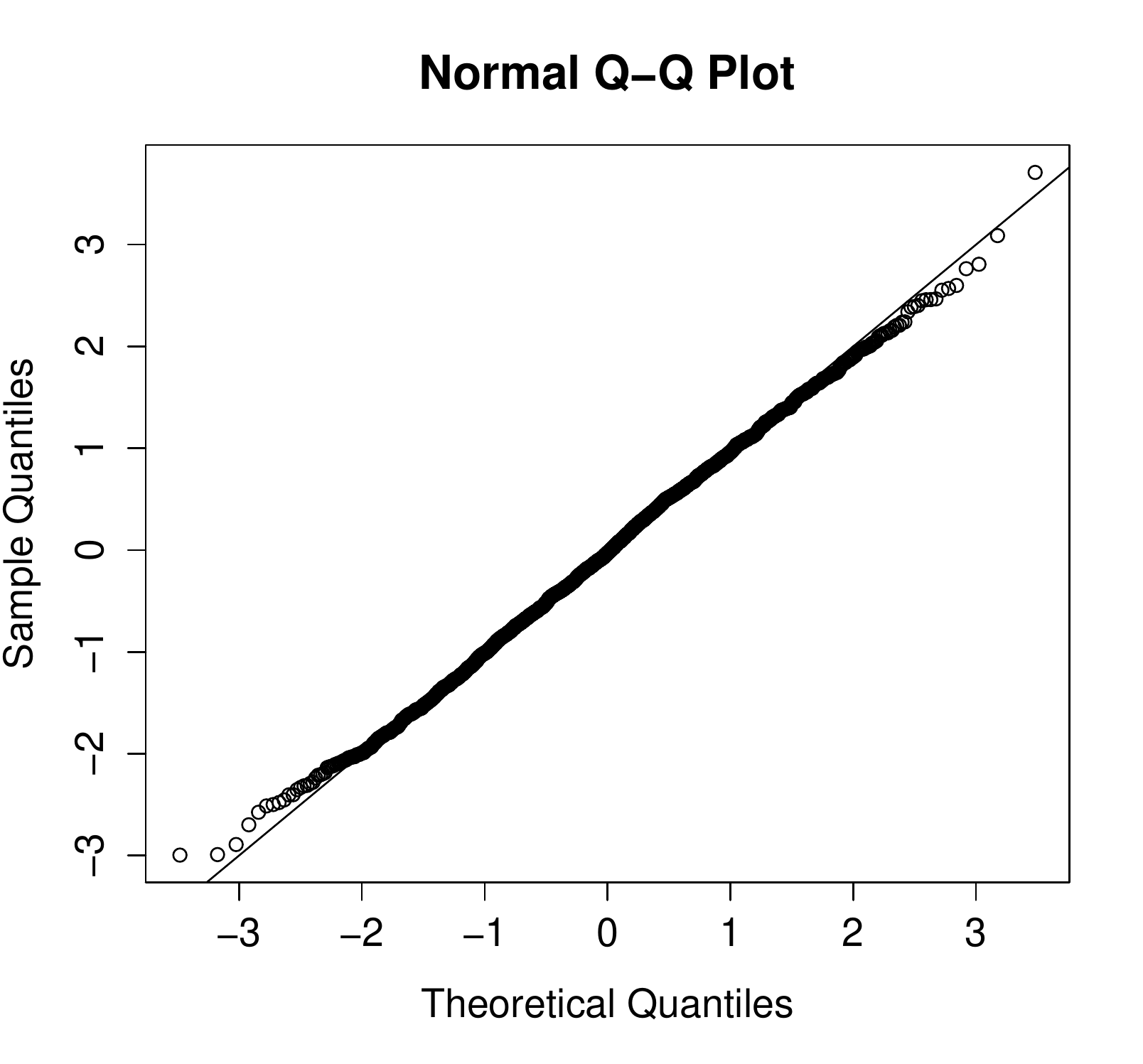} 
    \caption{$\omega=10, T=50$} 
    %\vspace{4ex}
  \end{subfigure}
 \begin{subfigure}[b]{0.32\linewidth}
    \centering
    \includegraphics[width=\linewidth]{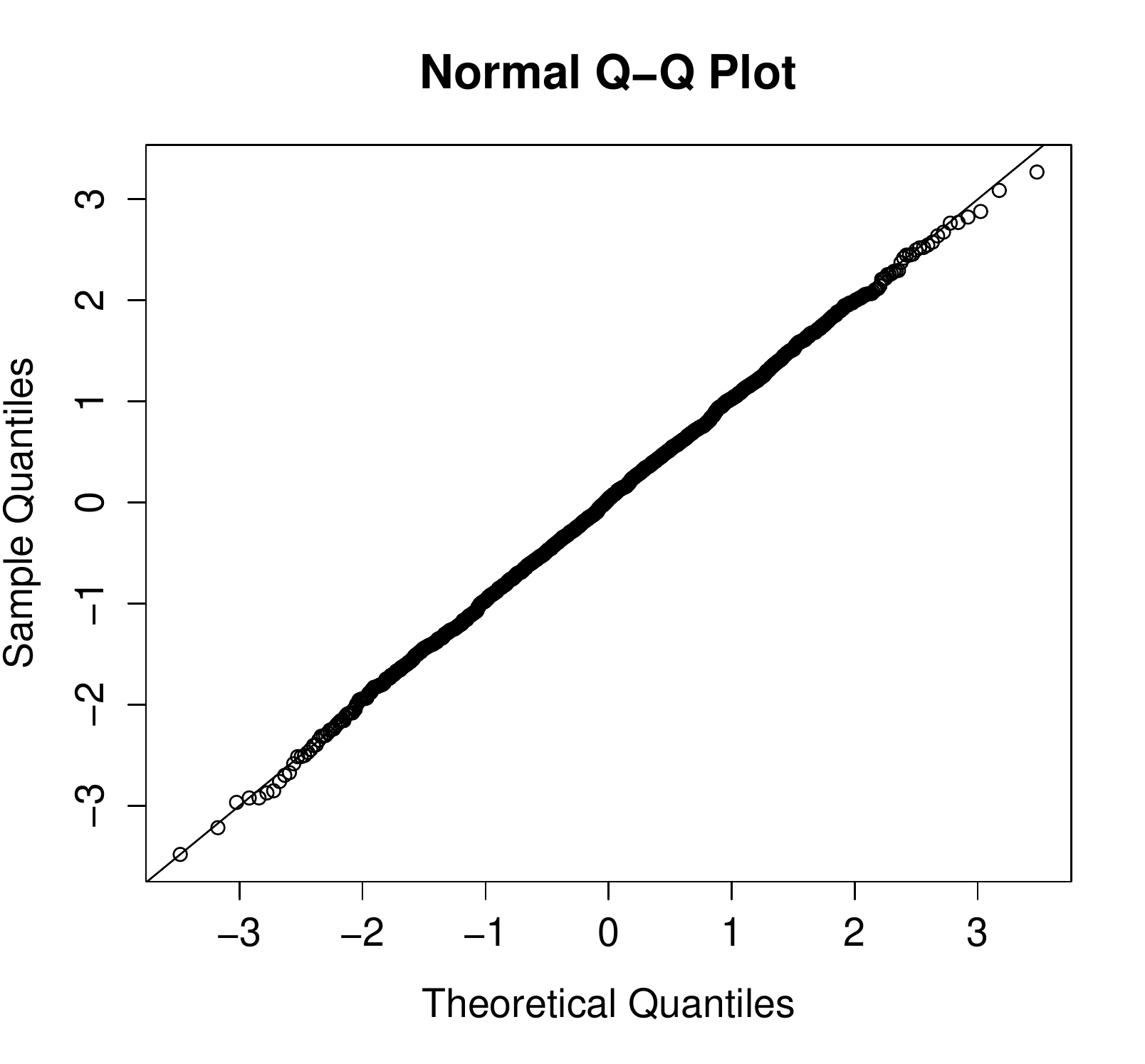} 
    \caption{$\omega=10, T=100$} 
    %\vspace{4ex}
  \end{subfigure}
  \caption{Normal QQ plots for the real part of the truncated Fourier transform of the Ornstein-Uhlenbeck type process driven by standard Brownian Motion for the frequencies $0, 0.1, 1 , 10$ (rows) and time horizons/maximum non-equidistant grid sizes $10/0.1, 50/0.05, 100/0.01$ (columns). The theoretical quantiles are coming from the (limiting) law described in Theorem \ref{thm:ApproxTFTDoubleLimitDistribution}. }\label{plot:QQCARNormal} 
\end{figure}

\begin{figure}[tp]    
  \begin{subfigure}[b]{0.32\linewidth}
    \centering
    \includegraphics[width=\linewidth]{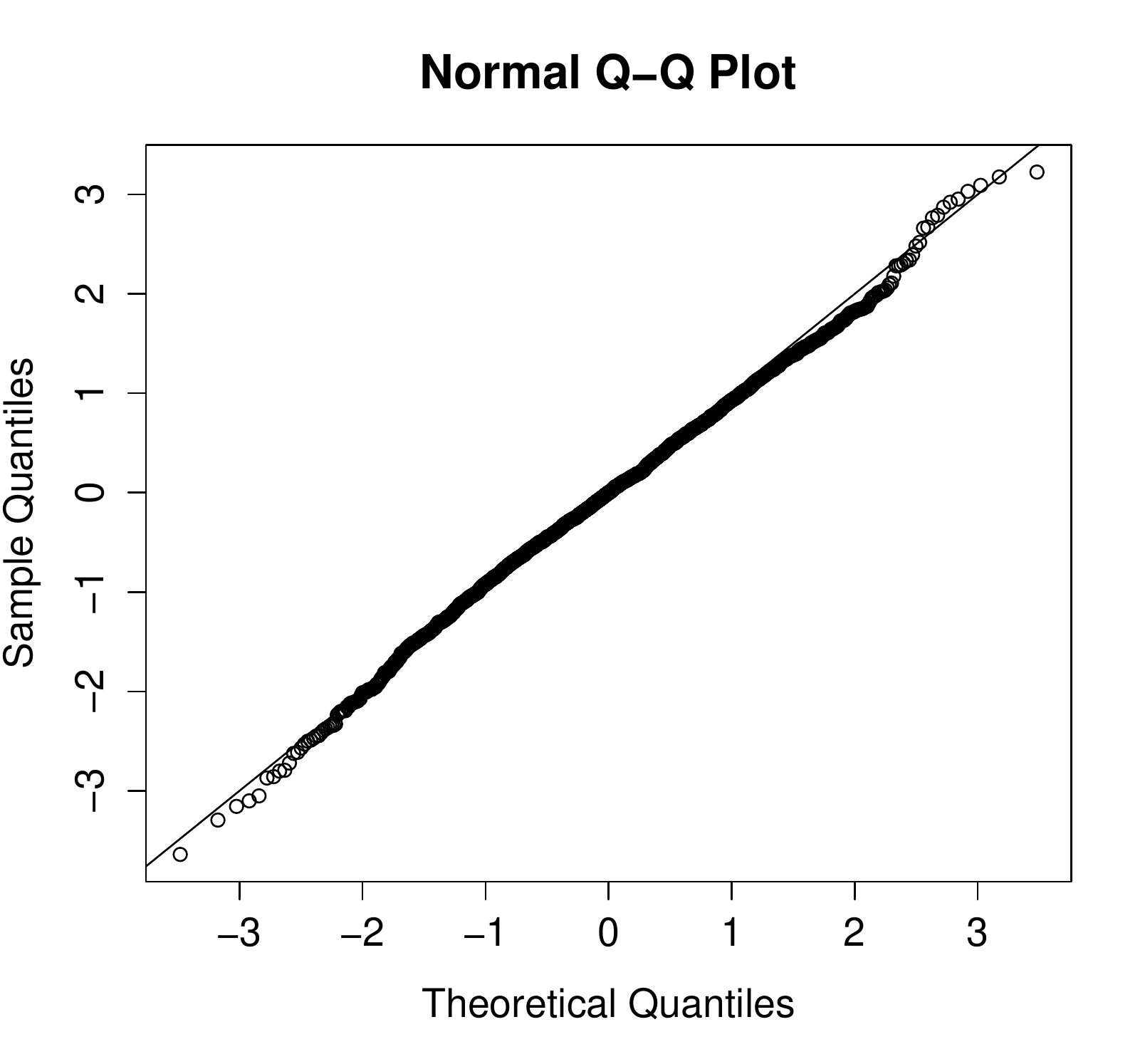} 
    \caption{ $\omega=0, T=10$} 
  
    %\vspace{4ex}
  \end{subfigure}%% 
 \begin{subfigure}[b]{0.32\linewidth}
    \centering
    \includegraphics[width=\linewidth]{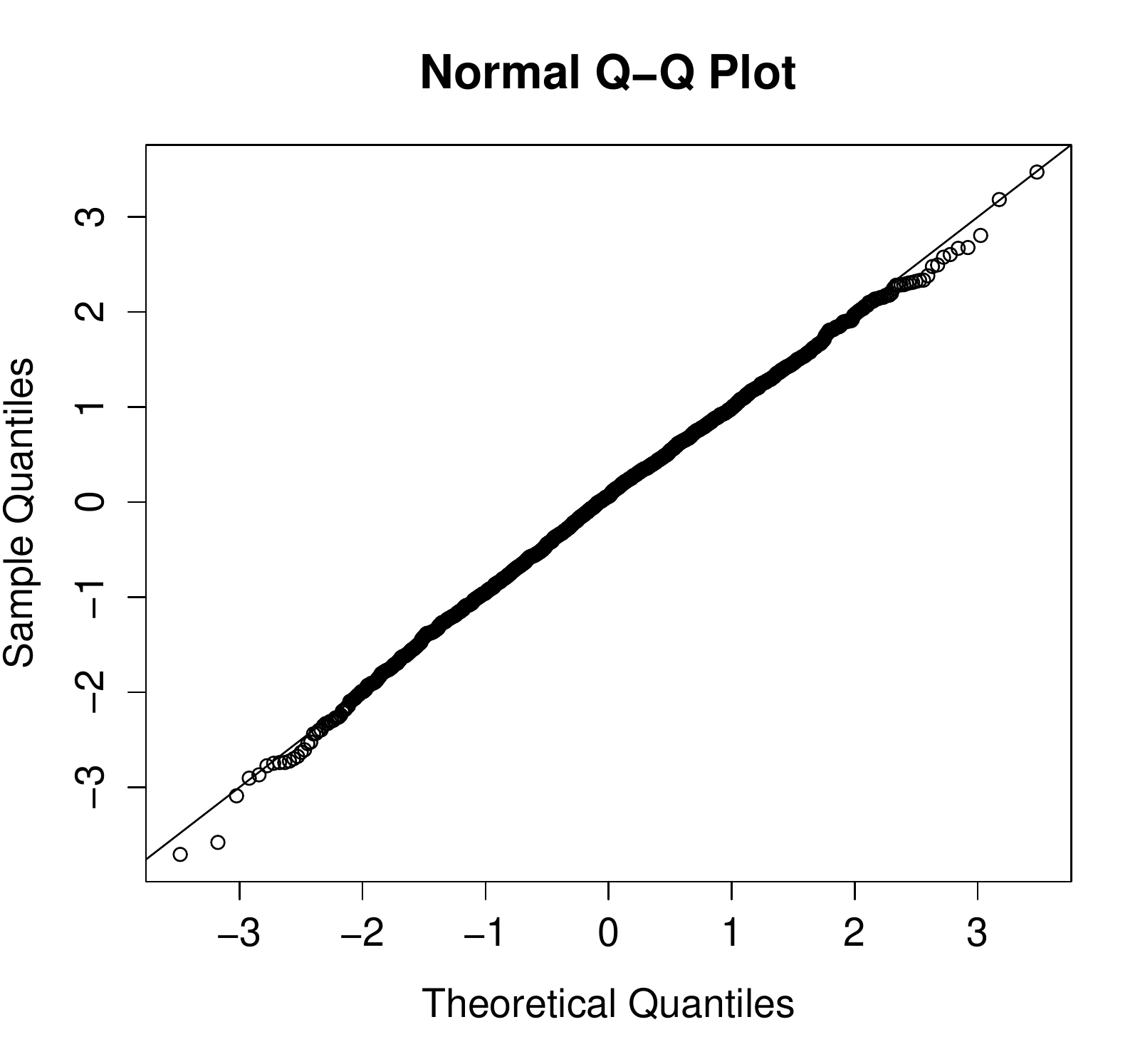} 
    \caption{$\omega=0, T=50$} 
    %\vspace{4ex}
  \end{subfigure}
 \begin{subfigure}[b]{0.32\linewidth}
    \centering
    \includegraphics[width=\linewidth]{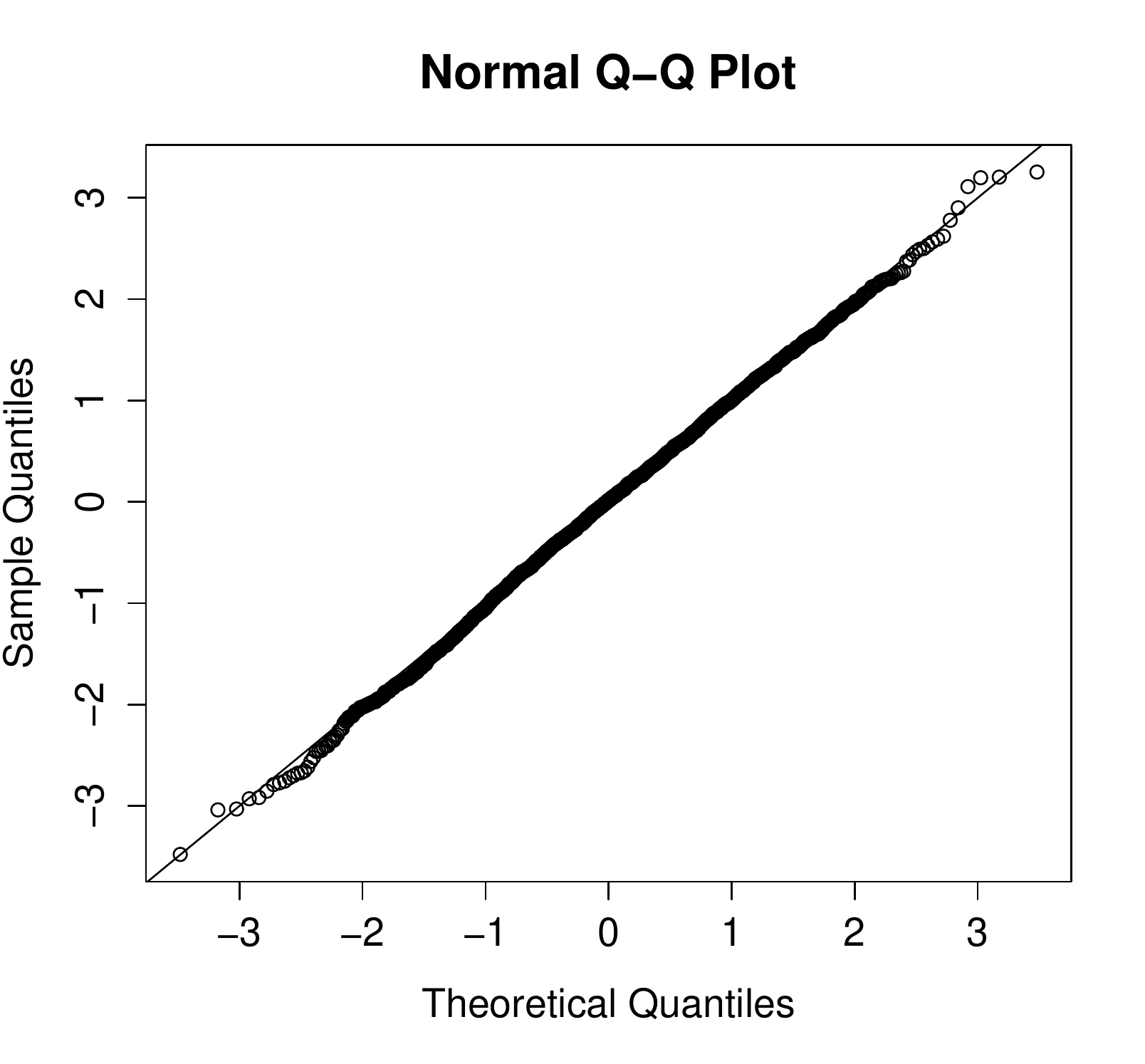} 
    \caption{$\omega=0, T=100$} 
    %\vspace{4ex}
  \end{subfigure}
  \begin{subfigure}[b]{0.32\linewidth}
    \centering
    \includegraphics[width=\linewidth]{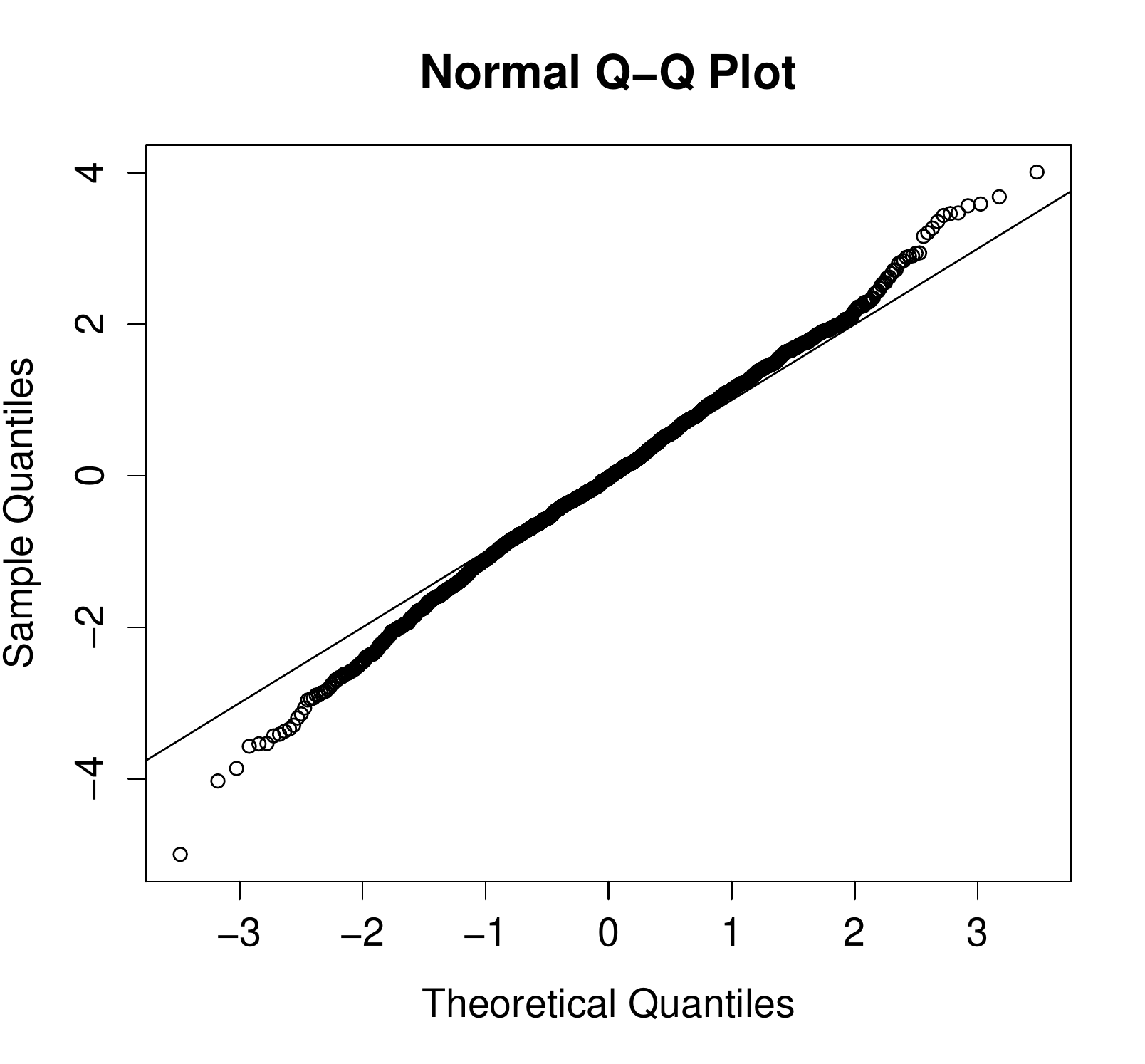} 
    \caption{$\omega=0.1, T=10$} 
    %\vspace{4ex}
  \end{subfigure}%% 
 \begin{subfigure}[b]{0.32\linewidth}
    \centering
    \includegraphics[width=\linewidth]{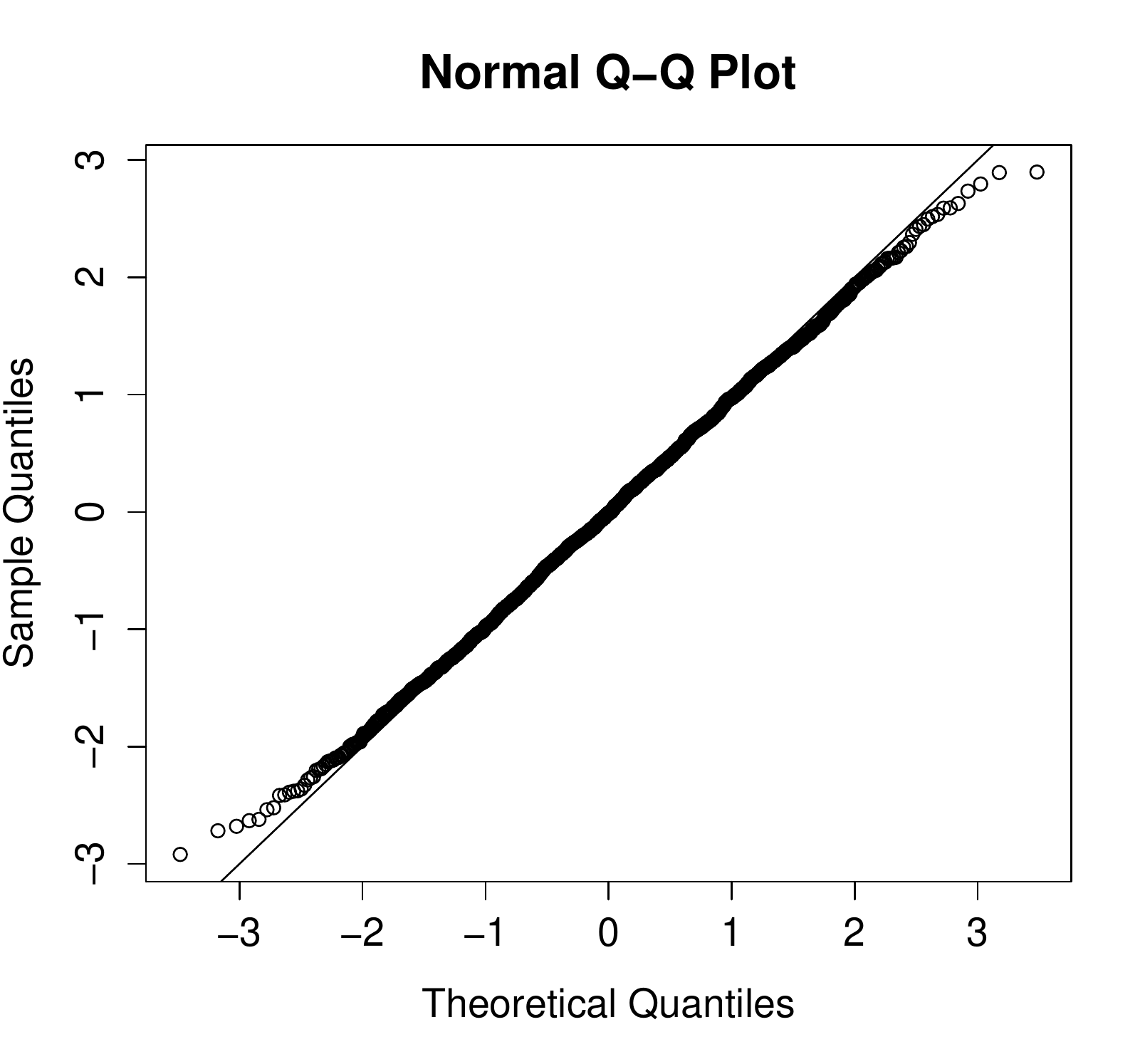} 
    \caption{$\omega=0.1, T=50$} 
    %\vspace{4ex}
  \end{subfigure}
 \begin{subfigure}[b]{0.32\linewidth}
    \centering
    \includegraphics[width=\linewidth]{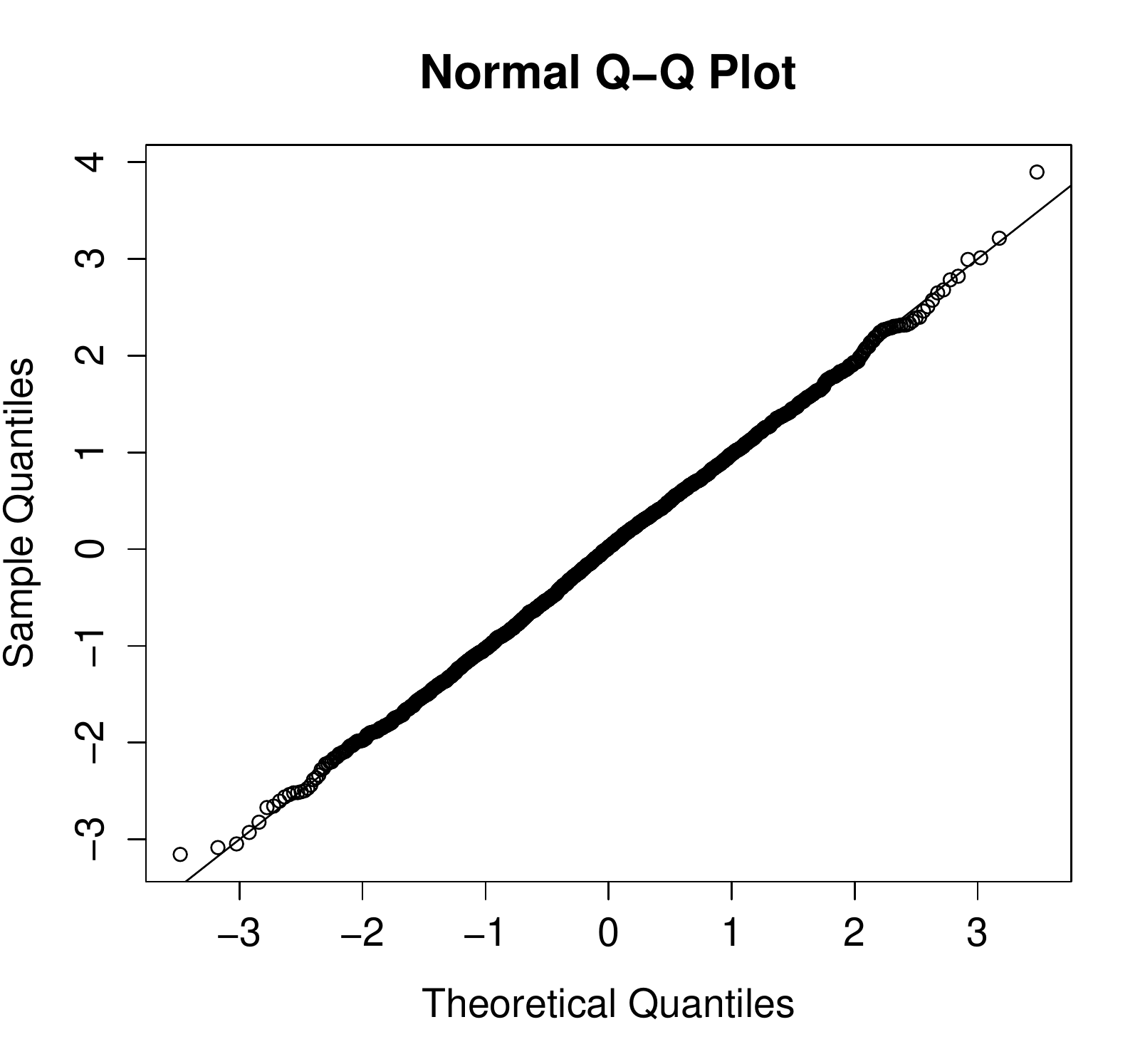} 
    \caption{$\omega=0.1, T=100$} 
    %\vspace{4ex}
  \end{subfigure}
  \begin{subfigure}[b]{0.32\linewidth}
    \centering
    \includegraphics[width=\linewidth]{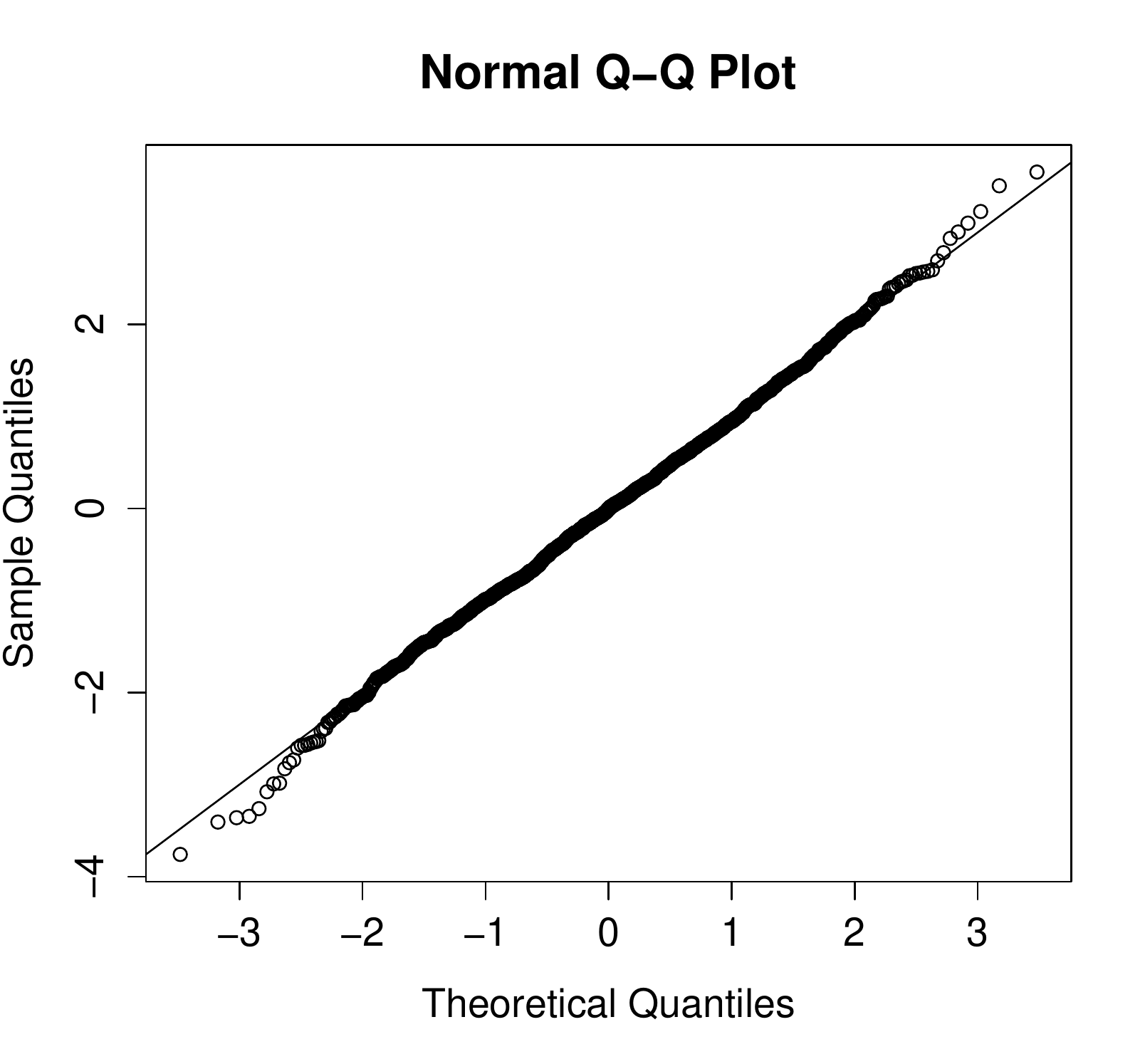} 
    \caption{$\omega=1, T=10$} 
    %\vspace{4ex}
  \end{subfigure}%% 
 \begin{subfigure}[b]{0.32\linewidth}
    \centering
    \includegraphics[width=\linewidth]{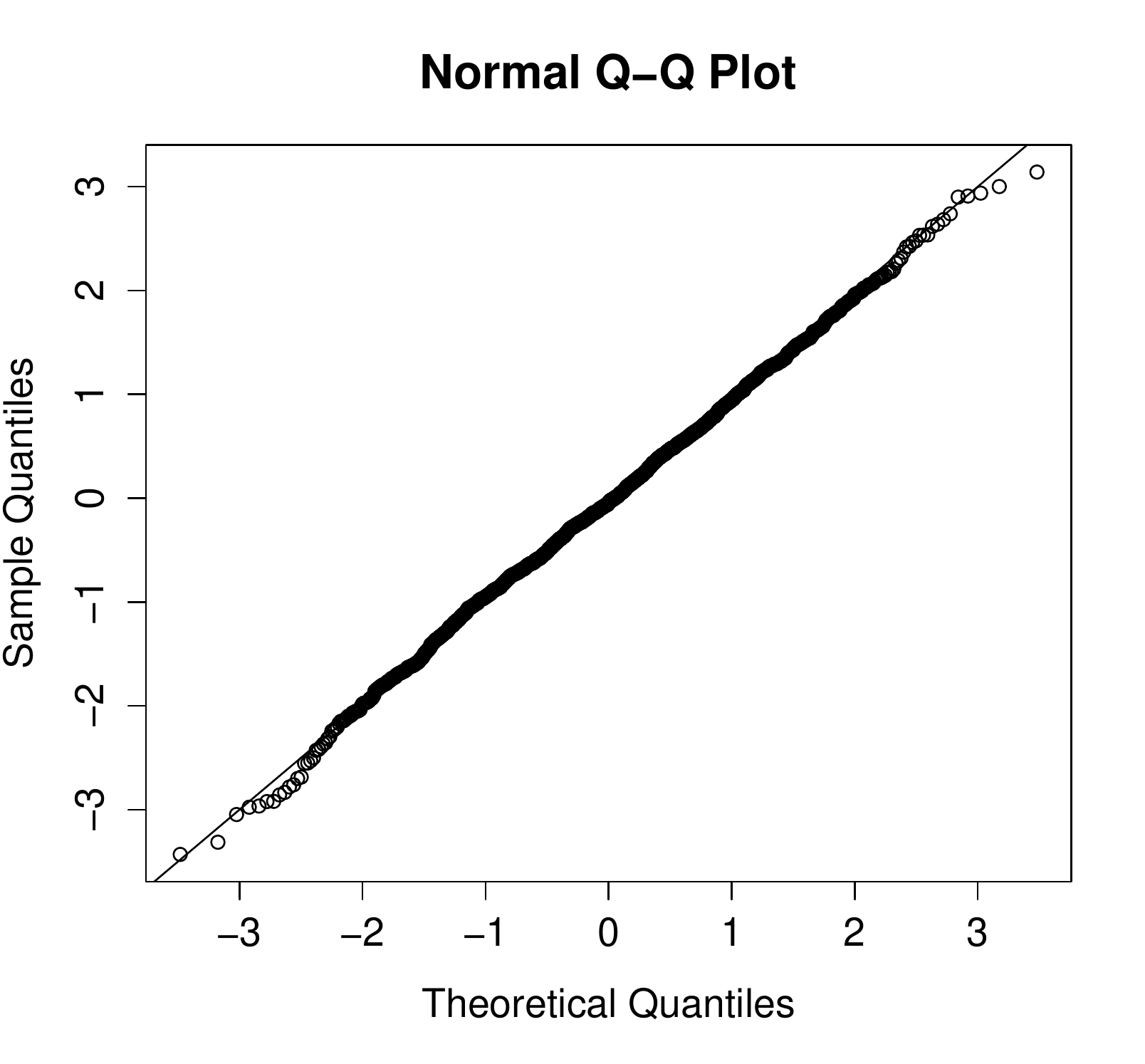} 
    \caption{$\omega=1, T=50$} 
    %\vspace{4ex}
  \end{subfigure}
 \begin{subfigure}[b]{0.32\linewidth}
    \centering
    \includegraphics[width=\linewidth]{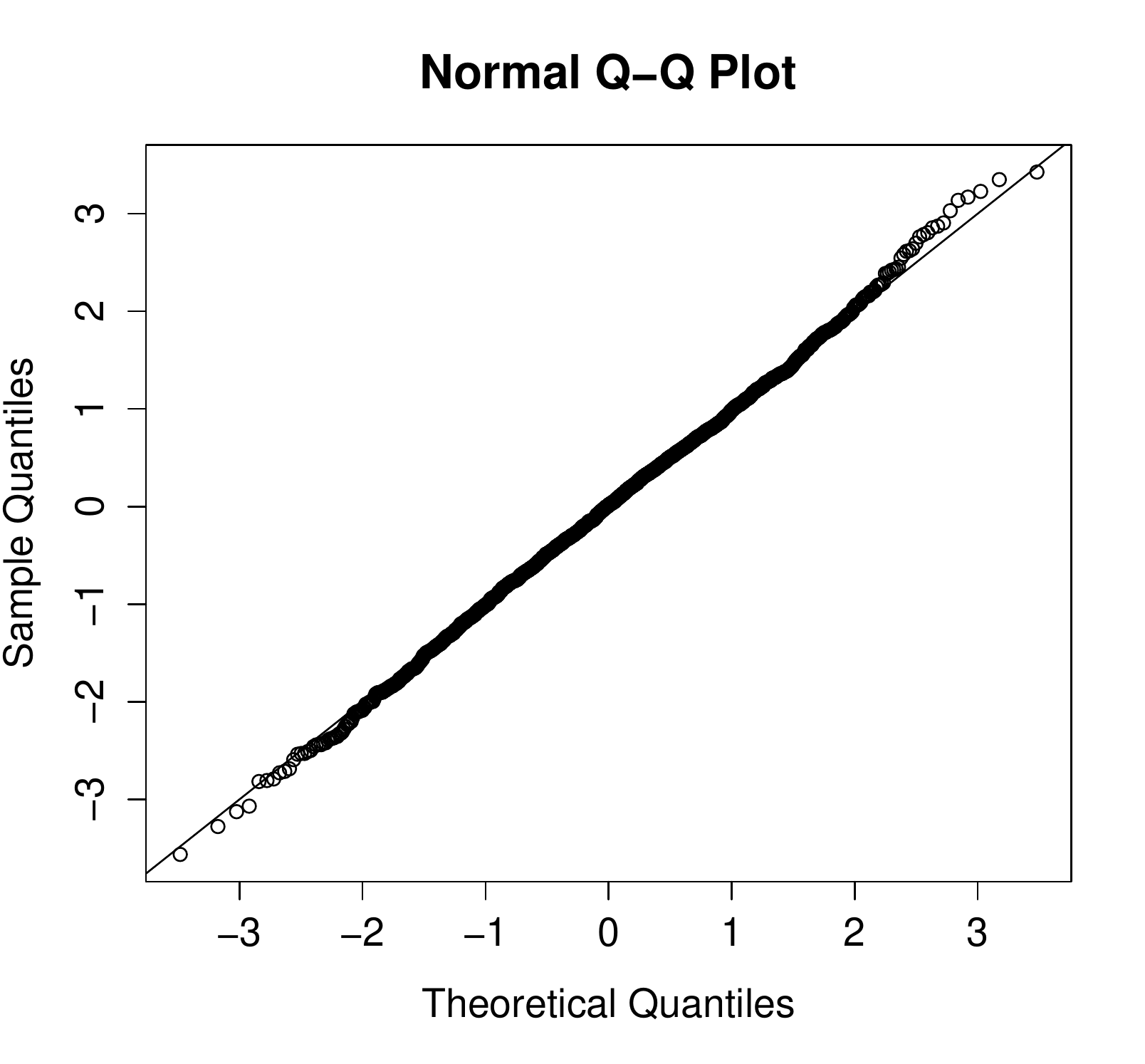} 
    \caption{$\omega=1, T=100$} 
    %\vspace{4ex}
  \end{subfigure}
  \begin{subfigure}[b]{0.32\linewidth}
    \centering
    \includegraphics[width=\linewidth]{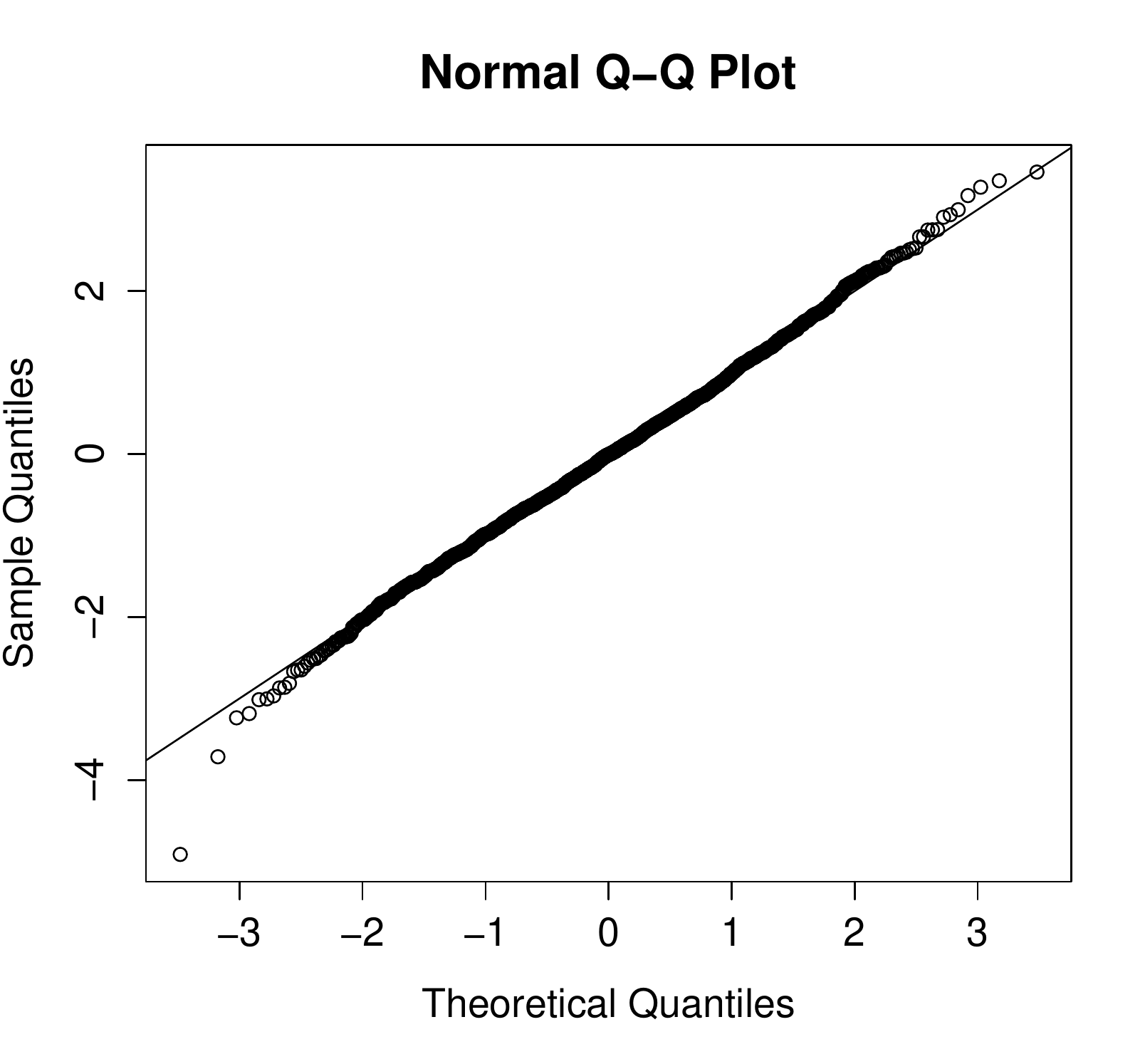} 
    \caption{$\omega=10, T=10$}  
    %\vspace{4ex}
  \end{subfigure}%% 
 \begin{subfigure}[b]{0.32\linewidth}
    \centering
    \includegraphics[width=\linewidth]{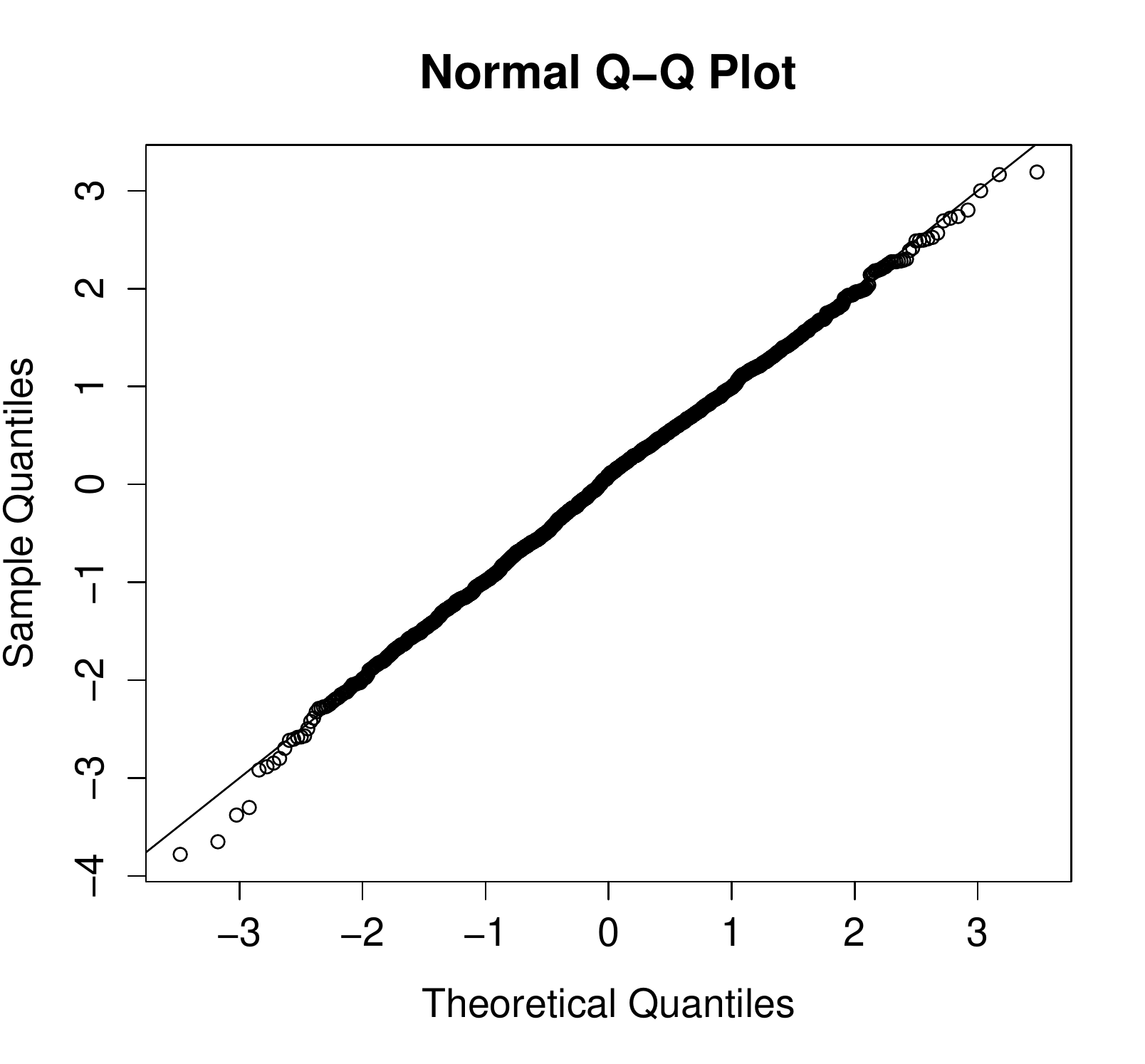} 
    \caption{$\omega=10, T=50$} 
    %\vspace{4ex}
  \end{subfigure}
 \begin{subfigure}[b]{0.32\linewidth}
    \centering
    \includegraphics[width=\linewidth]{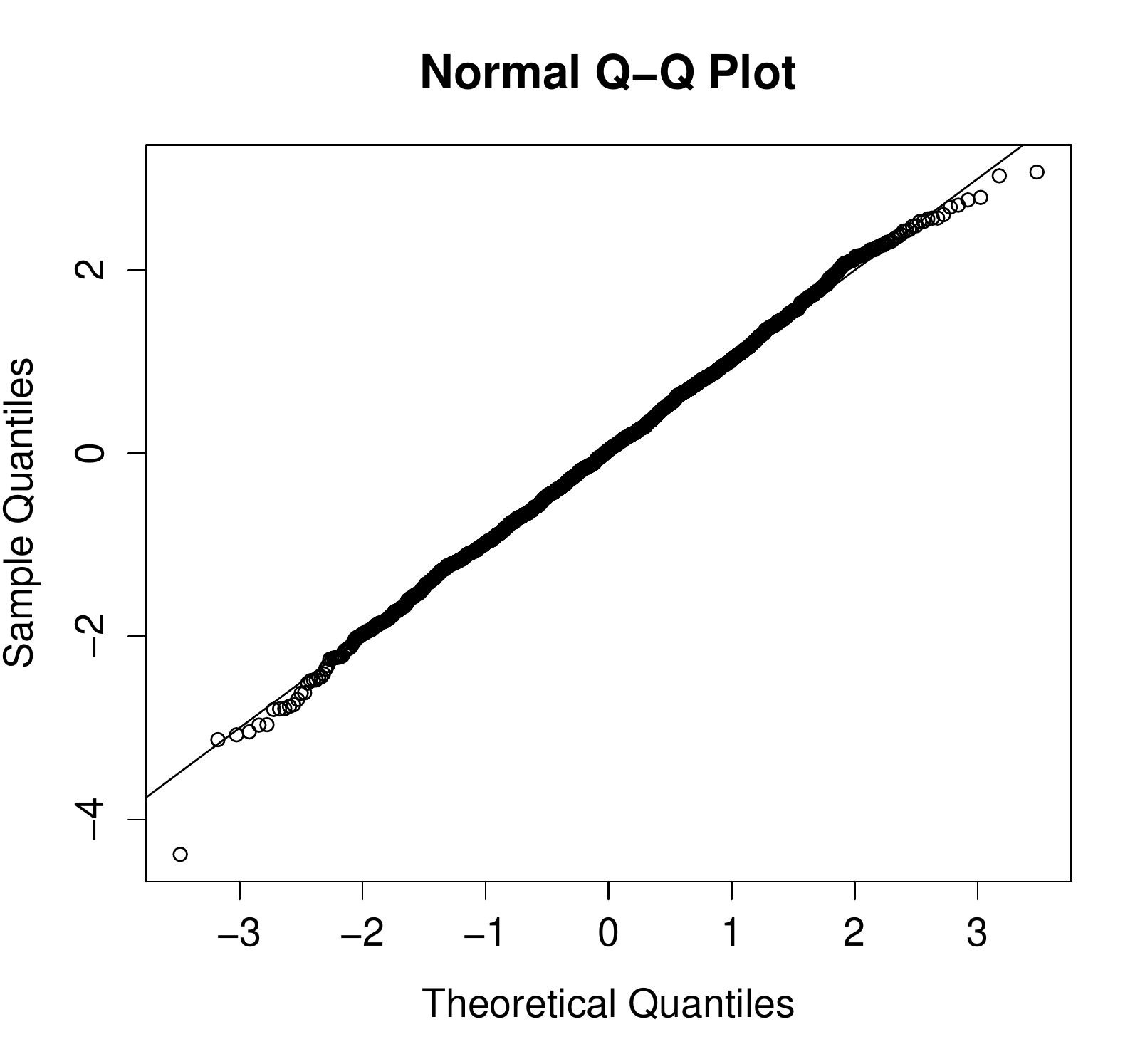} 
    \caption{$\omega=10, T=100$} 
    %\vspace{4ex}
  \end{subfigure}
  \caption{Normal QQ plots for the real part of the truncated Fourier transform of the Ornstein-Uhlenbeck type process driven by a Variance Gamma process for the frequencies $0, 0.1, 1 , 10$ (rows) and time horizons/maximum non-equidistant grid sizes $10/0.1, 50/0.05, 100/0.01$ (columns). The theoretical quantiles are coming from the (limiting) law described in Theorem \ref{thm:ApproxTFTDoubleLimitDistribution}. }\label{plot:QQCARVG} 
\end{figure}

\begin{figure}[tp]    
  \begin{subfigure}[b]{0.32\linewidth}
    \centering
    \includegraphics[width=\linewidth]{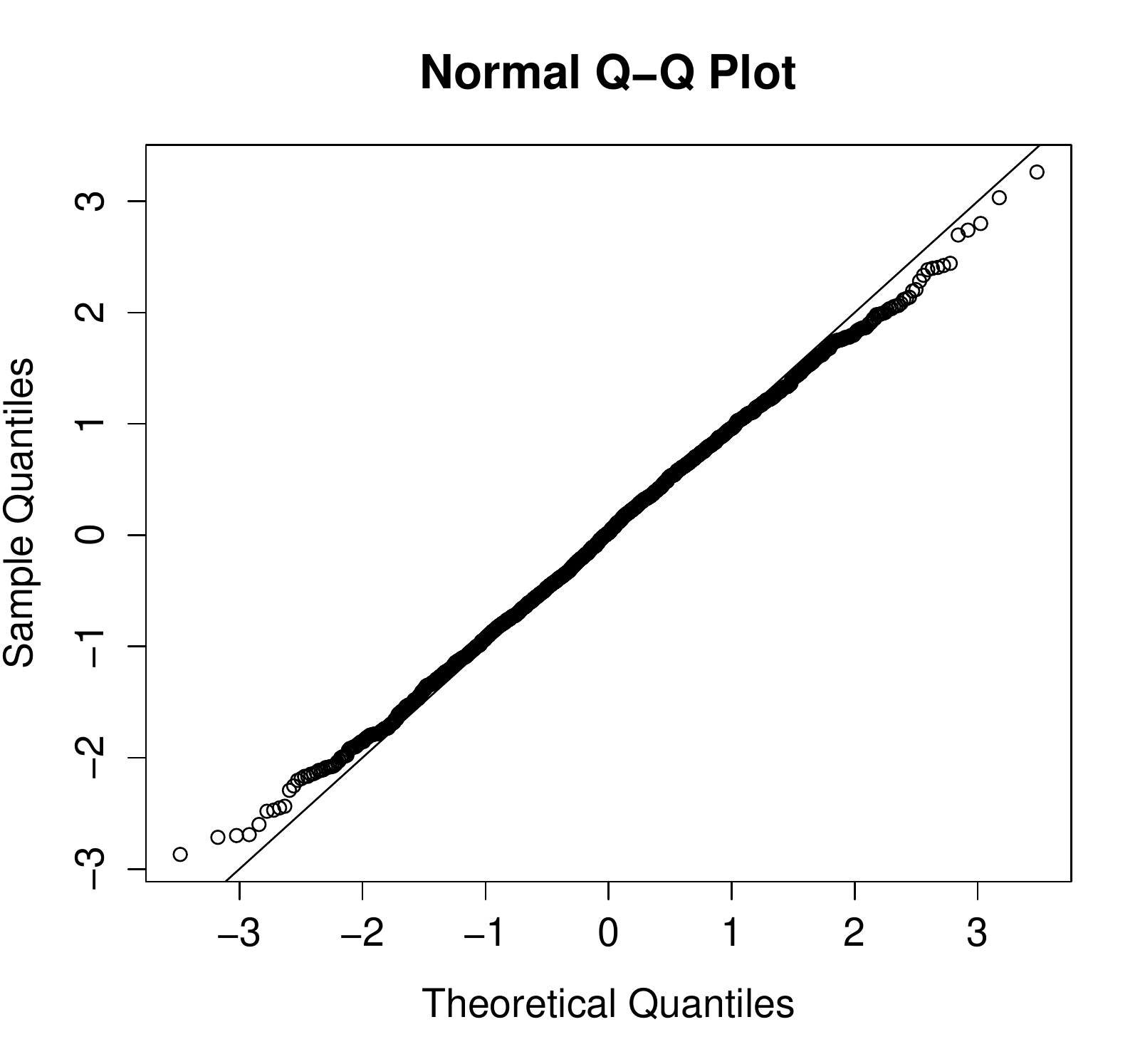} 
    \caption{ $\omega=0, T=10$} 
  
    %\vspace{4ex}
  \end{subfigure}%% 
 \begin{subfigure}[b]{0.32\linewidth}
    \centering
    \includegraphics[width=\linewidth]{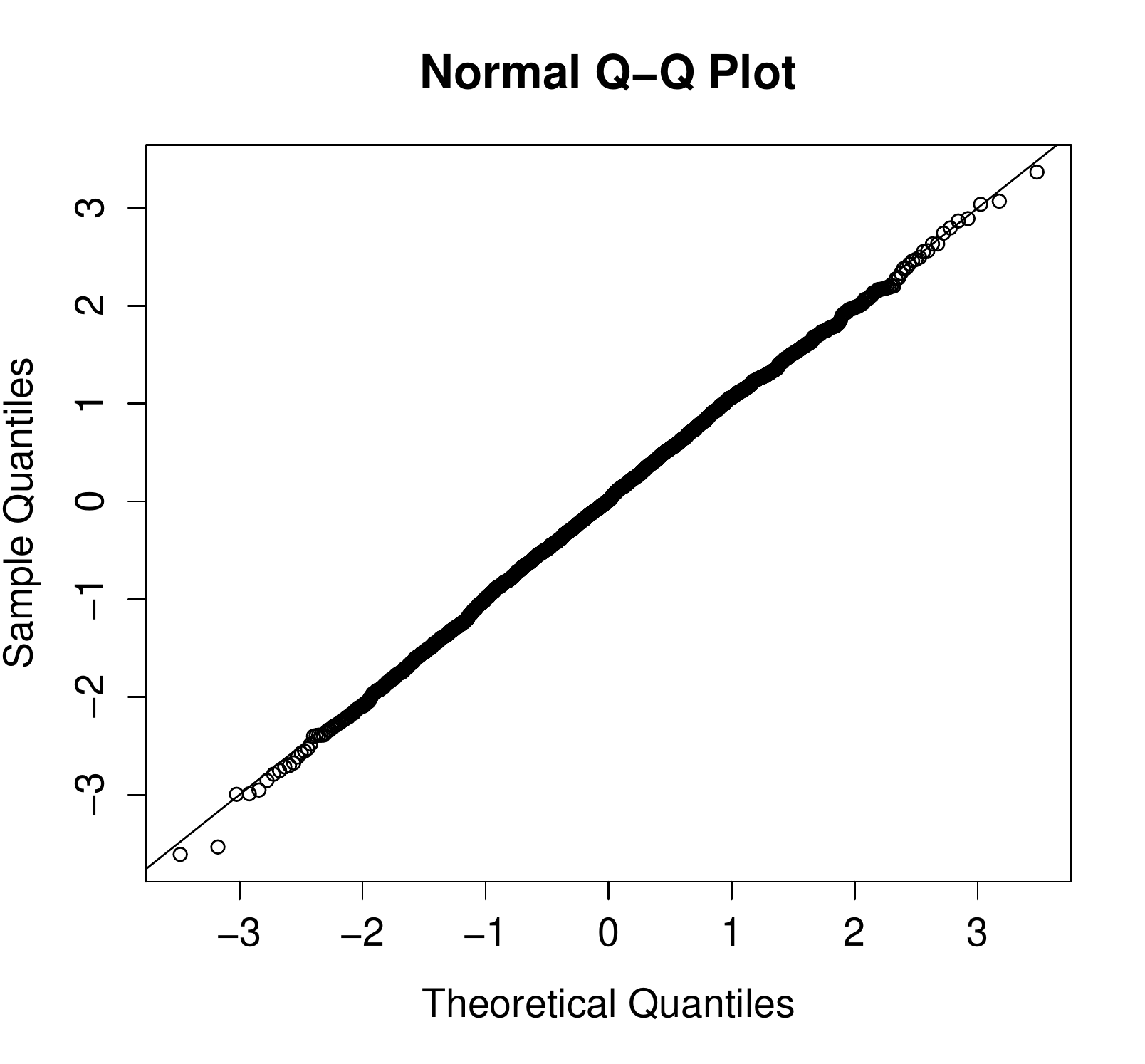} 
    \caption{$\omega=0, T=50$} 
    %\vspace{4ex}
  \end{subfigure}
 \begin{subfigure}[b]{0.32\linewidth}
    \centering
    \includegraphics[width=\linewidth]{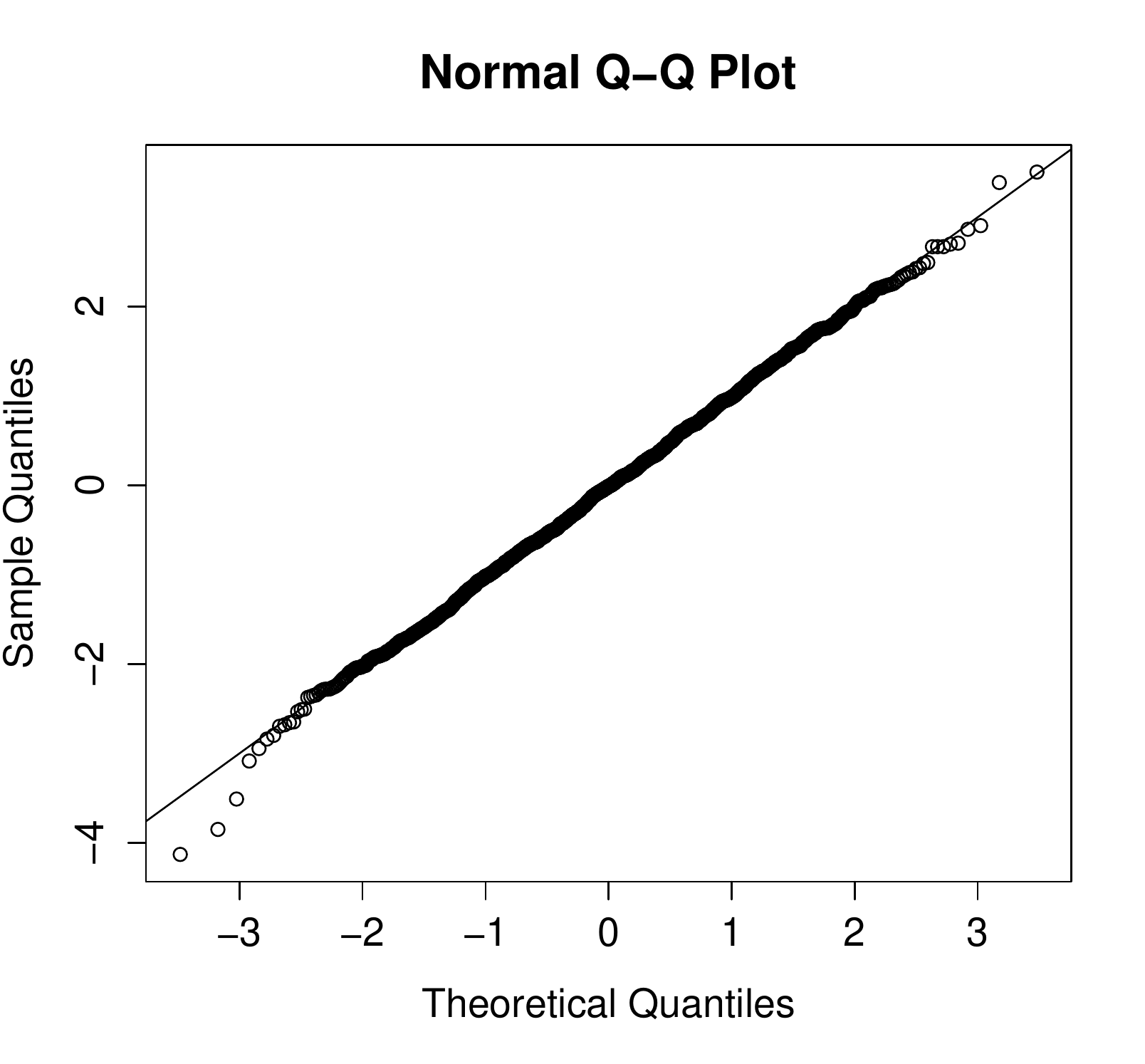} 
    \caption{$\omega=0, T=100$} 
    %\vspace{4ex}
  \end{subfigure}
  \begin{subfigure}[b]{0.32\linewidth}
    \centering
    \includegraphics[width=\linewidth]{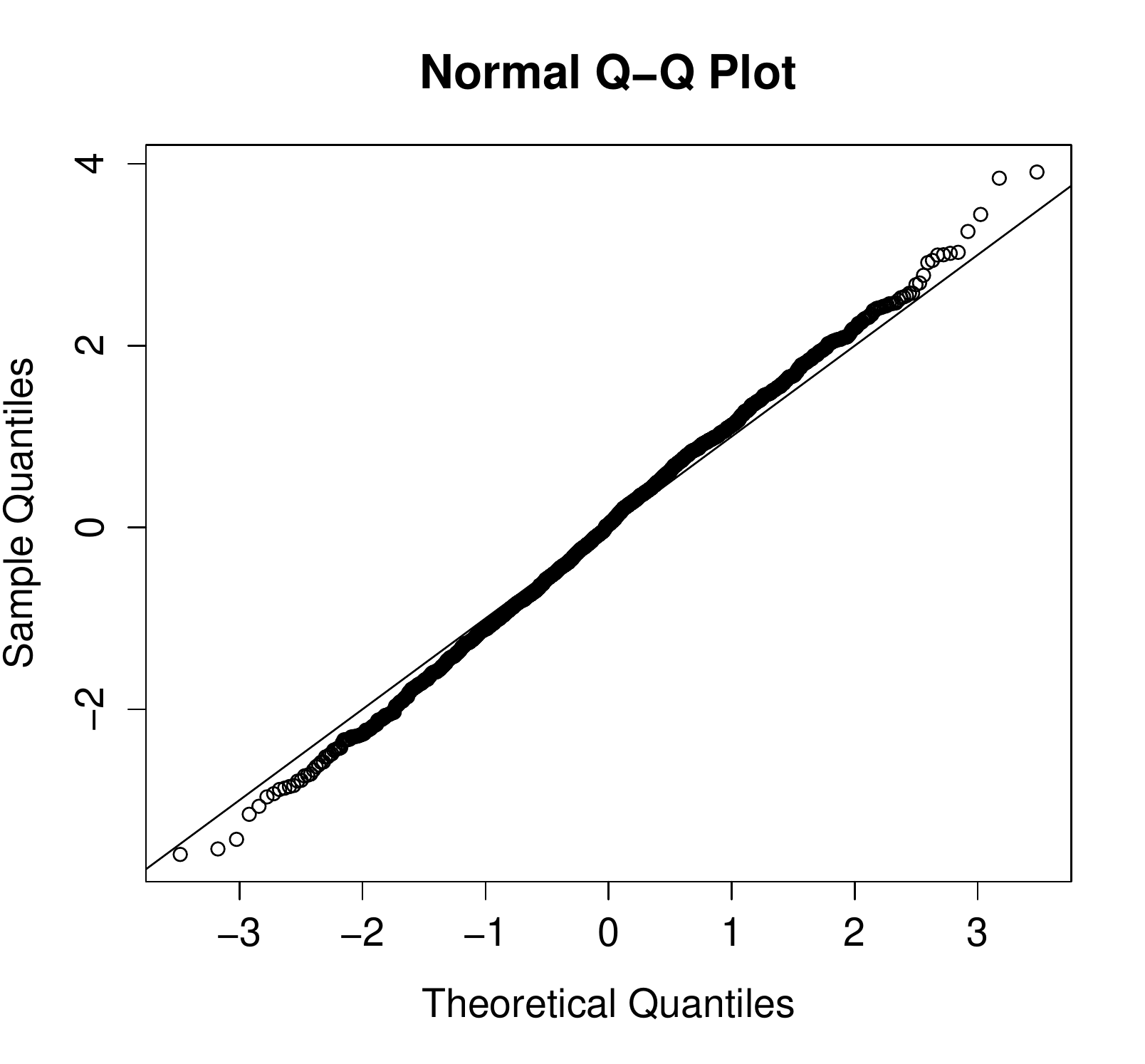} 
    \caption{$\omega=0.1, T=10$} 
    %\vspace{4ex}
  \end{subfigure}%% 
 \begin{subfigure}[b]{0.32\linewidth}
    \centering
    \includegraphics[width=\linewidth]{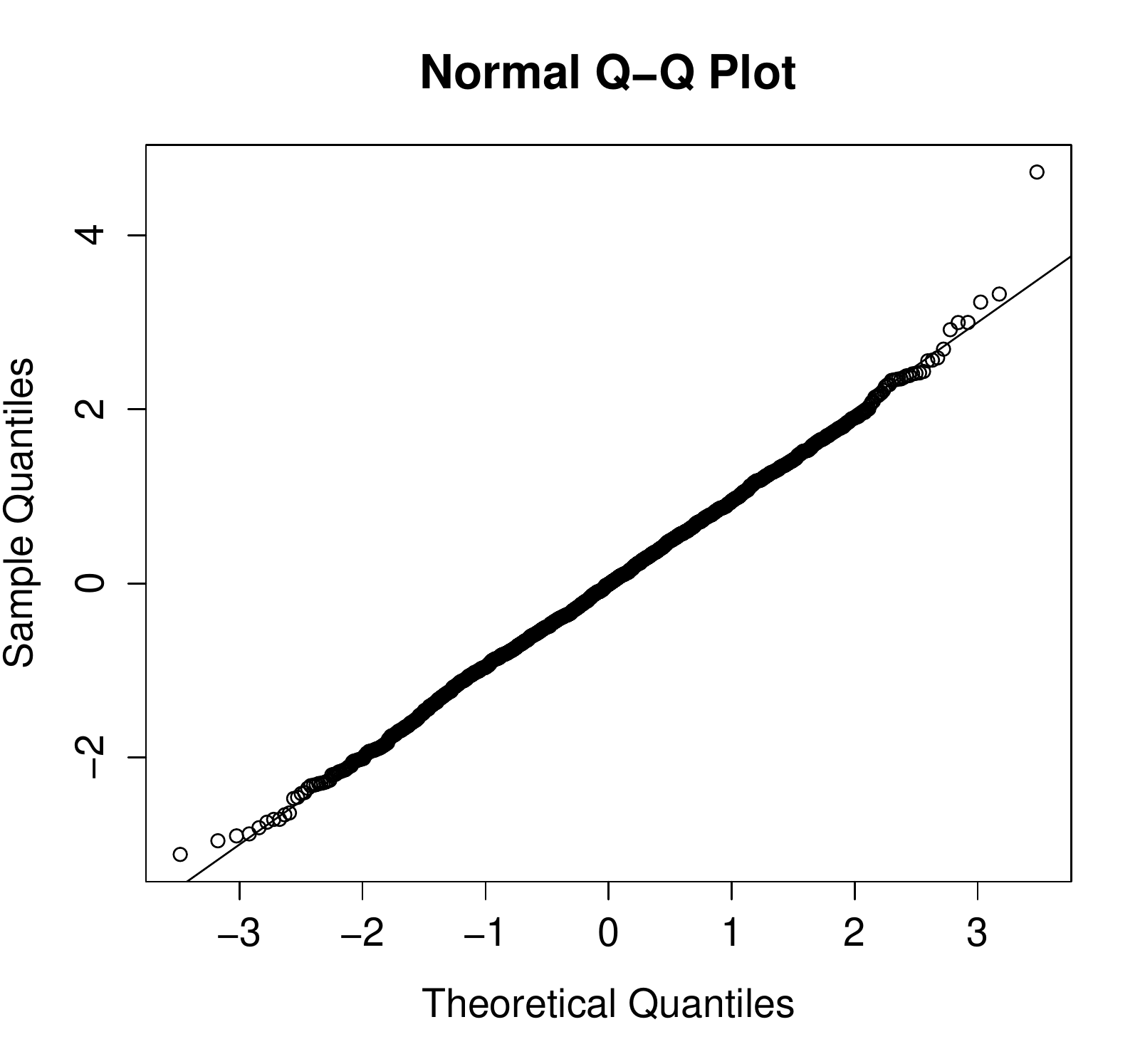} 
    \caption{$\omega=0.1, T=50$} 
    %\vspace{4ex}
  \end{subfigure}
 \begin{subfigure}[b]{0.32\linewidth}
    \centering
    \includegraphics[width=\linewidth]{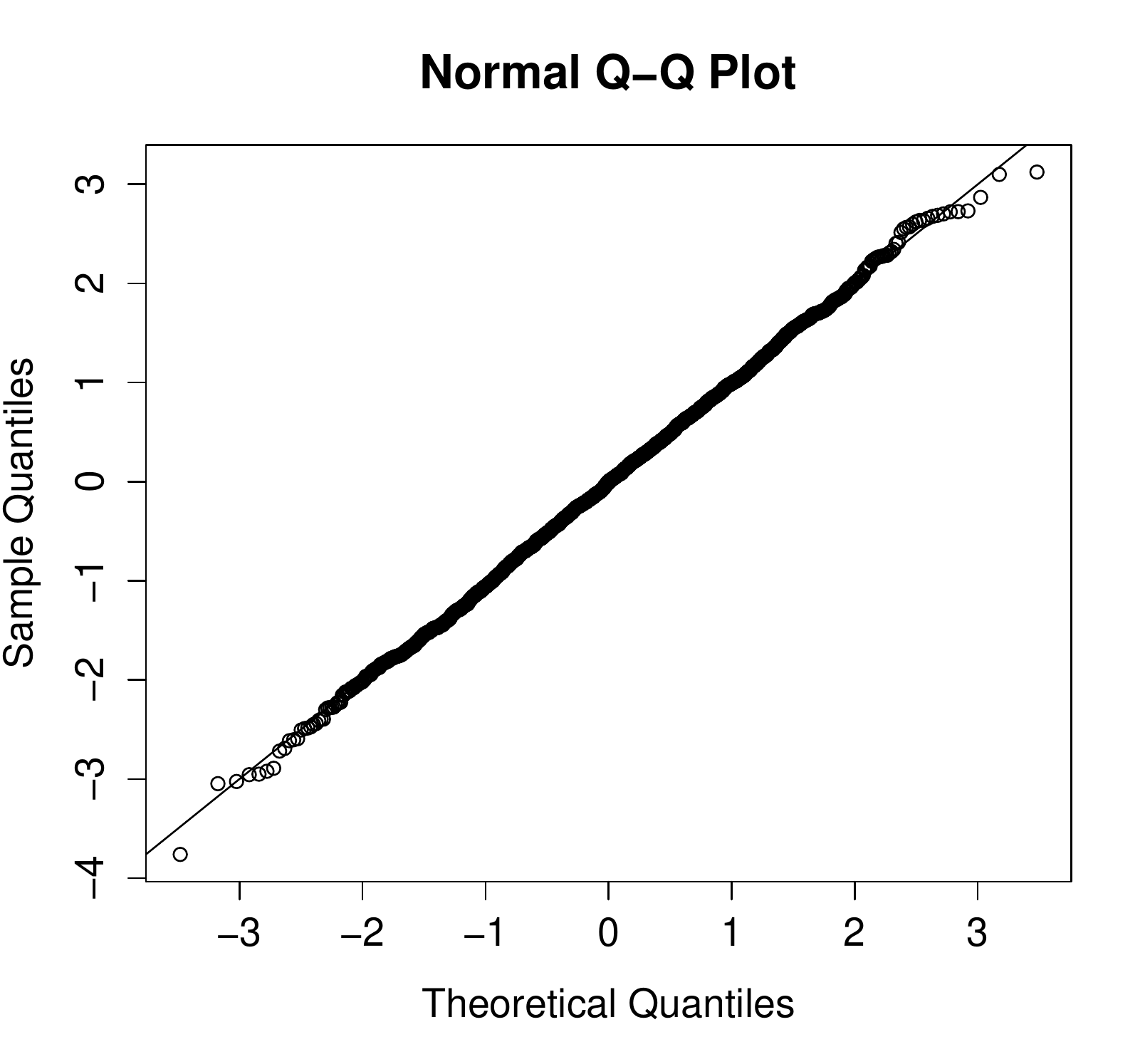} 
    \caption{$\omega=0.1, T=100$} 
    %\vspace{4ex}
  \end{subfigure}
  \begin{subfigure}[b]{0.32\linewidth}
    \centering
    \includegraphics[width=\linewidth]{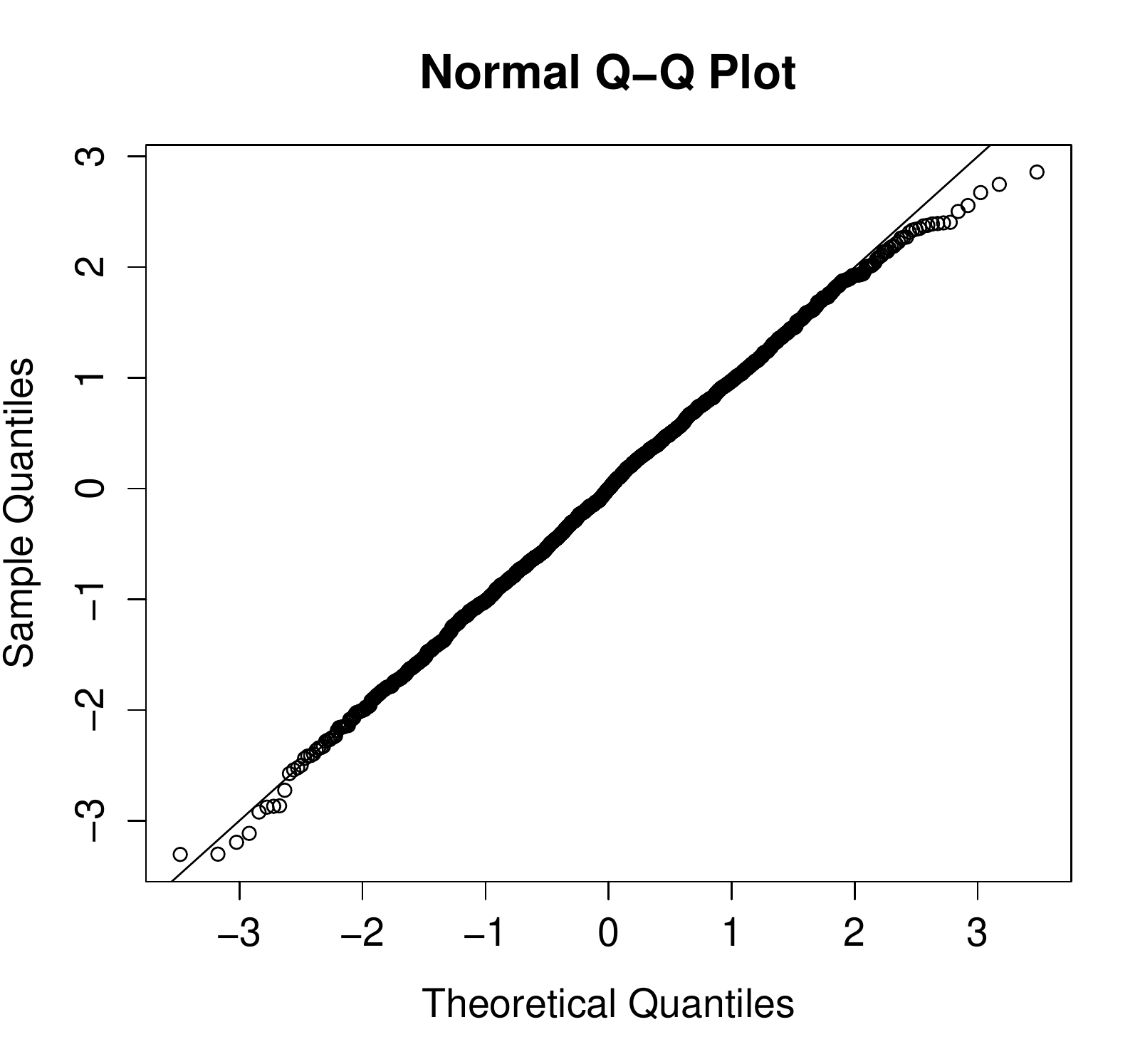} 
    \caption{$\omega=1, T=10$} 
    %\vspace{4ex}
  \end{subfigure}%% 
 \begin{subfigure}[b]{0.32\linewidth}
    \centering
    \includegraphics[width=\linewidth]{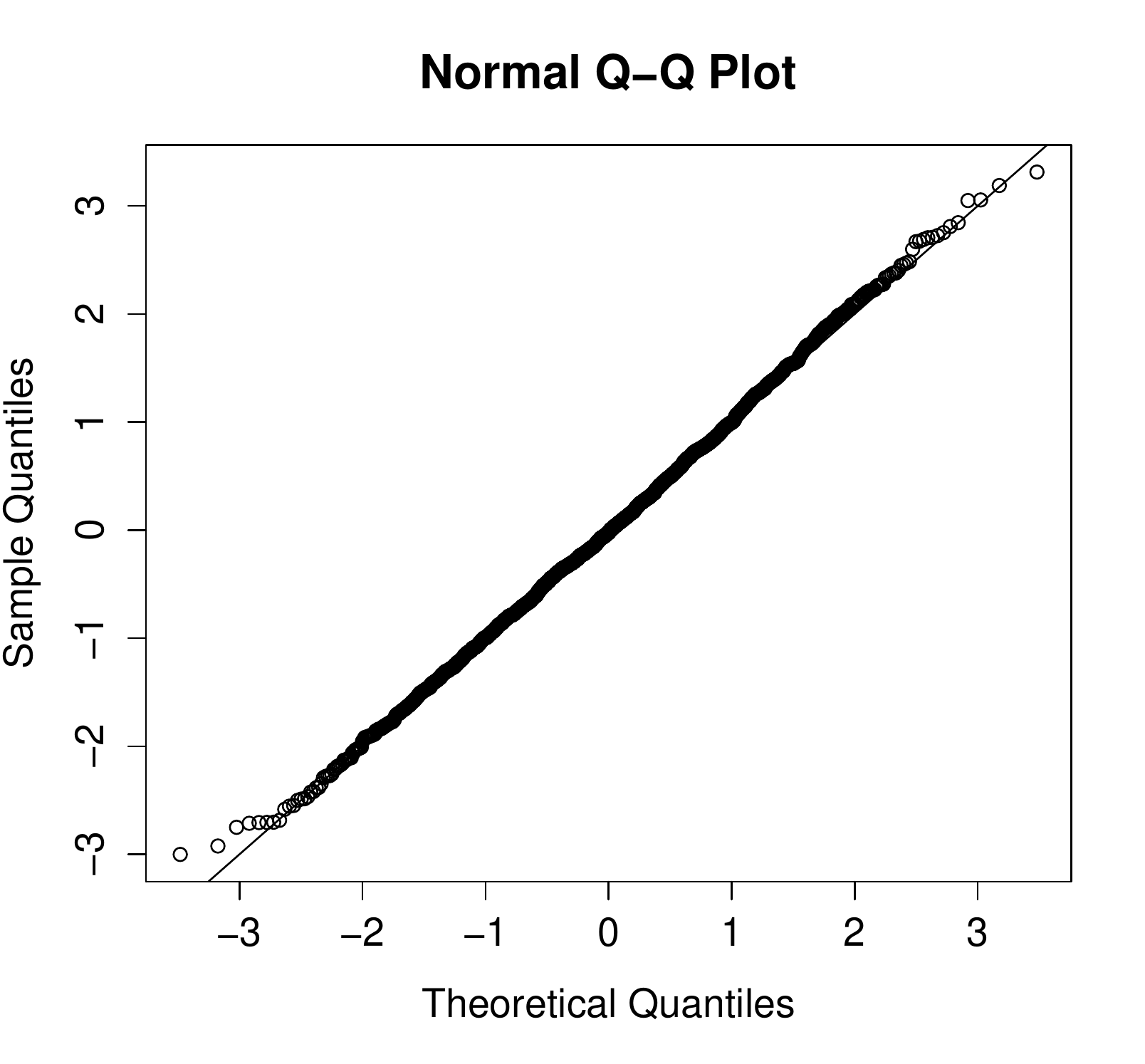} 
    \caption{$\omega=1, T=50$} 
    %\vspace{4ex}
  \end{subfigure}
 \begin{subfigure}[b]{0.32\linewidth}
    \centering
    \includegraphics[width=\linewidth]{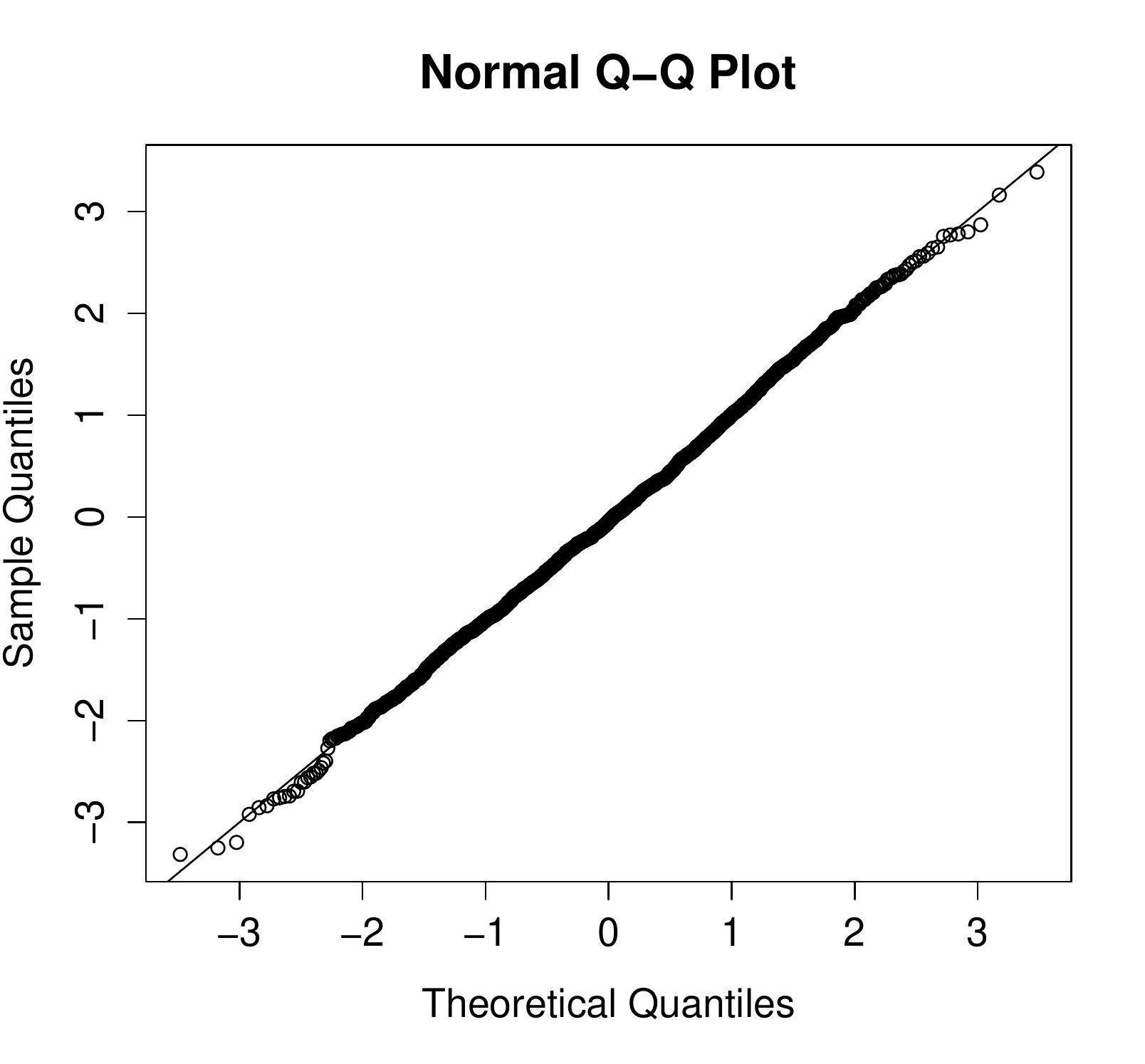} 
    \caption{$\omega=1, T=100$} 
    %\vspace{4ex}
  \end{subfigure}
  \begin{subfigure}[b]{0.32\linewidth}
    \centering
    \includegraphics[width=\linewidth]{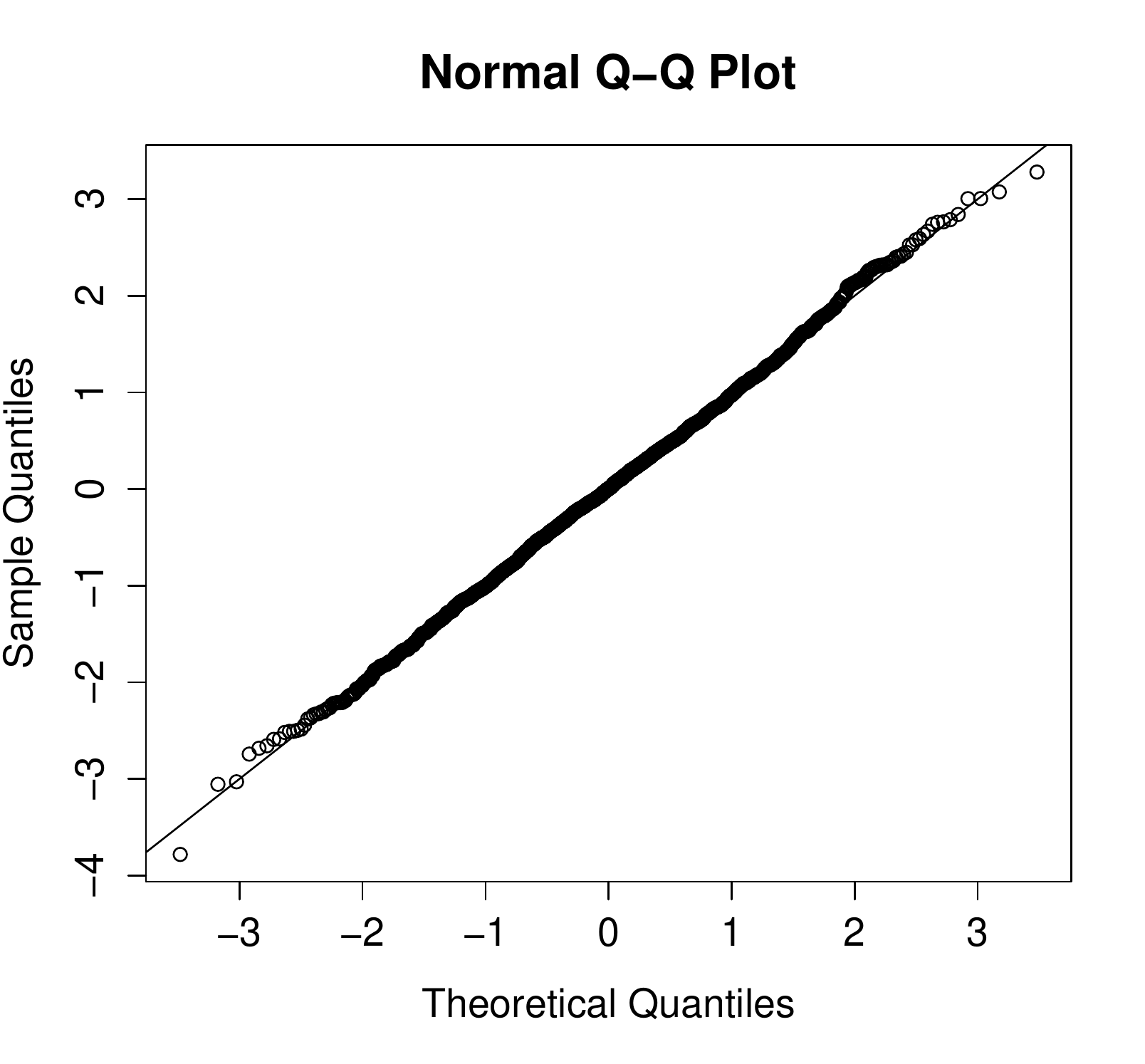} 
    \caption{$\omega=10, T=10$}  
    %\vspace{4ex}
  \end{subfigure}%% 
 \begin{subfigure}[b]{0.32\linewidth}
    \centering
    \includegraphics[width=\linewidth]{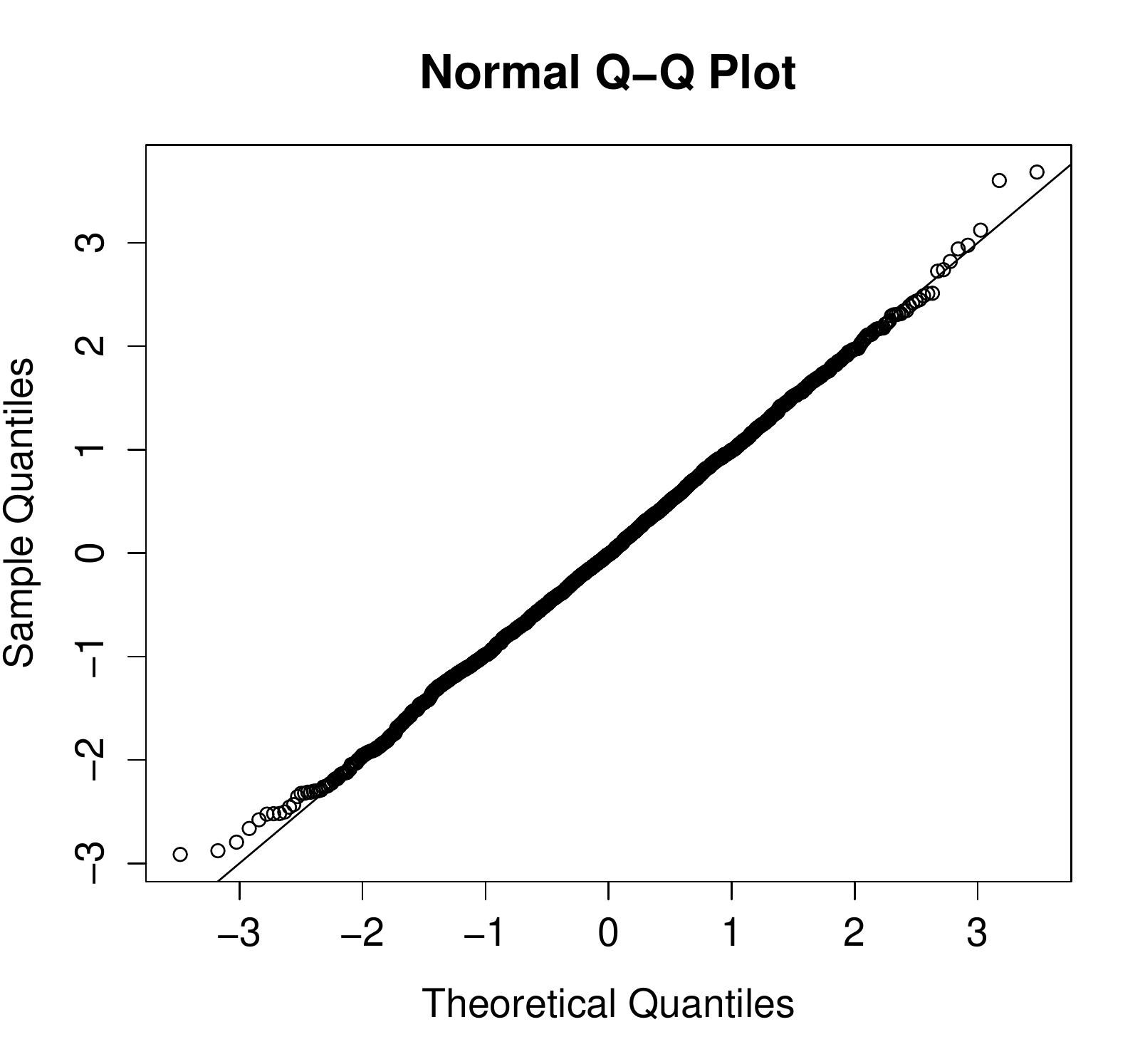} 
    \caption{$\omega=10, T=50$} 
    %\vspace{4ex}
  \end{subfigure}
 \begin{subfigure}[b]{0.32\linewidth}
    \centering
    \includegraphics[width=\linewidth]{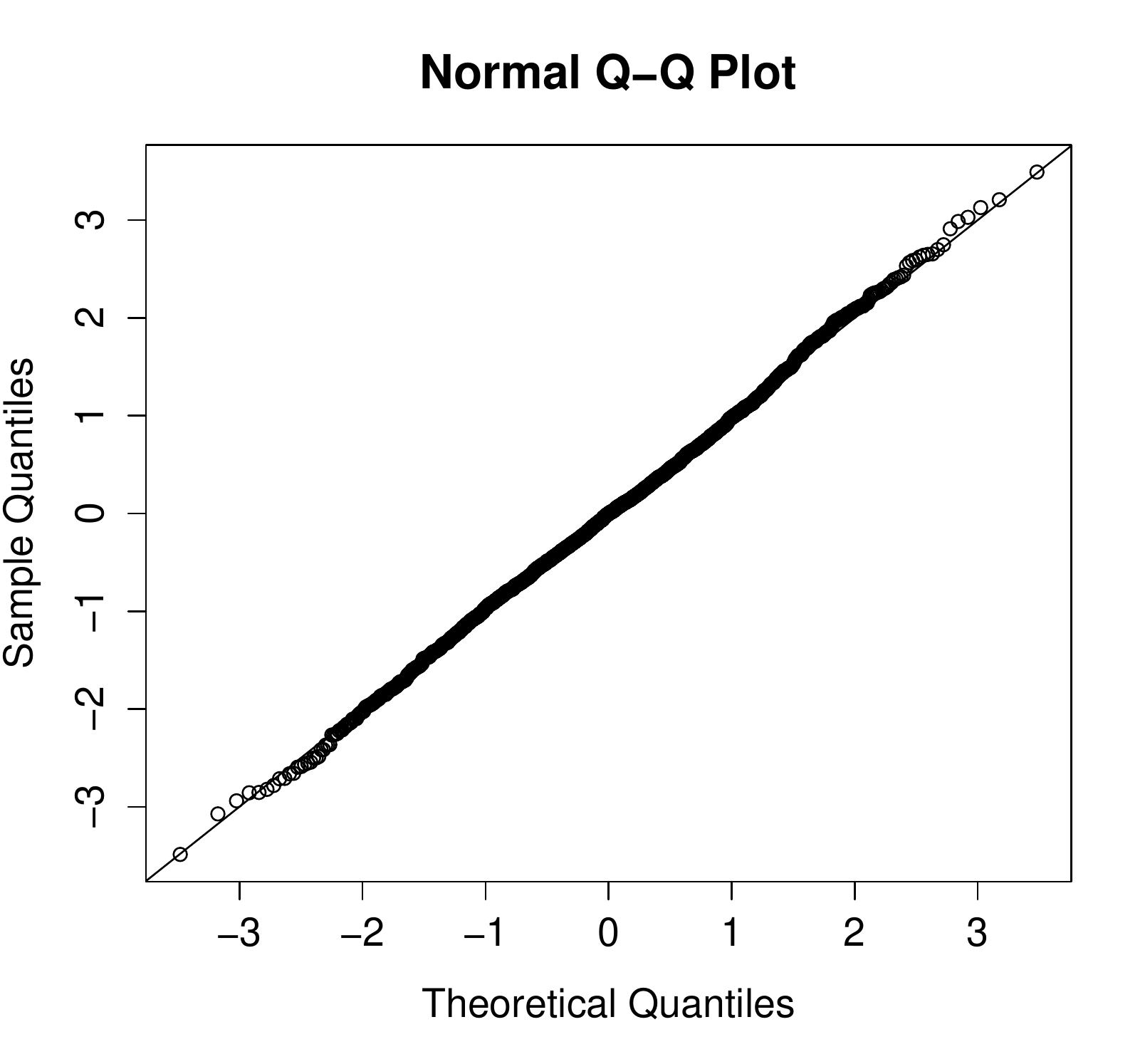} 
    \caption{$\omega=10, T=100$} 
    %\vspace{4ex}
  \end{subfigure}
  \caption{Normal QQ plots for the real part of the truncated Fourier transform of the Ornstein-Uhlenbeck type process driven by a two-sided Poisson process for the frequencies $0, 0.1, 1 , 10$ (rows) and time horizons/maximum non-equidistant grid sizes $10/0.1, 50/0.05, 100/0.01$ (columns). The theoretical quantiles are coming from the (limiting) law described in Theorem \ref{thm:ApproxTFTDoubleLimitDistribution}. }\label{plot:QQCARPois} 
\end{figure}

\begin{example}
 We consider the $\operatorname{CAR}(1)$ model. Then $\A=-a_1$ and 
 $$a(z):=z+a_1,\quad b(z)=b_0.$$
 So the spectral density is 
 $$ f(\omega)=\frac{\sigma^2}{2\pi}\left|\frac{b(i\omega)}{a(i\omega)}\right|^2=\frac{\sigma^2}{2\pi}\frac{b_0^2}{\omega^2+a_1^2}.$$
 For the simulations  we take $[b_0,a_1]=[1,2]$.

QQ-plots showing the results for 2000 simulated paths for the four different frequencies and three different combinations of time horizon and maximum grid width can be found in Figures \ref{plot:QQCARNormal}, \ref{plot:QQCARVG} and \ref{plot:QQCARPois} for the driving L\'evy process being a standard Brownian motion, a Variance Gamma and a two-sided Poisson process, respectively. 
 
\end{example}

\begin{figure}[tp]    
  \begin{subfigure}[b]{0.32\linewidth}
    \centering
    \includegraphics[width=\linewidth]{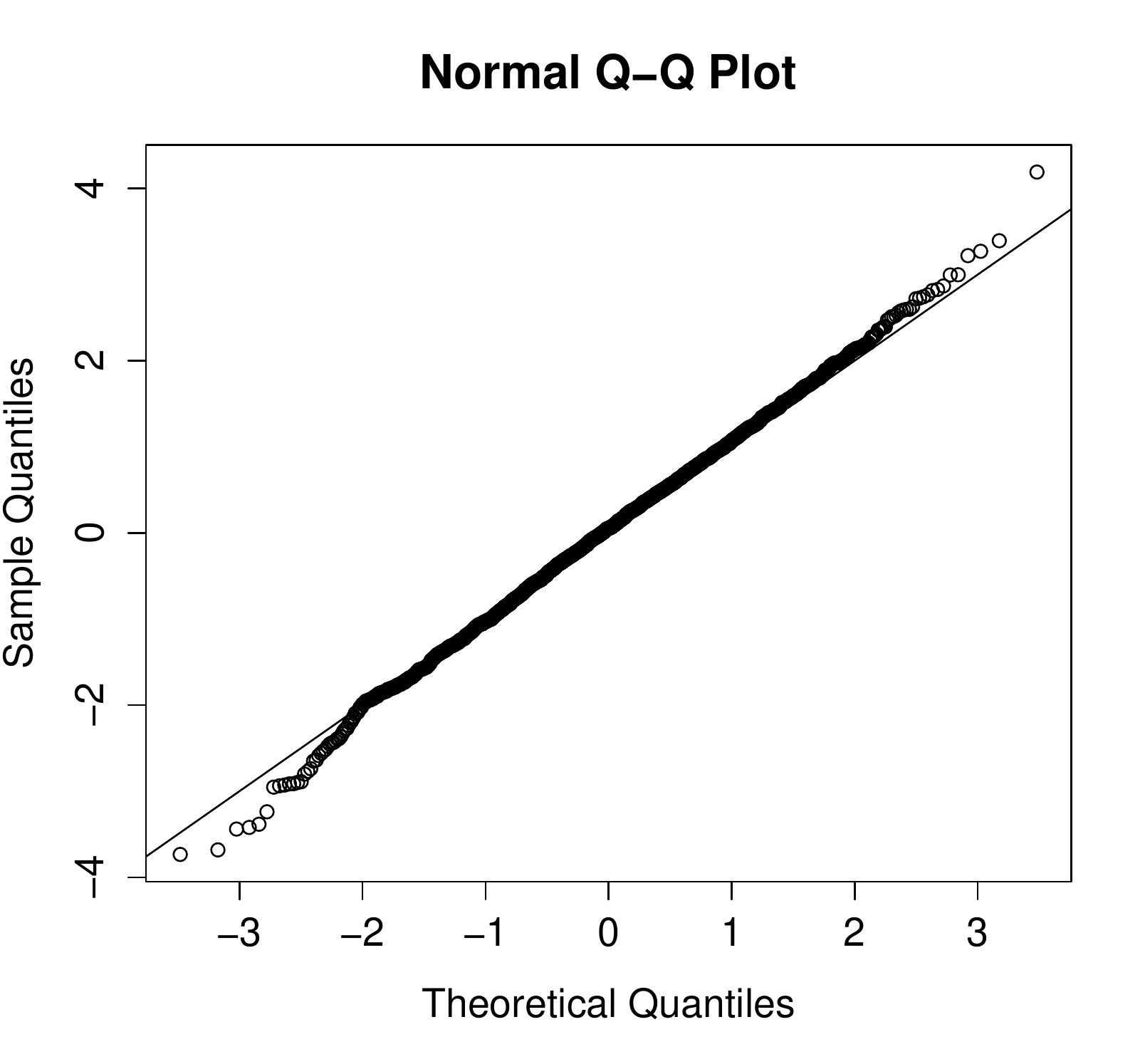} 
    \caption{ $\omega=0, T=10$} 
  
    %\vspace{4ex}
  \end{subfigure}%% 
 \begin{subfigure}[b]{0.32\linewidth}
    \centering
    \includegraphics[width=\linewidth]{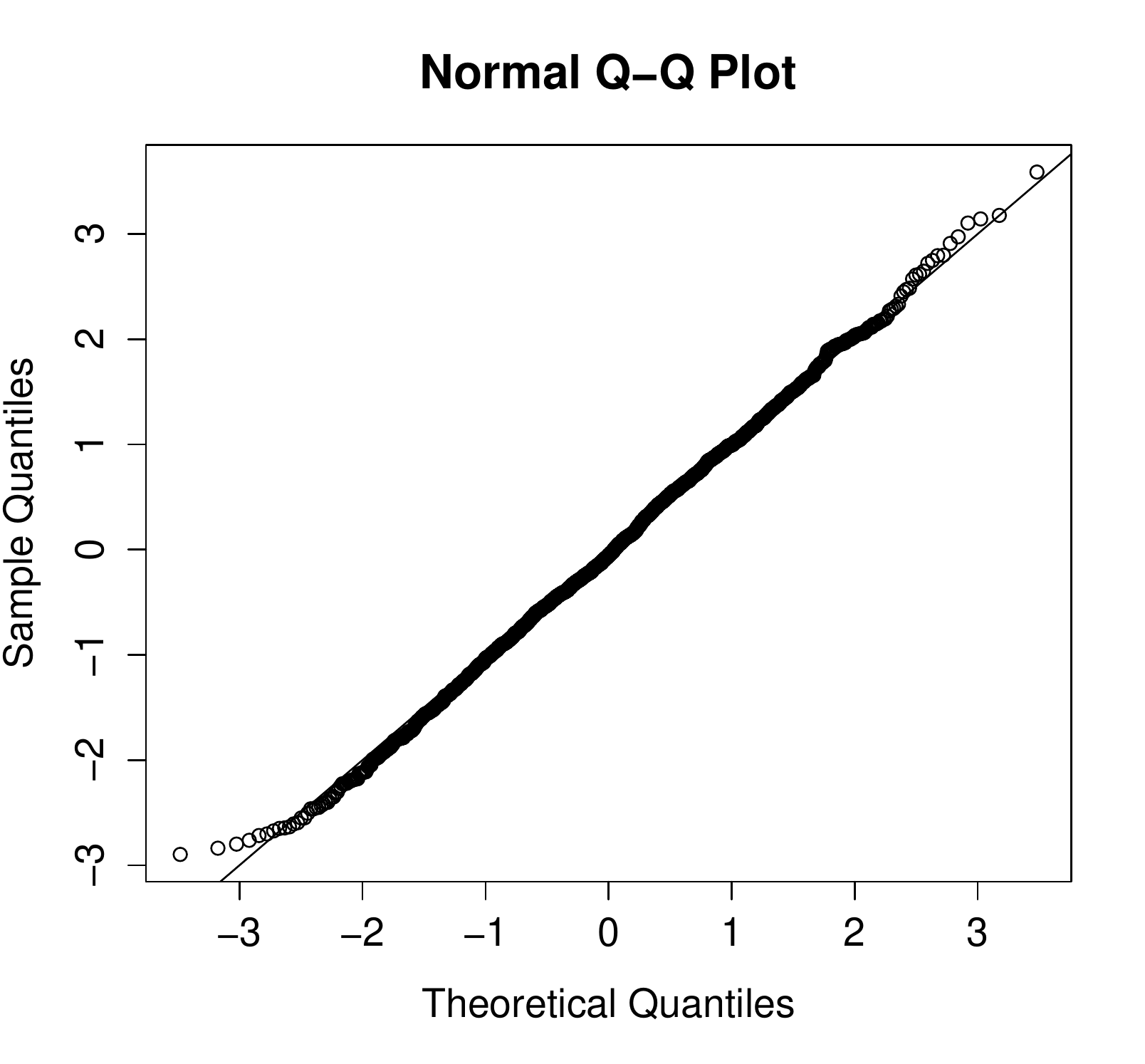} 
    \caption{$\omega=0, T=50$} 
    %\vspace{4ex}
  \end{subfigure}
 \begin{subfigure}[b]{0.32\linewidth}
    \centering
    \includegraphics[width=\linewidth]{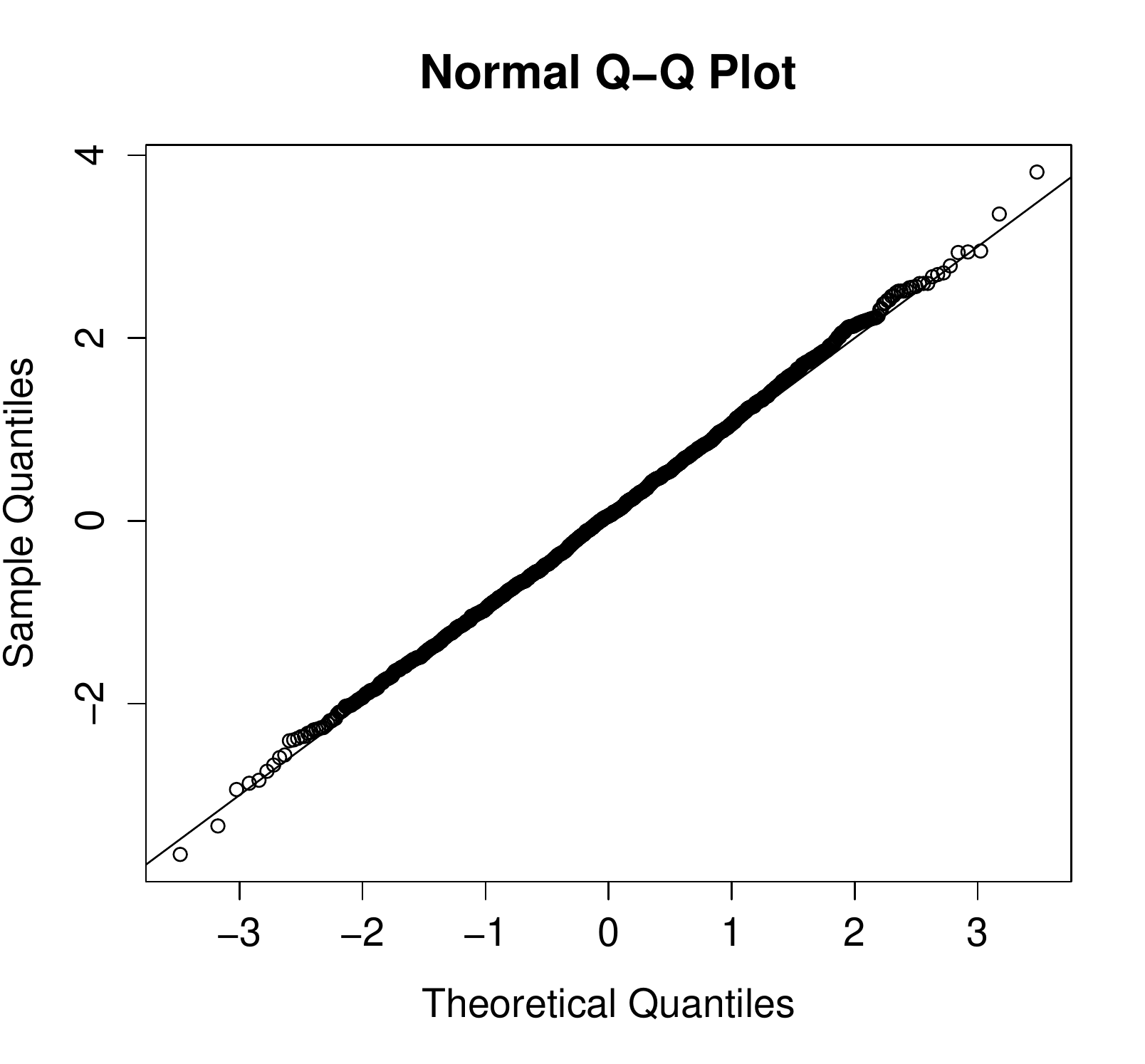} 
    \caption{$\omega=0, T=100$} 
    %\vspace{4ex}
  \end{subfigure}
  \begin{subfigure}[b]{0.32\linewidth}
    \centering
    \includegraphics[width=\linewidth]{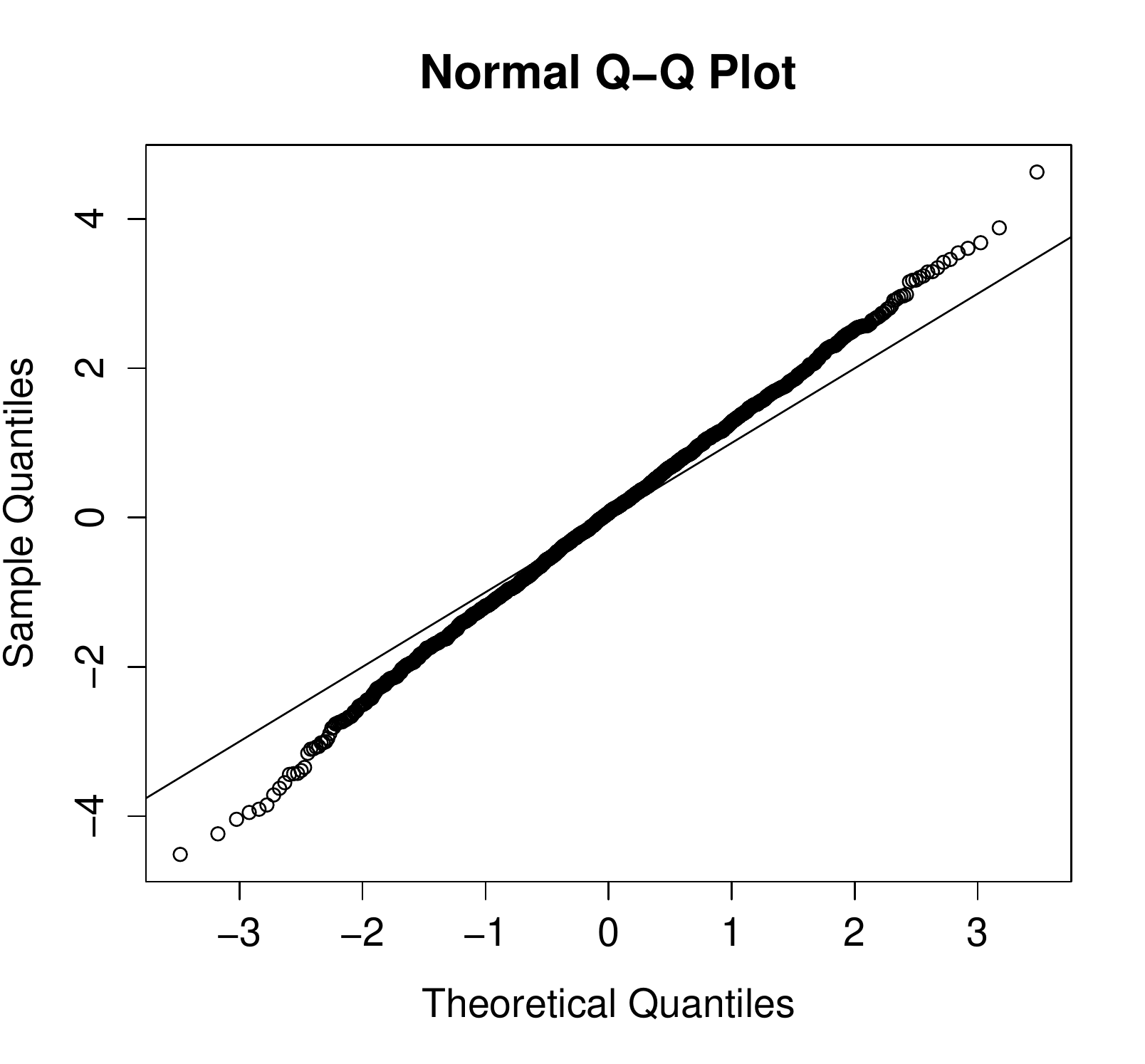} 
    \caption{$\omega=0.1, T=10$} 
    %\vspace{4ex}
  \end{subfigure}%% 
 \begin{subfigure}[b]{0.32\linewidth}
    \centering
    \includegraphics[width=\linewidth]{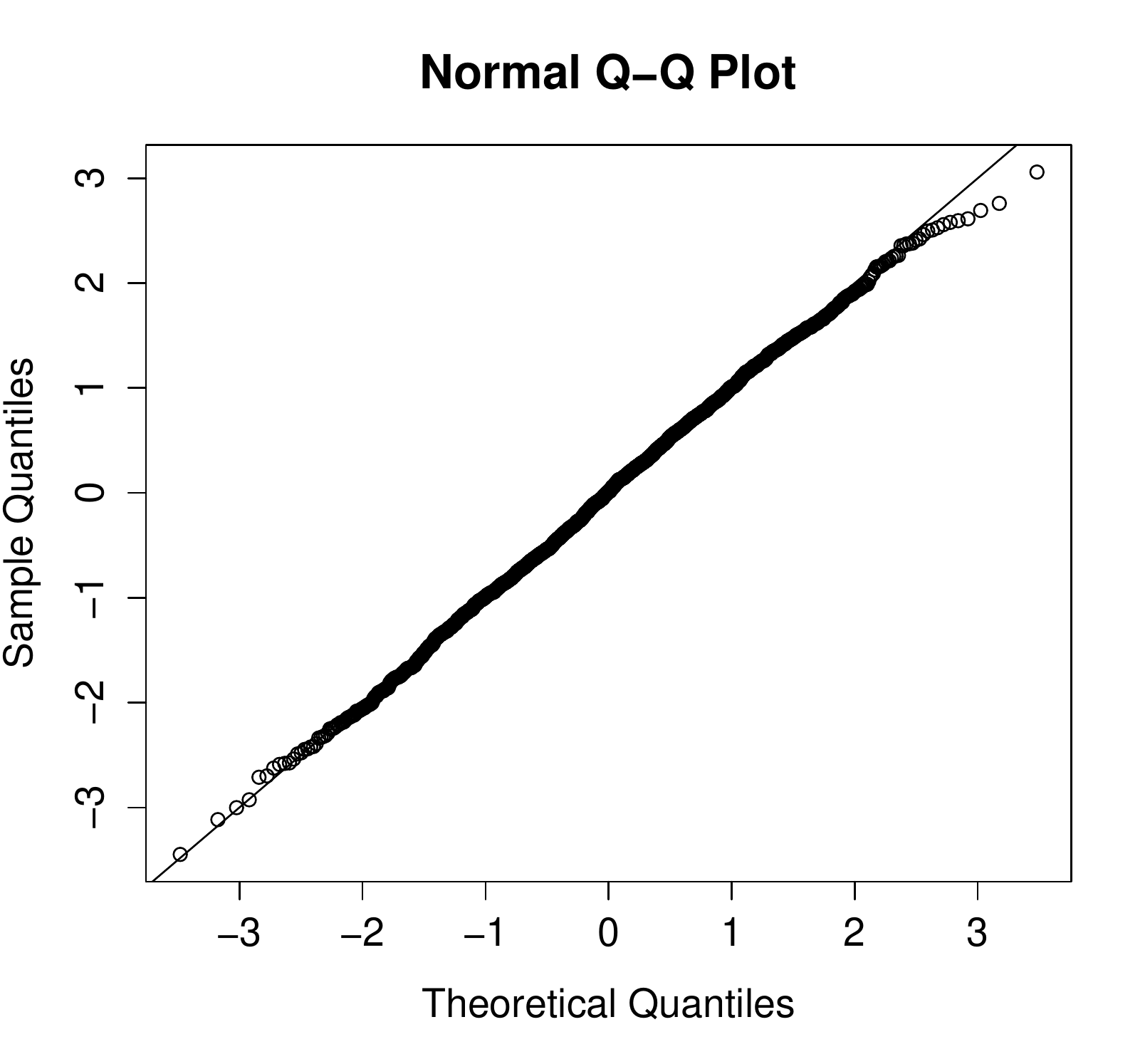} 
    \caption{$\omega=0.1, T=50$} 
    %\vspace{4ex}
  \end{subfigure}
 \begin{subfigure}[b]{0.32\linewidth}
    \centering
    \includegraphics[width=\linewidth]{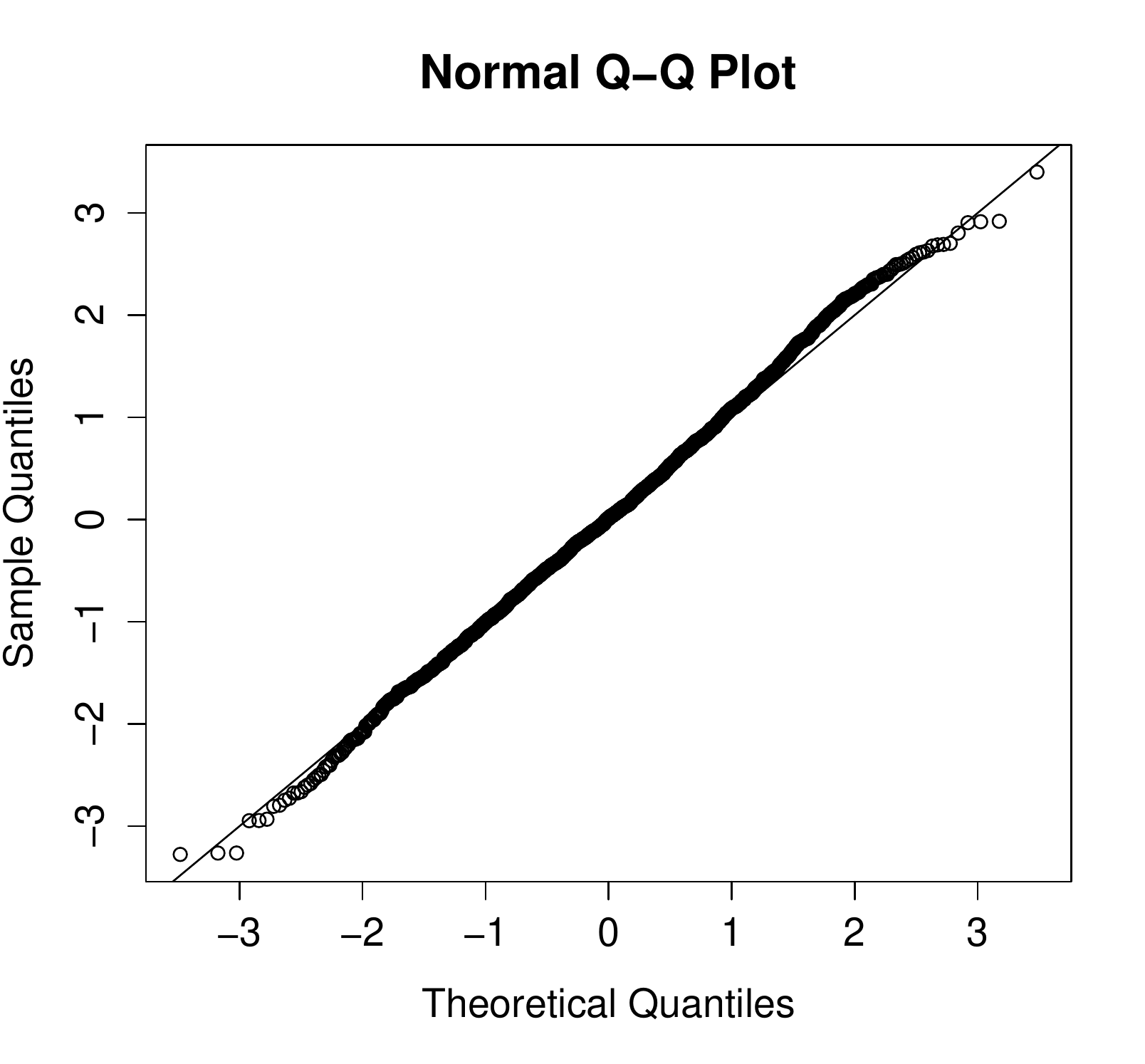} 
    \caption{$\omega=0.1, T=100$} 
    %\vspace{4ex}
  \end{subfigure}
  \begin{subfigure}[b]{0.32\linewidth}
    \centering
    \includegraphics[width=\linewidth]{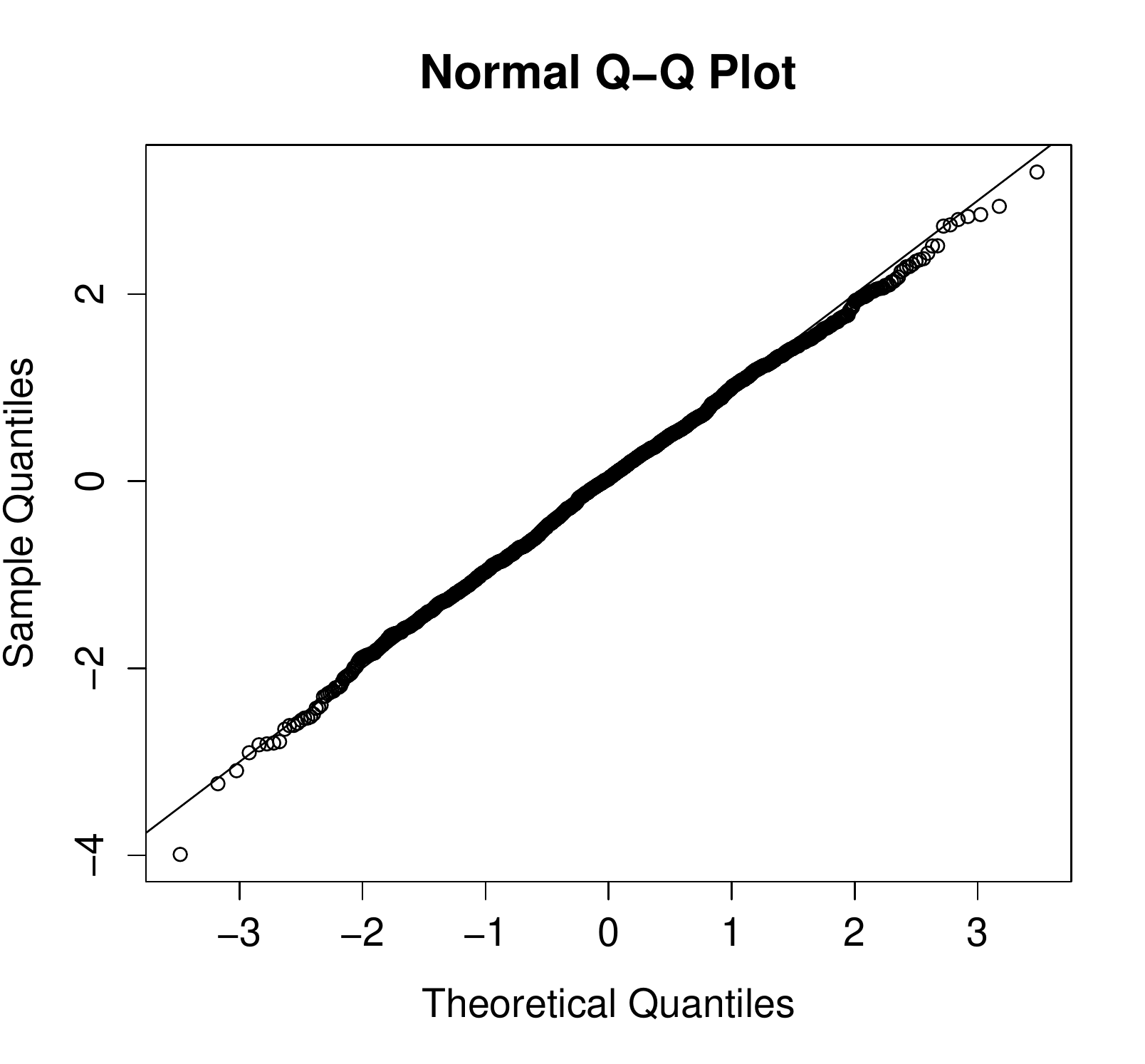} 
    \caption{$\omega=1, T=10$} 
    %\vspace{4ex}
  \end{subfigure}%% 
 \begin{subfigure}[b]{0.32\linewidth}
    \centering
    \includegraphics[width=\linewidth]{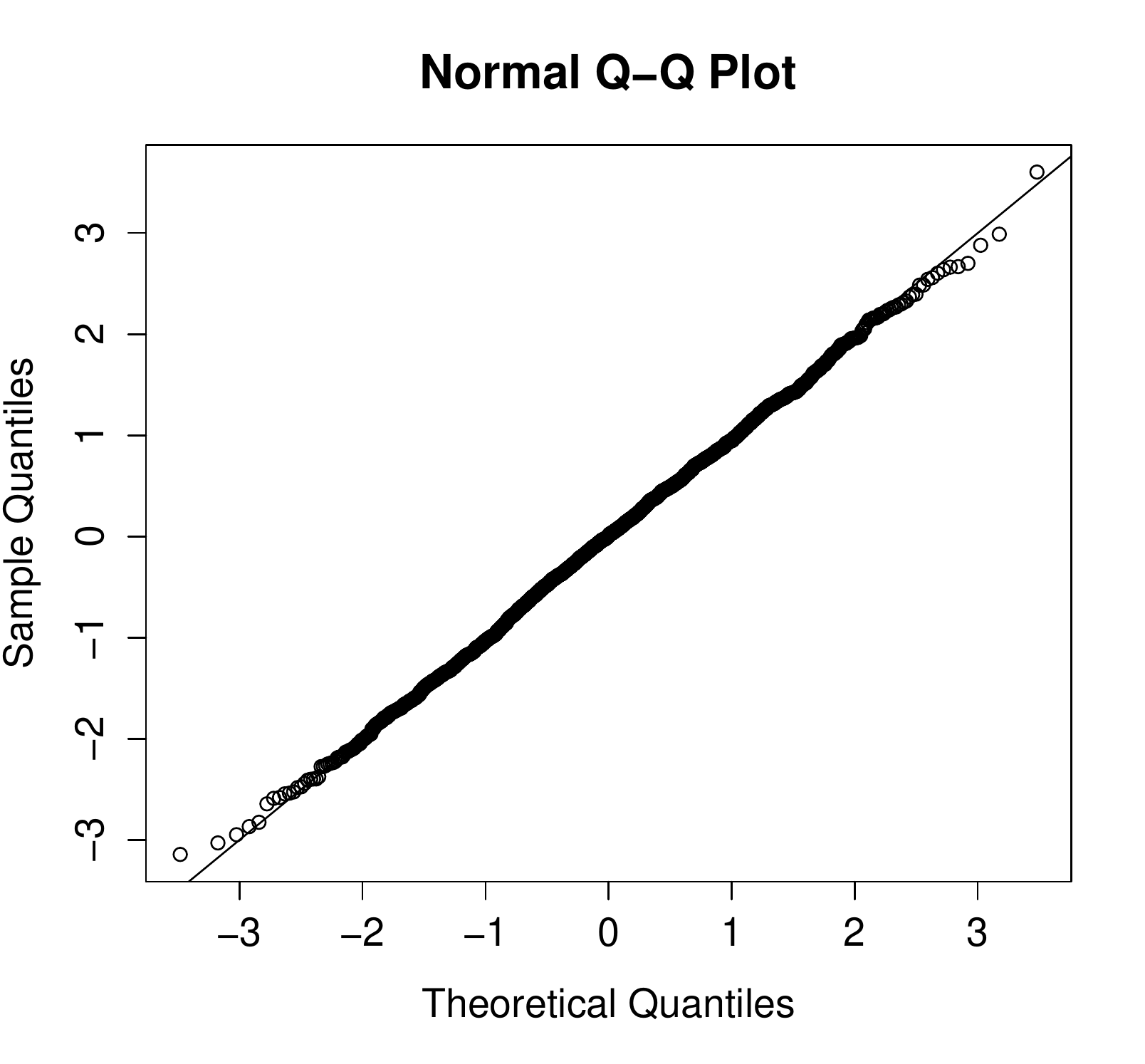} 
    \caption{$\omega=1, T=50$} 
    %\vspace{4ex}
  \end{subfigure}
 \begin{subfigure}[b]{0.32\linewidth}
    \centering
    \includegraphics[width=\linewidth]{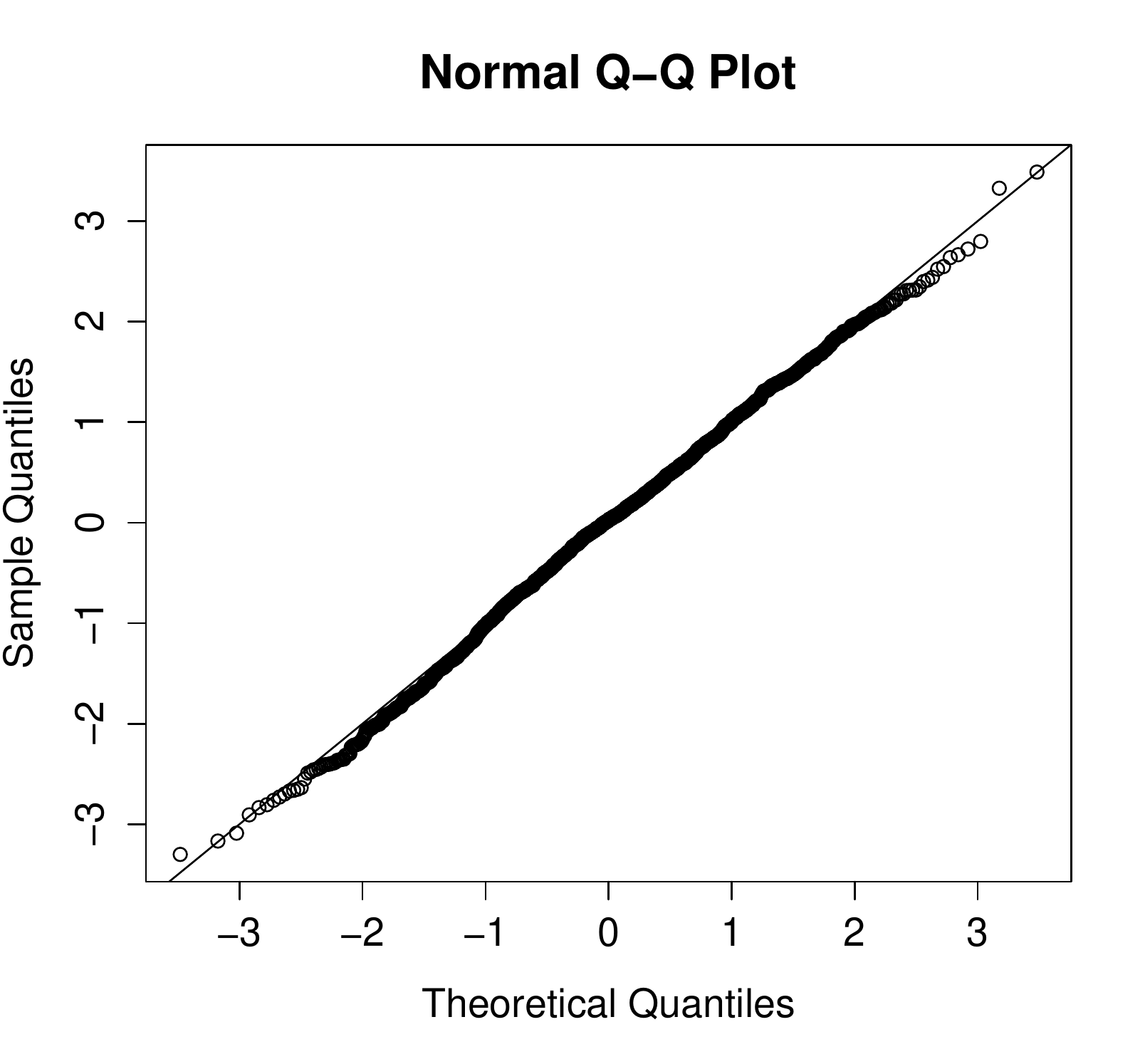} 
    \caption{$\omega=1, T=100$} 
    %\vspace{4ex}
  \end{subfigure}
  \begin{subfigure}[b]{0.32\linewidth}
    \centering
    \includegraphics[width=\linewidth]{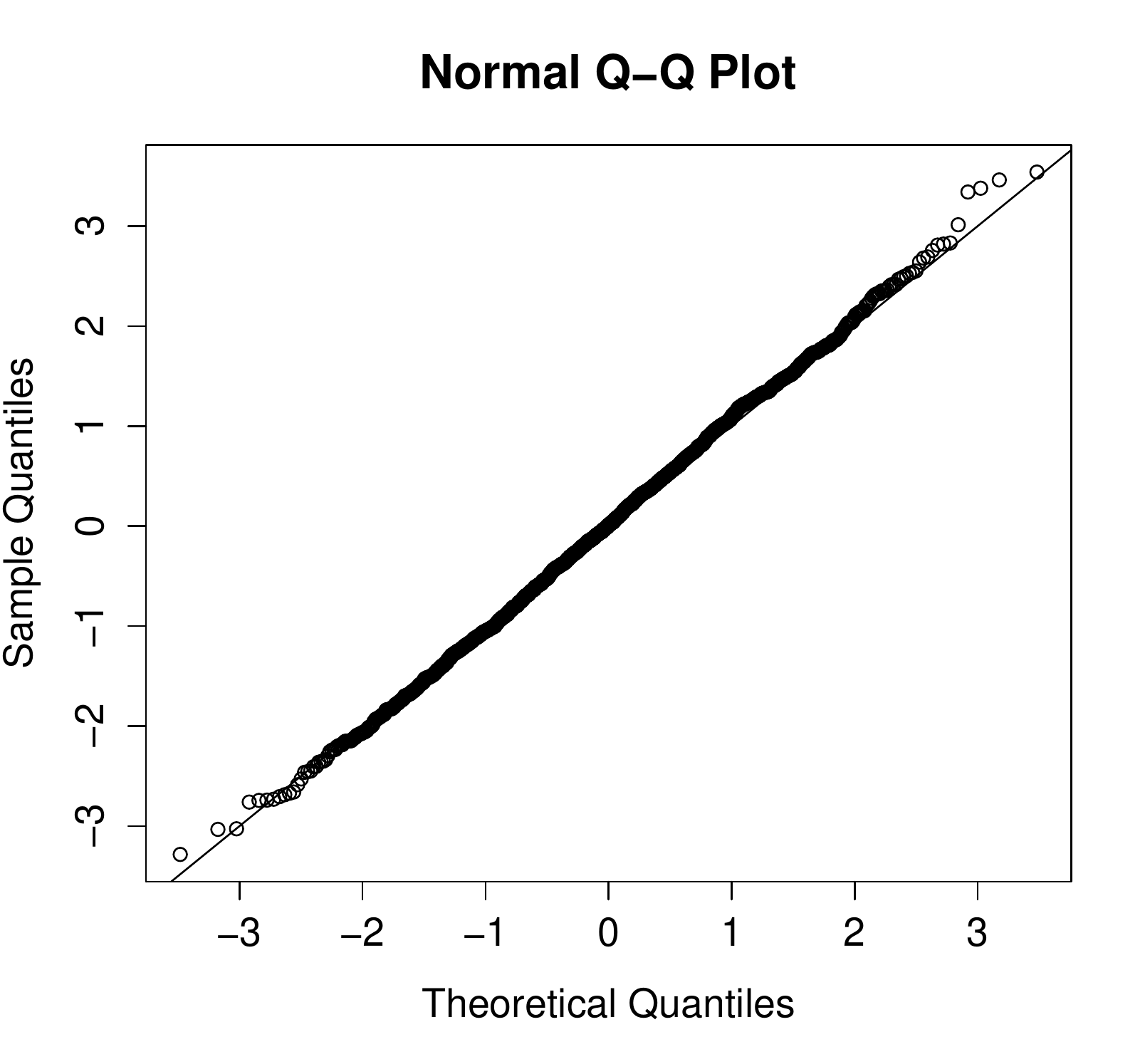} 
    \caption{$\omega=10, T=10$}  
    %\vspace{4ex}
  \end{subfigure}%% 
 \begin{subfigure}[b]{0.32\linewidth}
    \centering
    \includegraphics[width=\linewidth]{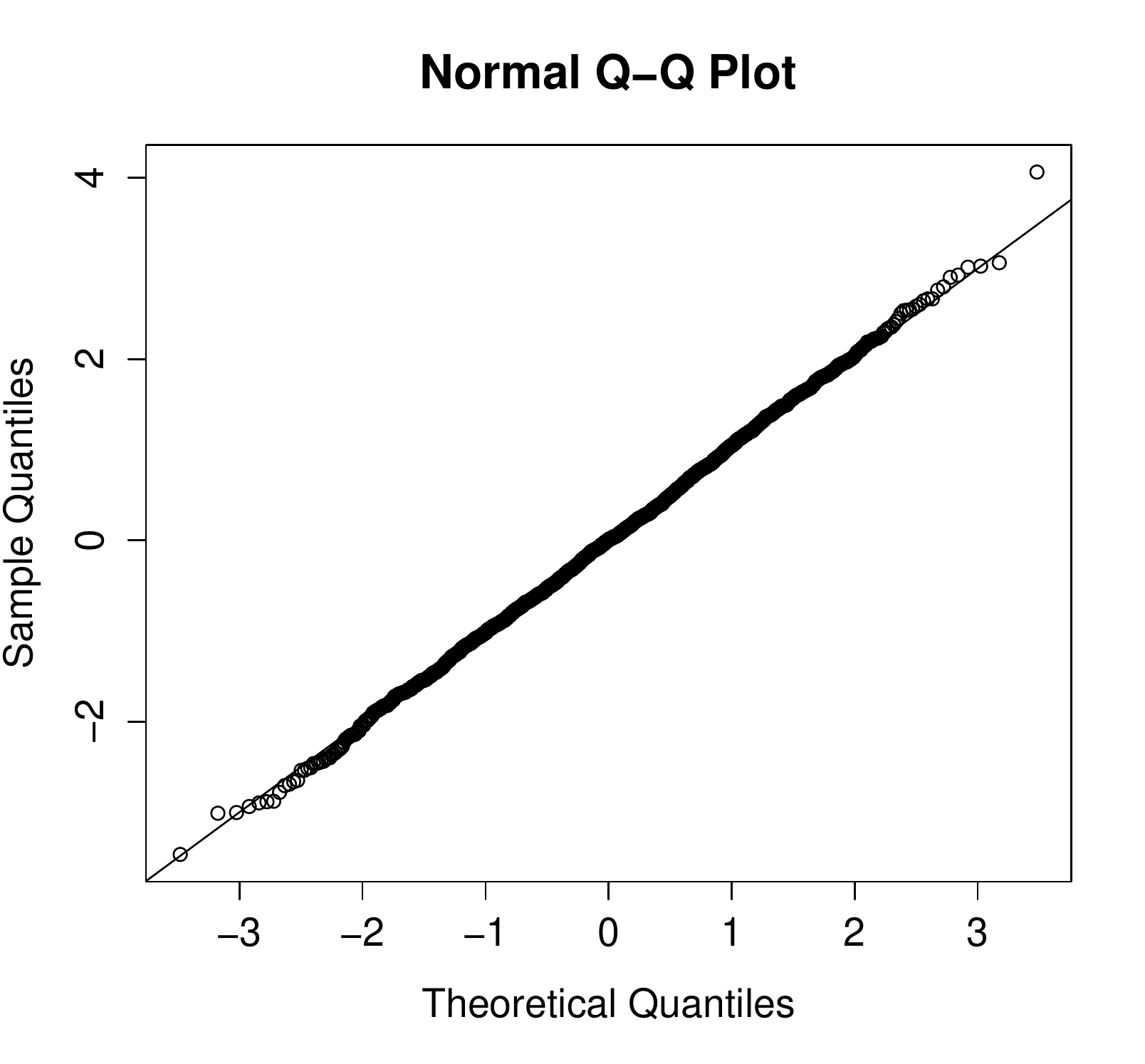} 
    \caption{$\omega=10, T=50$} 
    %\vspace{4ex}
  \end{subfigure}
 \begin{subfigure}[b]{0.32\linewidth}
    \centering
    \includegraphics[width=\linewidth]{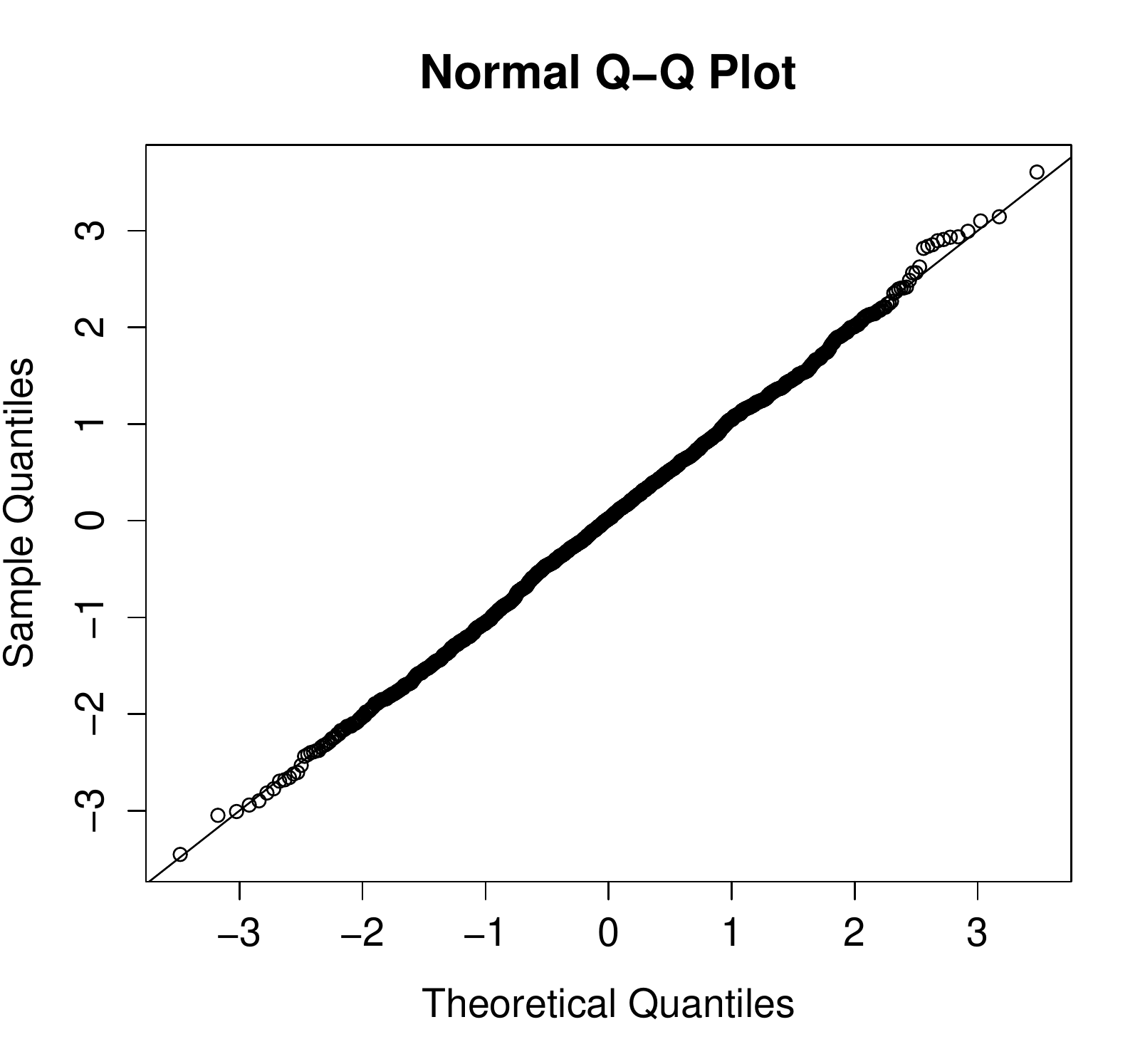} 
    \caption{$\omega=10, T=100$} 
    %\vspace{4ex}
  \end{subfigure}
  \caption{Normal QQ plots for the real part of the truncated Fourier transform of the simulated CARMA(2,1) processes driven by standard Brownian Motion for the frequencies $0, 0.1, 1 , 10$ (rows) and time horizons/maximum non-equidistant grid sizes $10/0.1, 50/0.05, 100/0.01$ (columns). The theoretical quantiles are coming from the (limiting) law described in Theorem \ref{thm:ApproxTFTDoubleLimitDistribution}. }\label{plot:QQCARMANormal} 
\end{figure}

\begin{figure}[tp]    
  \begin{subfigure}[b]{0.32\linewidth}
    \centering
    \includegraphics[width=\linewidth]{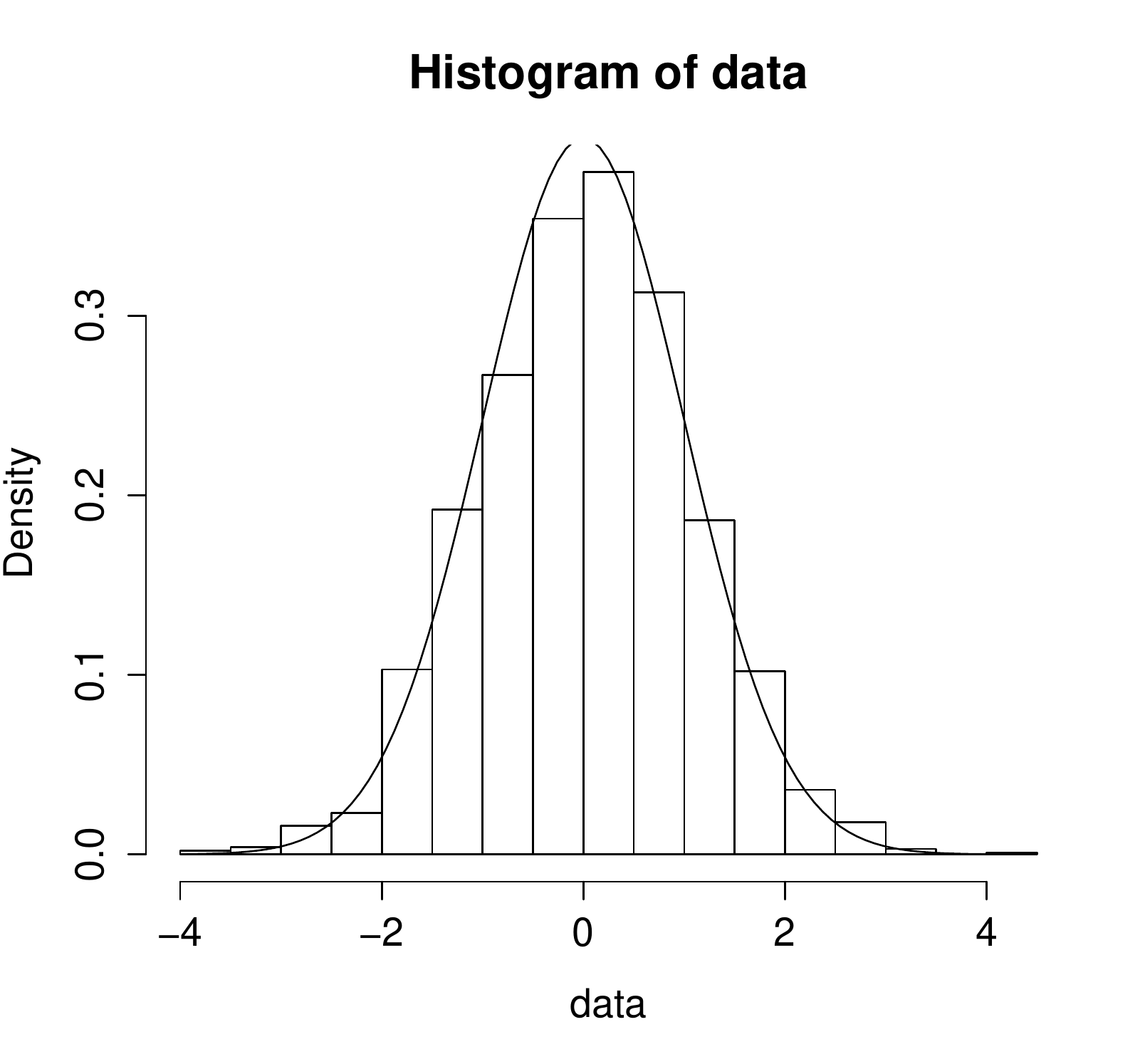} 
    \caption{ $\omega=0, T=10$} 
  
    %\vspace{4ex}
  \end{subfigure}%% 
 \begin{subfigure}[b]{0.32\linewidth}
    \centering
    \includegraphics[width=\linewidth]{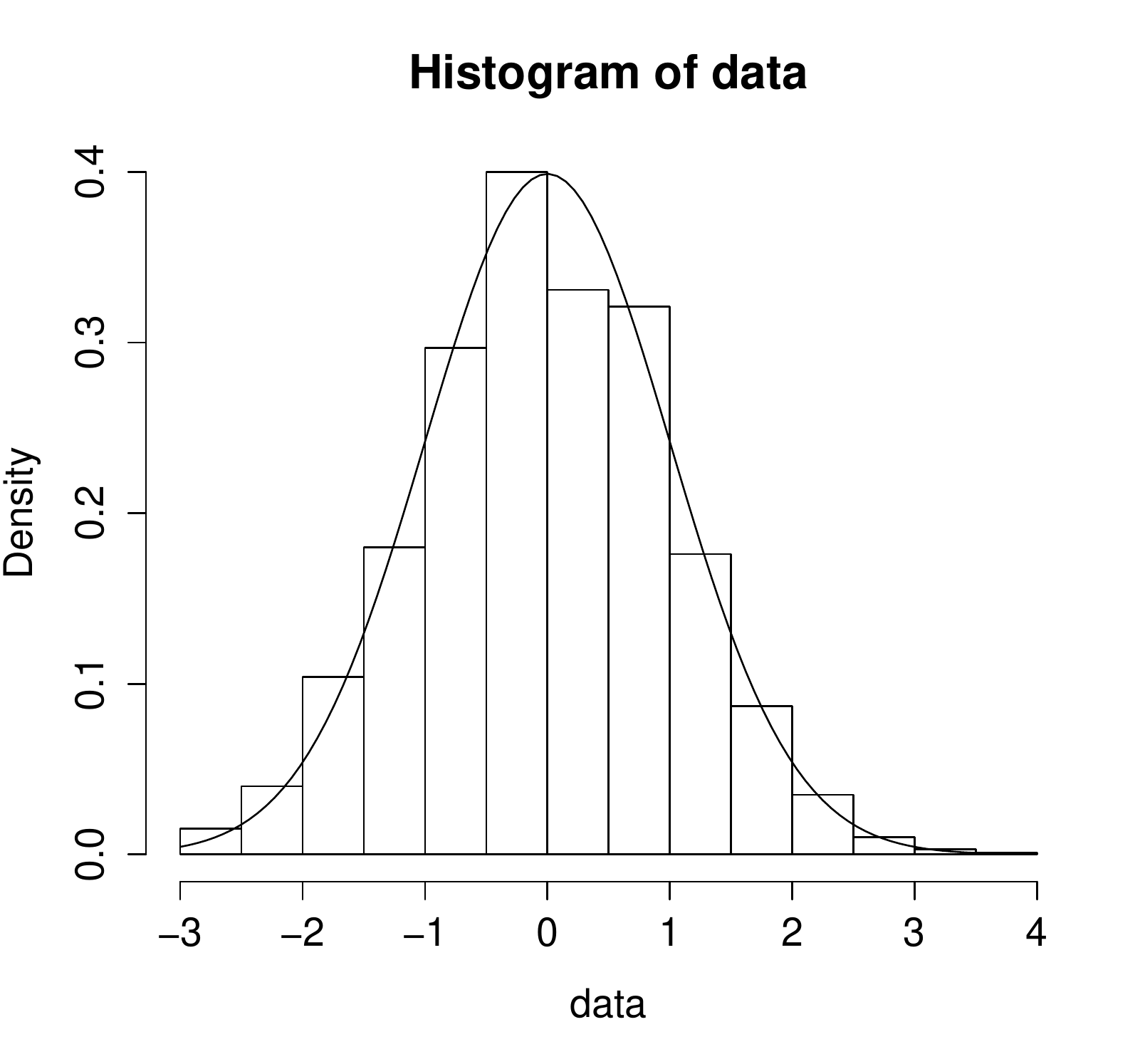} 
    \caption{$\omega=0, T=50$} 
    %\vspace{4ex}
  \end{subfigure}
 \begin{subfigure}[b]{0.32\linewidth}
    \centering
    \includegraphics[width=\linewidth]{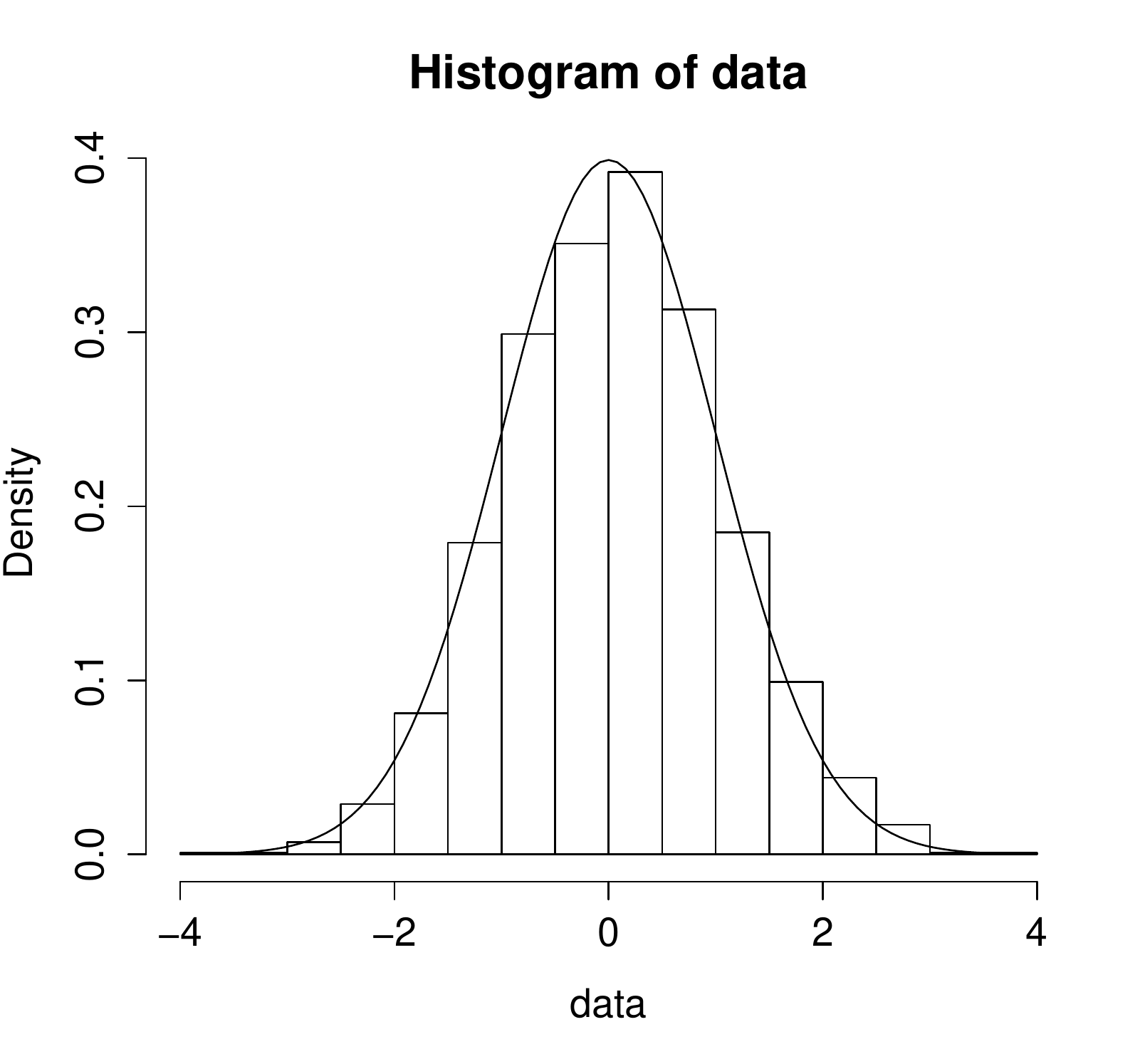} 
    \caption{$\omega=0, T=100$} 
    %\vspace{4ex}
  \end{subfigure}
  \begin{subfigure}[b]{0.32\linewidth}
    \centering
    \includegraphics[width=\linewidth]{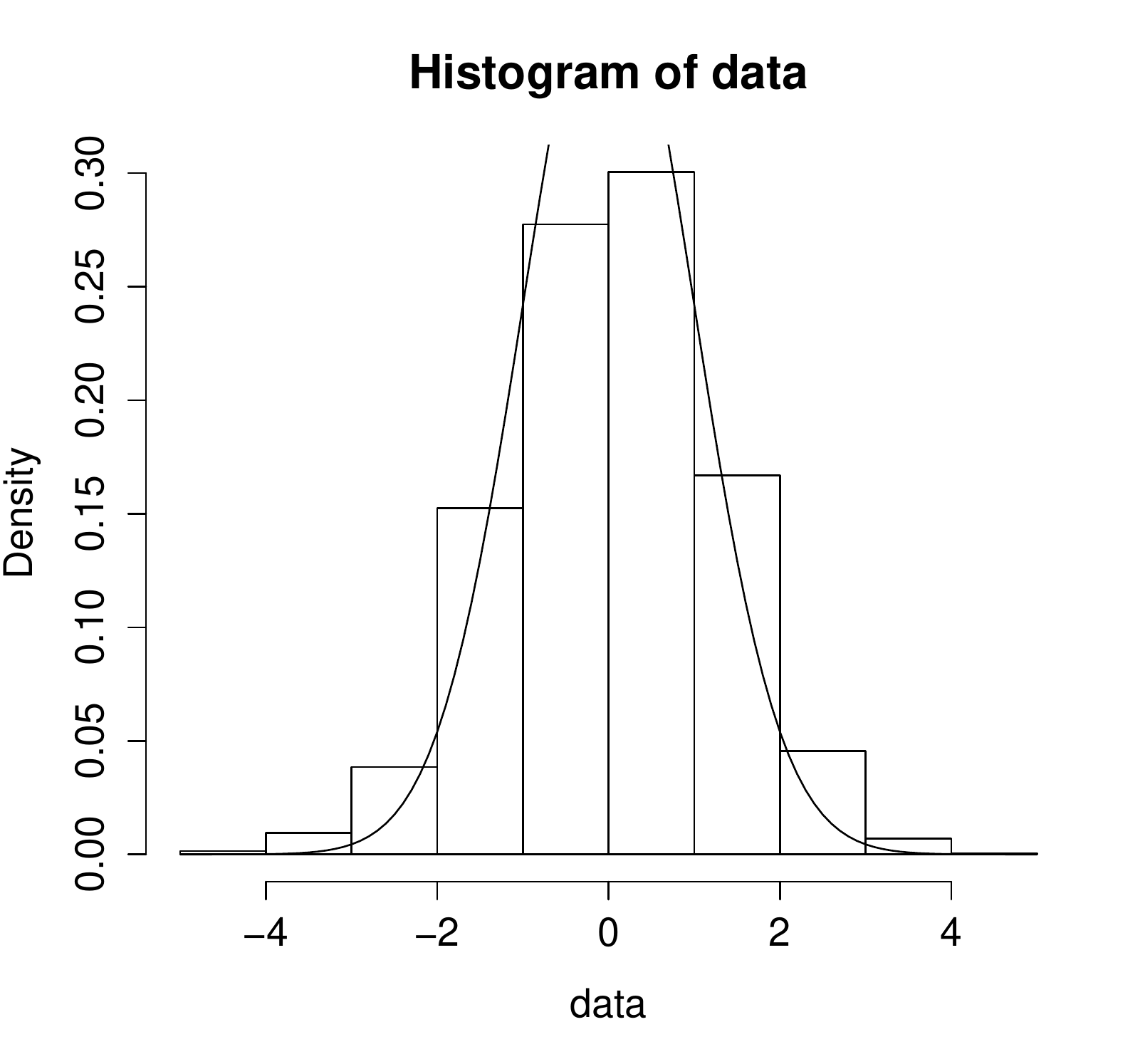} 
    \caption{$\omega=0.1, T=10$} 
    %\vspace{4ex}
  \end{subfigure}%% 
 \begin{subfigure}[b]{0.32\linewidth}
    \centering
    \includegraphics[width=\linewidth]{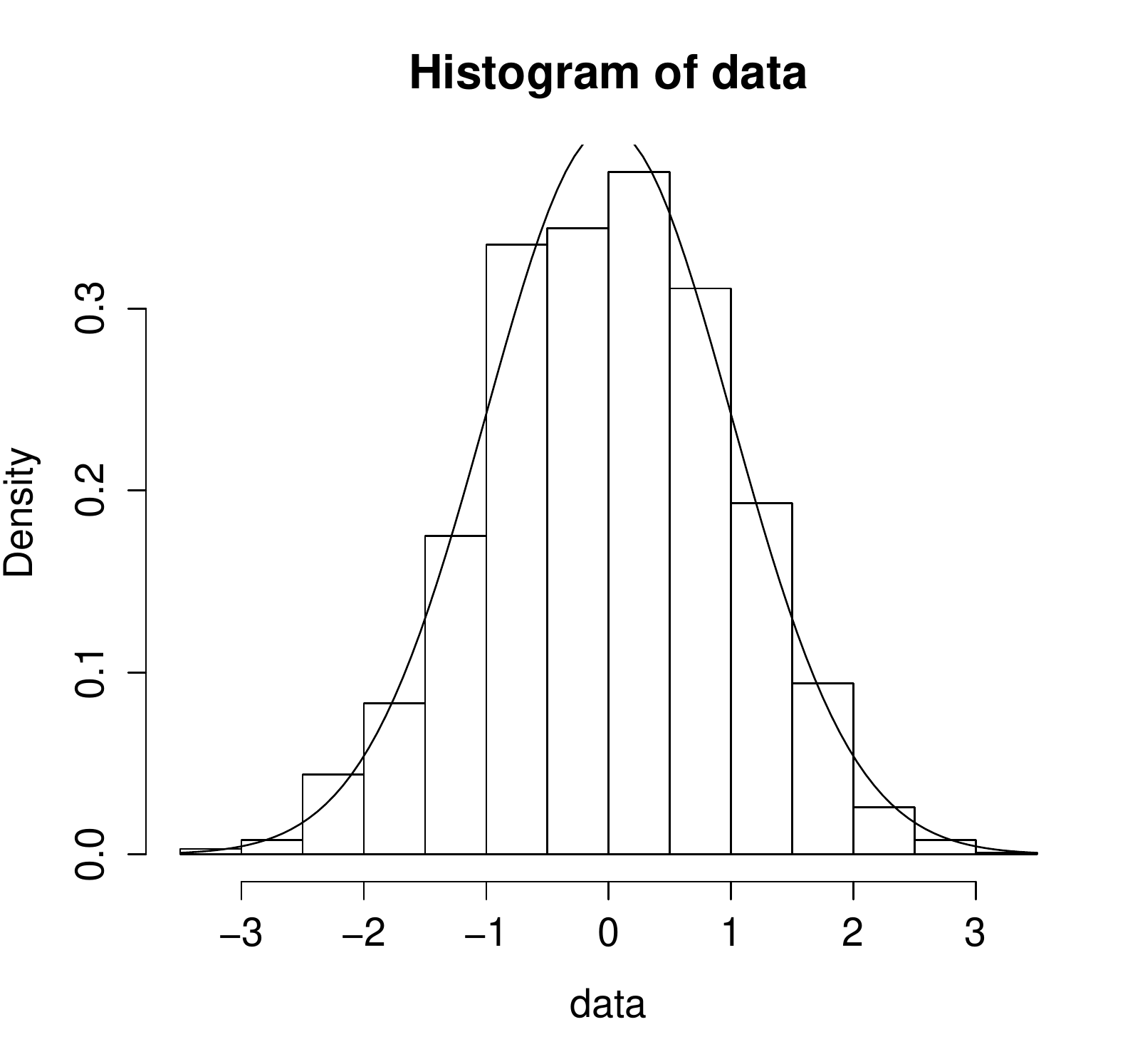} 
    \caption{$\omega=0.1, T=50$} 
    %\vspace{4ex}
  \end{subfigure}
 \begin{subfigure}[b]{0.32\linewidth}
    \centering
    \includegraphics[width=\linewidth]{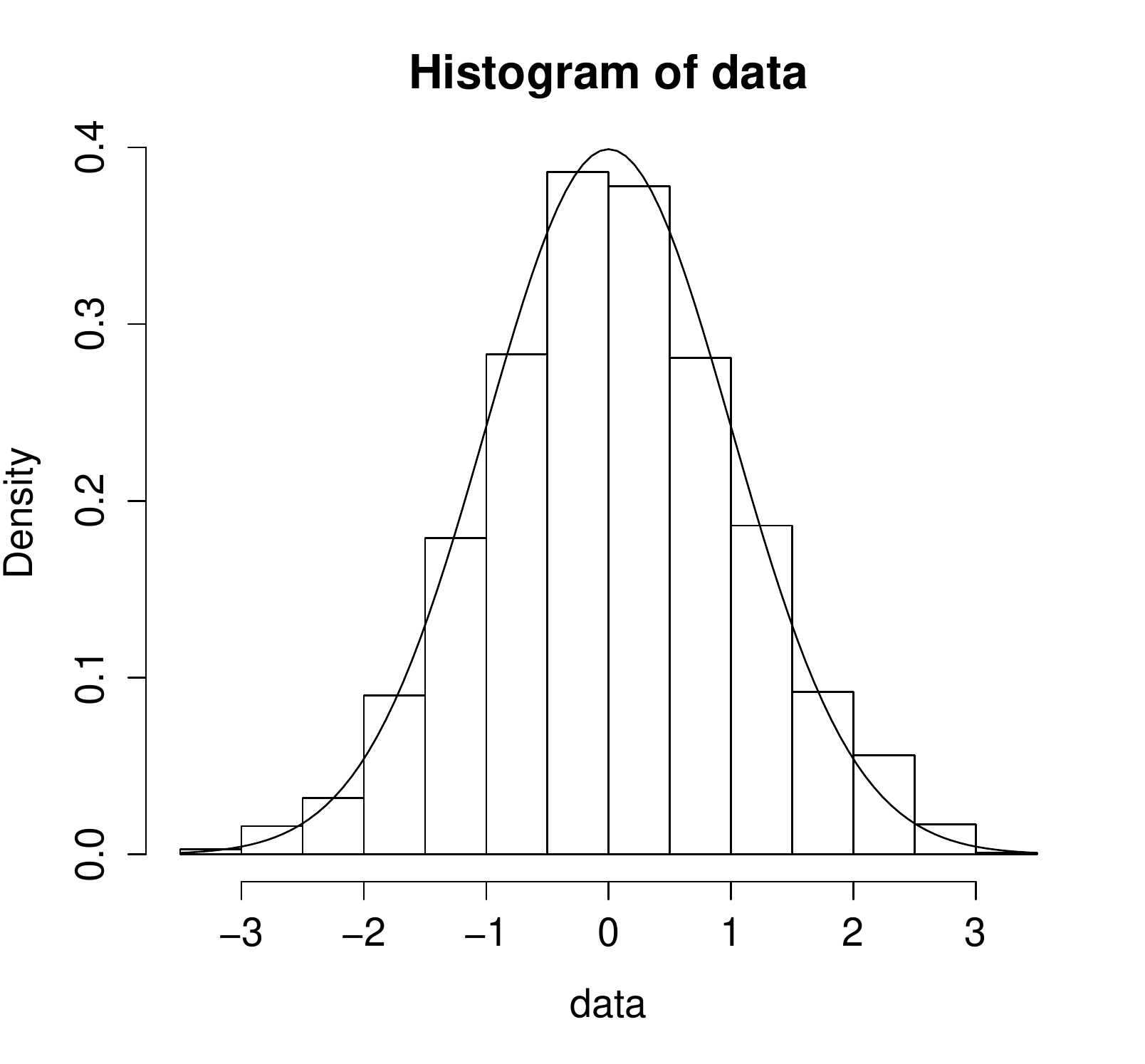} 
    \caption{$\omega=0.1, T=100$} 
    %\vspace{4ex}
  \end{subfigure}
  \begin{subfigure}[b]{0.32\linewidth}
    \centering
    \includegraphics[width=\linewidth]{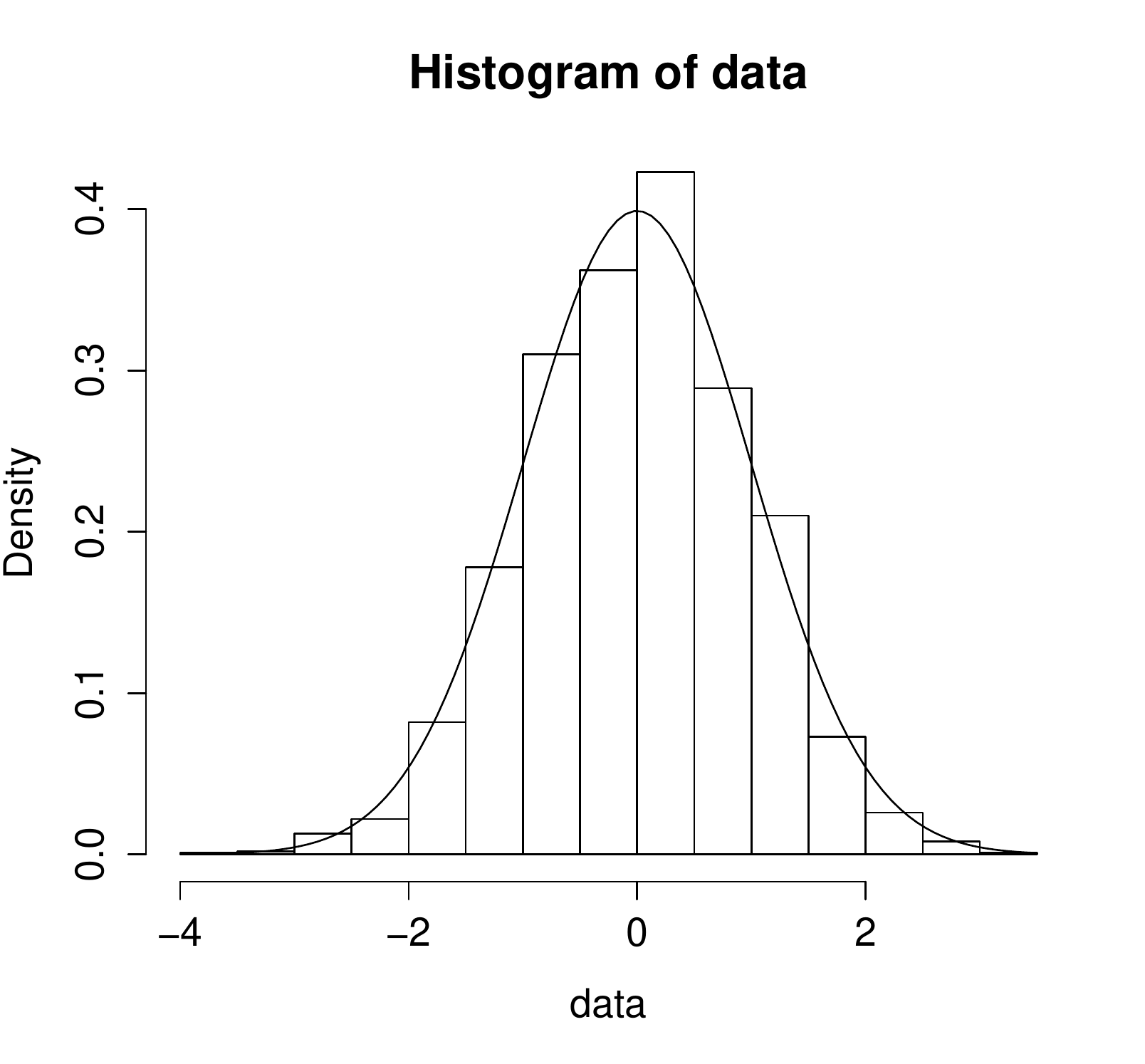} 
    \caption{$\omega=1, T=10$} 
    %\vspace{4ex}
  \end{subfigure}%% 
 \begin{subfigure}[b]{0.32\linewidth}
    \centering
    \includegraphics[width=\linewidth]{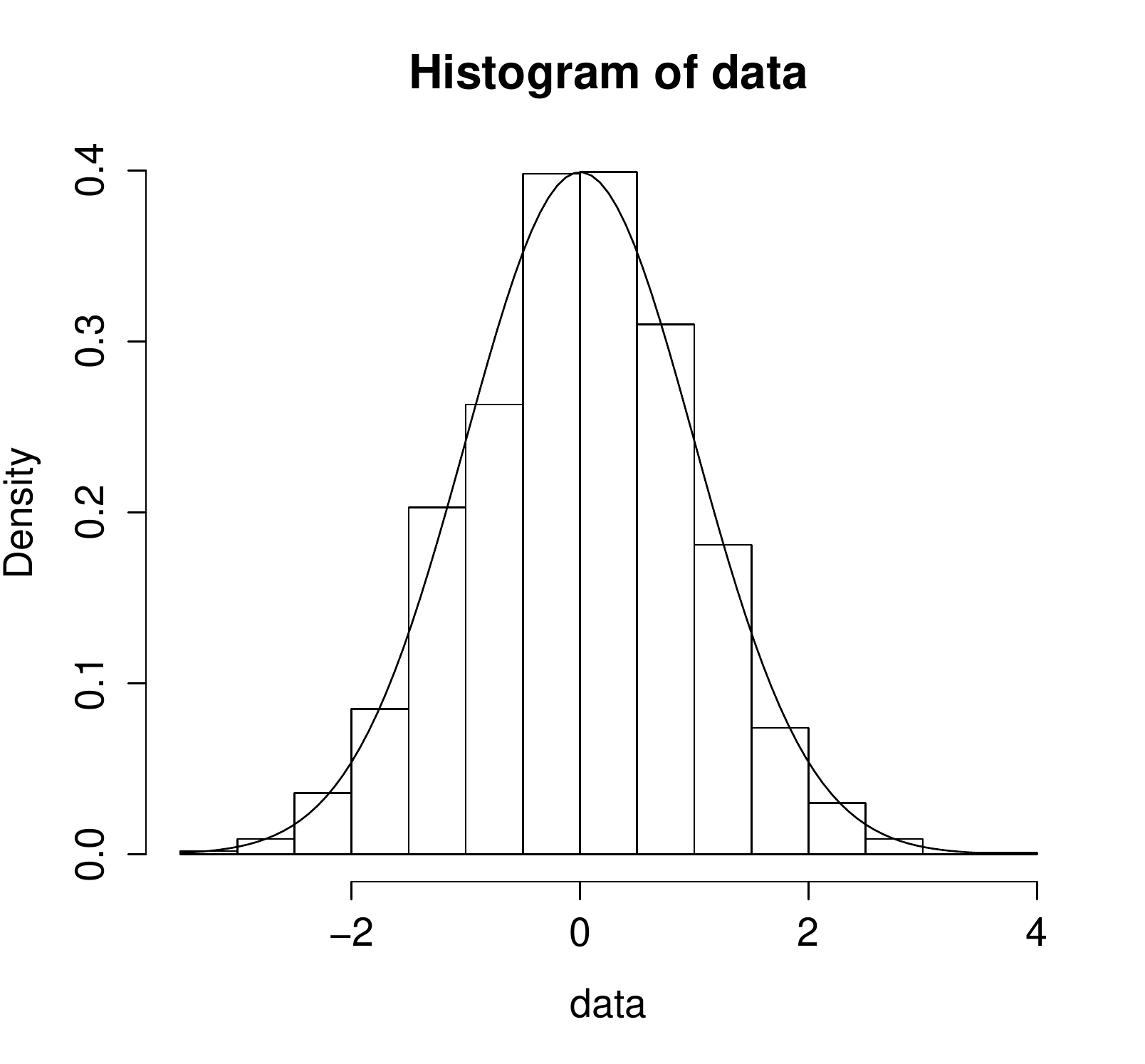} 
    \caption{$\omega=1, T=50$} 
    %\vspace{4ex}
  \end{subfigure}
 \begin{subfigure}[b]{0.32\linewidth}
    \centering
    \includegraphics[width=\linewidth]{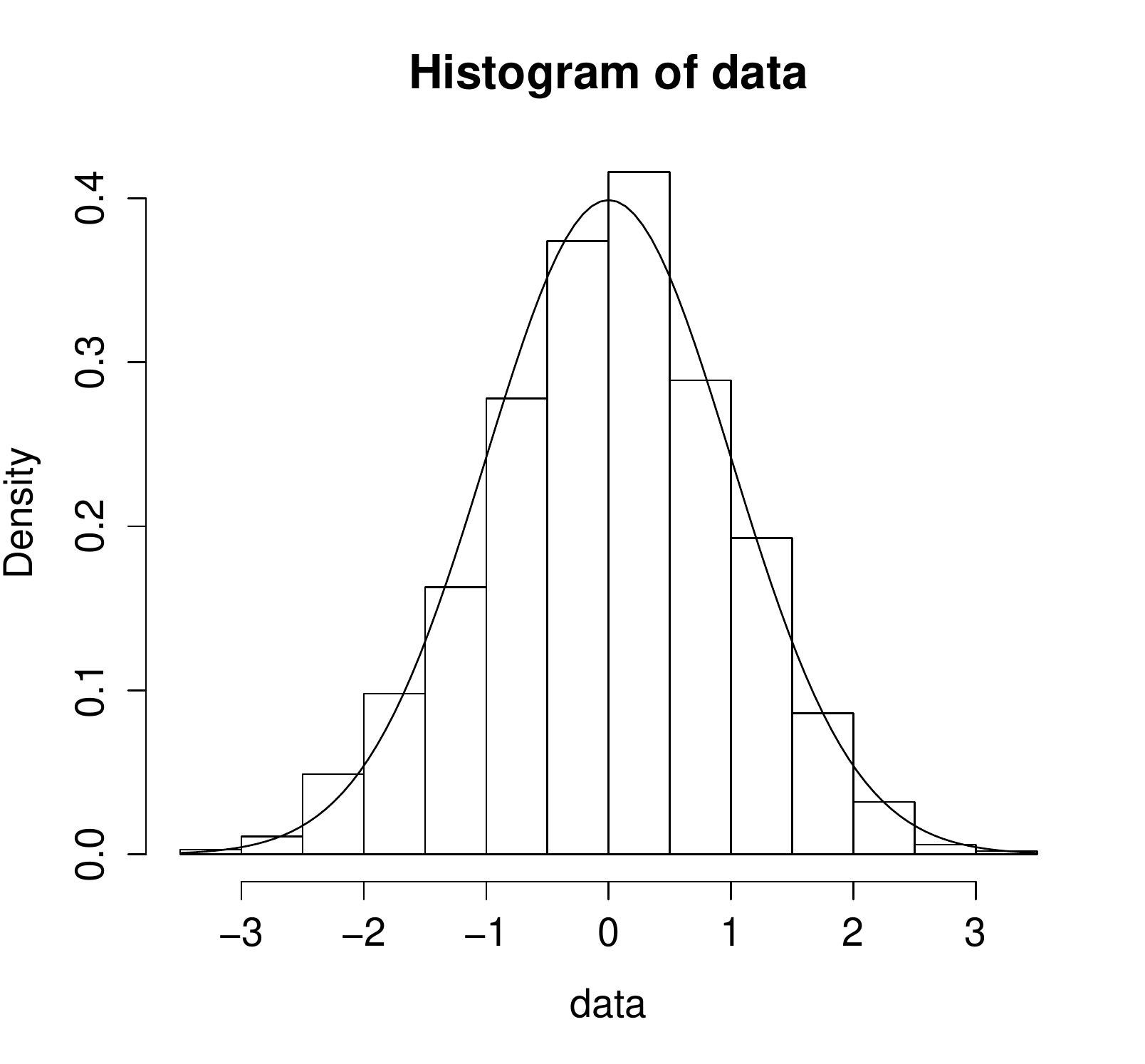} 
    \caption{$\omega=1, T=100$} 
    %\vspace{4ex}
  \end{subfigure}
  \begin{subfigure}[b]{0.32\linewidth}
    \centering
    \includegraphics[width=\linewidth]{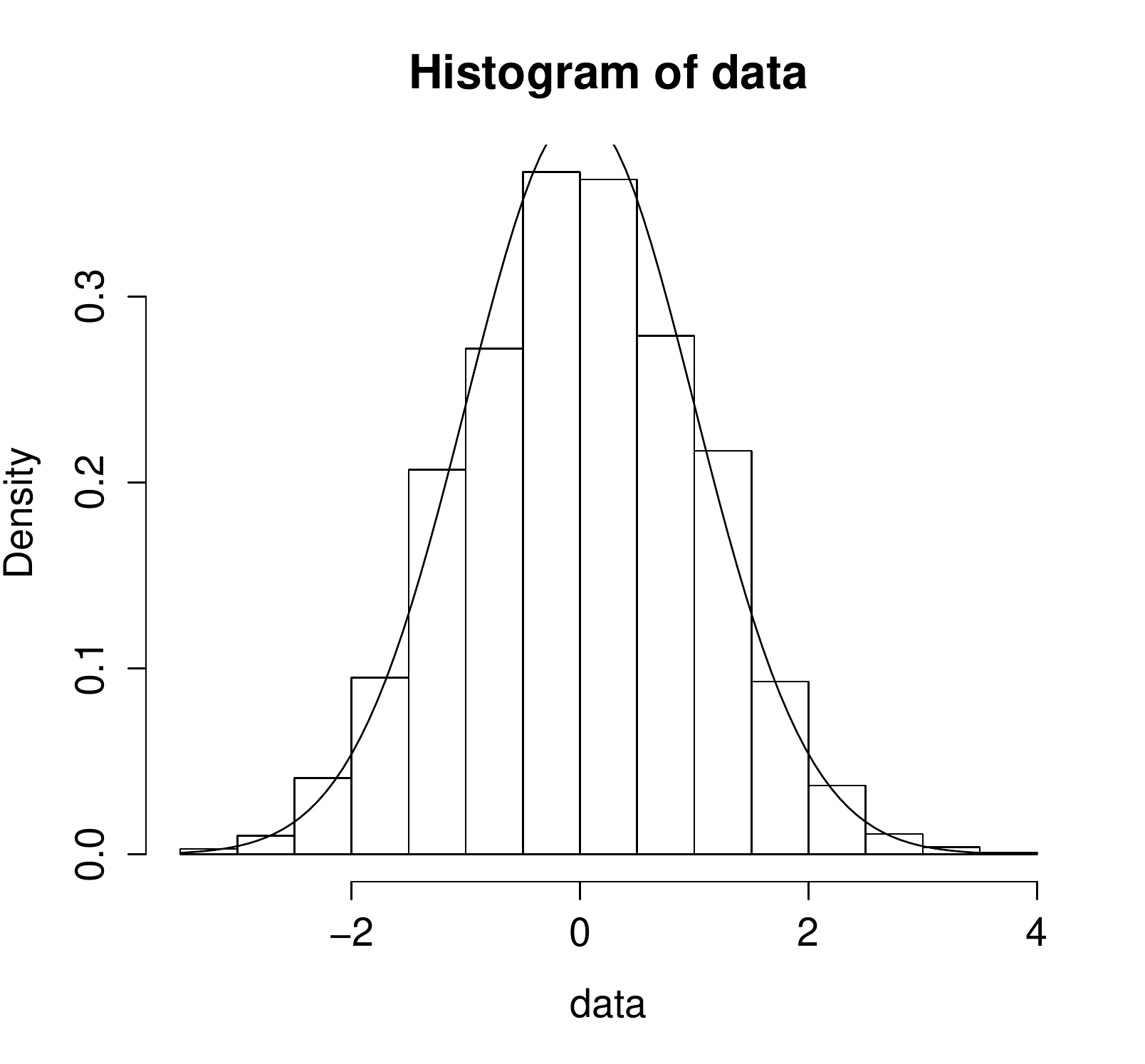} 
    \caption{$\omega=10, T=10$}  
    %\vspace{4ex}
  \end{subfigure}%% 
 \begin{subfigure}[b]{0.32\linewidth}
    \centering
    \includegraphics[width=\linewidth]{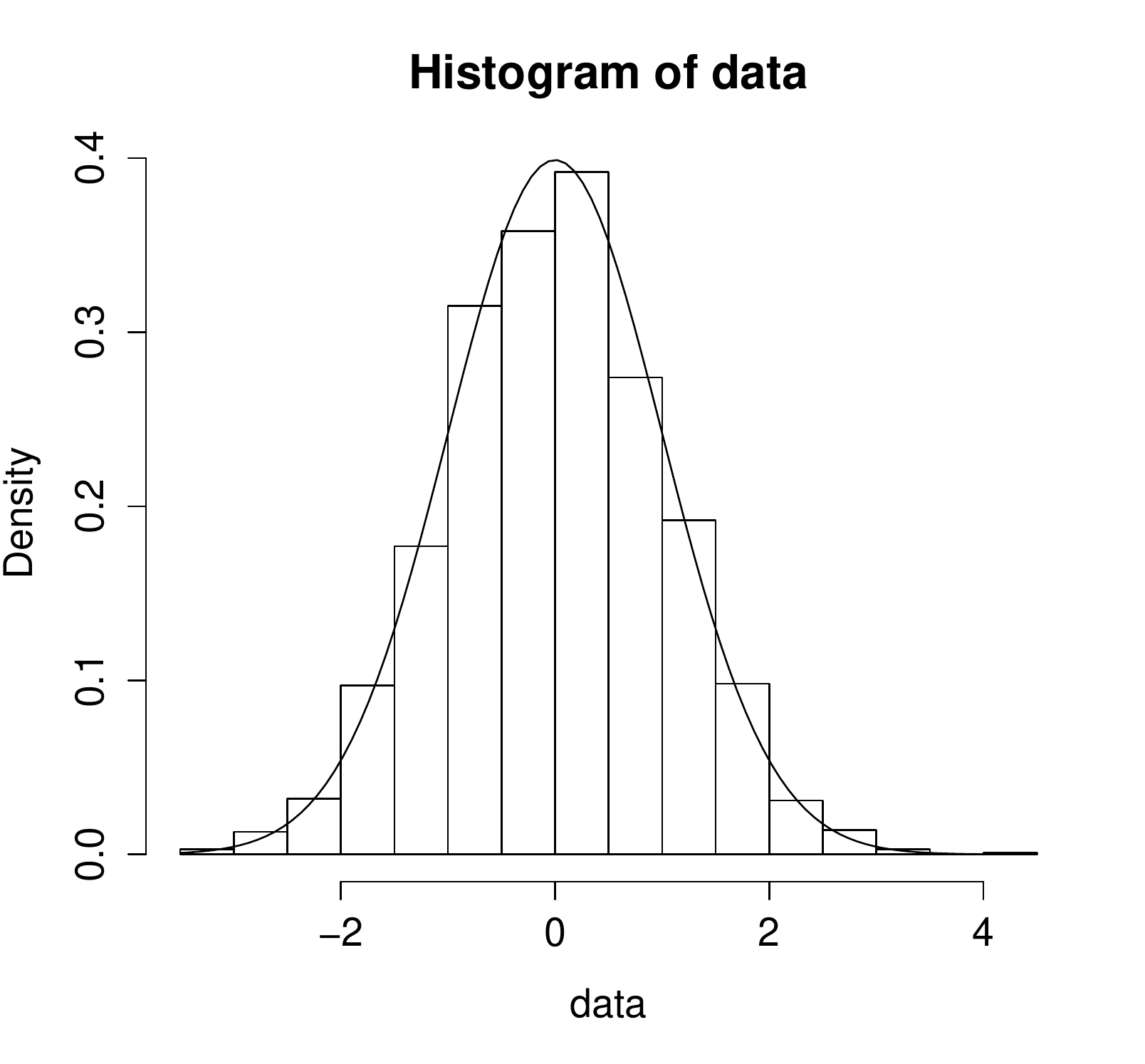} 
    \caption{$\omega=10, T=50$} 
    %\vspace{4ex}
  \end{subfigure}
 \begin{subfigure}[b]{0.32\linewidth}
    \centering
    \includegraphics[width=\linewidth]{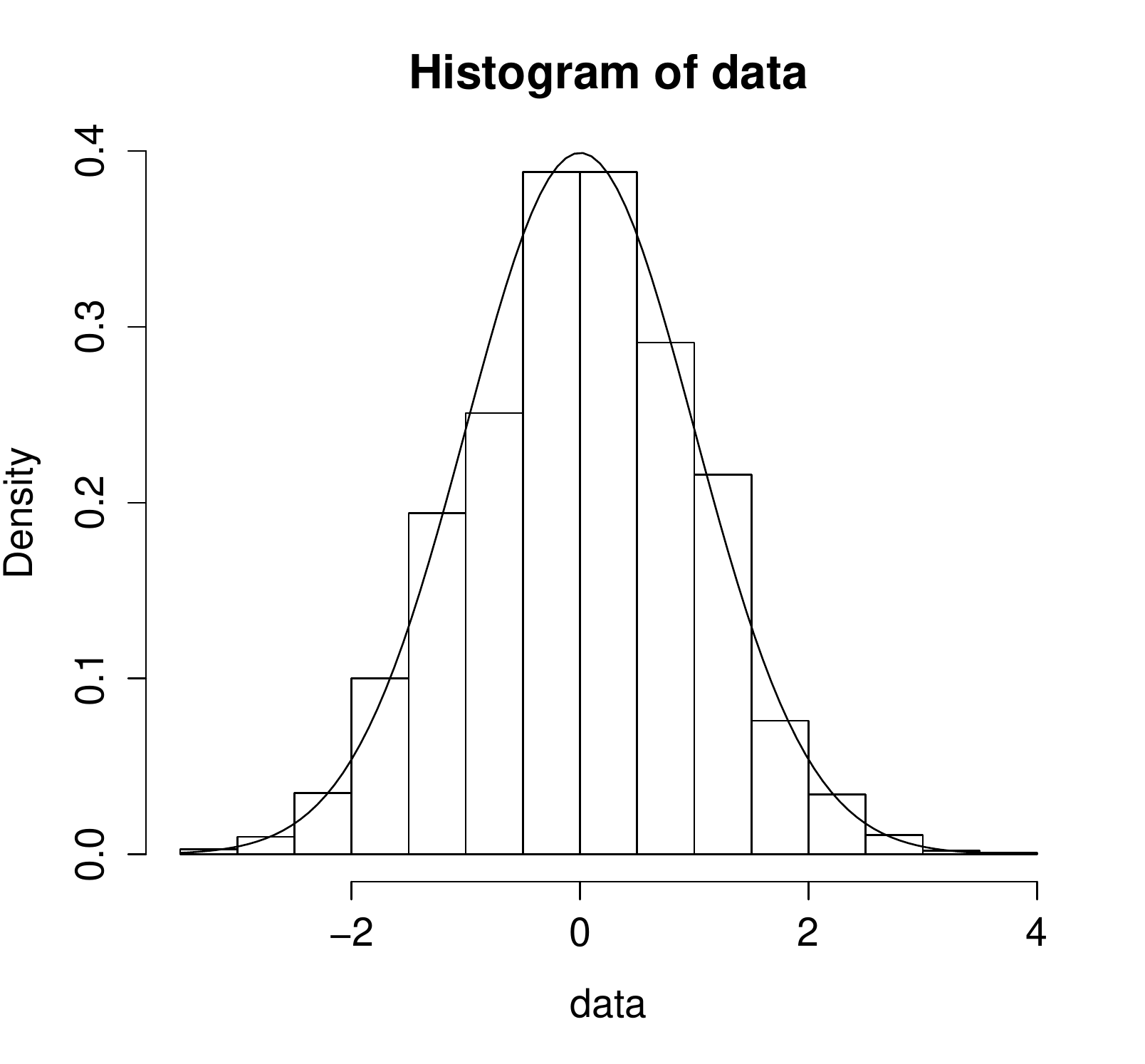} 
    \caption{$\omega=10, T=100$} 
    %\vspace{4ex}
  \end{subfigure}
  \caption{Histograms and limiting density for the real part of the truncated Fourier transform of the simulated CARMA(2,1) processes driven by standard Brownian Motion for the frequencies $0, 0.1, 1 , 10$ (rows) and time horizons/maximum non-equidistant grid sizes $10/0.1, 50/0.05, 100/0.01$ (columns)}\label{plot:HistCARMANormal} 
\end{figure}

\begin{figure}[tp]    
  \begin{subfigure}[b]{0.32\linewidth}
    \centering
    \includegraphics[width=\linewidth]{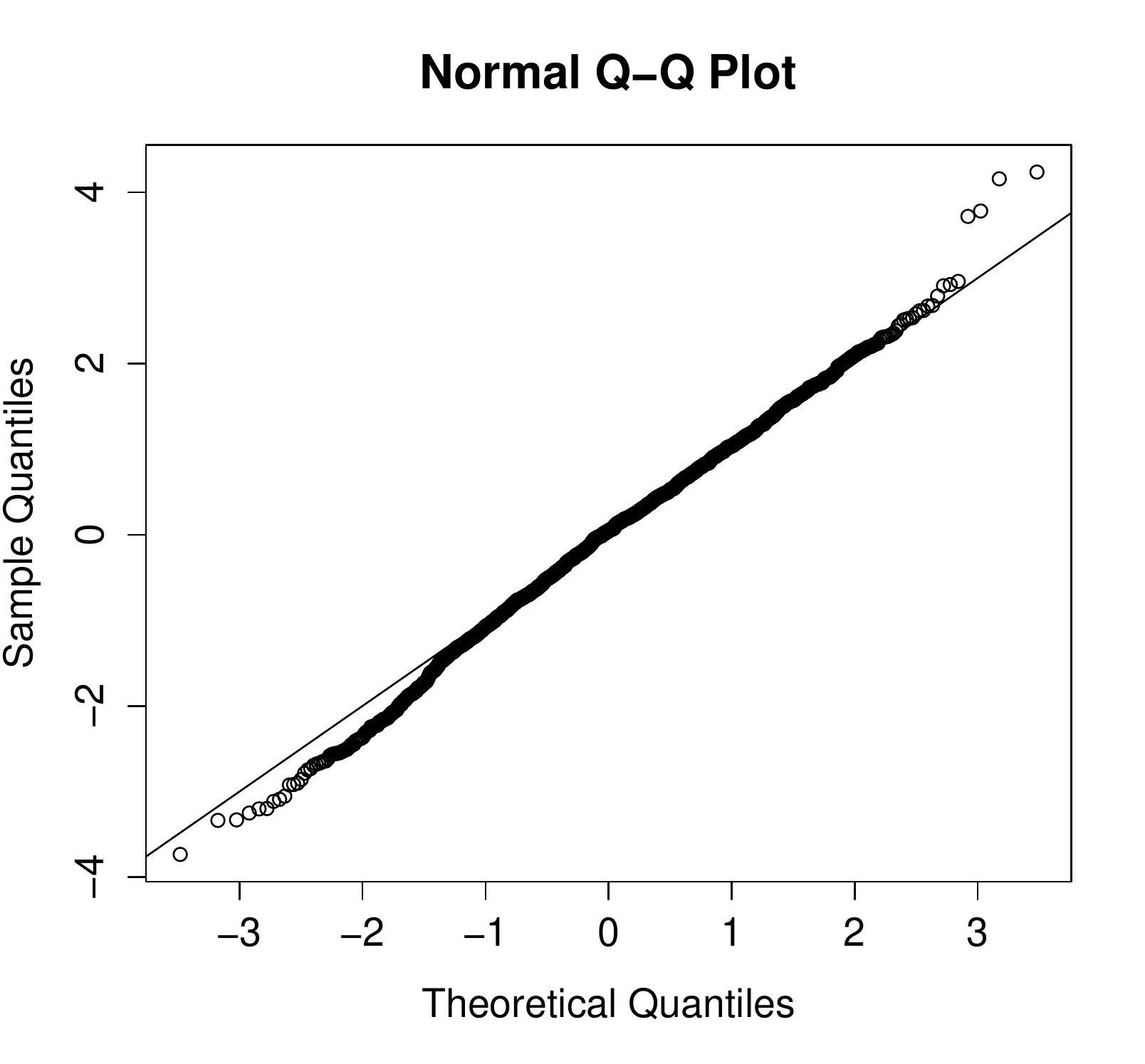} 
    \caption{ $\omega=0, T=10$} 
  
    %\vspace{4ex}
  \end{subfigure}%% 
 \begin{subfigure}[b]{0.32\linewidth}
    \centering
    \includegraphics[width=\linewidth]{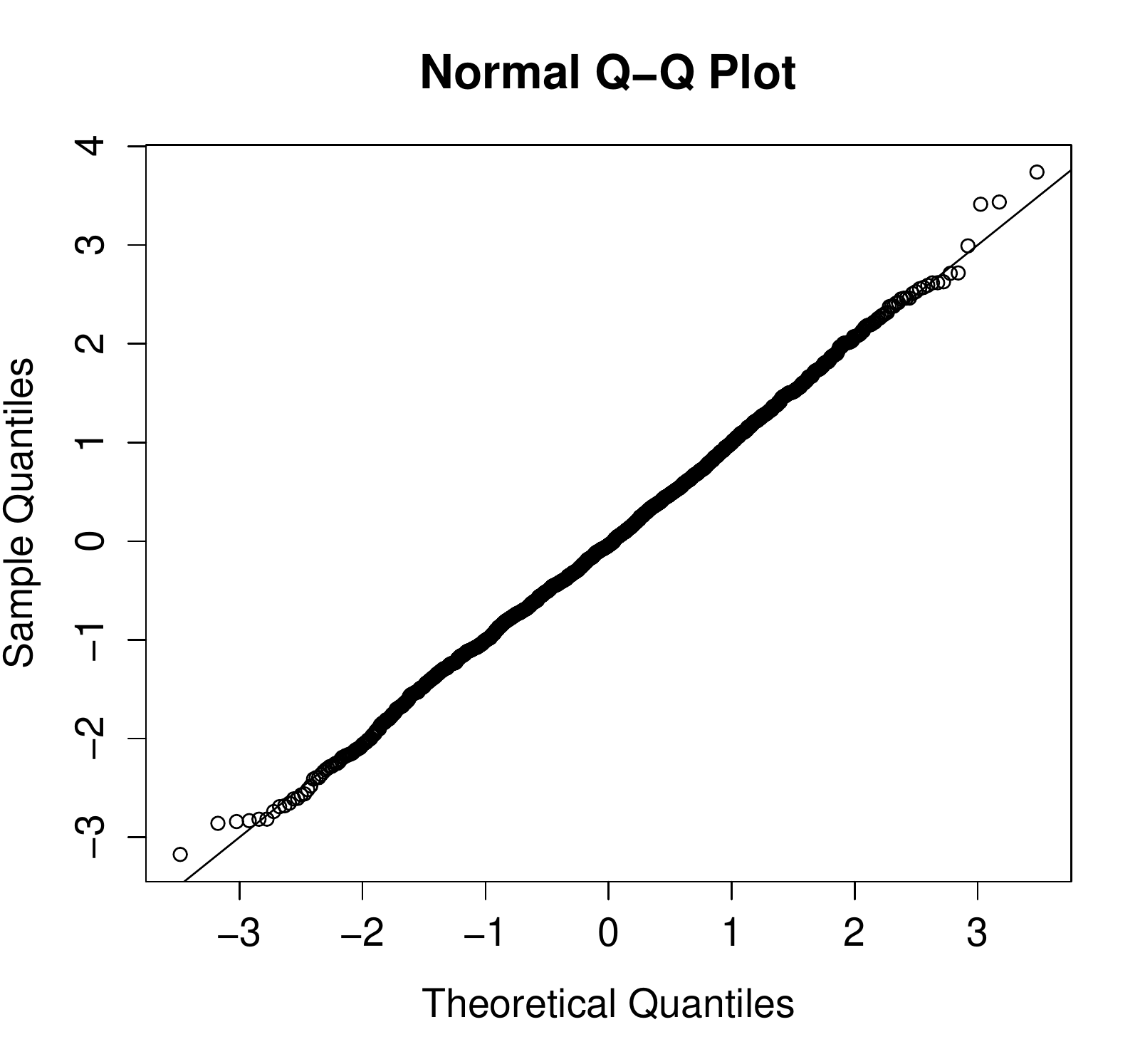} 
    \caption{$\omega=0, T=50$} 
    %\vspace{4ex}
  \end{subfigure}
 \begin{subfigure}[b]{0.32\linewidth}
    \centering
    \includegraphics[width=\linewidth]{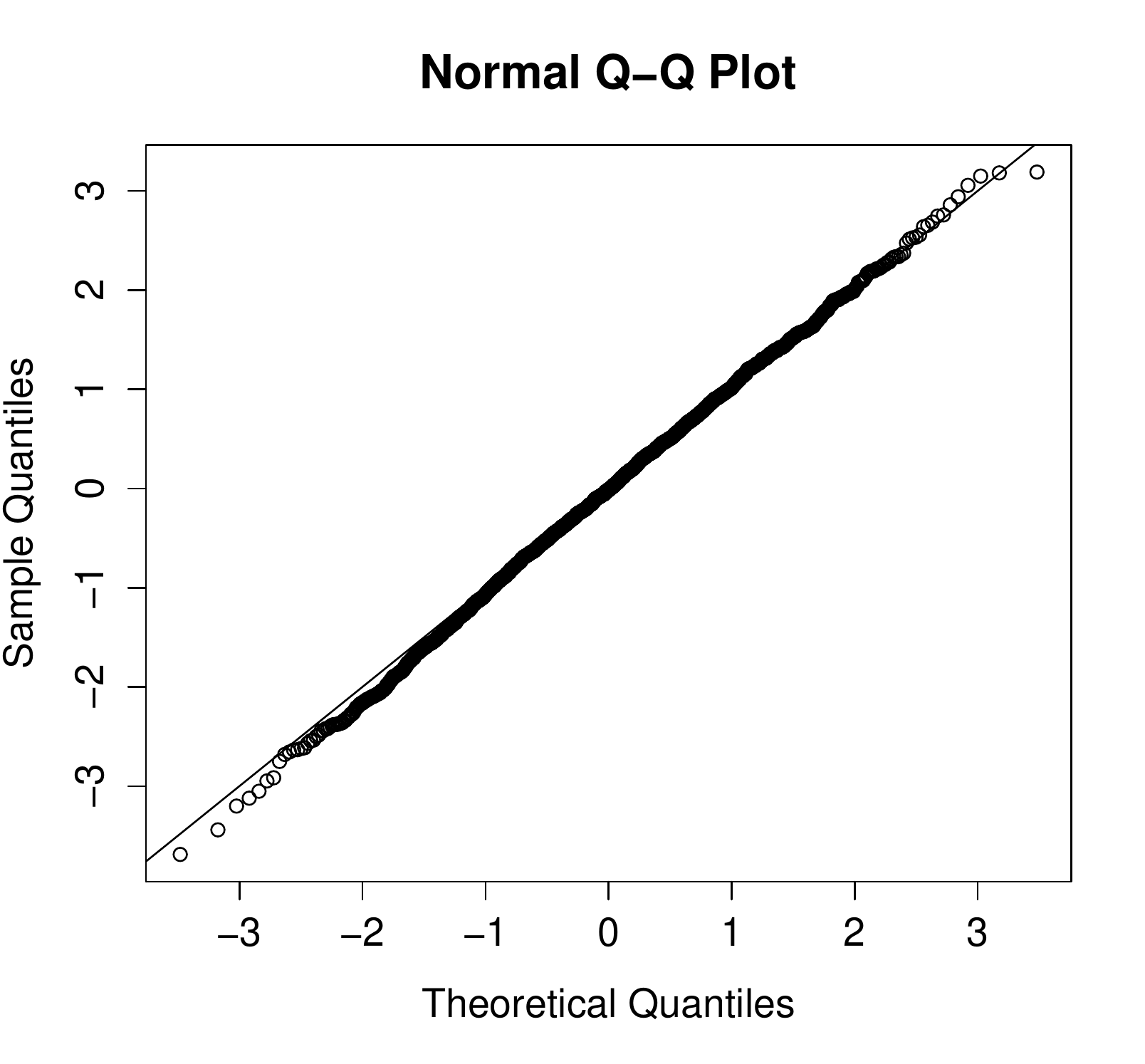} 
    \caption{$\omega=0, T=100$} 
    %\vspace{4ex}
  \end{subfigure}
  \begin{subfigure}[b]{0.32\linewidth}
    \centering
    \includegraphics[width=\linewidth]{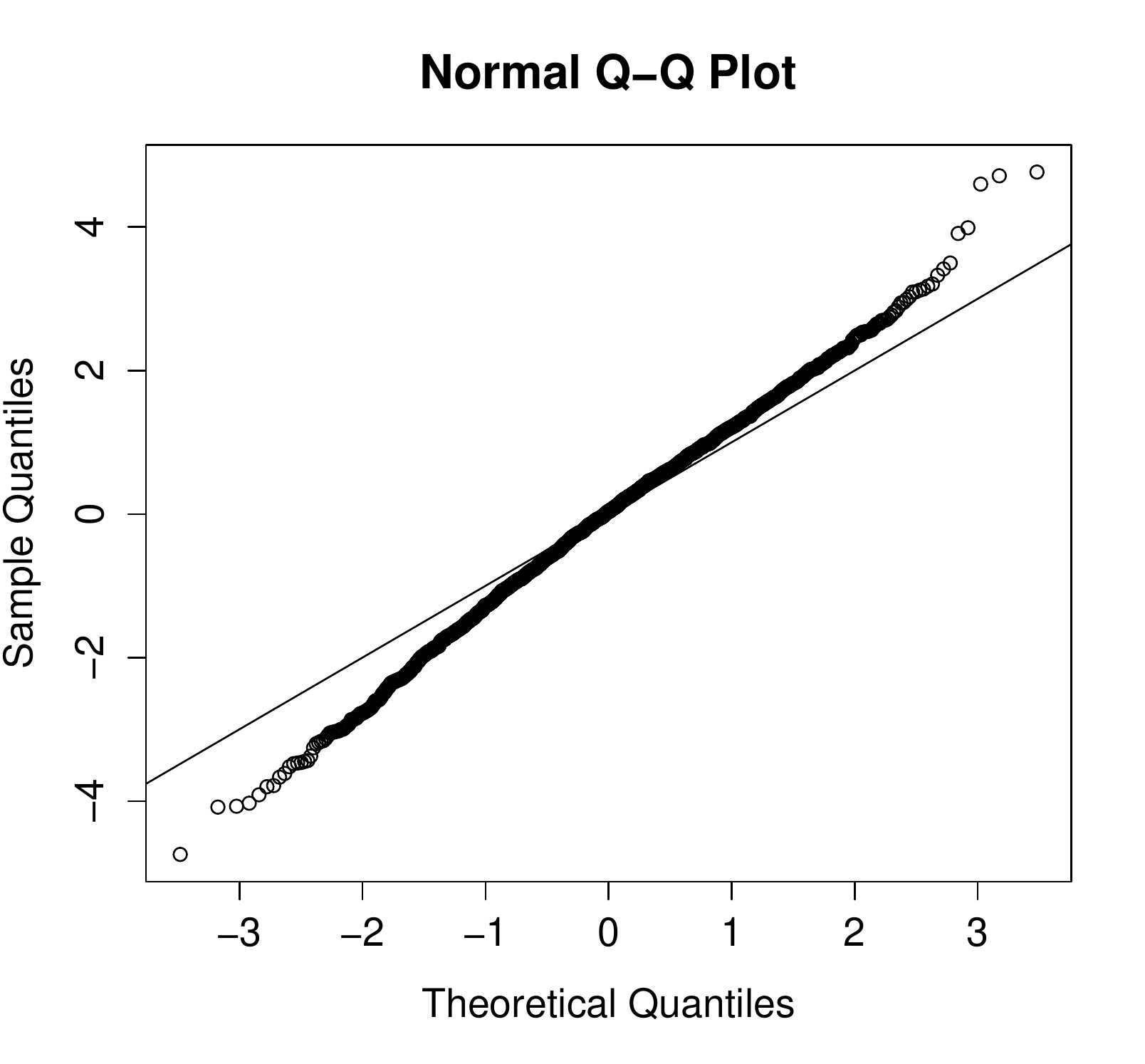} 
    \caption{$\omega=0.1, T=10$} 
    %\vspace{4ex}
  \end{subfigure}%% 
 \begin{subfigure}[b]{0.32\linewidth}
    \centering
    \includegraphics[width=\linewidth]{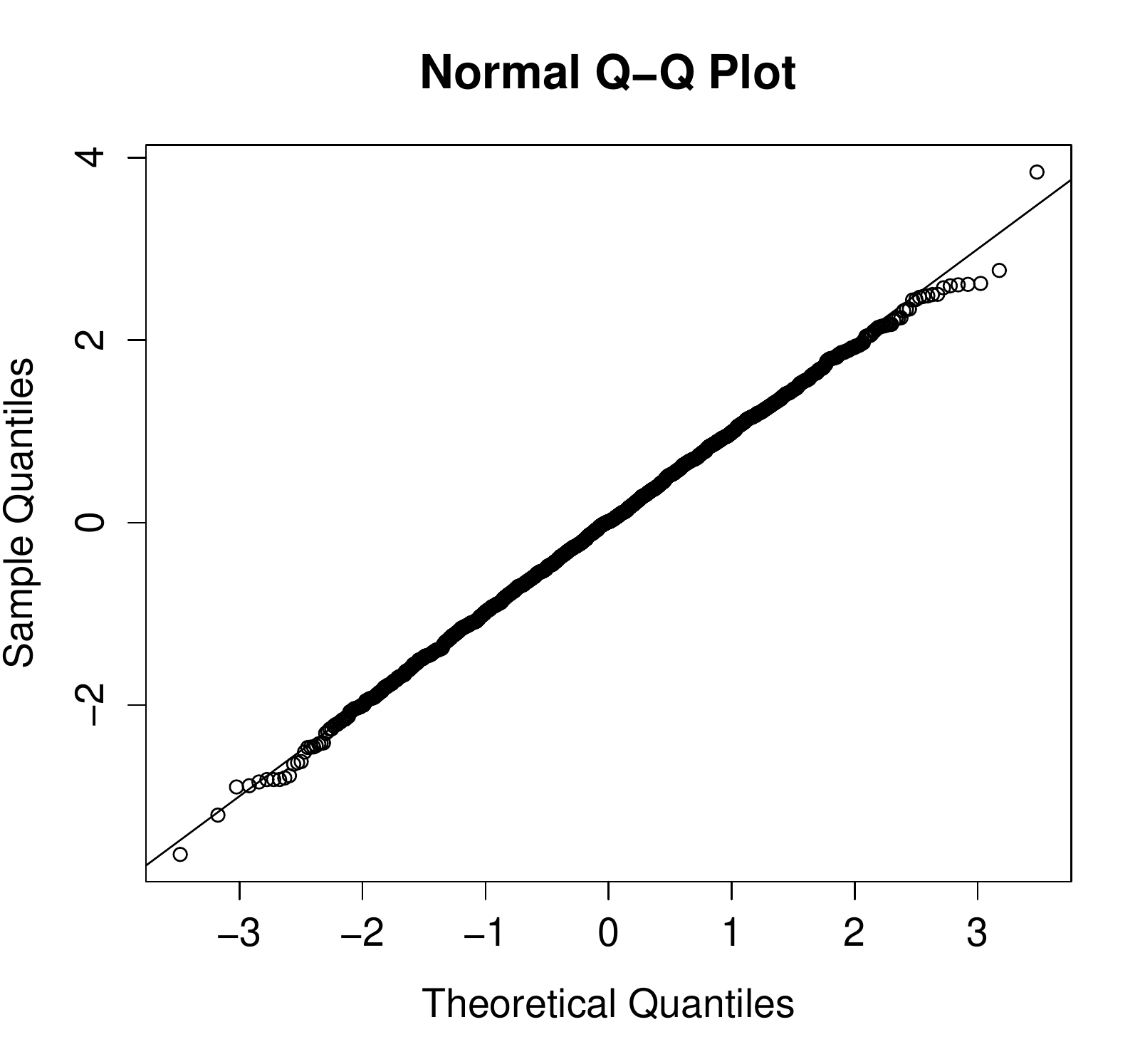} 
    \caption{$\omega=0.1, T=50$} 
    %\vspace{4ex}
  \end{subfigure}
 \begin{subfigure}[b]{0.32\linewidth}
    \centering
    \includegraphics[width=\linewidth]{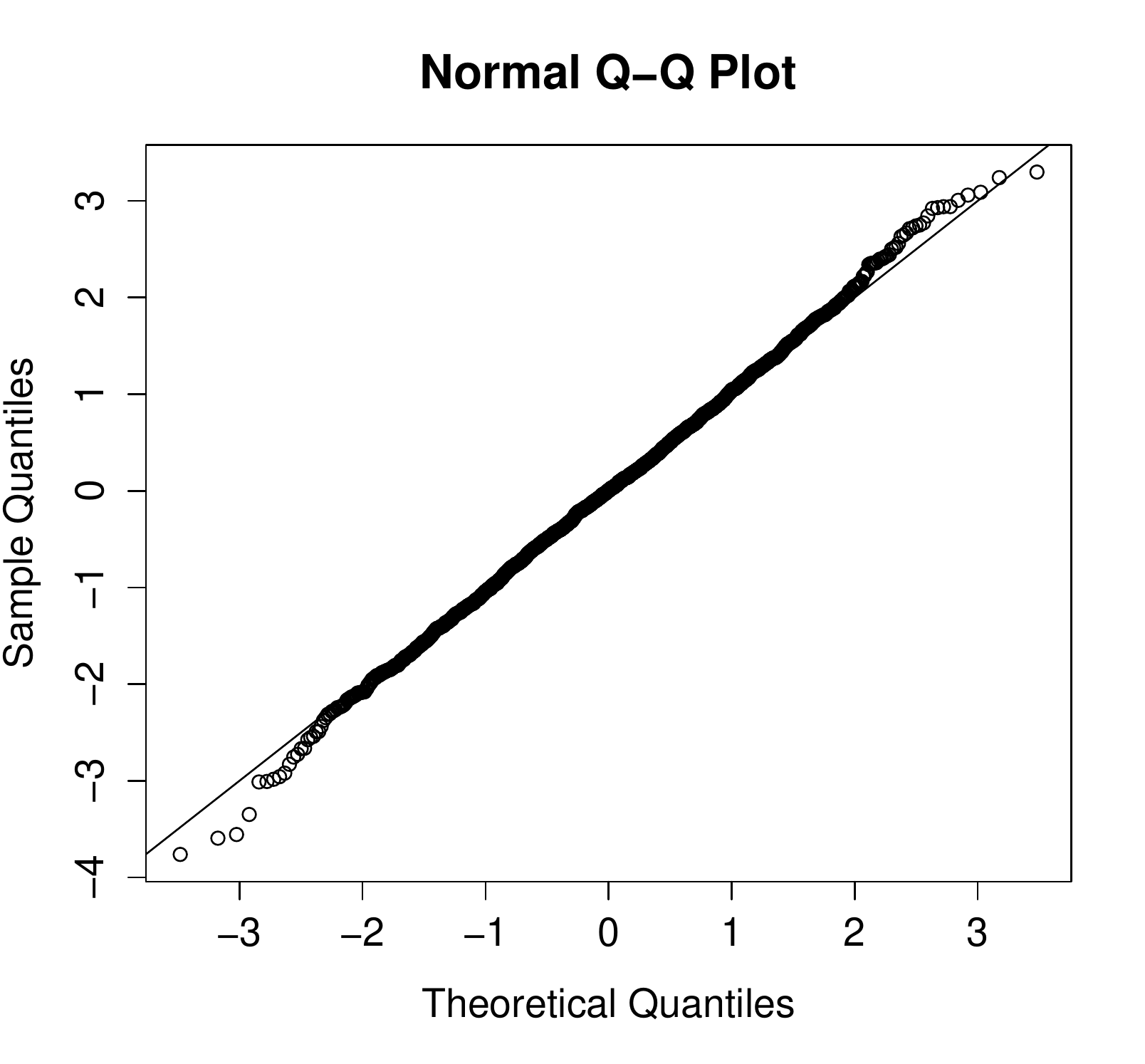} 
    \caption{$\omega=0.1, T=100$} 
    %\vspace{4ex}
  \end{subfigure}
  \begin{subfigure}[b]{0.32\linewidth}
    \centering
    \includegraphics[width=\linewidth]{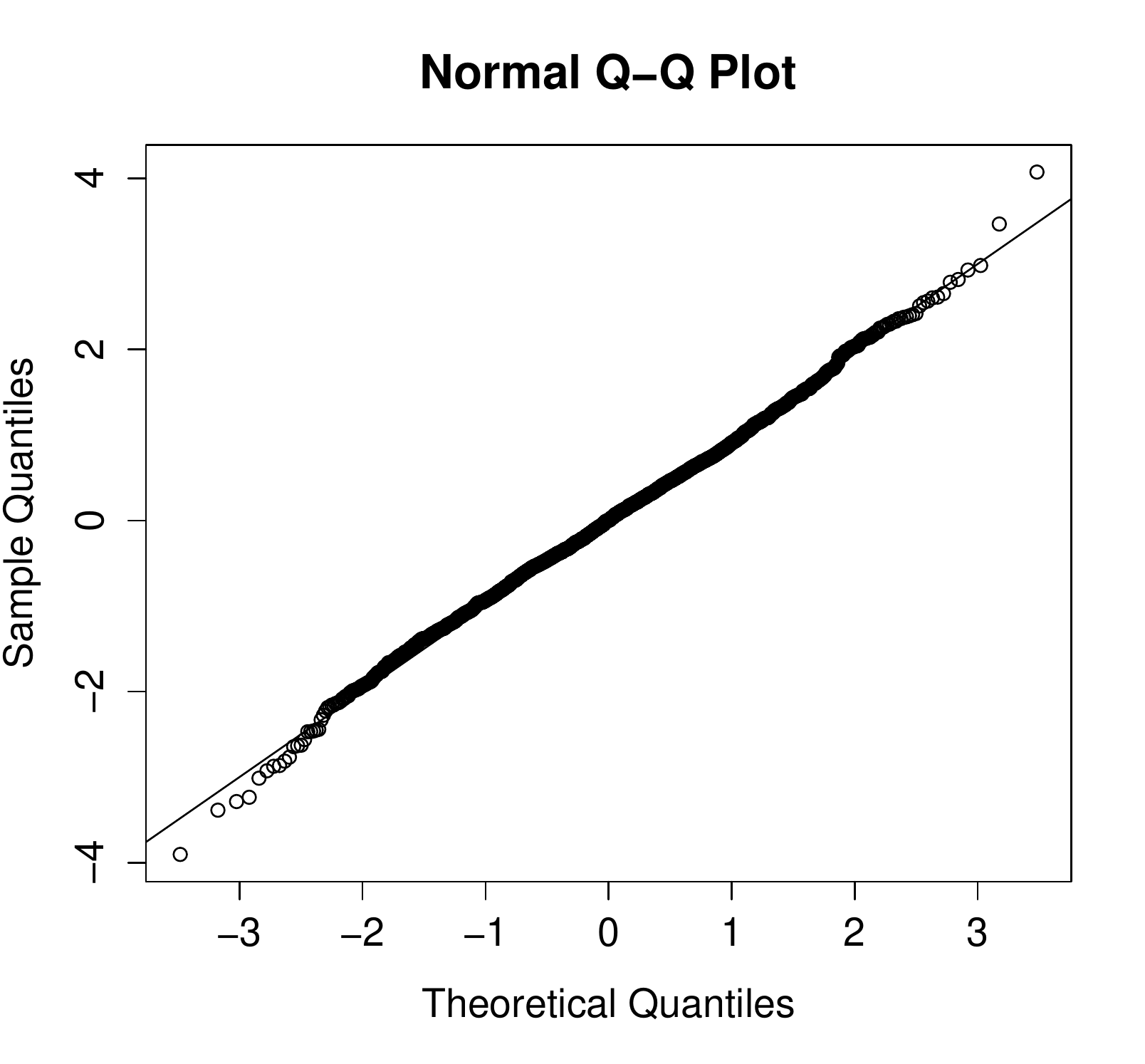} 
    \caption{$\omega=1, T=10$} 
    %\vspace{4ex}
  \end{subfigure}%% 
 \begin{subfigure}[b]{0.32\linewidth}
    \centering
    \includegraphics[width=\linewidth]{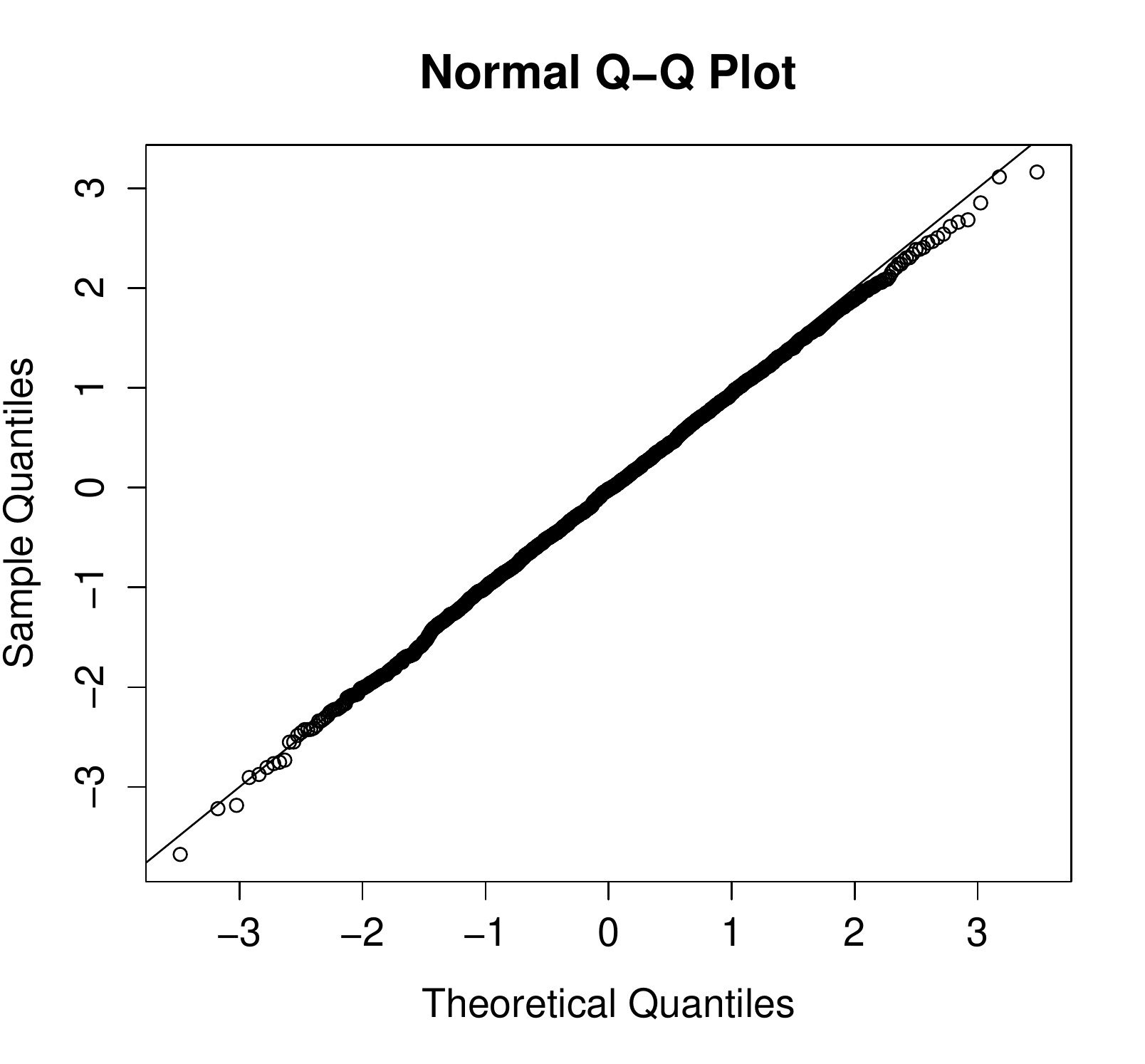} 
    \caption{$\omega=1, T=50$} 
    %\vspace{4ex}
  \end{subfigure}
 \begin{subfigure}[b]{0.32\linewidth}
    \centering
    \includegraphics[width=\linewidth]{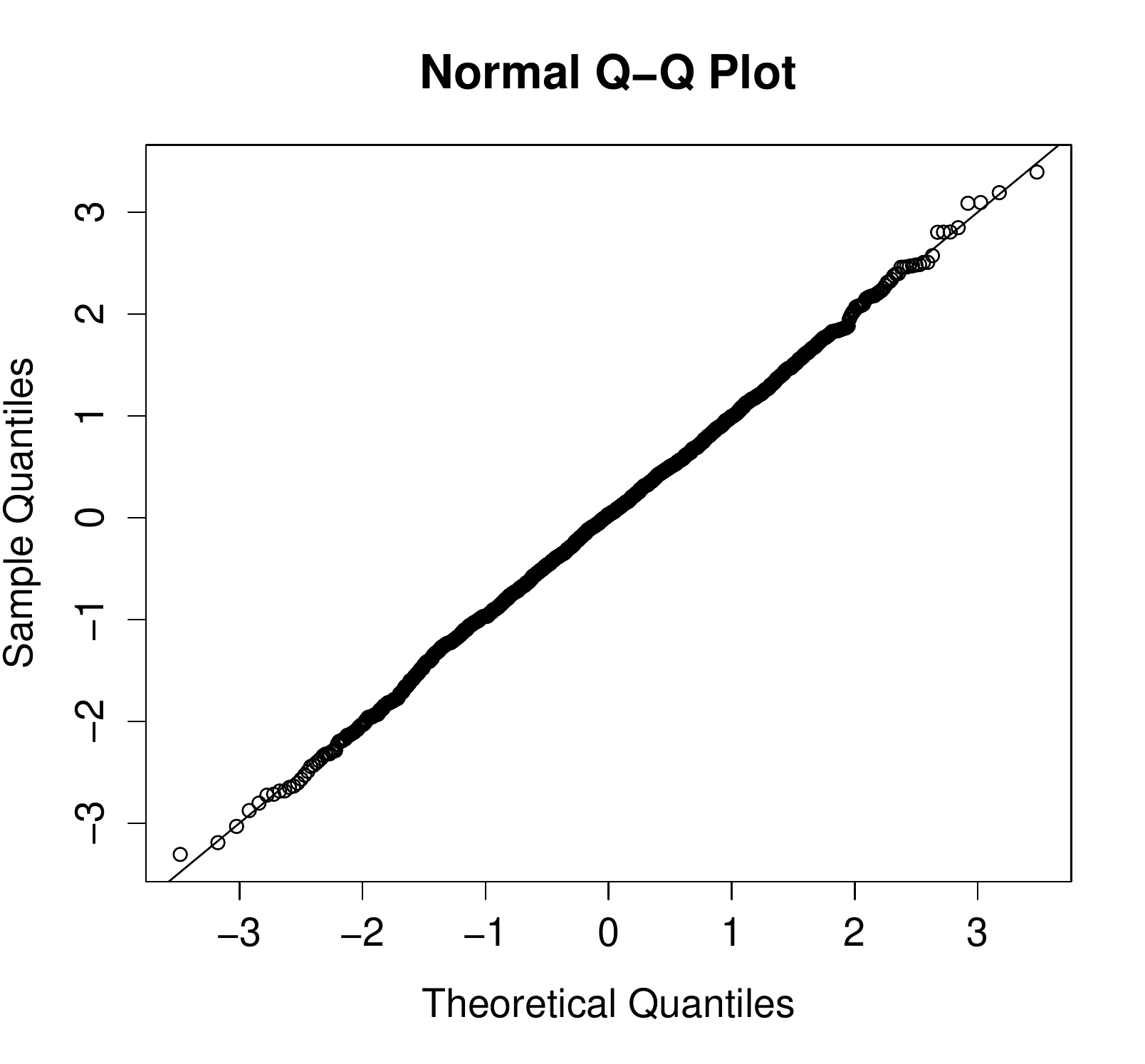} 
    \caption{$\omega=1, T=100$} 
    %\vspace{4ex}
  \end{subfigure}
  \begin{subfigure}[b]{0.32\linewidth}
    \centering
    \includegraphics[width=\linewidth]{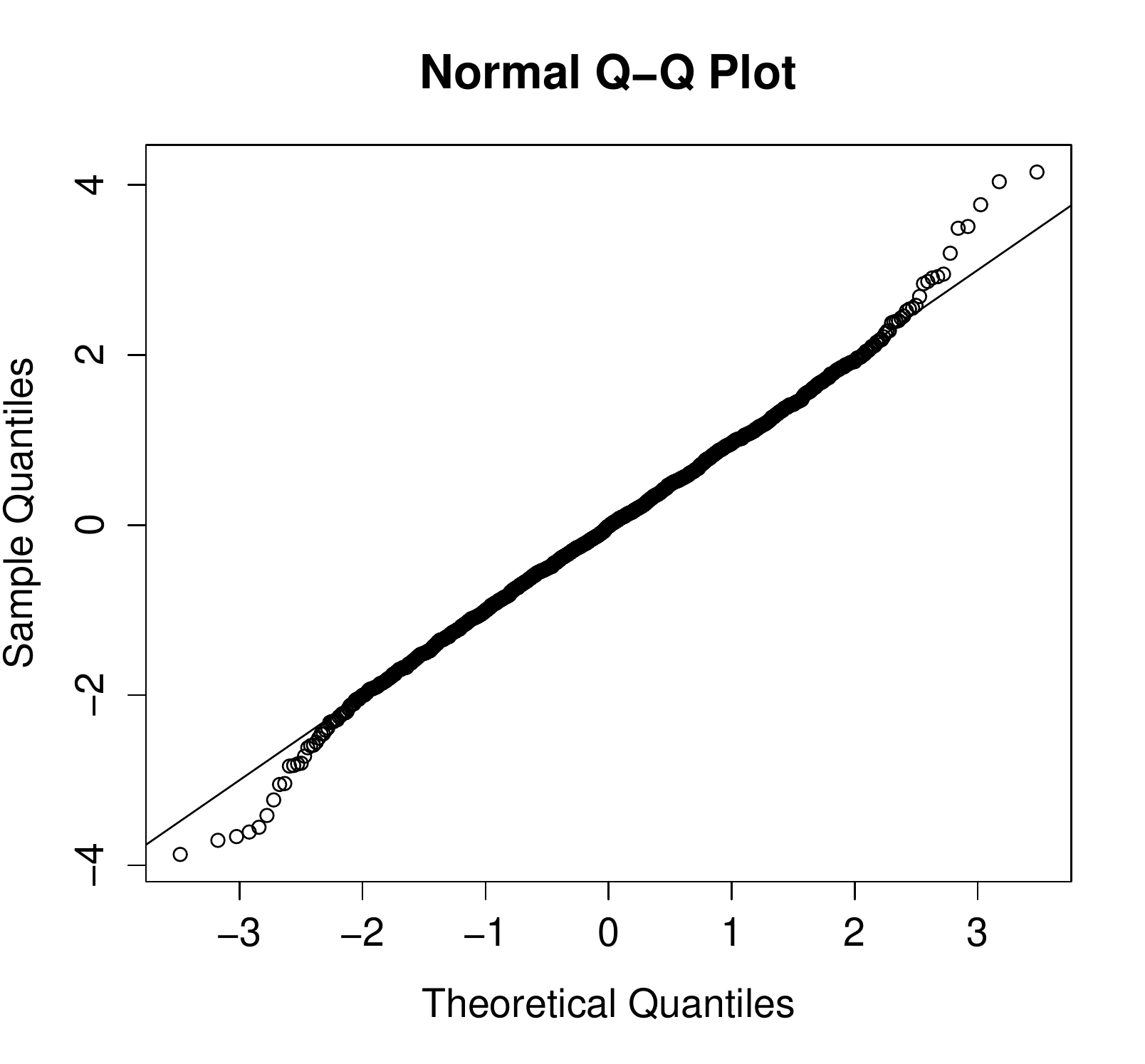} 
    \caption{$\omega=10, T=10$}  
    %\vspace{4ex}
  \end{subfigure}%% 
 \begin{subfigure}[b]{0.32\linewidth}
    \centering
    \includegraphics[width=\linewidth]{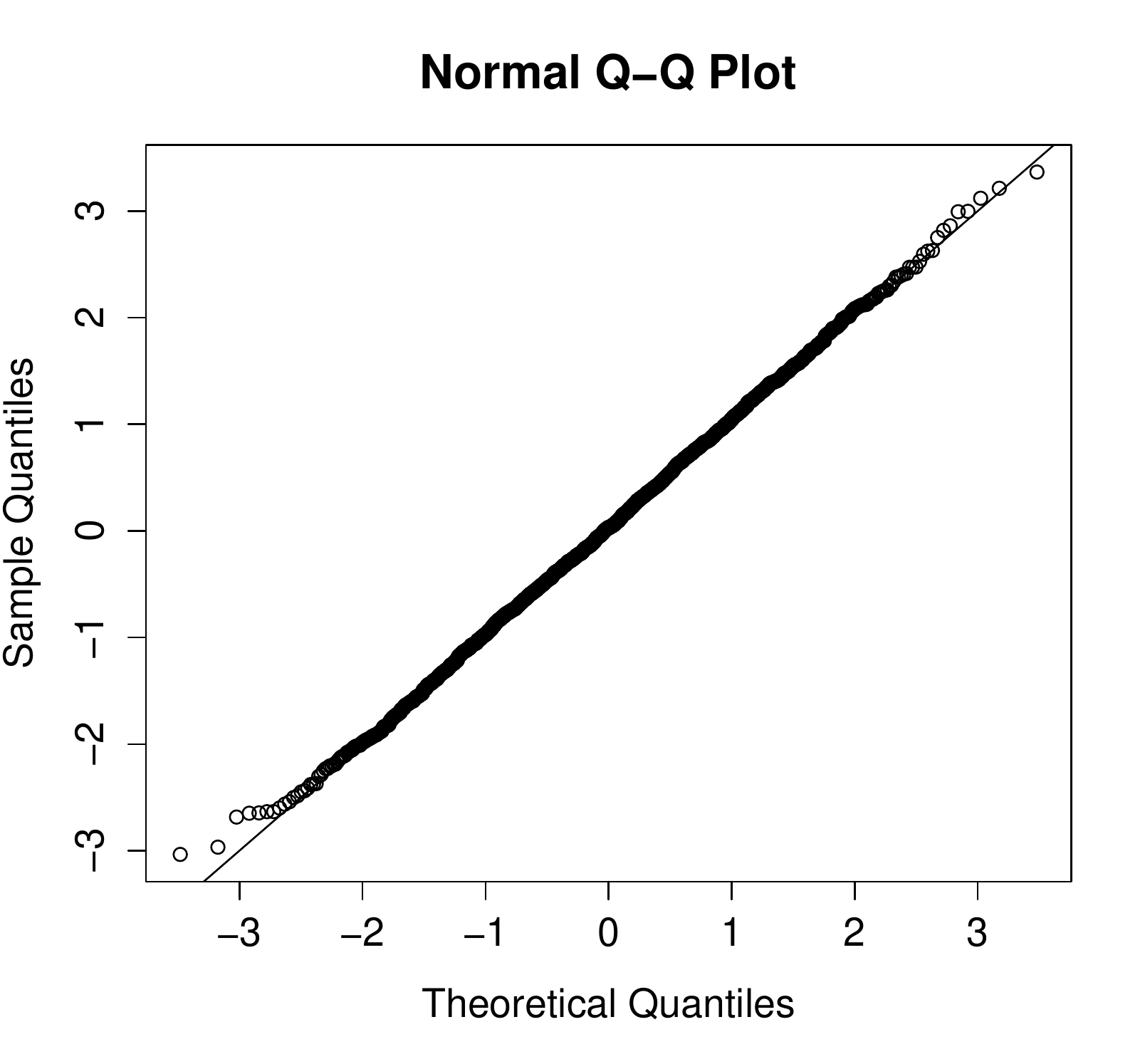} 
    \caption{$\omega=10, T=50$} 
    %\vspace{4ex}
  \end{subfigure}
 \begin{subfigure}[b]{0.32\linewidth}
    \centering
    \includegraphics[width=\linewidth]{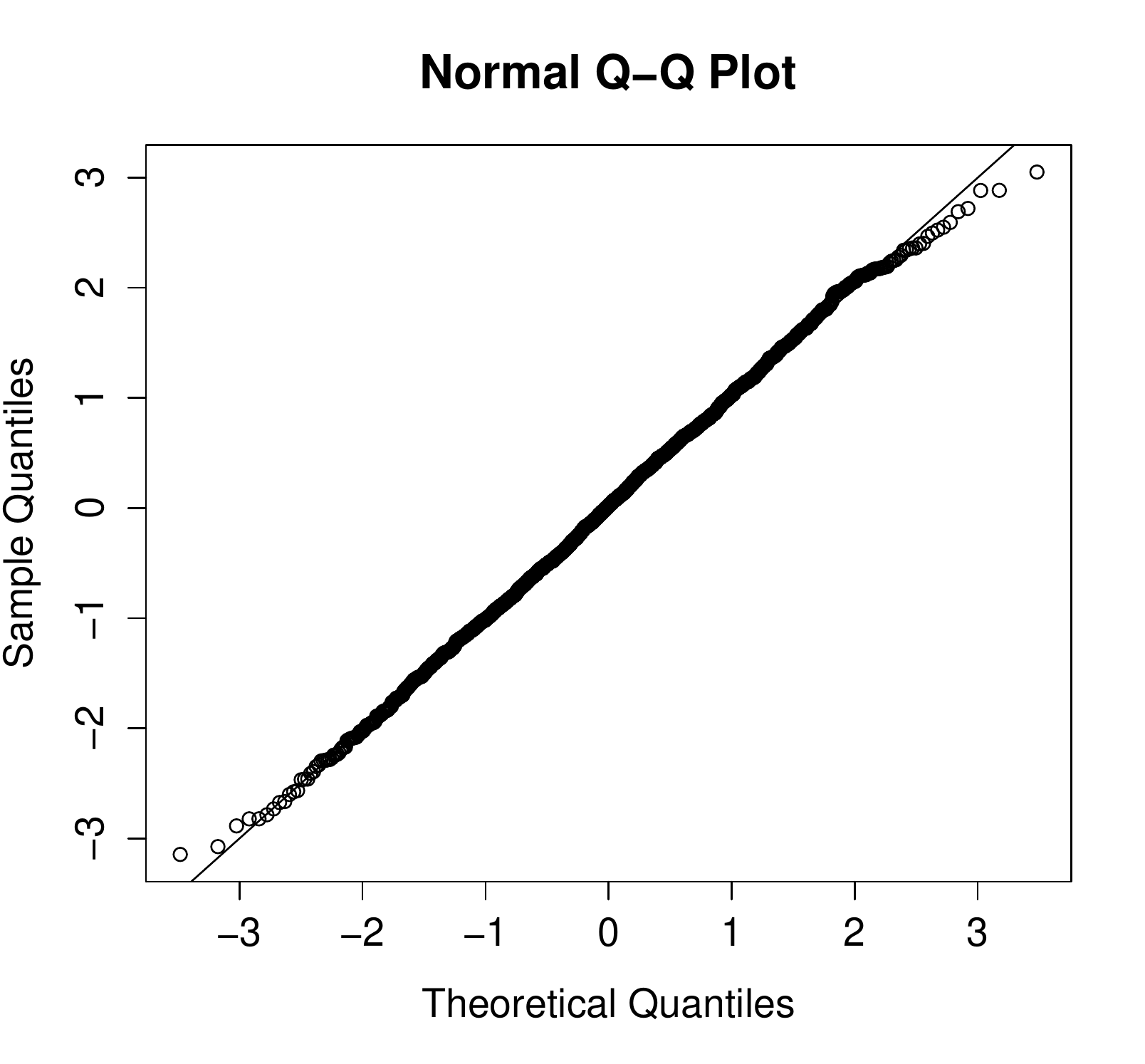} 
    \caption{$\omega=10, T=100$} 
    %\vspace{4ex}
  \end{subfigure}
  \caption{Normal QQ plots for the real part of the truncated Fourier transform of the simulated CARMA(2,1) processes driven by a Variance Gamma process for the frequencies $0, 0.1, 1 , 10$ (rows) and time horizons/maximum non-equidistant grid sizes $10/0.1, 50/0.05, 100/0.01$ (columns). The theoretical quantiles are coming from the (limiting) law described in Theorem \ref{thm:ApproxTFTDoubleLimitDistribution}. }\label{plot:QQCARMAVG} 
\end{figure}

\begin{figure}[tp]    
  \begin{subfigure}[b]{0.32\linewidth}
    \centering
    \includegraphics[width=\linewidth]{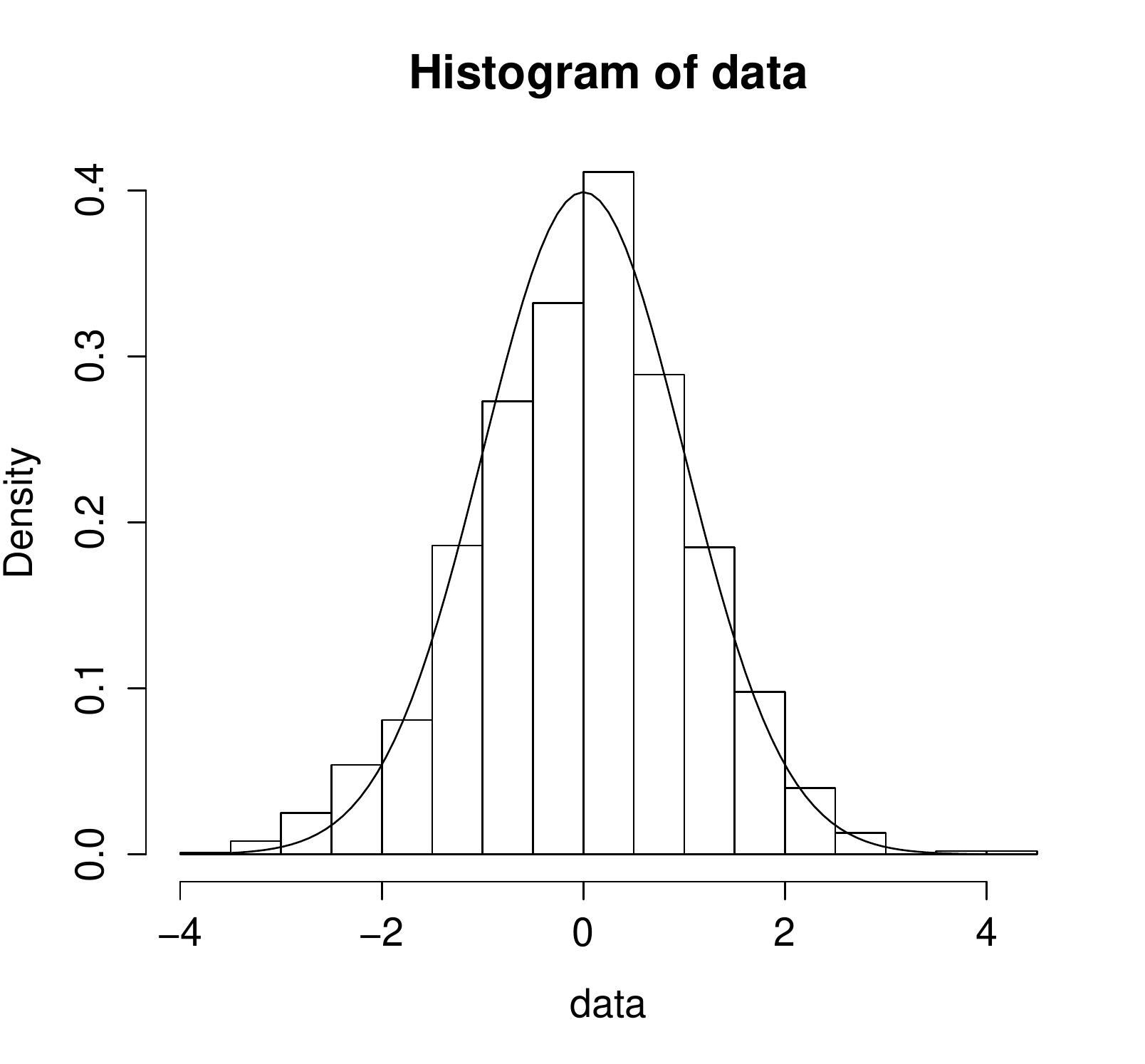} 
    \caption{ $\omega=0, T=10$} 
  
    %\vspace{4ex}
  \end{subfigure}%% 
 \begin{subfigure}[b]{0.32\linewidth}
    \centering
    \includegraphics[width=\linewidth]{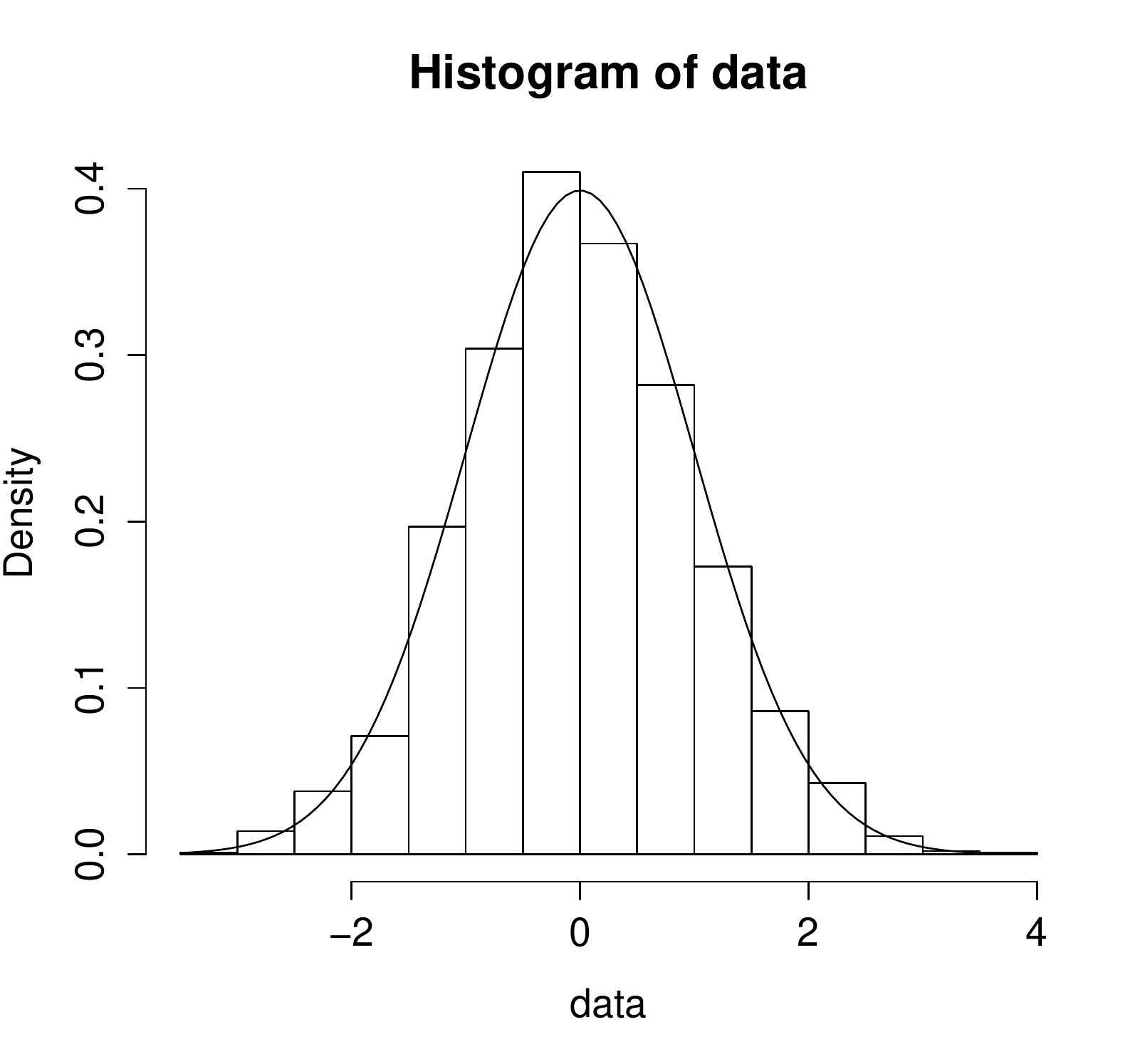} 
    \caption{$\omega=0, T=50$} 
    %\vspace{4ex}
  \end{subfigure}
 \begin{subfigure}[b]{0.32\linewidth}
    \centering
    \includegraphics[width=\linewidth]{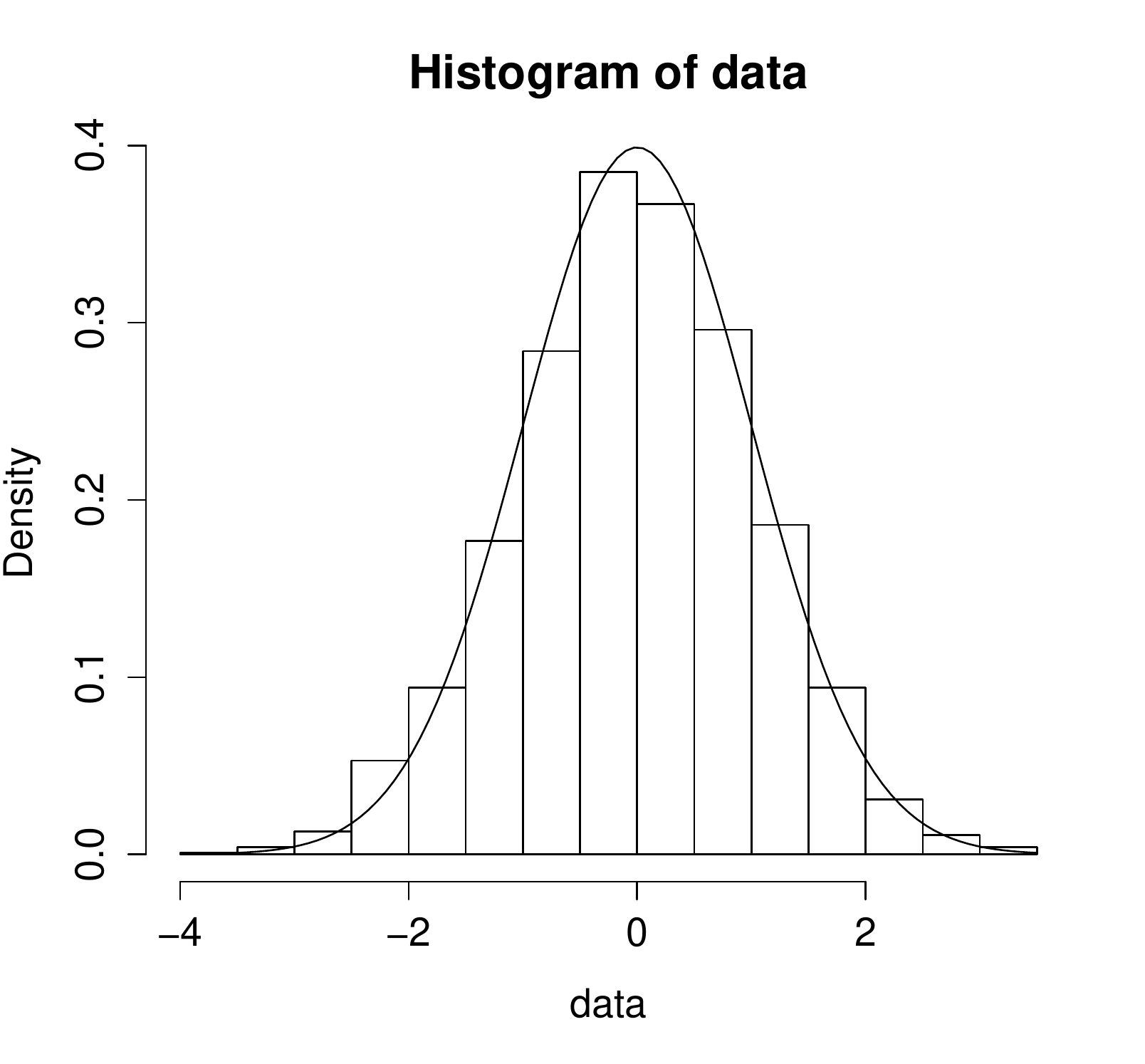} 
    \caption{$\omega=0, T=100$} 
    %\vspace{4ex}
  \end{subfigure}
  \begin{subfigure}[b]{0.32\linewidth}
    \centering
    \includegraphics[width=\linewidth]{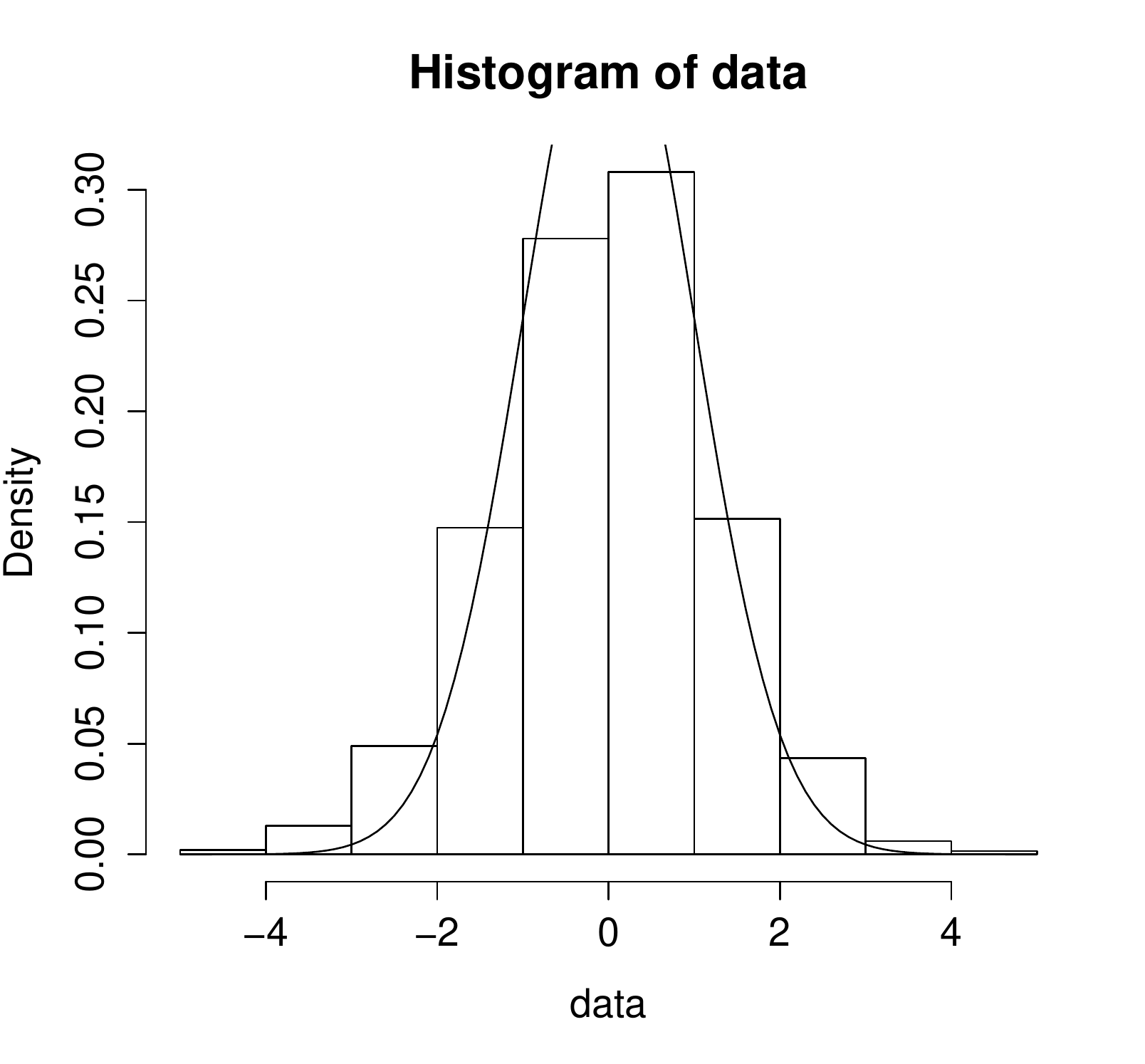} 
    \caption{$\omega=0.1, T=10$} 
    %\vspace{4ex}
  \end{subfigure}%% 
 \begin{subfigure}[b]{0.32\linewidth}
    \centering
    \includegraphics[width=\linewidth]{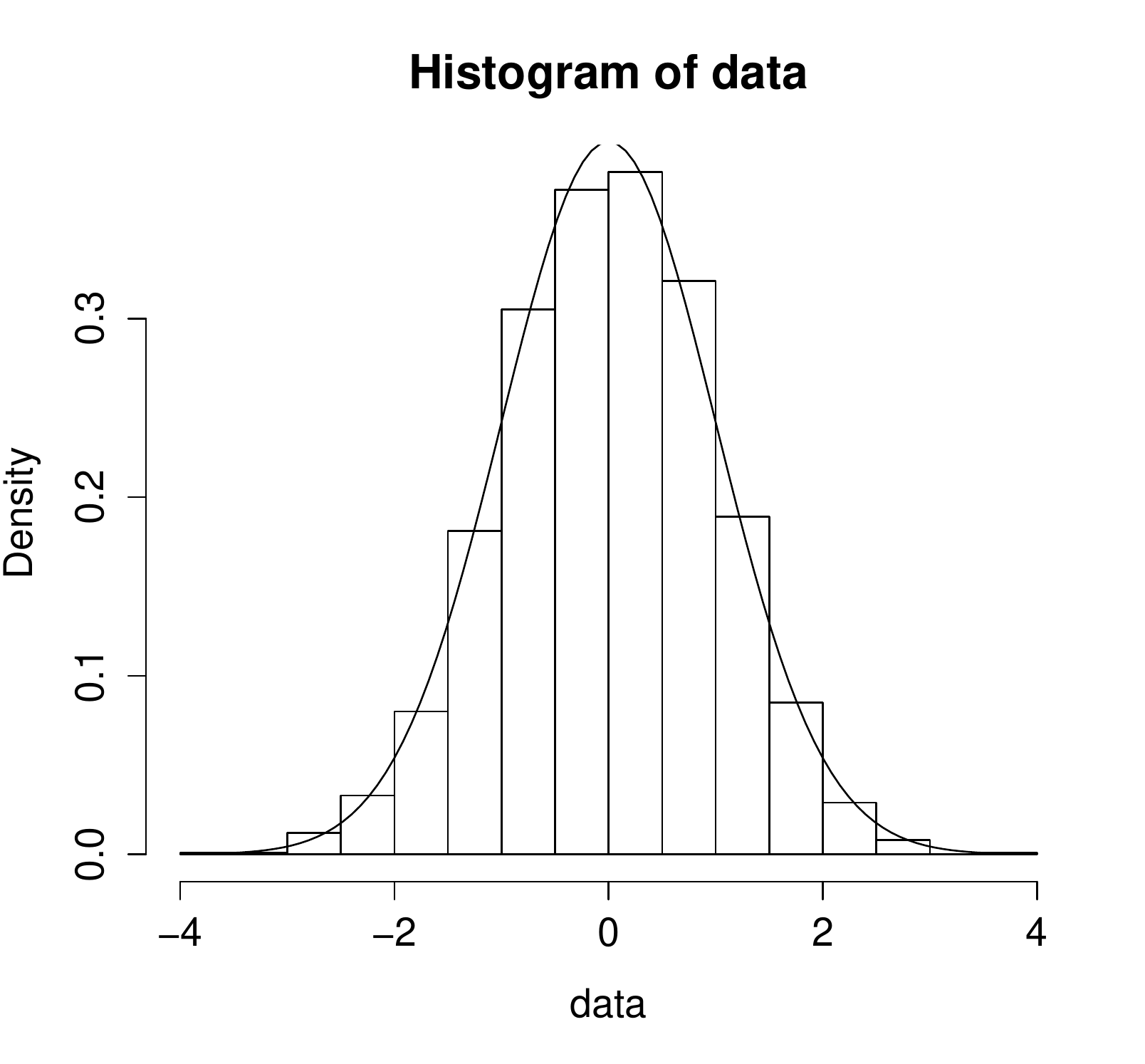} 
    \caption{$\omega=0.1, T=50$} 
    %\vspace{4ex}
  \end{subfigure}
 \begin{subfigure}[b]{0.32\linewidth}
    \centering
    \includegraphics[width=\linewidth]{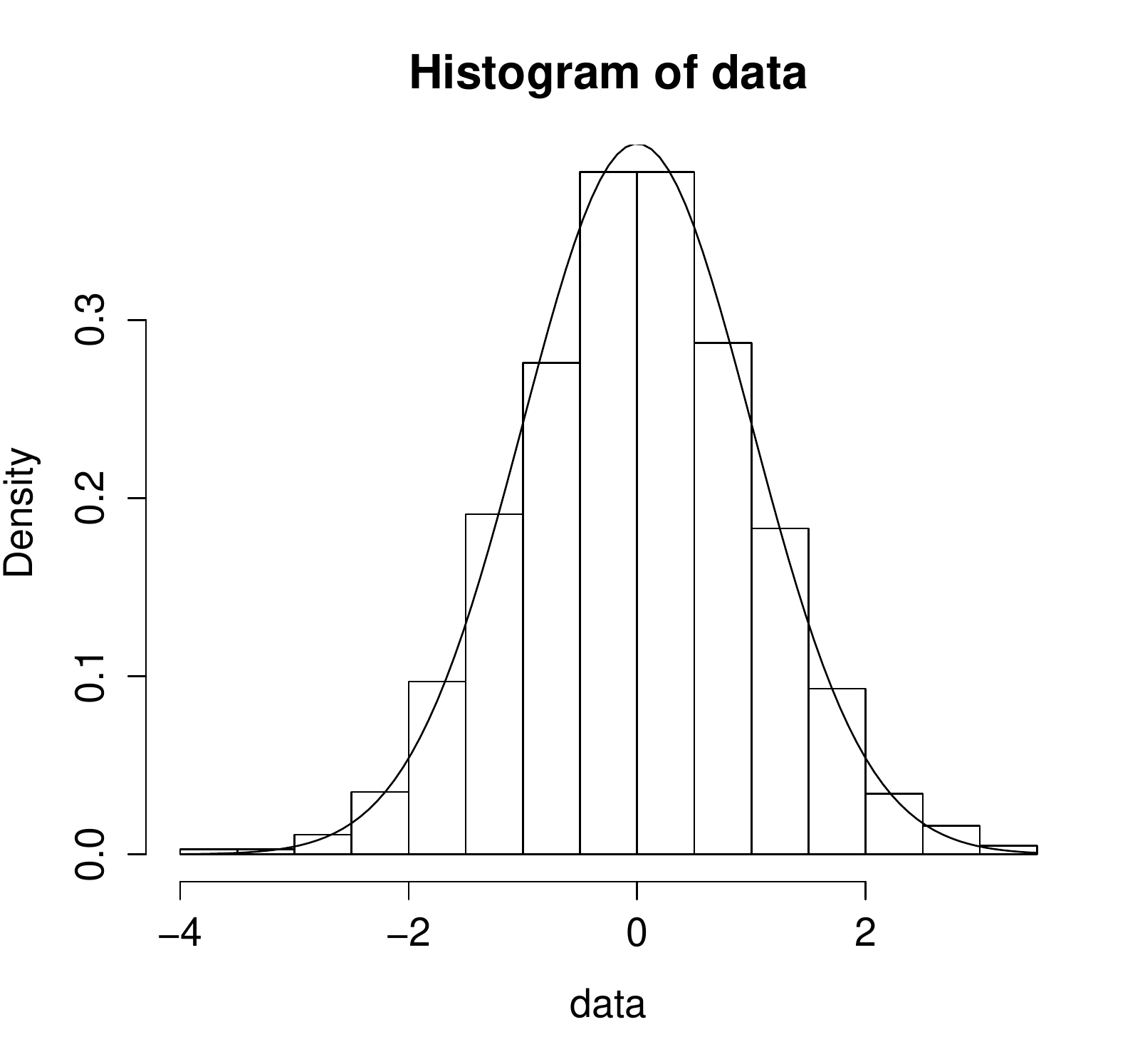} 
    \caption{$\omega=0.1, T=100$} 
    %\vspace{4ex}
  \end{subfigure}
  \begin{subfigure}[b]{0.32\linewidth}
    \centering
    \includegraphics[width=\linewidth]{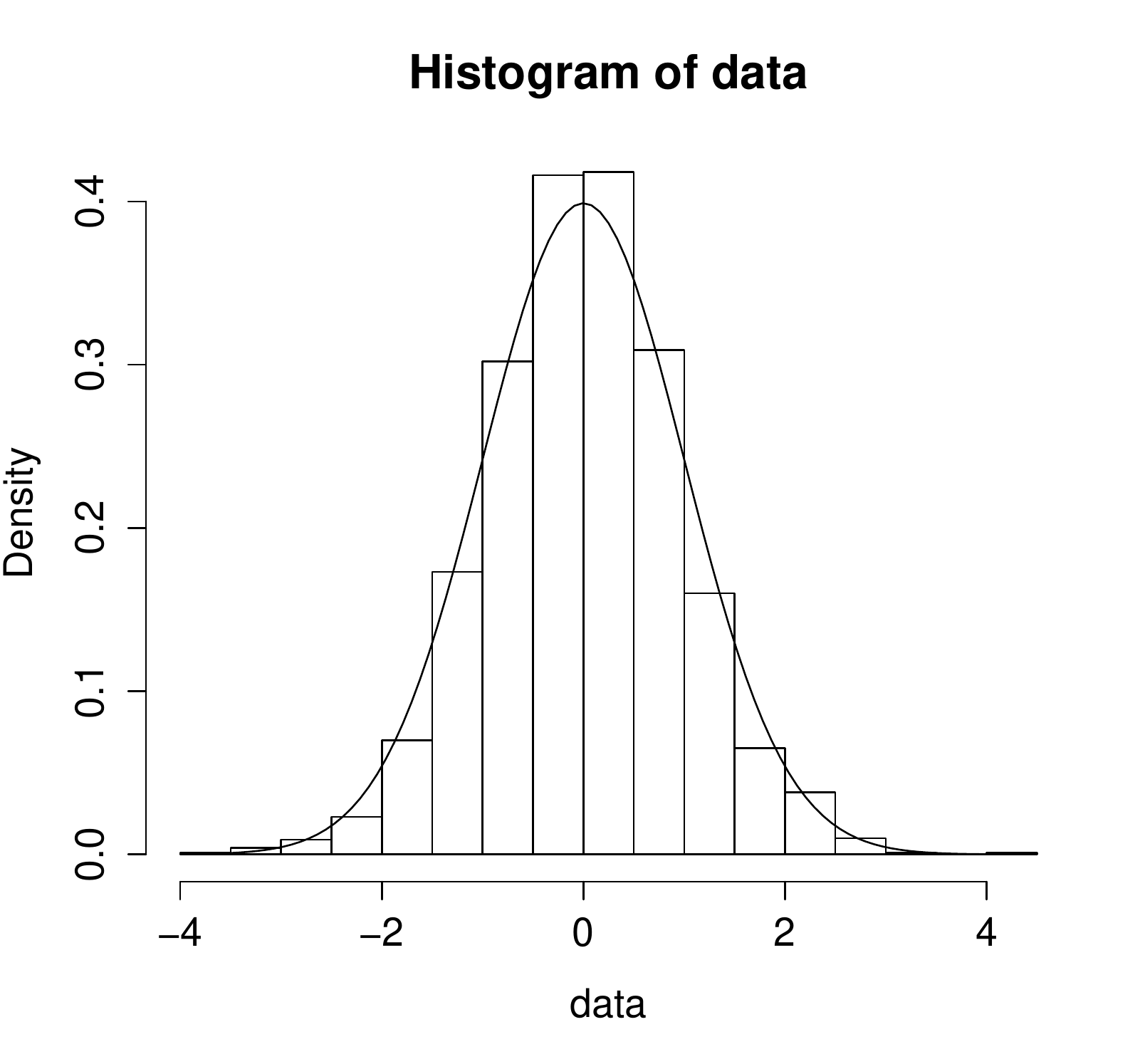} 
    \caption{$\omega=1, T=10$} 
    %\vspace{4ex}
  \end{subfigure}%% 
 \begin{subfigure}[b]{0.32\linewidth}
    \centering
    \includegraphics[width=\linewidth]{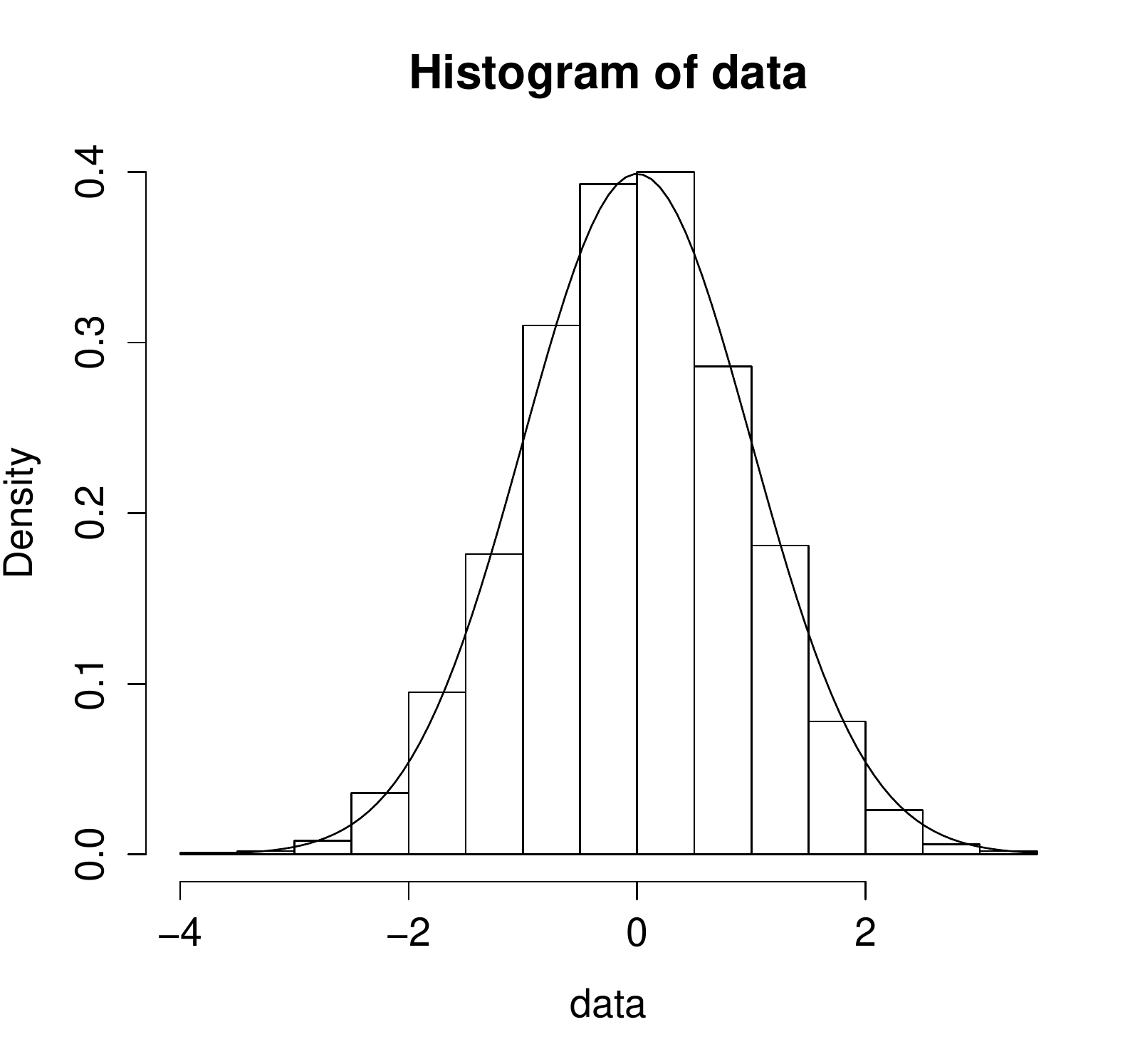} 
    \caption{$\omega=1, T=50$} 
    %\vspace{4ex}
  \end{subfigure}
 \begin{subfigure}[b]{0.32\linewidth}
    \centering
    \includegraphics[width=\linewidth]{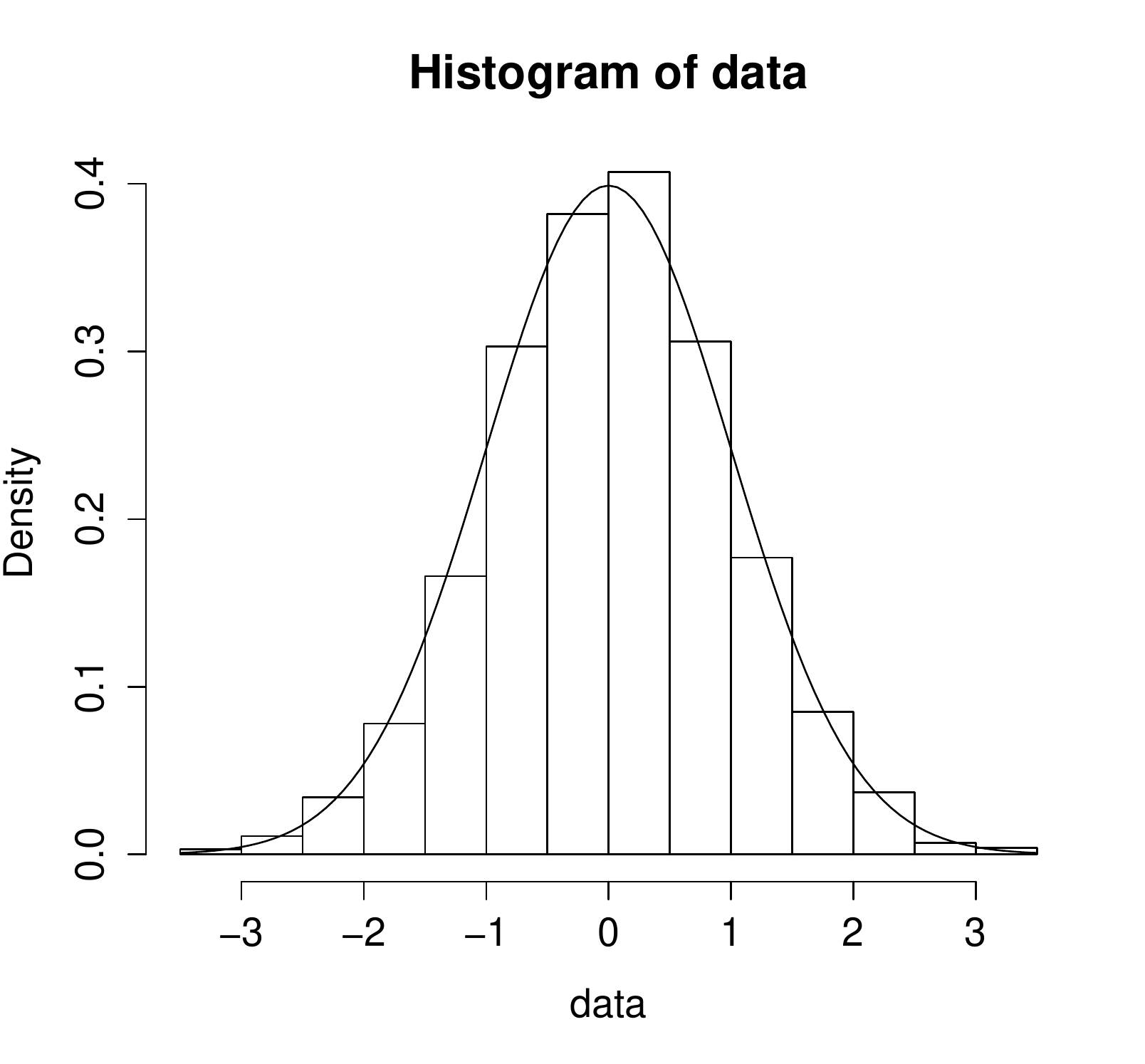} 
    \caption{$\omega=1, T=100$} 
    %\vspace{4ex}
  \end{subfigure}
  \begin{subfigure}[b]{0.32\linewidth}
    \centering
    \includegraphics[width=\linewidth]{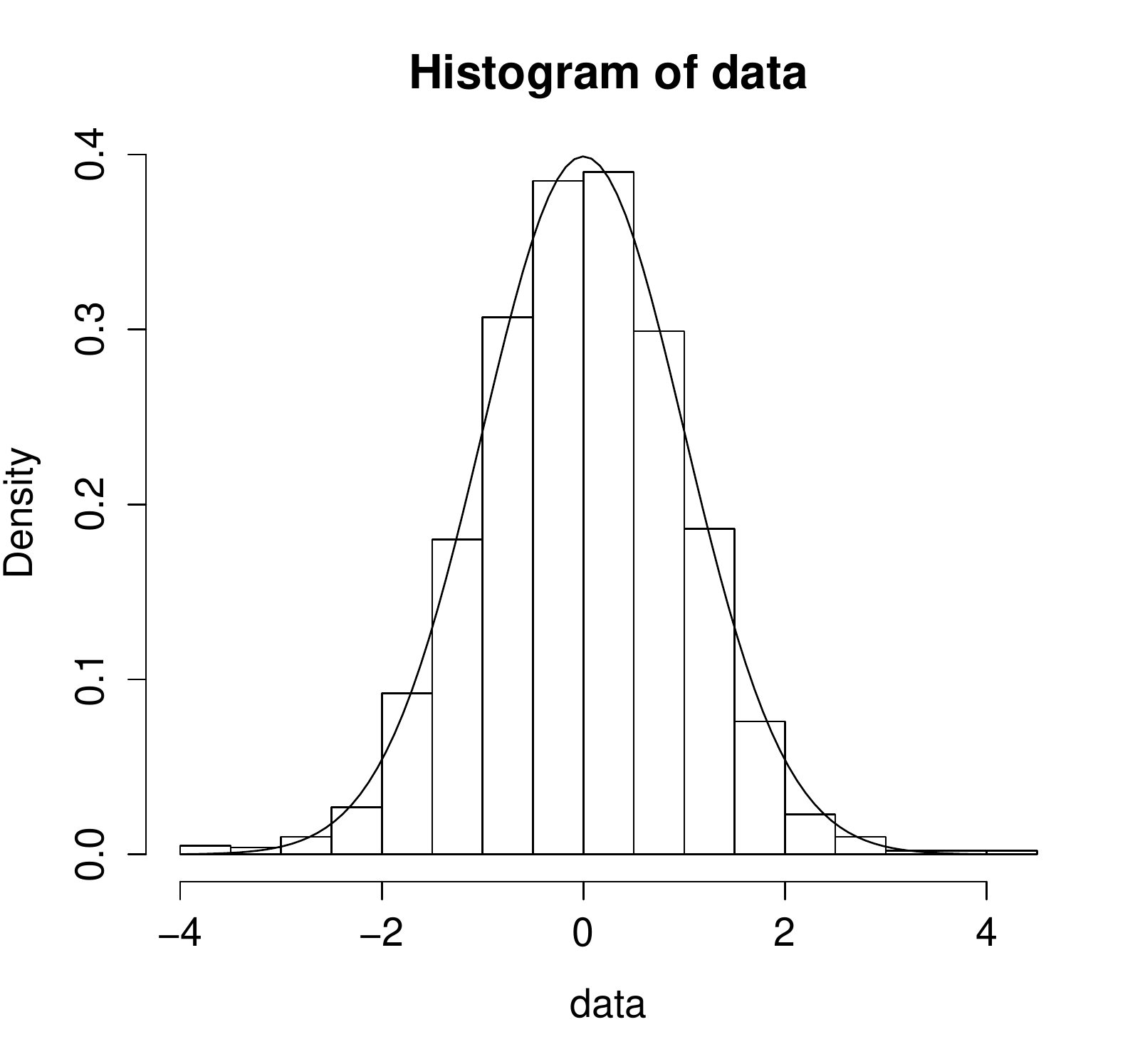} 
    \caption{$\omega=10, T=10$}  
    %\vspace{4ex}
  \end{subfigure}%% 
 \begin{subfigure}[b]{0.32\linewidth}
    \centering
    \includegraphics[width=\linewidth]{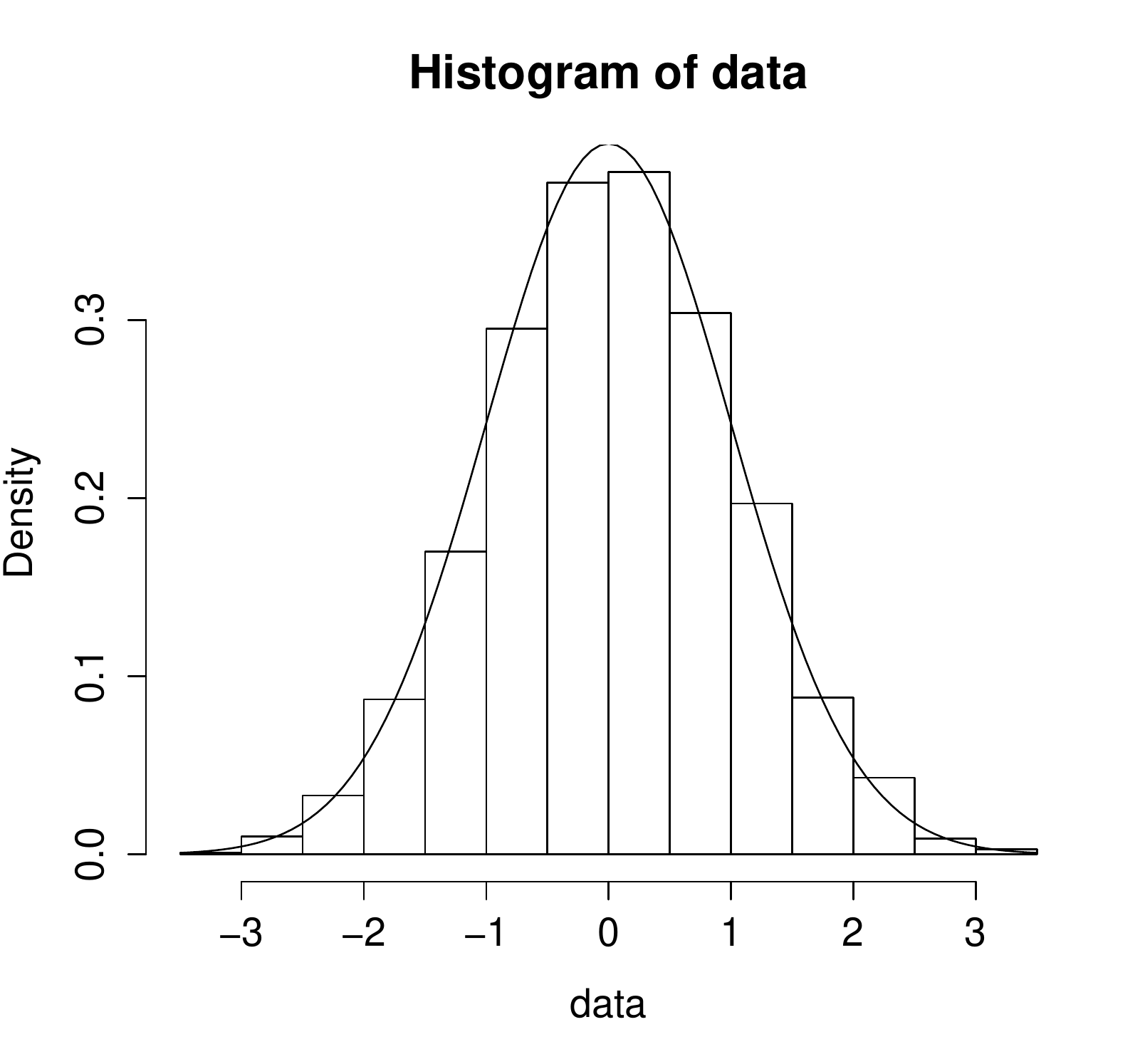} 
    \caption{$\omega=10, T=50$} 
    %\vspace{4ex}
  \end{subfigure}
 \begin{subfigure}[b]{0.32\linewidth}
    \centering
    \includegraphics[width=\linewidth]{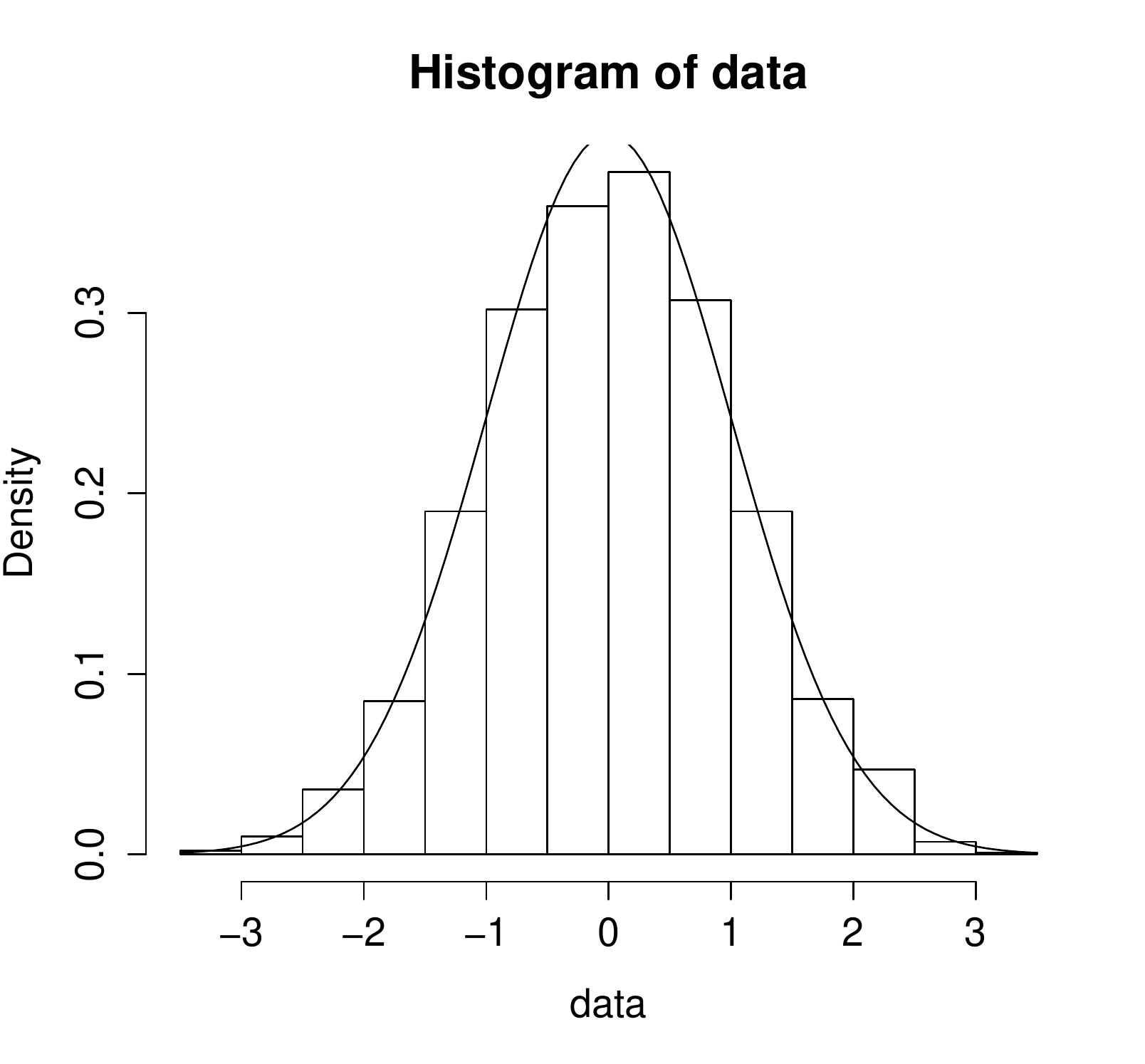} 
    \caption{$\omega=10, T=100$} 
    %\vspace{4ex}
  \end{subfigure}
  \caption{Histograms and limiting density for the real part of the truncated Fourier transform of the simulated CARMA(2,1) processes driven by a Variance Gamma process for the frequencies $0, 0.1, 1 , 10$ (rows) and time horizons/maximum non-equidistant grid sizes $10/0.1, 50/0.05, 100/0.01$ (columns)}\label{plot:HistCARMAVG} 
\end{figure}

\begin{figure}[tp]    
  \begin{subfigure}[b]{0.32\linewidth}
    \centering
    \includegraphics[width=\linewidth]{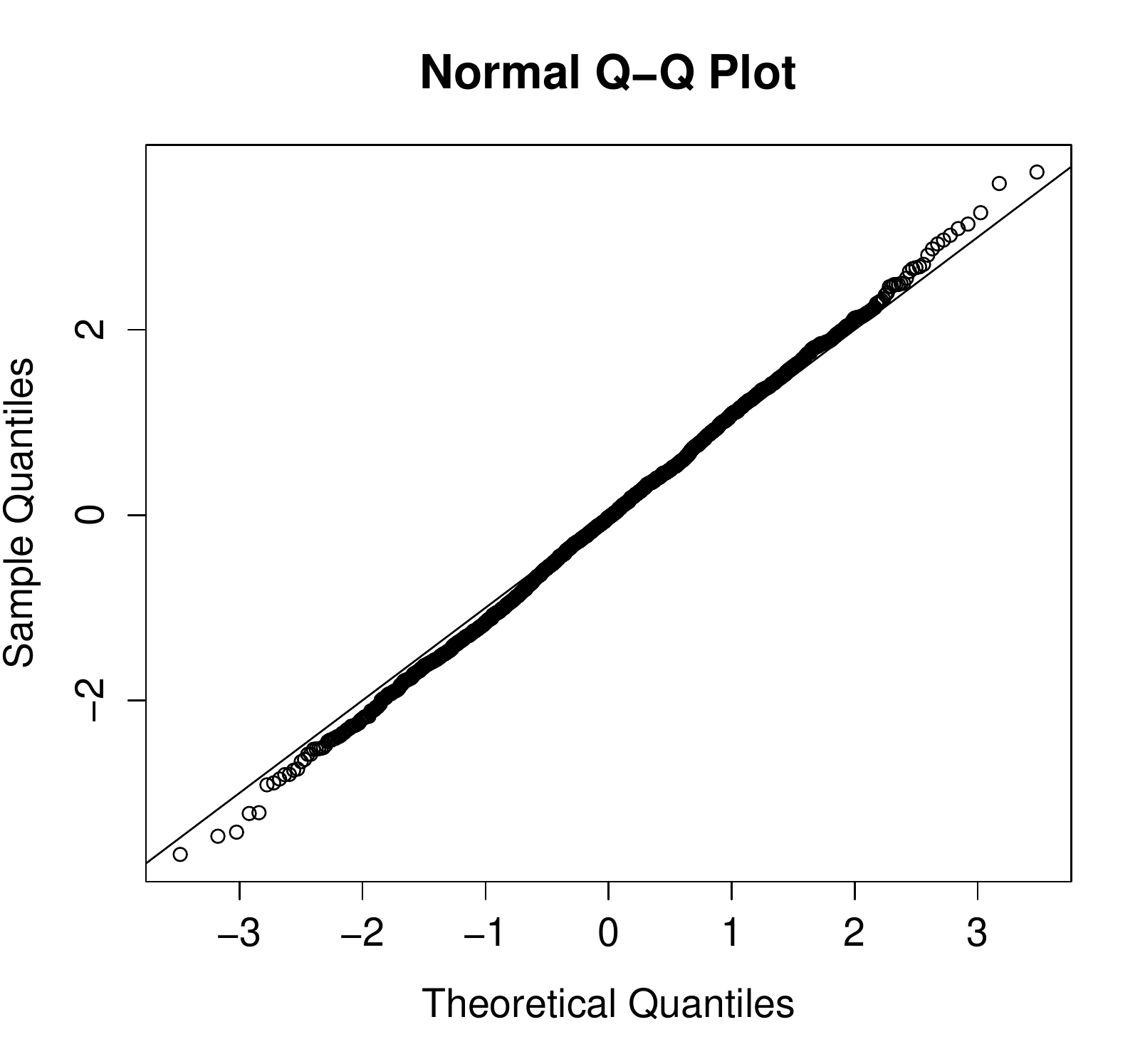} 
    \caption{ $\omega=0, T=10$} 
  
    %\vspace{4ex}
  \end{subfigure}%% 
 \begin{subfigure}[b]{0.32\linewidth}
    \centering
    \includegraphics[width=\linewidth]{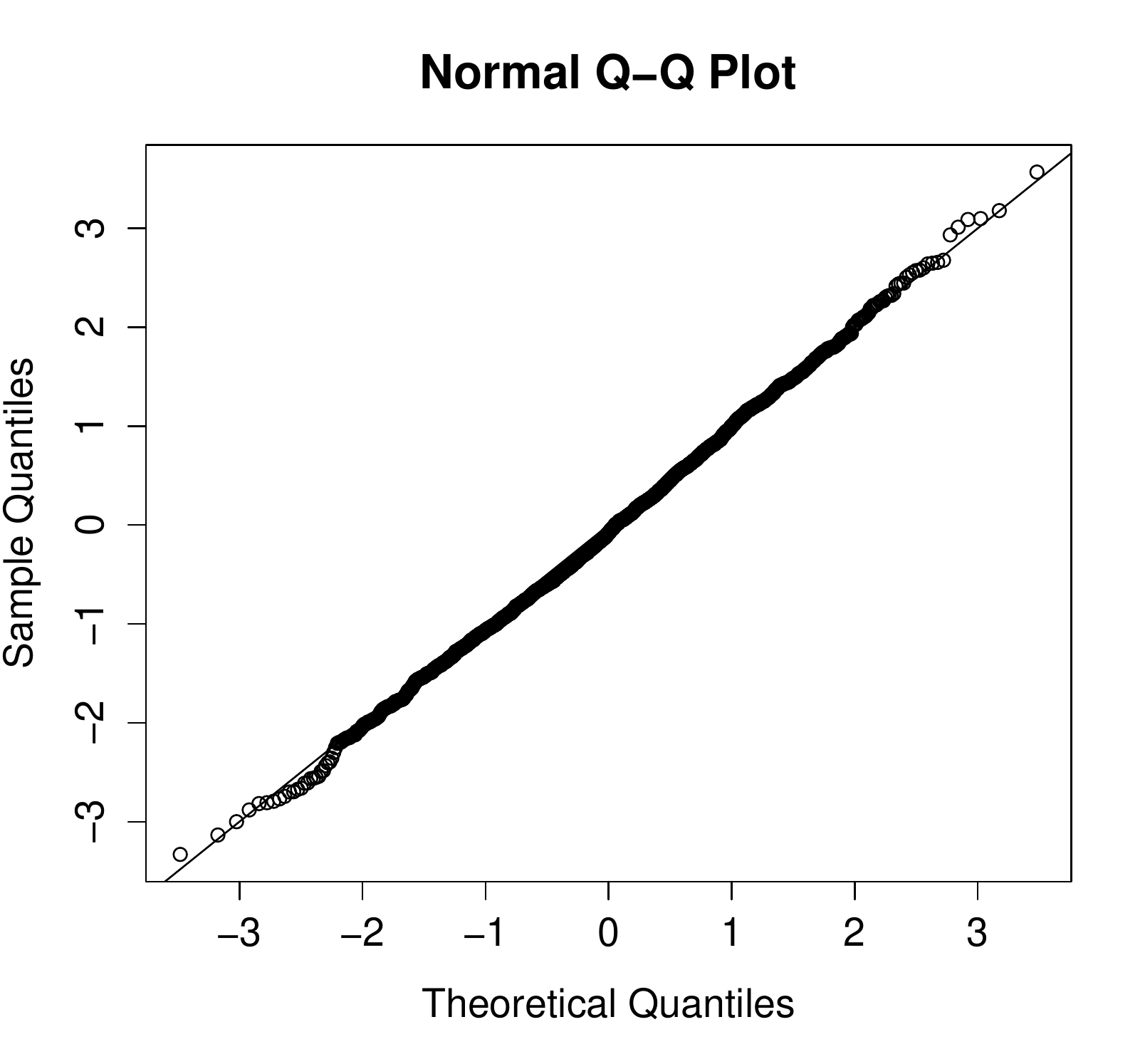} 
    \caption{$\omega=0, T=50$} 
    %\vspace{4ex}
  \end{subfigure}
 \begin{subfigure}[b]{0.32\linewidth}
    \centering
    \includegraphics[width=\linewidth]{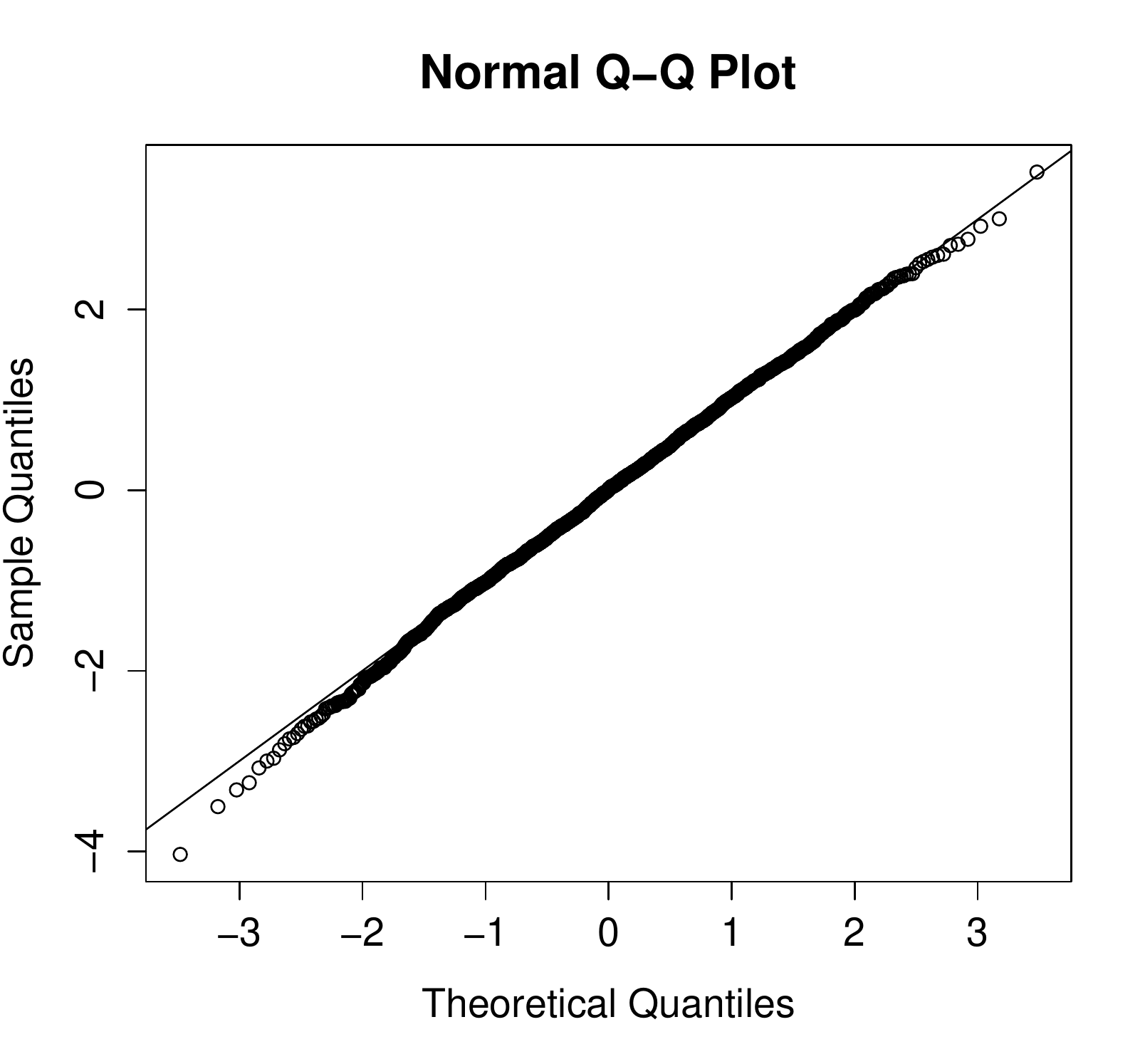} 
    \caption{$\omega=0, T=100$} 
    %\vspace{4ex}
  \end{subfigure}
  \begin{subfigure}[b]{0.32\linewidth}
    \centering
    \includegraphics[width=\linewidth]{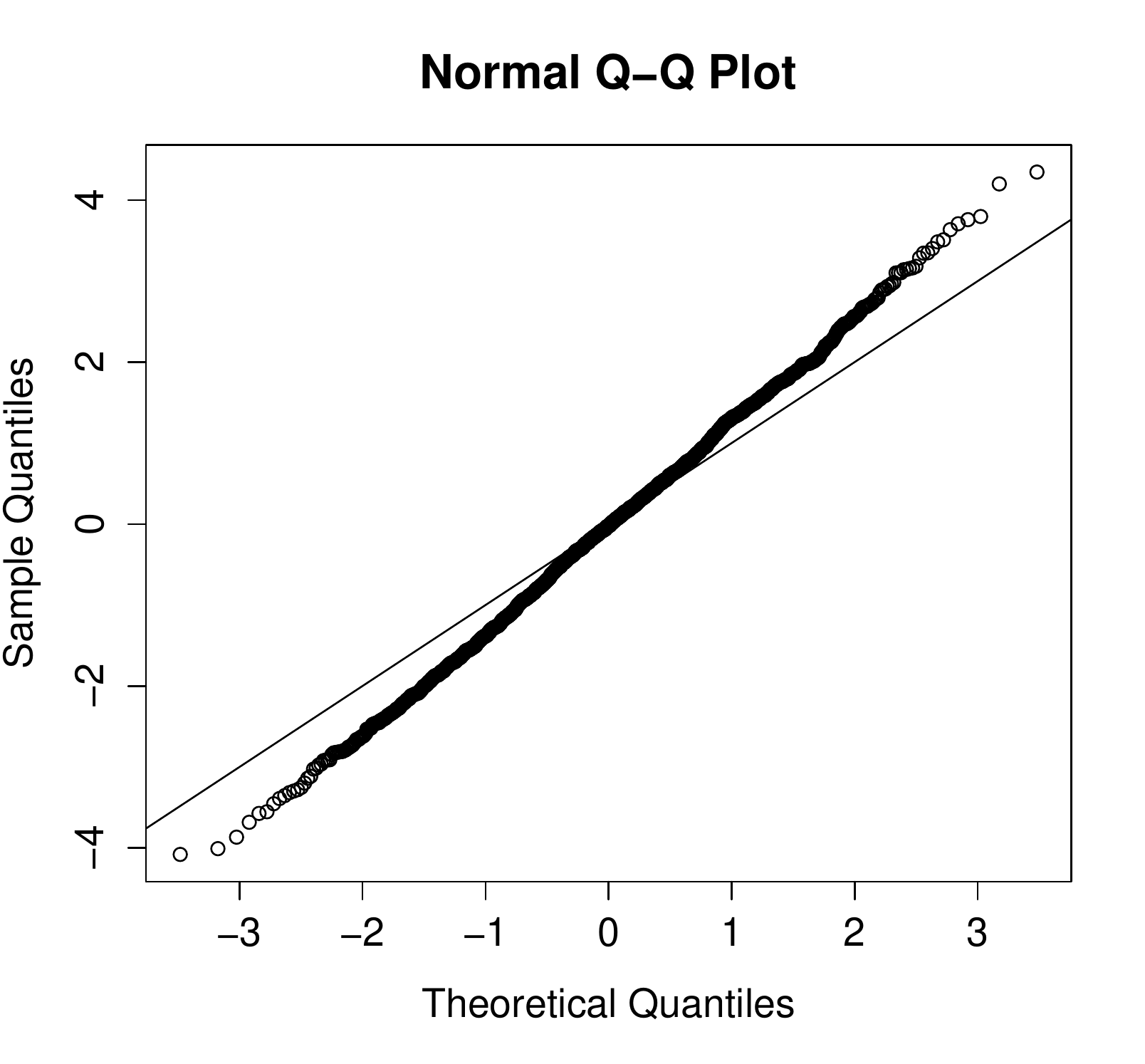} 
    \caption{$\omega=0.1, T=10$} 
    %\vspace{4ex}
  \end{subfigure}%% 
 \begin{subfigure}[b]{0.32\linewidth}
    \centering
    \includegraphics[width=\linewidth]{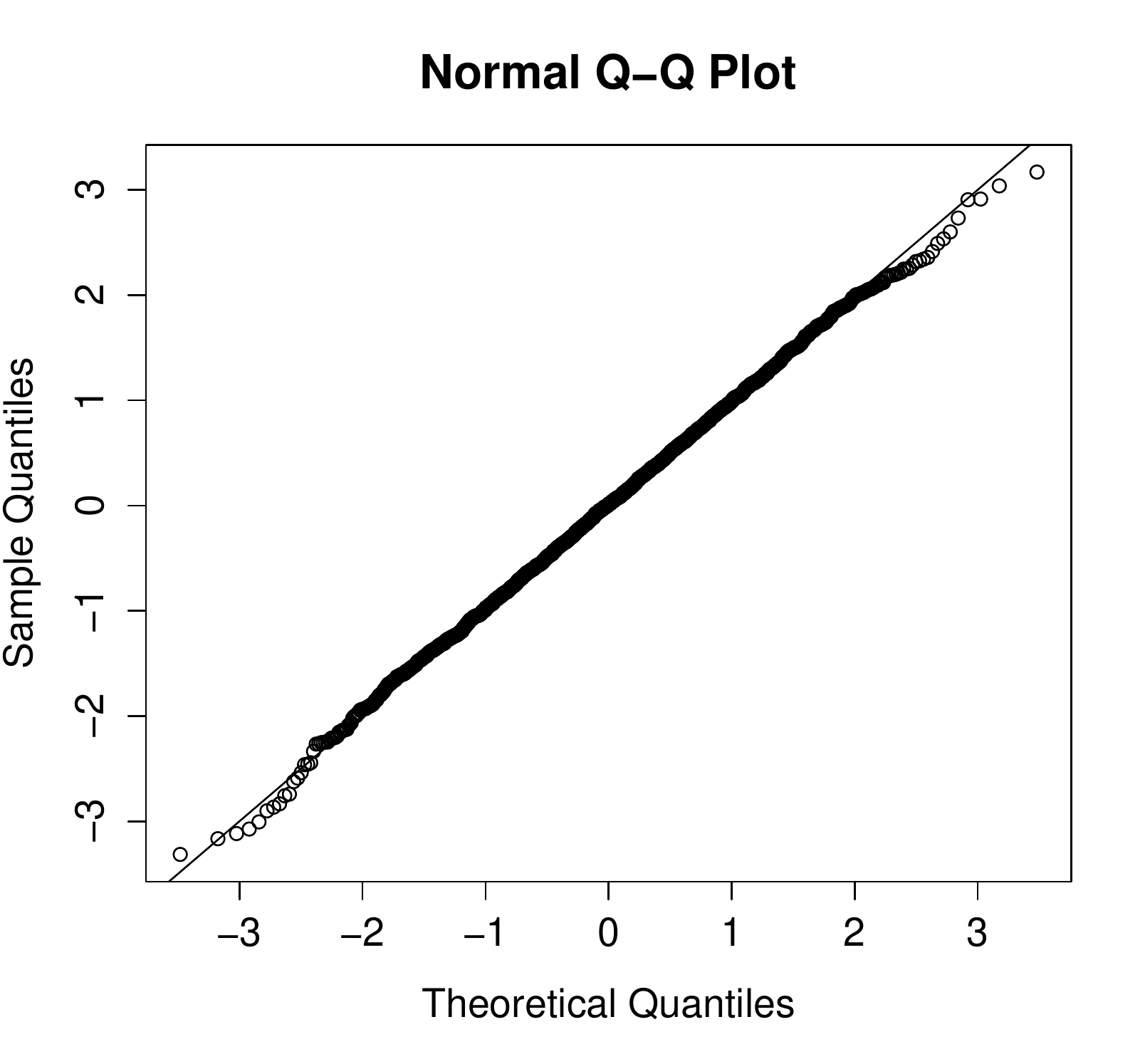} 
    \caption{$\omega=0.1, T=50$} 
    %\vspace{4ex}
  \end{subfigure}
 \begin{subfigure}[b]{0.32\linewidth}
    \centering
    \includegraphics[width=\linewidth]{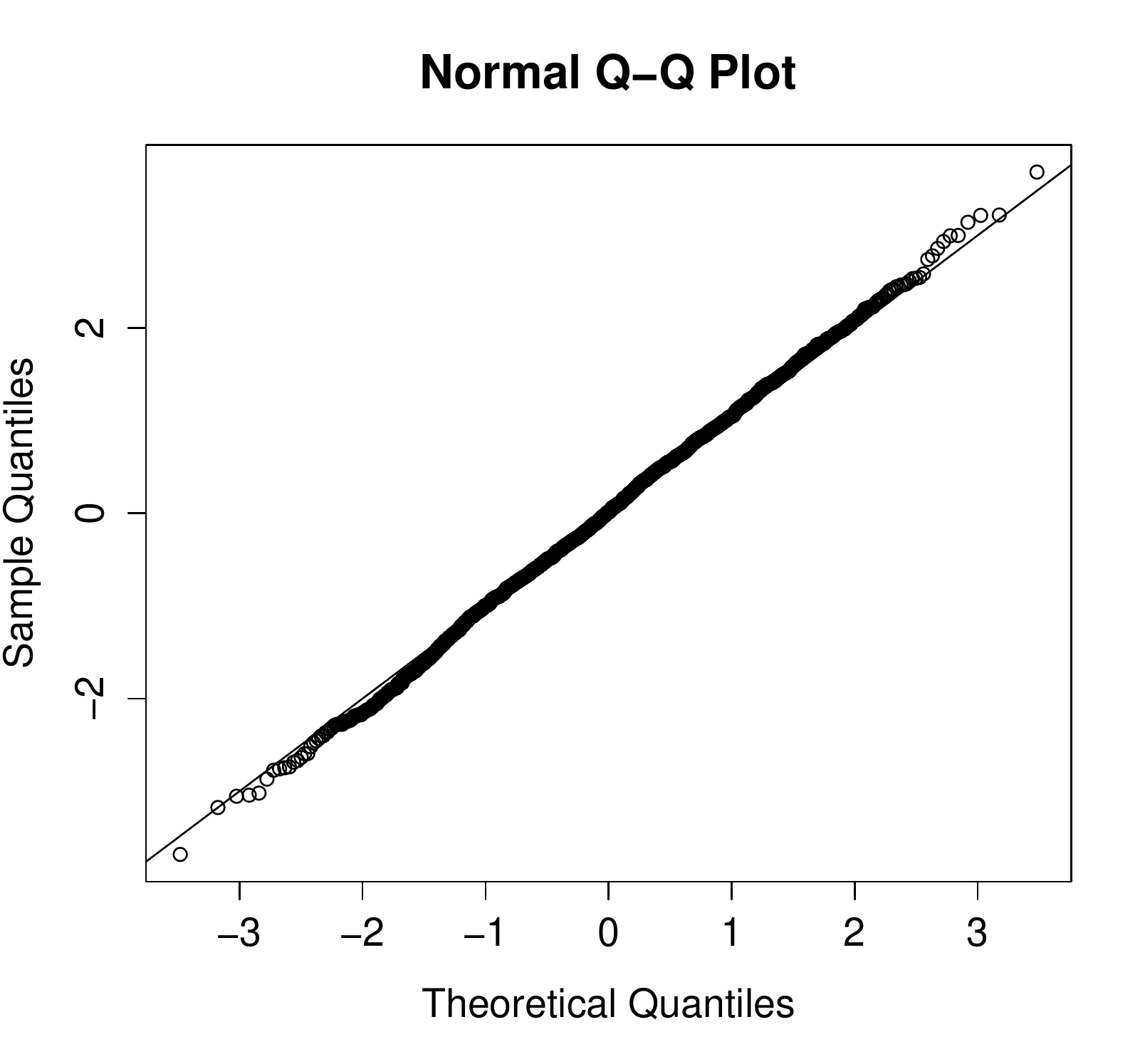} 
    \caption{$\omega=0.1, T=100$} 
    %\vspace{4ex}
  \end{subfigure}
  \begin{subfigure}[b]{0.32\linewidth}
    \centering
    \includegraphics[width=\linewidth]{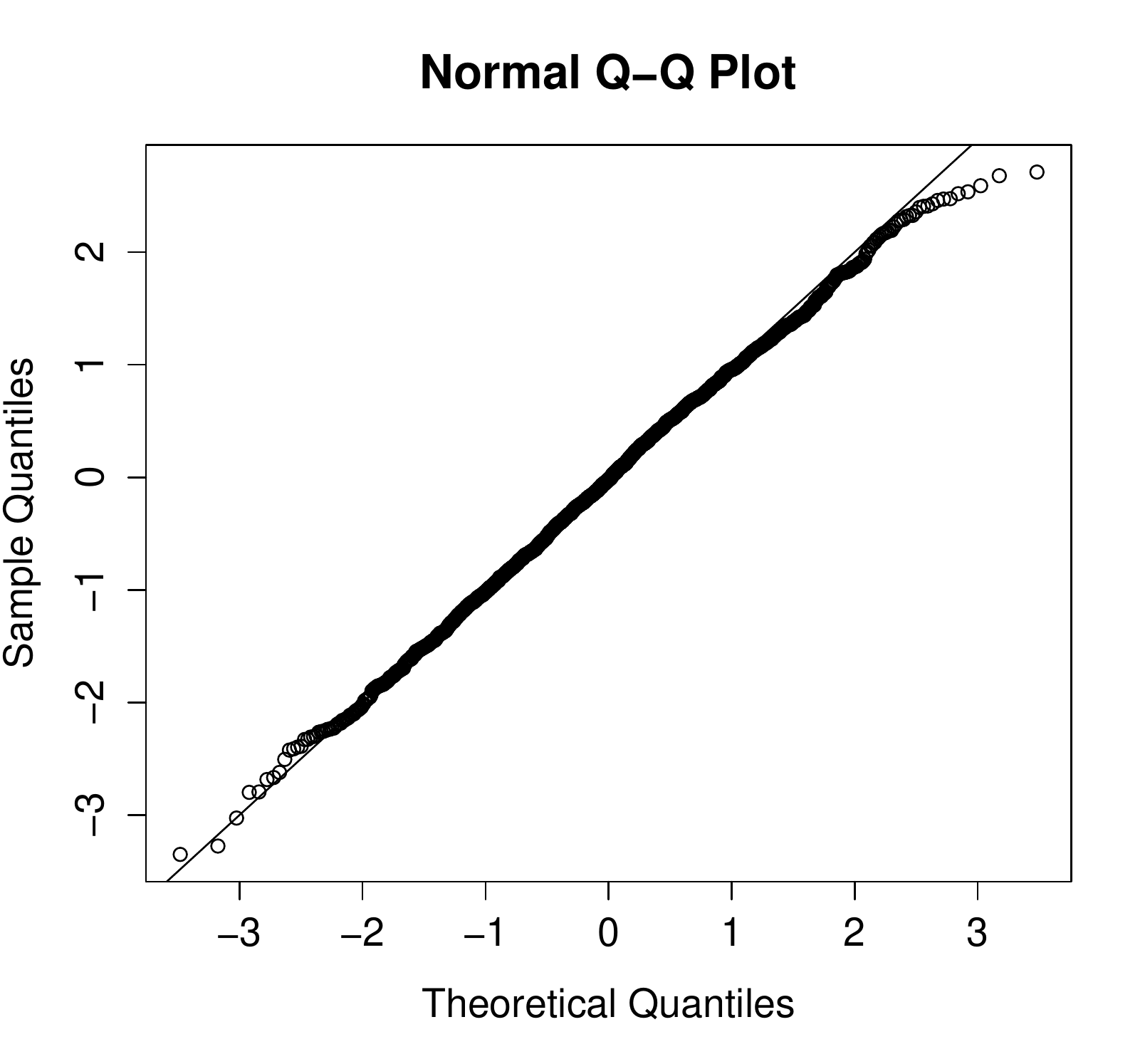} 
    \caption{$\omega=1, T=10$} 
    %\vspace{4ex}
  \end{subfigure}%% 
 \begin{subfigure}[b]{0.32\linewidth}
    \centering
    \includegraphics[width=\linewidth]{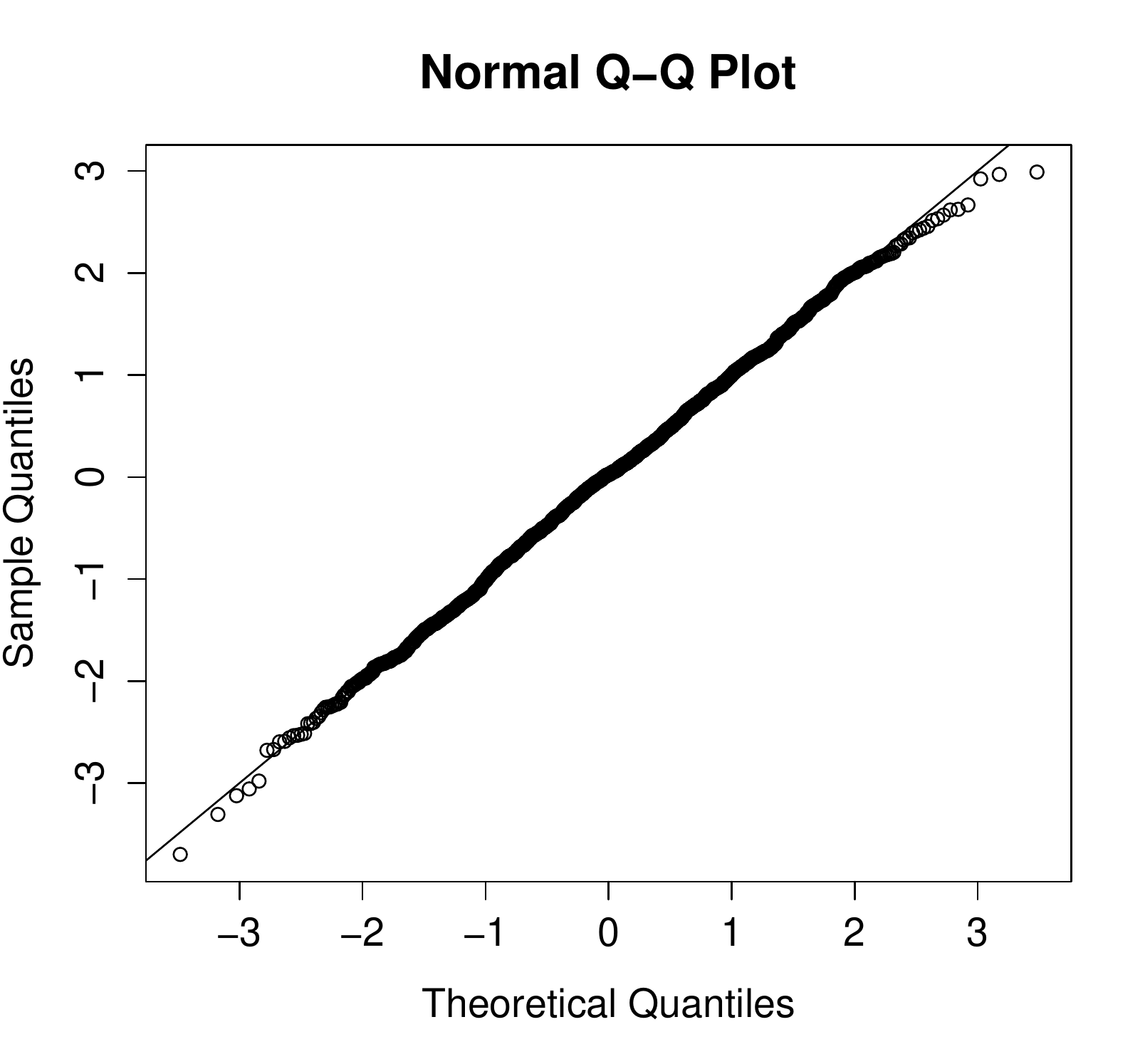} 
    \caption{$\omega=1, T=50$} 
    %\vspace{4ex}
  \end{subfigure}
 \begin{subfigure}[b]{0.32\linewidth}
    \centering
    \includegraphics[width=\linewidth]{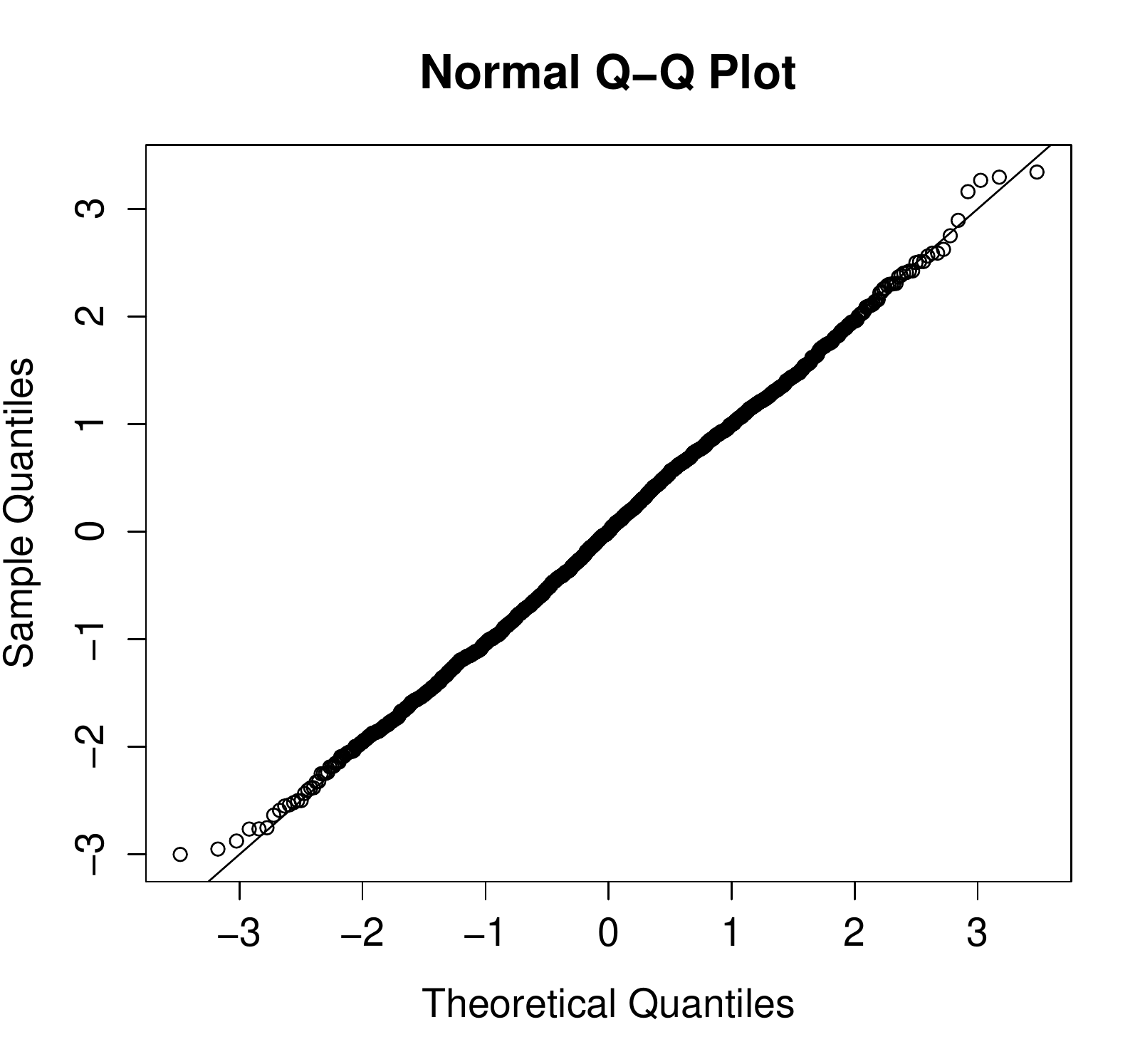} 
    \caption{$\omega=1, T=100$} 
    %\vspace{4ex}
  \end{subfigure}
  \begin{subfigure}[b]{0.32\linewidth}
    \centering
    \includegraphics[width=\linewidth]{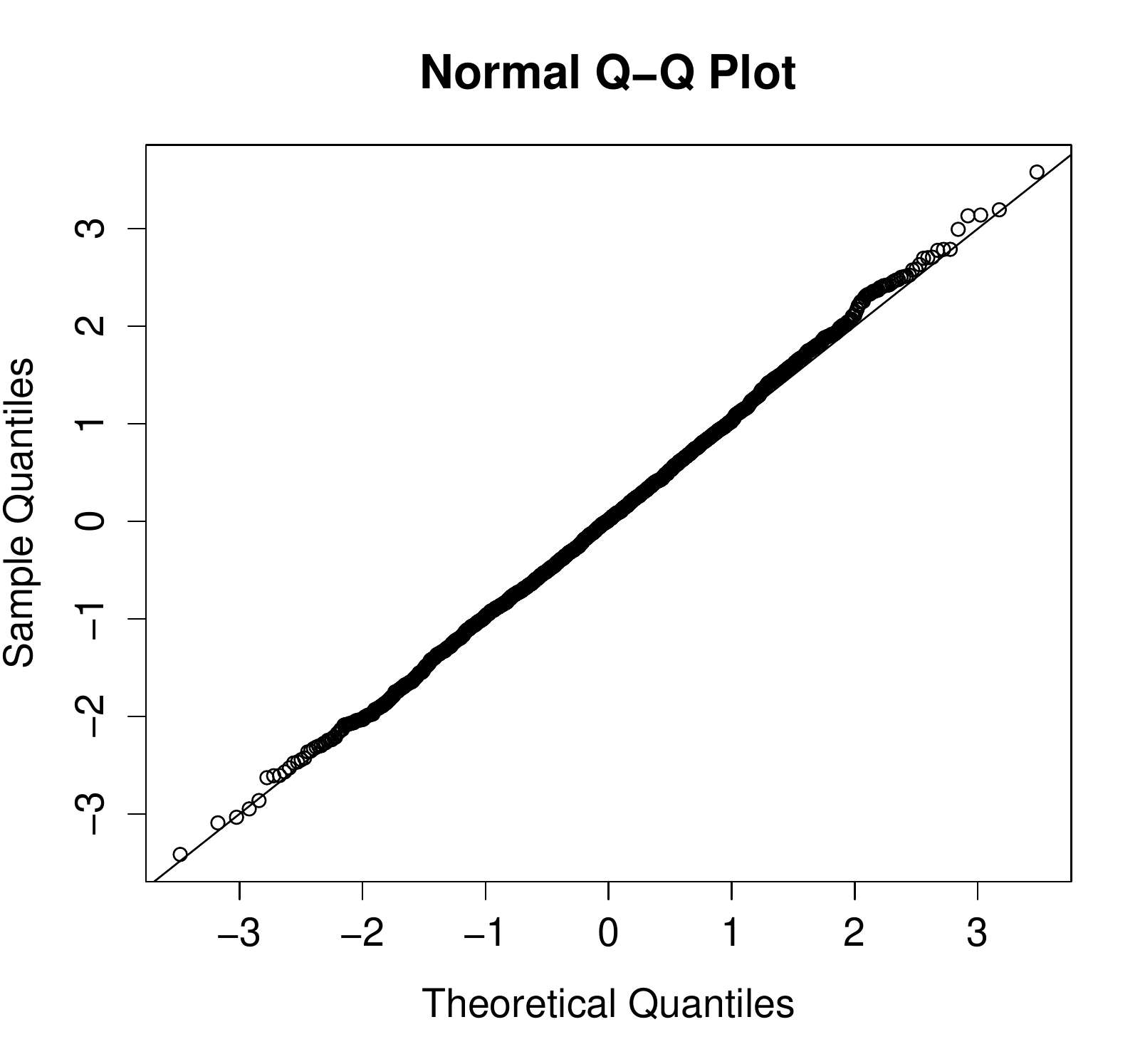} 
    \caption{$\omega=10, T=10$}  
    %\vspace{4ex}
  \end{subfigure}%% 
 \begin{subfigure}[b]{0.32\linewidth}
    \centering
    \includegraphics[width=\linewidth]{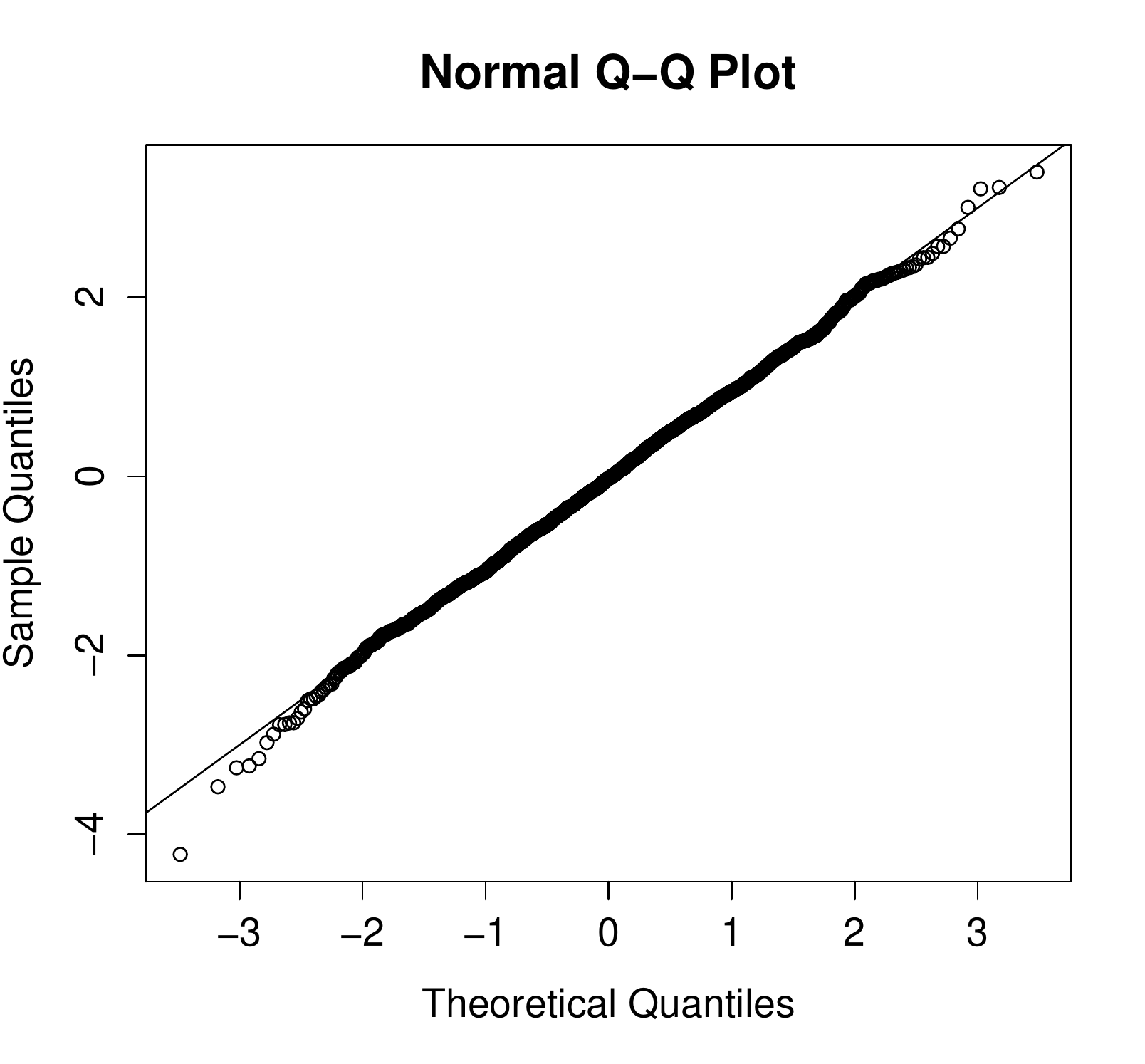} 
    \caption{$\omega=10, T=50$} 
    %\vspace{4ex}
  \end{subfigure}
 \begin{subfigure}[b]{0.32\linewidth}
    \centering
    \includegraphics[width=\linewidth]{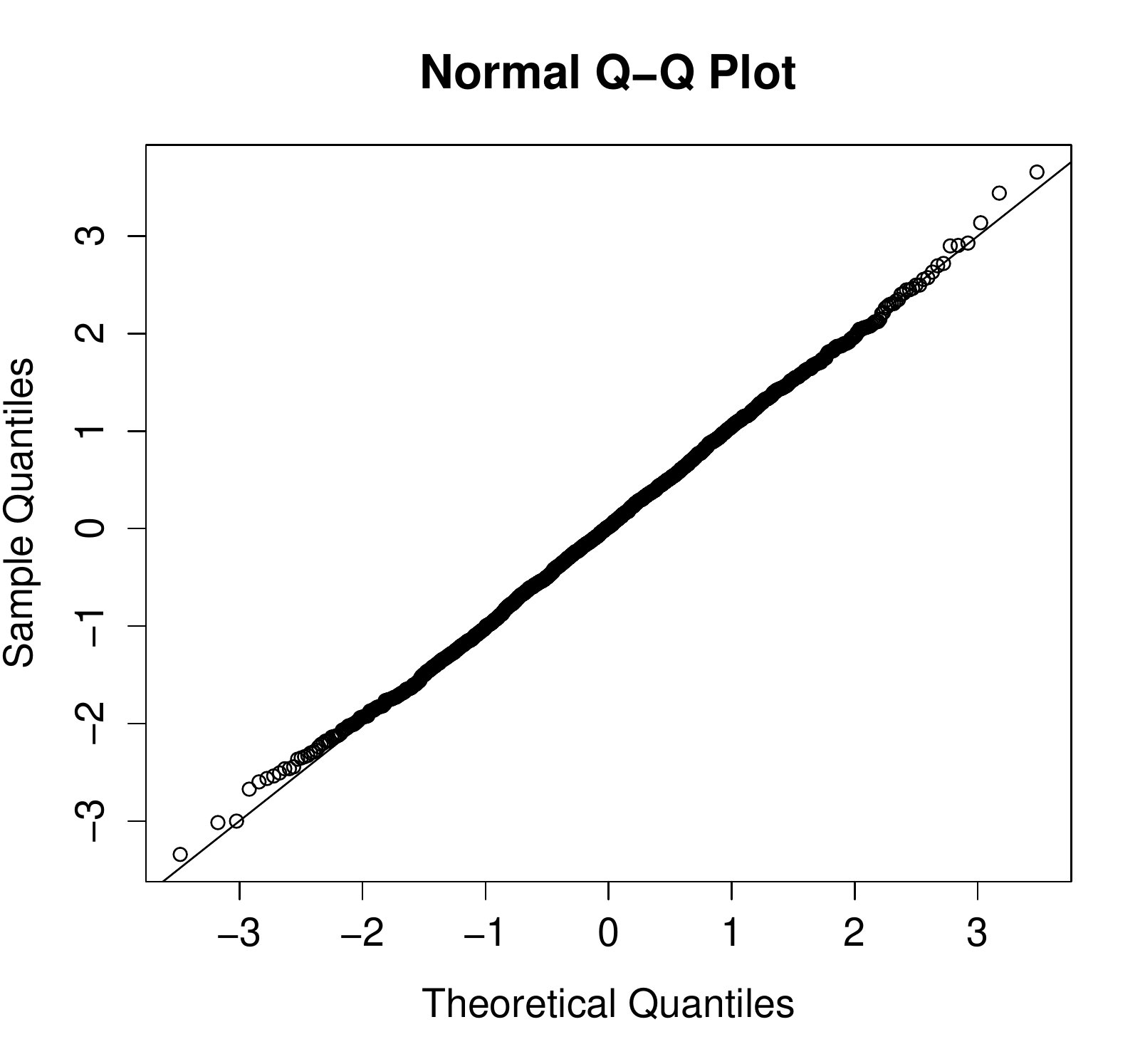} 
    \caption{$\omega=10, T=100$} 
    %\vspace{4ex}
  \end{subfigure}
  \caption{Normal QQ plots for the real part of the truncated Fourier transform of the simulated CARMA(2,1) processes driven by a two sided Poisson process for the frequencies $0, 0.1, 1 , 10$ (rows) and time horizons/maximum non-equidistant grid sizes $10/0.1, 50/0.05, 100/0.01$ (columns). The theoretical quantiles are coming from the (limiting) law described in Theorem \ref{thm:ApproxTFTDoubleLimitDistribution}. }\label{plot:QQCARMAPois} 
\end{figure}

\begin{figure}[tp]    
  \begin{subfigure}[b]{0.32\linewidth}
    \centering
    \includegraphics[width=\linewidth]{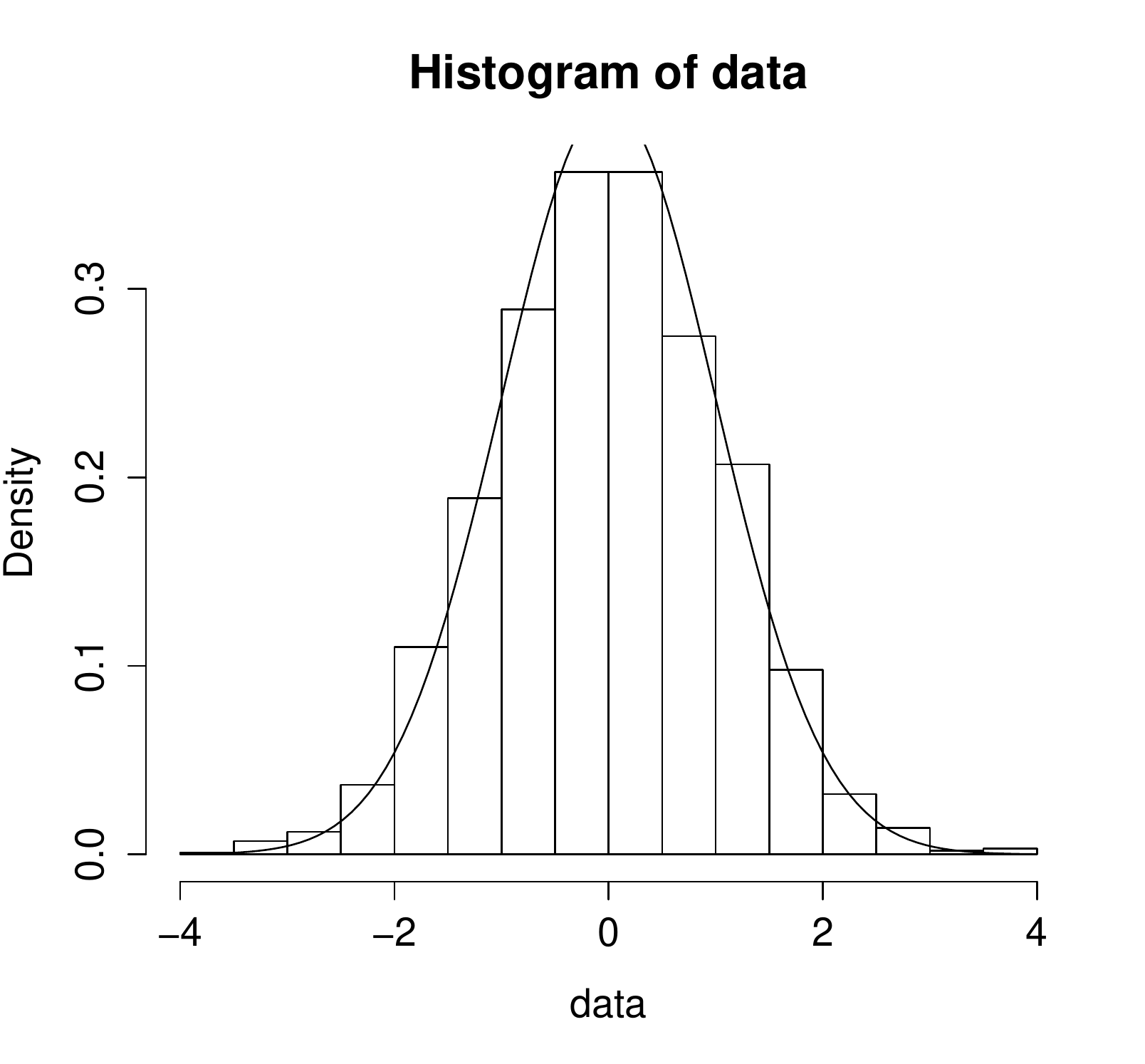} 
    \caption{ $\omega=0, T=10$} 
  
    %\vspace{4ex}
  \end{subfigure}%% 
 \begin{subfigure}[b]{0.32\linewidth}
    \centering
    \includegraphics[width=\linewidth]{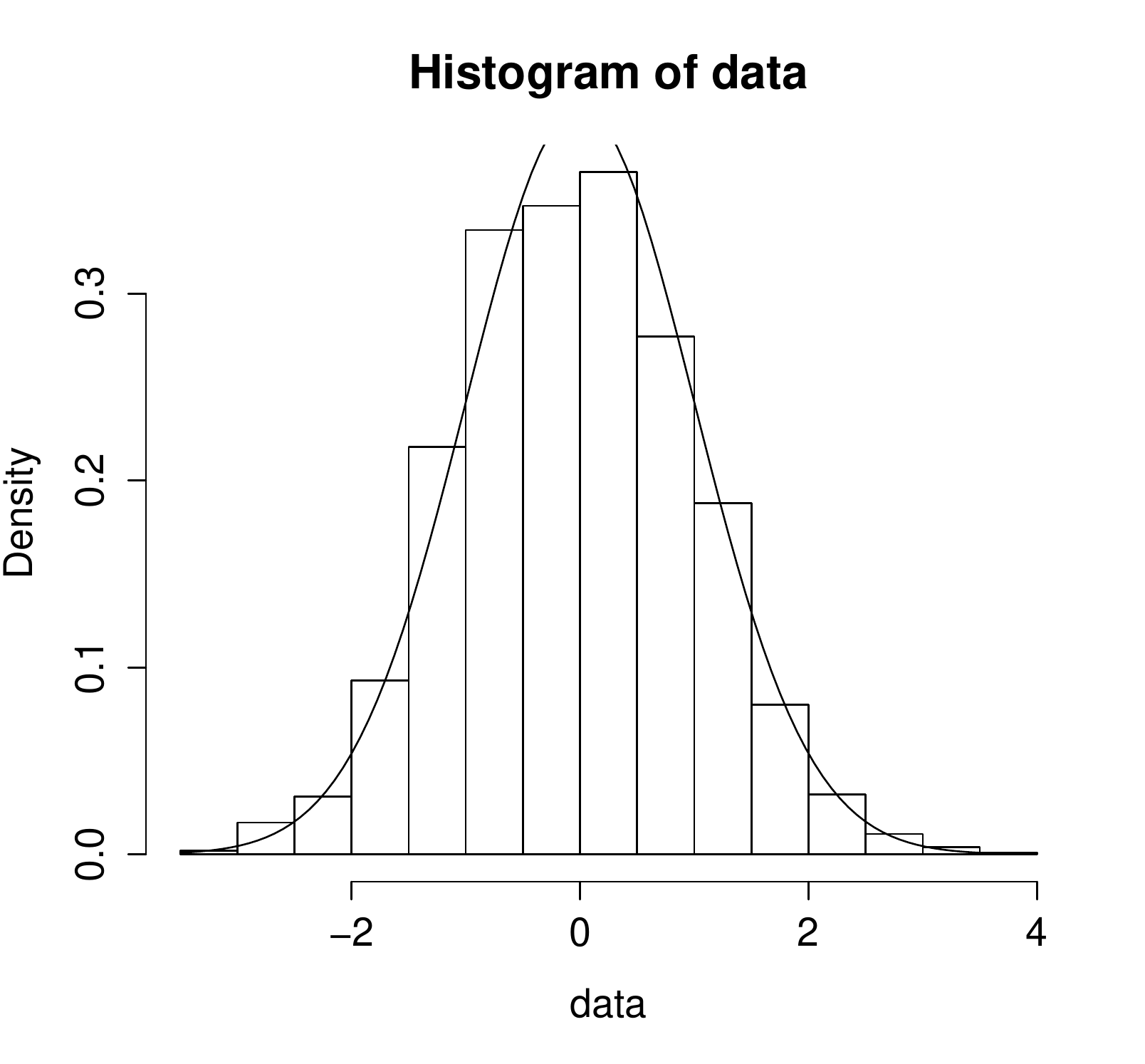} 
    \caption{$\omega=0, T=50$} 
    %\vspace{4ex}
  \end{subfigure}
 \begin{subfigure}[b]{0.32\linewidth}
    \centering
    \includegraphics[width=\linewidth]{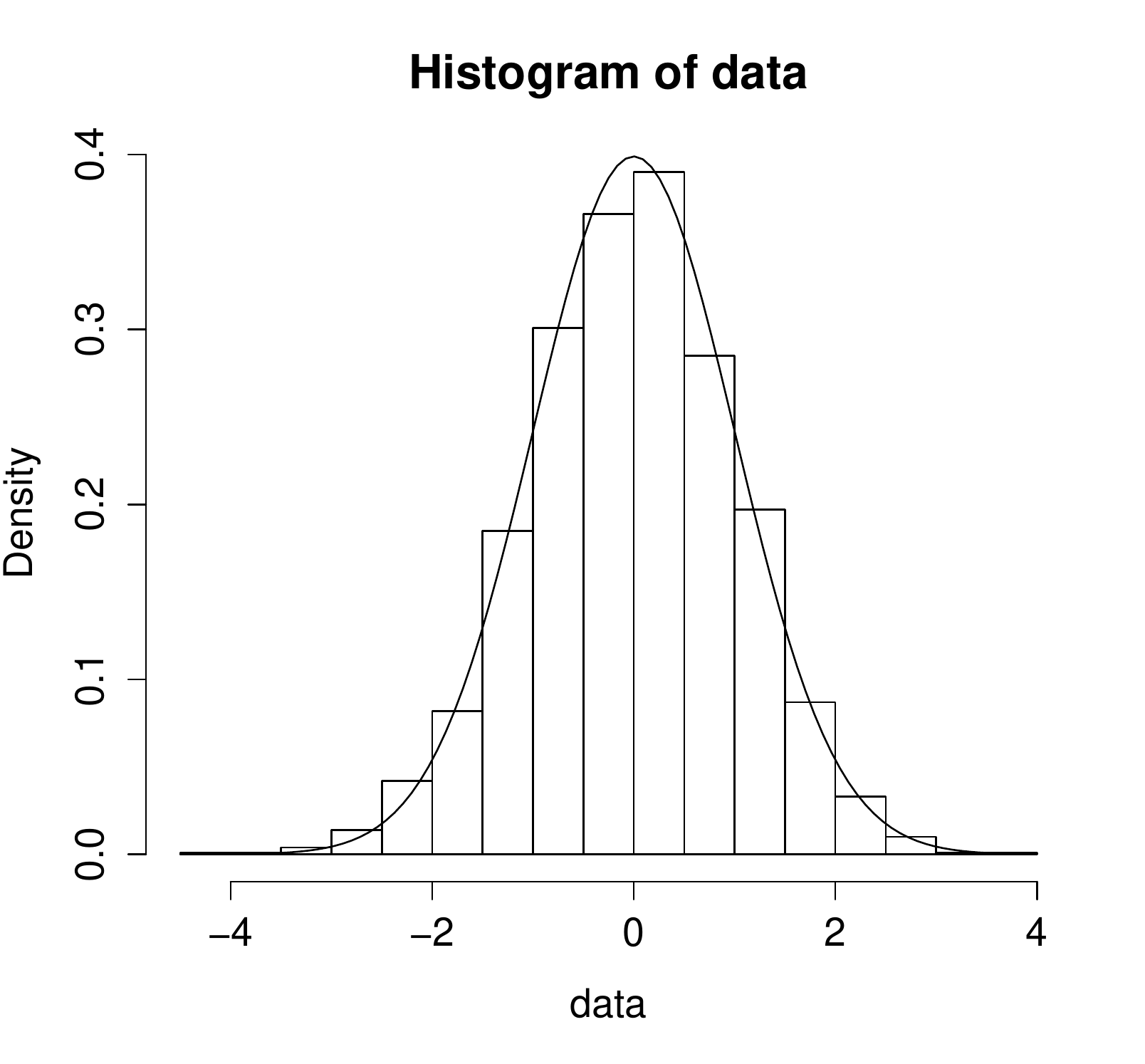} 
    \caption{$\omega=0, T=100$} 
    %\vspace{4ex}
  \end{subfigure}
  \begin{subfigure}[b]{0.32\linewidth}
    \centering
    \includegraphics[width=\linewidth]{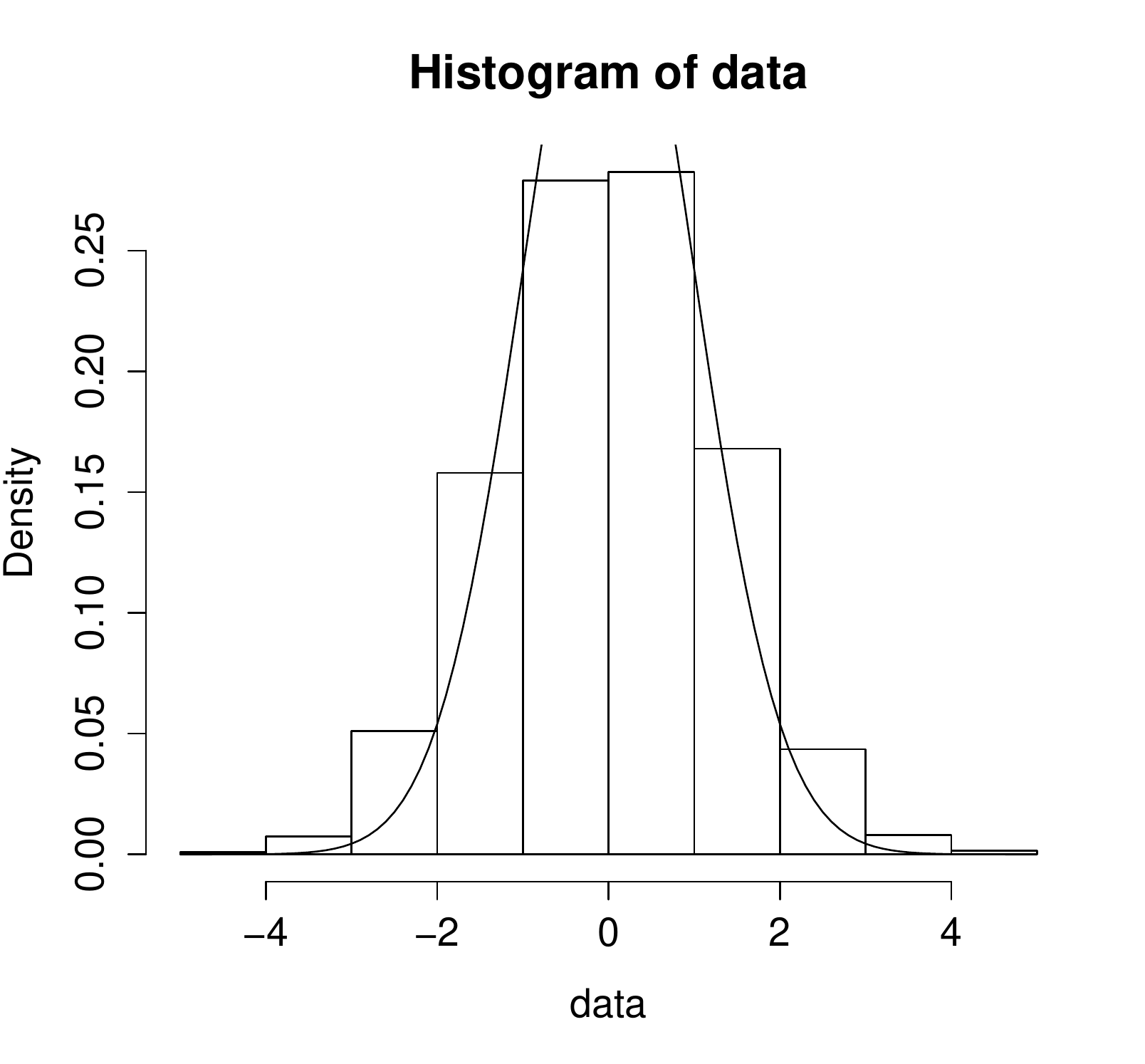} 
    \caption{$\omega=0.1, T=10$} 
    %\vspace{4ex}
  \end{subfigure}%% 
 \begin{subfigure}[b]{0.32\linewidth}
    \centering
    \includegraphics[width=\linewidth]{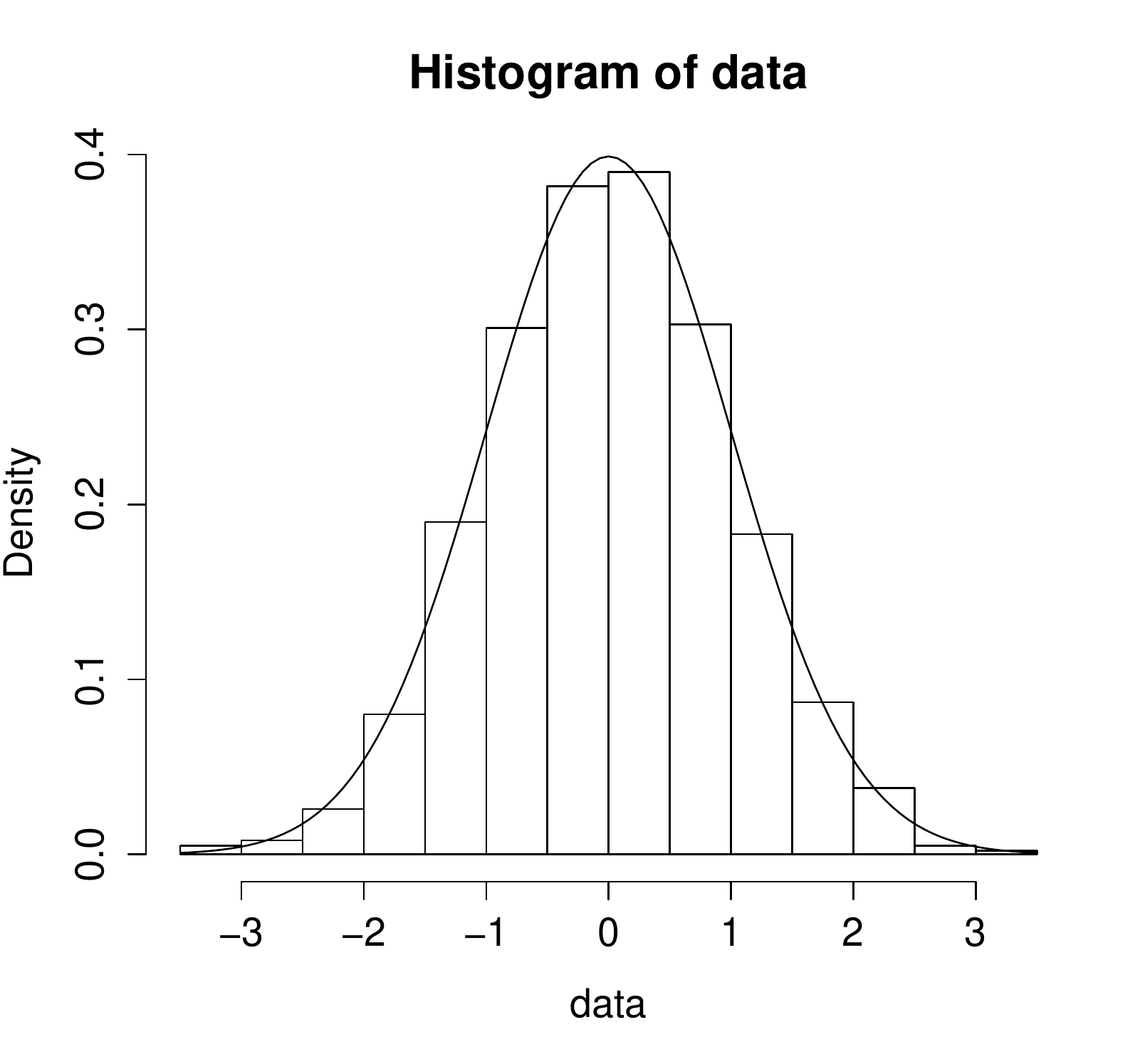} 
    \caption{$\omega=0.1, T=50$} 
    %\vspace{4ex}
  \end{subfigure}
 \begin{subfigure}[b]{0.32\linewidth}
    \centering
    \includegraphics[width=\linewidth]{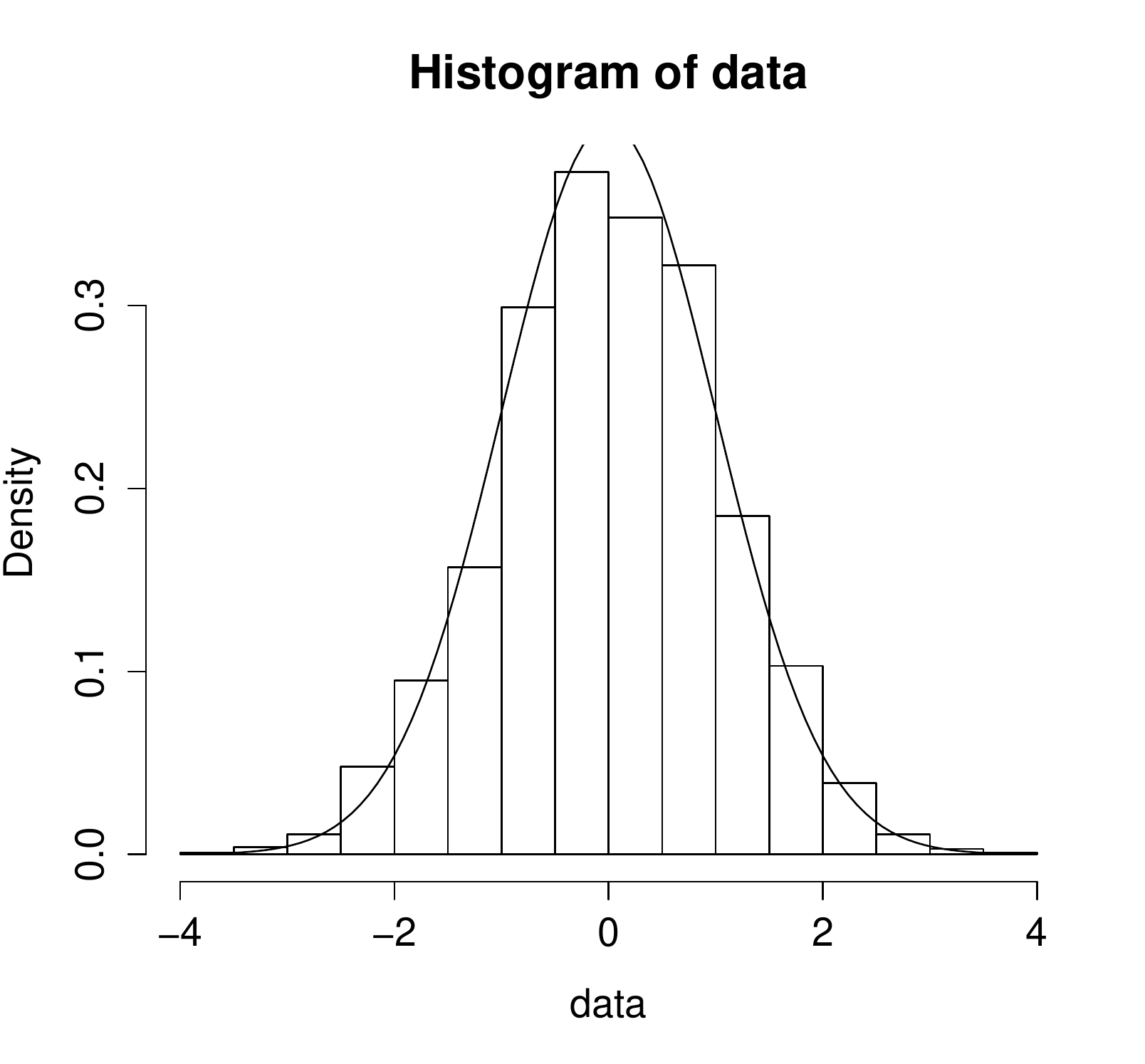} 
    \caption{$\omega=0.1, T=100$} 
    %\vspace{4ex}
  \end{subfigure}
  \begin{subfigure}[b]{0.32\linewidth}
    \centering
    \includegraphics[width=\linewidth]{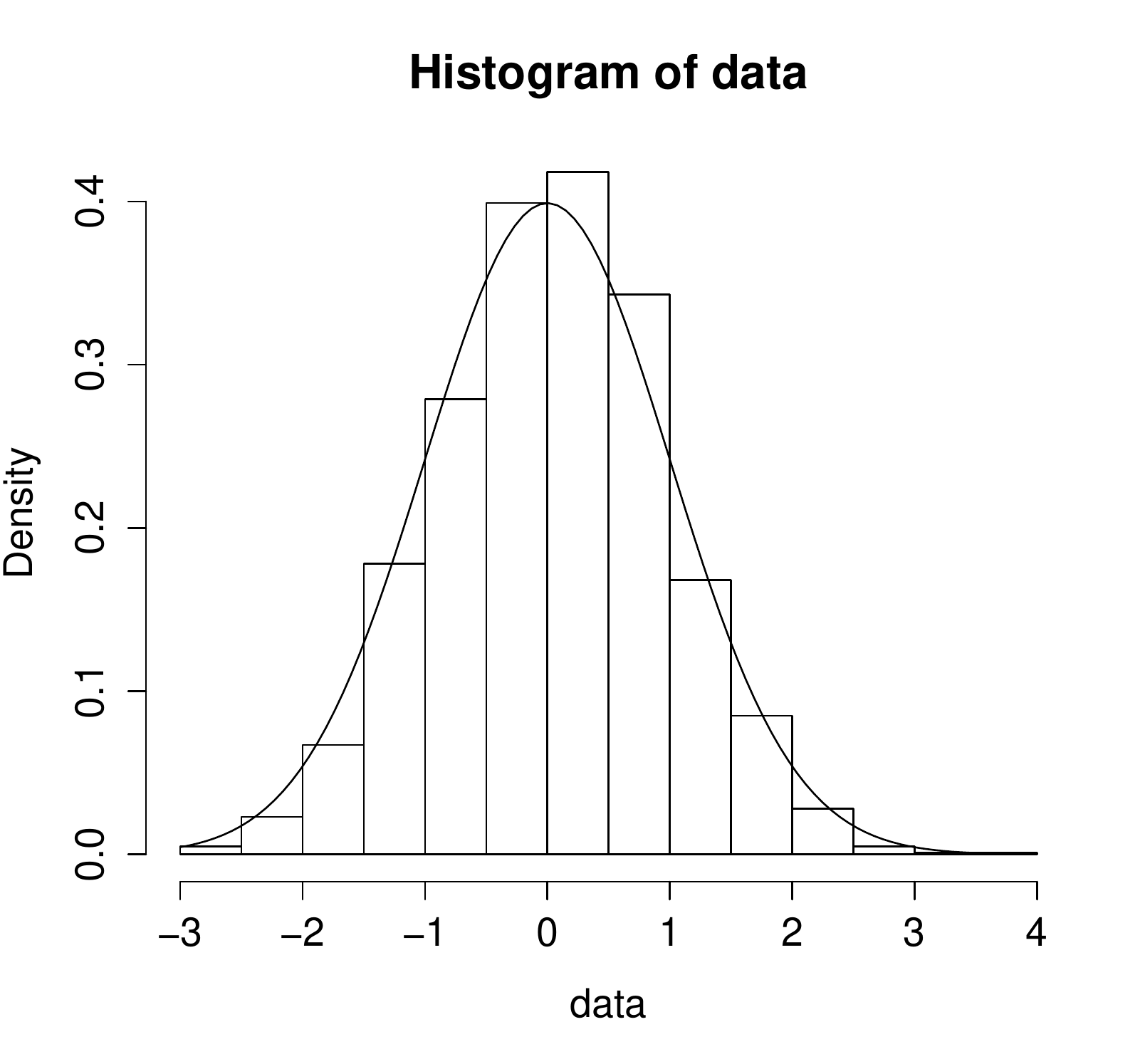} 
    \caption{$\omega=1, T=10$} 
    %\vspace{4ex}
  \end{subfigure}%% 
 \begin{subfigure}[b]{0.32\linewidth}
    \centering
    \includegraphics[width=\linewidth]{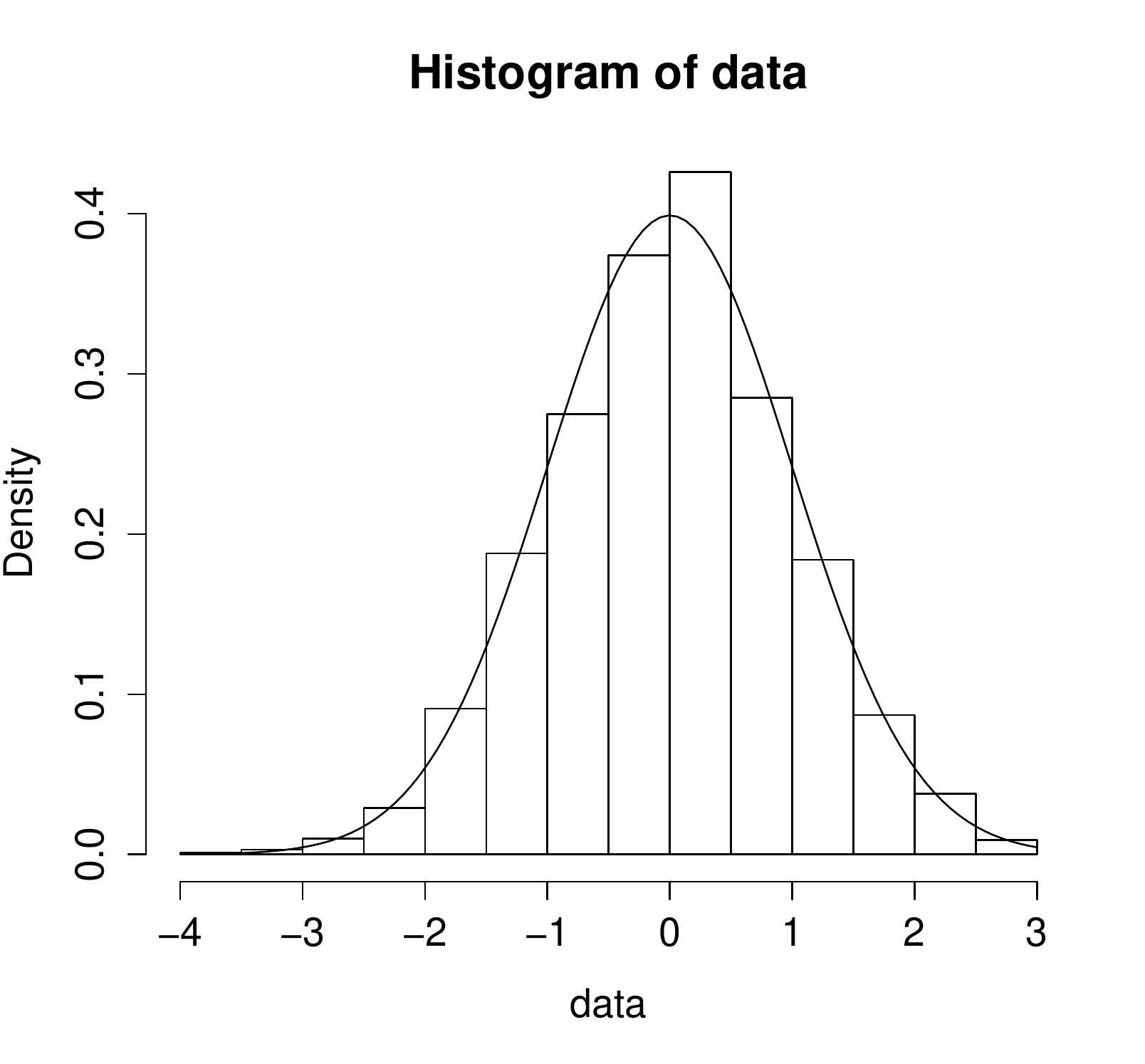} 
    \caption{$\omega=1, T=50$} 
    %\vspace{4ex}
  \end{subfigure}
 \begin{subfigure}[b]{0.32\linewidth}
    \centering
    \includegraphics[width=\linewidth]{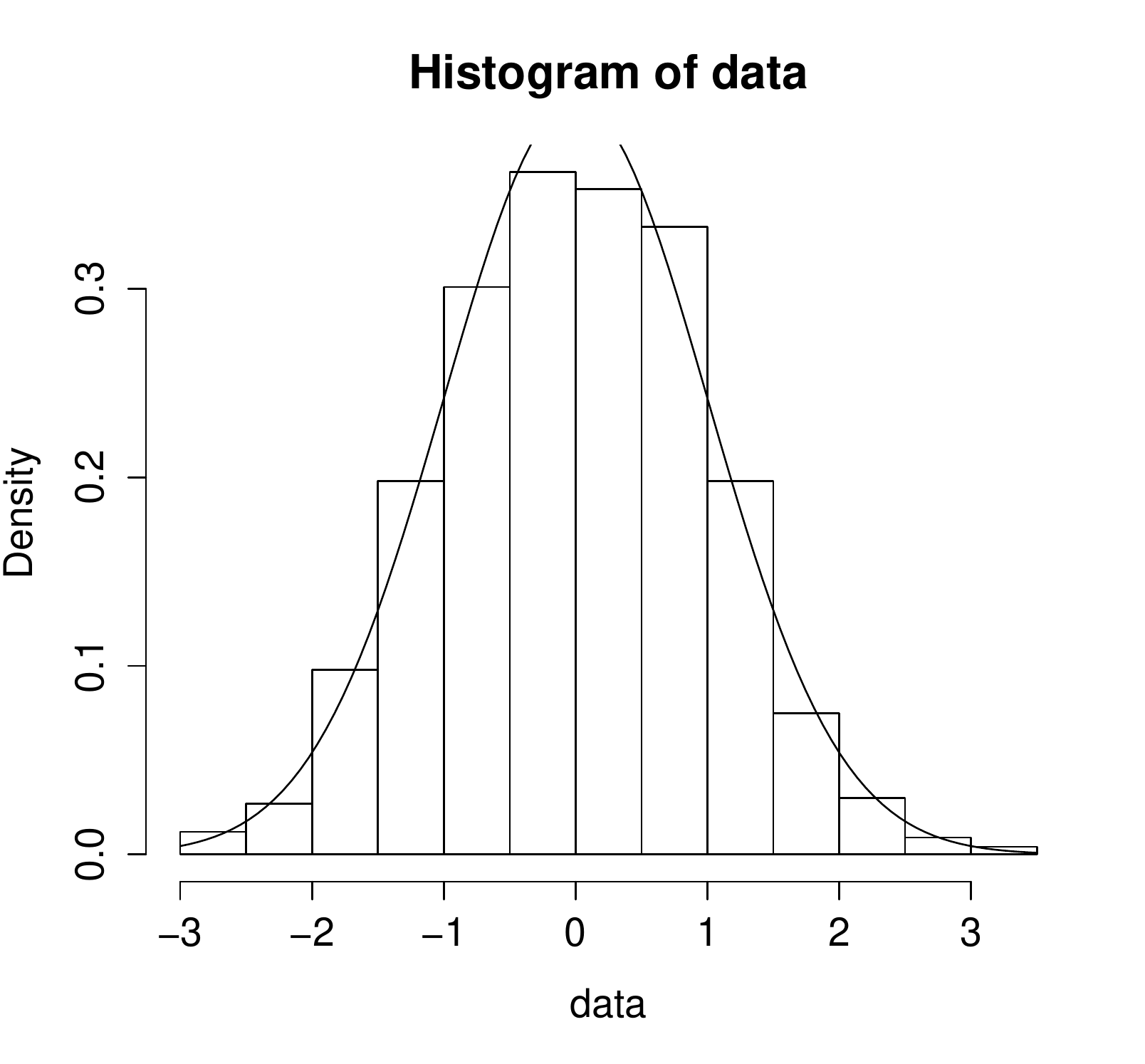} 
    \caption{$\omega=1, T=100$} 
    %\vspace{4ex}
  \end{subfigure}
  \begin{subfigure}[b]{0.32\linewidth}
    \centering
    \includegraphics[width=\linewidth]{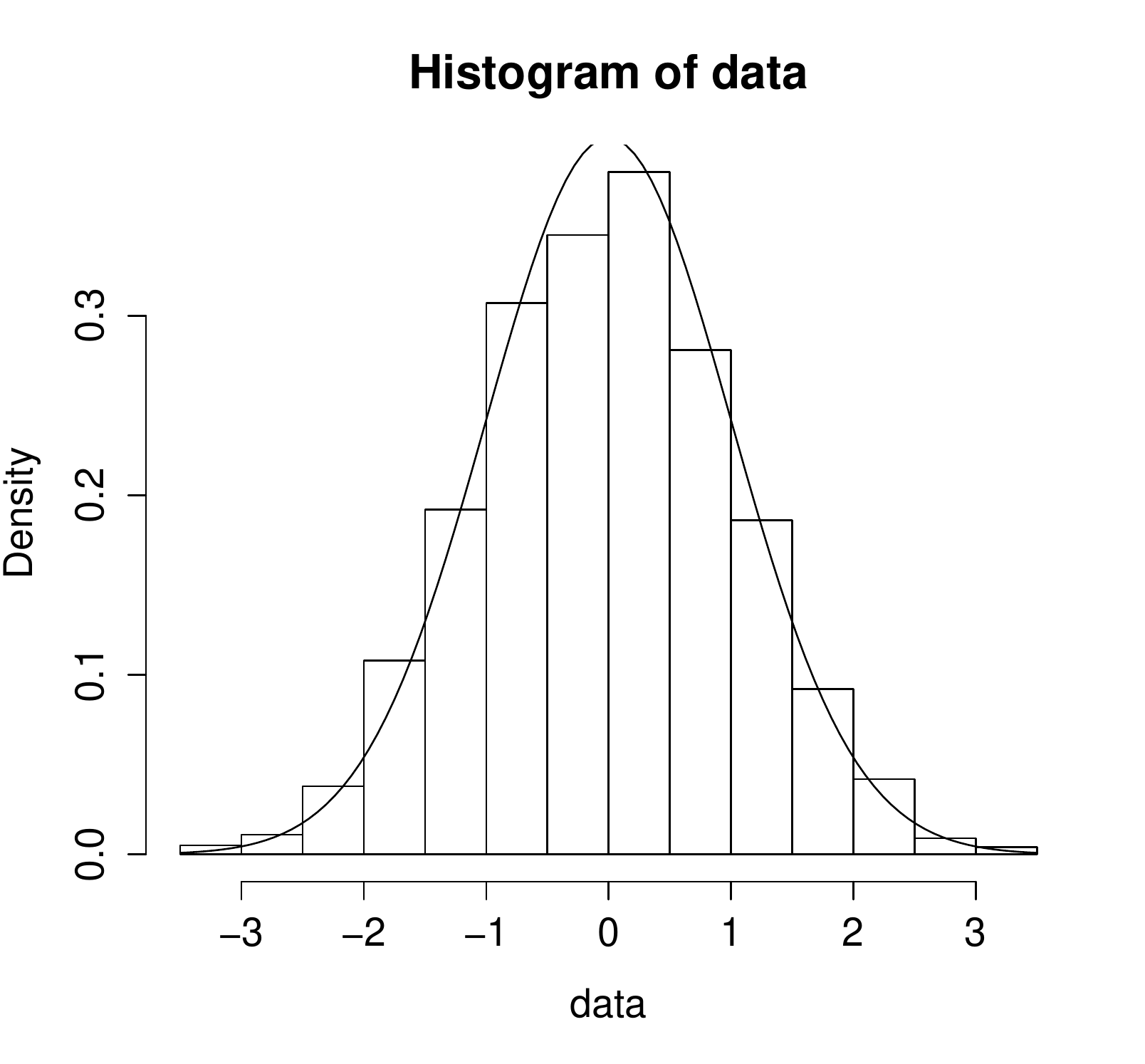} 
    \caption{$\omega=10, T=10$}  
    %\vspace{4ex}
  \end{subfigure}%% 
 \begin{subfigure}[b]{0.32\linewidth}
    \centering
    \includegraphics[width=\linewidth]{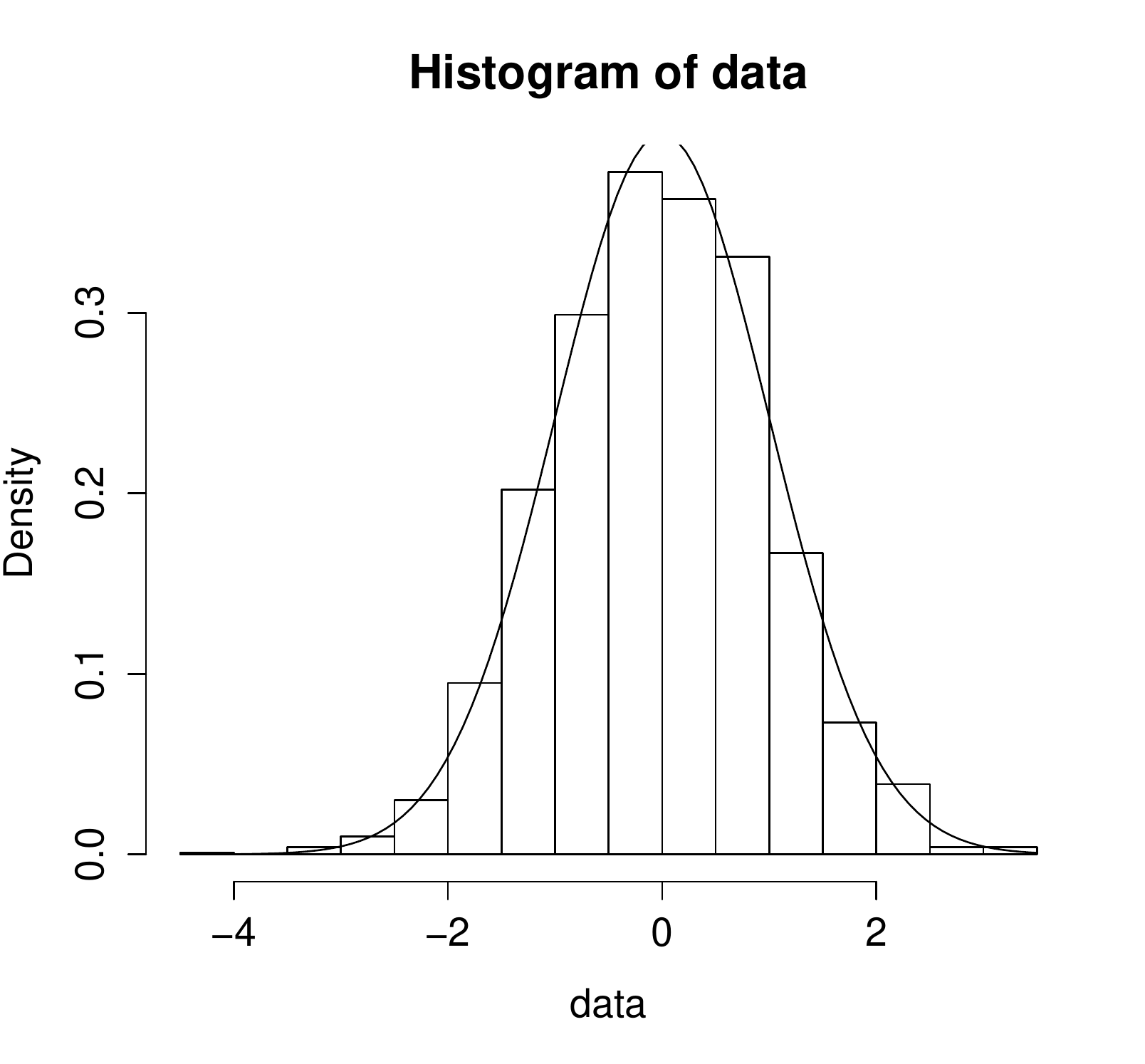} 
    \caption{$\omega=10, T=50$} 
    %\vspace{4ex}
  \end{subfigure}
 \begin{subfigure}[b]{0.32\linewidth}
    \centering
    \includegraphics[width=\linewidth]{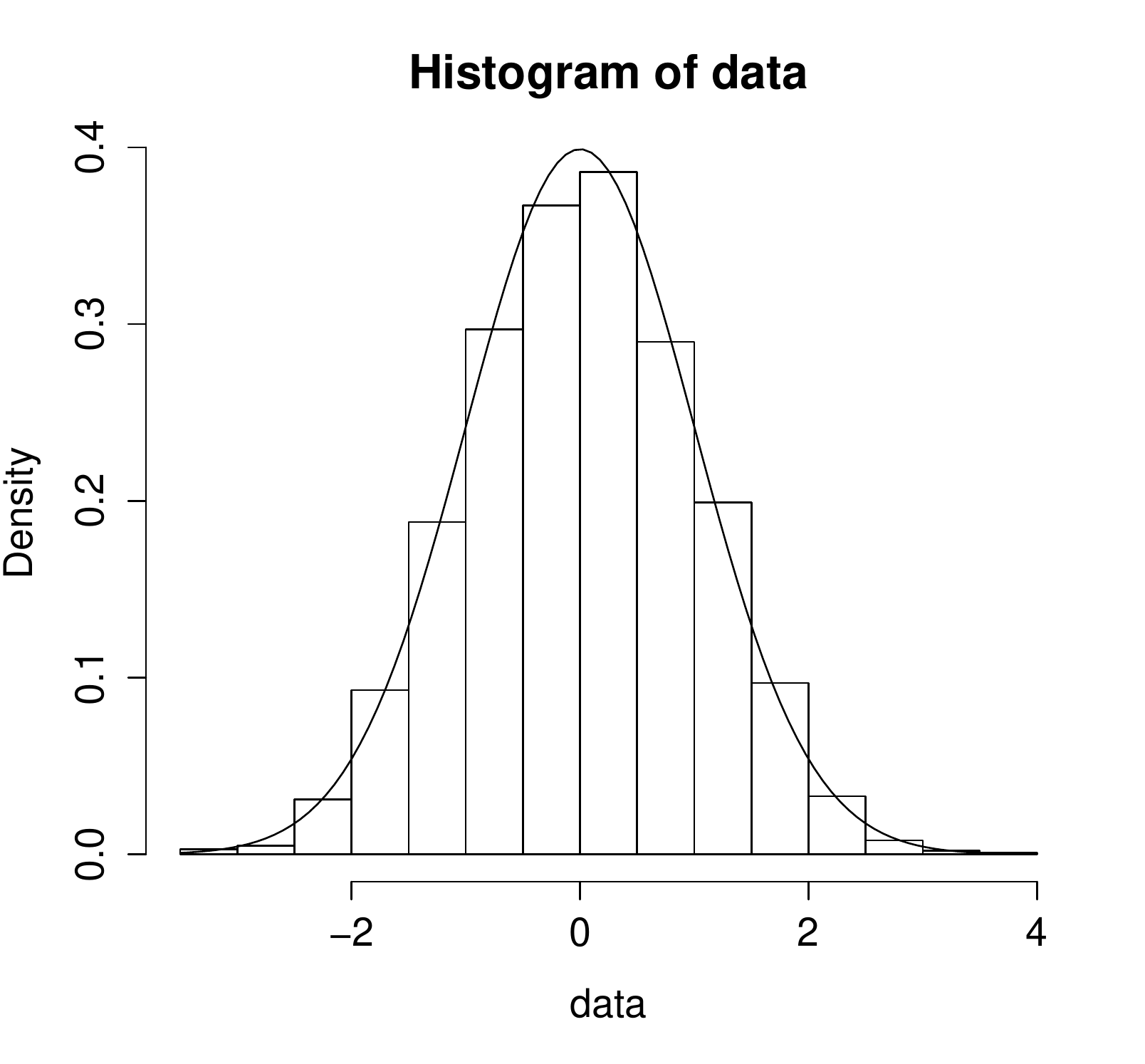} 
    \caption{$\omega=10, T=100$} 
    %\vspace{4ex}
  \end{subfigure}
  \caption{Histograms and limiting density for the real part of the truncated Fourier transform of the simulated CARMA(2,1) processes driven by a two-sided Poisson process for the frequencies $0, 0.1, 1 , 10$ (rows) and time horizons/maximum non-equidistant grid sizes $10/0.1, 50/0.05, 100/0.01$ (columns)}\label{plot:HistCARMAPois} 
\end{figure}

\begin{example}
We consider the $\operatorname{CARMA}(2,1)$ model. We have
$$ 
\A:=\[ 
\begin{matrix}
0 &1 \\
-a_2& -a_1\\
\end{matrix}
 \], \quad \b:=\[ 
\begin{matrix}
b_0  \\
1\\
\end{matrix}
 \],\quad \e:=\[ 
\begin{matrix}
0  \\
1 \\
\end{matrix}
 \],\quad \X_t:=\[ 
\begin{matrix}
X(t)  \\
X^{(1)}(t) \\
\end{matrix}
 \]
$$
The autoregressive and moving-average polynomials are of the form
$$
a(z)=z^2+a_1z+a_2,\quad
b(z)=z+b_0.
$$
We have
$$
\frac{b(i\omega)}{a(i\omega)}=\frac{ i\omega+b_0}{(i\omega)^2+(i\omega)a_1+a_2},\quad 
f(\omega)=\frac{\sigma^2}{2\pi}\left|\frac{b(i\omega)}{a(i\omega)}\right|^2=\frac{\sigma^2}{2\pi}\frac{b_0^2+\omega^2}{\omega^4+(a_1^2-2a_2)\omega^2+a_2^2}
$$
 For the simulation procedure we take $[b_0,b_1,a_1,a_2]=[1,1,1,2]$.

QQ-plots showing the results for 2000 simulated paths for the four different frequencies and three different combinations of time horizon and maximum grid width can be found in Figures \ref{plot:QQCARMANormal}, \ref{plot:QQCARMAVG} and \ref{plot:QQCARMAPois} for the driving L\'evy process being a standard Brownian motion, a Variance Gamma and a two-sided Poisson process, respectively. Likewise, Figures \ref{plot:HistCARMANormal}, \ref{plot:HistCARMAVG} and \ref{plot:HistCARMAPois} show corresponding histograms.

\end{example}


The simulation results seem to indicate the following. 

In Figure \ref{plot:QQCARNormal} we notice at first a pretty good fit of the empirical quantiles from the simulations with the theoretical ones of the asymptotic distribution across all time horizons and frequencies. Looking more carefully, the fit in the tails clearly improves when the time horizon/fineness of the grid increases, but it is never  bad. For the longest time horizon and finest grid the fit is clearly very good. Of course, it should not be forgotten that in this case the distribution of the (trapezoidal approximation of the) truncated Fourier transform is always exactly Gaussian and not only asymptotically. When looking across the non-zero frequencies one notes that for the shortest time horizon the quantiles for the smallest frequency $0.1$ appear to lie on a line which is somewhat different from the line of the theoretical quantiles. This indicates that the quantiles of the simulated paths come from a normal distribution, but one with a different variance then the asymptotic one. It is no surprise that this occurs for the lowest frequency and the smallest time interval, as for low frequencies one observes -- regardless of the fineness of the sampling -- the fewest full cycles over a time interval of fixed length. For this combination of time horizon and frequency we see only one full cycle.

Turning to Figure \ref{plot:QQCARVG}, we first notice that the fit in the tails improves again clearly with  increasing $(T, 1/h_{\max})$. Especially, for the highest $(T, 1/h_{\max})$ one sees that the fit in the tails is a bit worse now for a driving Variance Gamma process compared with the driving Brownian motion in Figure \ref{plot:QQCARNormal}. Of course, now the simulated values are indeed only asymptotically following a Gaussian distribution. Looking at the different non-zero frequencies one again sees that the fit improves for the higher frequencies. Most notably for the lowest frequency one sees for the smallest $T=10$ again that the points do seem to lie on a straight line in the normal QQ-plot, but one with a different slope than for the theoretical quantiles. Hence, the variance is clearly different from the asmpytotic one. Obviously, this effect is now more pronounced than in the case of the driving Brownian motion. 

Moving on to the case of the driving process being a two-sided Poisson process in Figure  \ref{plot:QQCARPois} we first of all note that again the fit in particular in the tails clearly improves with increasing $(T, 1/h_{\max})$. For the highest $(T, 1/h_{\max})$ the simulated and theoretical asymptotic quantiles agree again extremely well. Again it is certainly a bit worse than in the case of a driving Brownian motion, but it seems to be very similar to the Variance Gamma case, although maybe for frequency $0$ the agreement of the quantiles is slightly worse. Turning to the behaviour across non-zero frequencies, we see again that the quantiles are closer for higher frequencies and that for the smallest non-zero frequency and time horizon the empirical quantiles seem to be in line with a normal distribution with a somewhat different variance compared to the asymptotic one. The size of this effect seems to be rather similar to the Variance Gamma case. Interestingly, also at frequency $0$ the QQ-plot seems to indicate for $T=10$ that the empirical quantiles are close to the ones of a normal distribution with a slightly different variance than the asymptotic one.

To summarize the simulation study in the CAR(1)/OU-type case we can clearly conclude that the asymptotic distribution result approximates the finite-sample  distribution of the trapezoidal approximation of the truncated Fourier transform in our simulations very well and that the convergence to the asymptotic distribution is fast. For small frequencies, especially when one has about one full cycle or less over the time horizon considered, one has to be careful, as  then the distribution tends to be somewhat different from the asymptotic one for good reasons. This effect seems to be more pronounced when one considers a L\'evy process with jumps compared to a Brownian motion. In general the quality of the approximation of the simulated quantiles by the asymptotic ones is somewhat better in the case of a Brownian motion than in a pure jump process. Comparing the driving jump processes, the finite activity rather discrete two-sided Poisson process with the infinite activity Variance Gamma process, we do not see any significant differences. It should be noted that both jump processes are, however, light-tailed in the sense that they have exponential moments. It would not be surprising if this picture changes when considering a really heavily  tailed driving L\'evy process. Note that our theoretical results are valid also in rather heavily-tailed cases. For the asymptotic normality of the (trapezoidal approximation of the) truncated Fourier transform we only needed finite second moments. 

Turning to the simulations of CARMA(2,1) processes, most of the findings of the CAR(1)/\allowbreak OU case remain valid, so we only point out the differences. In the case of a driving Brownian motion, depicted in Figure \ref{plot:QQCARMANormal}, the only difference seems to be that for $T=10$ and $\omega=0.1$ the empirical quantiles are now appearing to lie on a line farther away from the theoretical quantiles which implies that in the CARMA(2,1) case the variance in the simulations is clearly farther away from the asymptotic one than in the OU case. The same applies for the Variance Gamma case of Figure \ref{plot:QQCARMAVG} and the two-sided Poisson case of Figure \ref{plot:QQCARMAVG}. On top of the QQ plots we now also provide histograms in Figures \ref{plot:HistCARMANormal}, \ref{plot:HistCARMAVG} and \ref{plot:HistCARMAPois}, respectively, together with plots of the limiting normal density. To us  it seems very hard to see the convergence to normality with increasing $T$ in the histograms, which reflects the fact that it is essentially  the tails which need to converge and they are much clearer visible in the QQ plots than in histograms. It is also not easy to see in them that for $\omega=0.1,\,T=10$ the variance of the simulated values is different from the asymptotic theoretical one. The only thing one notices is that for $\omega=0.1,\,T=10$ the histogram routine of R tends to use very different bins than in all the other cases. Note that all histograms were obtained using the default parameters of the \texttt{hist} function in R, so the binning was done by the standard automatic selection to give ``nice'' histograms. Hence, from our simulations of CARMA(2,1) processes we can conclude that the orders of the CARMA processes and the particular autoregressive and moving average parameters  appear not to really matter for the (qualitative) behaviour of the (trapezoidal approximation of the) truncated Fourier transform.

\section{Conclusion and Outlook}
We have obtained an asymptotic normality result for the (trapezoidal approximation of the) truncated Fourier transform under essentially minimal assumptions (i.e. second moments) and seen via a simulation study that this result approximates the finite sample behaviour very well, unless the frequency is too low compared to the length of the considered time interval. This suggests clearly that it should be very promising to develop statistical inference techniques for non-equidistantly sampled CARMA processes by considering continuous observation techniques and using numerical approximation schemes to compute the quantities of interest based on the observed non-equidistant data. The appropriate set-up to get convergence and asymptotic distribution results is to send the time horizon to infinity and to send at the same time the maximum distance of observation time points to zero.

Based on our results in this paper it seems natural to locally smooth the  trapezoidal approximation of the truncated Fourier transform to get consistent estimators of the spectral density and to use it in a Whittle type estimator for the AR and MA parameters. Considering this is beyond the scope of the present paper.

\section*{Acknowledgements}
The authors gratefully acknowledge the support of Deutsche
Forschungsgemeinschaft (DFG) by research grant STE 2005/1-2.

The authors would like to cordially thank W\l{}odzimierz Fechner for careful reading of the manuscript and his helpful remarks.

%It is worth to underlying that the truncated Fourier transform itself is not a good estimator of the spectral density. In order to perform a reasonable estimation procedure one should investigate the smoothed version.

\end{document}